\documentclass[12pt]{article}

\usepackage{amsmath}
\usepackage{amssymb}
\usepackage{amsfonts}
\usepackage{a4}
\usepackage{graphicx}
\usepackage[abs]{overpic}
\usepackage{caption}
\usepackage{wrapfig}
\usepackage{float}
\usepackage{graphicx}
\usepackage{multicol}

\newtheorem{theorem}{Theorem}[section]
\newtheorem{lemma}[theorem]{Lemma}

\newtheorem{corollary}[theorem]{Corollary}

\newtheorem{rem}[theorem]{Remark}

\newcommand{\Proof}{\par\noindent{\em Proof. }}
\newcommand{\eop}{\nopagebreak\hspace*{\fill}$\Box$\smallskip}

\newcommand{\N}{\Bbb N}
\newcommand{\Z}{\Bbb Z}
\newcommand{\R}{\Bbb R}

\def\id{\mathbf{id}}

\def\eps{\varepsilon}

\def\e{\mathbf{e}}

\def\dist{\operatorname{dist}}
\def\Xint#1{\mathchoice
   {\XXint\displaystyle\textstyle{#1}}%
   {\XXint\textstyle\scriptstyle{#1}}%
   {\XXint\scriptstyle\scriptscriptstyle{#1}}%
   {\XXint\scriptscriptstyle\scriptscriptstyle{#1}}%
   \!\int}
\def\XXint#1#2#3{{\setbox0=\hbox{$#1{#2#3}{\int}$}
     \vcenter{\hbox{$#2#3$}}\kern-.5\wd0}}

\def\dashint{\Xint-}

\usepackage{color}

\numberwithin{equation}{section}

\begin{document}

\begin{center}
\begin{Large}
{\bf {A Korn-Poincar\'e-type inequality for special functions of bounded deformation}}
\end{Large}
\end{center}

\begin{center}
\begin{large}
Manuel Friedrich\footnote{Faculty of Mathematics, University of Vienna, 
Oskar-Morgenstern-Platz 1, 1090 Vienna, Austria. {\tt manuel.friedrich@univie.ac.at}}
\end{large}
\end{center}

\begin{center}
\today
\end{center}
\bigskip

\begin{abstract}
We present a Korn-Poincar\'e-type inequality in a planar setting which is in the spirit of the Poincar\'e inequality in SBV due to De Giorgi, Carriero, Leaci (see \cite{DeGiorgiCarrieroLeaci:1989}). We show that for each function in SBD$^2$ with small jump set one can find a modification which differs from the original  displacement field only on   a finite union of rectangles with small measure such that the distance of the modification from a suitable infinitesimal rigid motion can be controlled by an appropriate combination of the elastic and the surface energy. In particular, the result can be used to obtain compactness estimates for functions of bounded deformation.
\end{abstract}
\bigskip

\begin{small}

\end{small}

\tableofcontents

\section{Introduction}

The propagation of crack has been studied in the realm of linearized elasticity since the seminal work of Griffith (see \cite{Griffith:1921}) and led to a lot of realistic applications in engineering. Also from a mathematical point of view the theory is well developed (see \cite{Ambrosio-Coscia-Dal Maso:1997, Bellettini-Coscia-DalMaso:98}) and adopted in many recent works in applied analysis (see e.g. \cite{Bellettini-Coscia-DalMaso:98, Bourdin-Francfort-Marigo:2008, Chambolle:2003, Chambolle:2004, Focardi-Iurlano:13, Iurlano:13, SchmidtFraternaliOrtiz:2009}). A natural framework for the investigation of fracture models in a geometrically linear setting is given by the space of functions of bounded deformation, denoted by $BD(\Omega,\R^d)$, which consists of all functions $u \in L^1(\Omega,\R^d)$ whose  symmetrized distributional derivative $Eu := \frac{1}{2}((Du)^T + Du)$ is a finite $\R^{d \times d}_{\rm sym}$-valued Radon measure. To study problems in fracture mechanics with variational methods Francfort and Marigo \cite{Francfort-Marigo:1998} have introduced energy functionals which are essentially of the form
\begin{align}\label{rig-eq: general energy}
\int_\Omega |e(u)|^2 \,dx + {\cal H}^{d-1}(J_u),
\end{align}
where $u \in SBD^2(\Omega,\R^d)$. (For the definition and properties of this space we refer to Section \ref{rig-sec: sub, bd} below.) These so-called Griffith functionals comprise elastic bulk contributions for the unfractured regions of the body represented by the linear elastic strain $e(u) := \frac{1}{2}(\nabla u^T + \nabla u)$ and surface terms that assign energy contributions on the crack paths comparable to the size of the crack ${\cal H}^{d-1}(J_u)$. 

A major difficulty of these problems is given by the fact that there is no control over the skew-symmetric part of the distributional derivative. Indeed, it would be desirable that uniform bounds  on \eqref{rig-eq: general energy} induce estimates for the strain  $\nabla u$ or the function $u$ itself  and that accordingly suitable compactness results can be derived. However, simple examples (see e.g. \cite{Ambrosio-Coscia-Dal Maso:1997}) show that such properties cannot be inferred in general as the behavior of small pieces being almost or completely separated from the bulk part may not be controlled. On the one hand, this observation particularly implies that SBD is not contained in SBV. (For the definition and properties of SBV we refer to \cite{Ambrosio-Fusco-Pallara:2000}.) On the other hand, it leads to the natural question if certain estimates still hold (1) up to a small exceptional set or (2) after passing to a slightly modified displacement field. The goal of the work at hand is to show that indeed the distance of the function $u$ from an infinitesimal rigid motion can be estimated by an appropriate combination of the energy terms given in \eqref{rig-eq: general energy}. 

The starting point of our analysis is the classical Korn-Poincar\'e inequality in BD (see \cite{Temam:85}) stating that there is a constant $C(\Omega)$ depending only on the domain $\Omega \subset \R^d$ such that
\begin{align}\label{rig-eq: help-korn}
\Vert u - Pu \Vert_{L^{\frac{d}{d-1}}(\Omega)} \le C(\Omega) |Eu|(\Omega)
\end{align}
for all $u \in BD(\Omega,\R^d)$, where $P$ is a linear projection onto the space of infinitesimal rigid motions and $|Eu|$ denotes the total variation of the symmetrized distributional derivative. This is a remarkable result in consideration of the fact that corresponding estimates also involving the  absolutely continuous part of the derivative $\nabla u$ do not hold since  Korn's inequality fails in BD (cf. \cite{ContiFaracoMaggi:2005}).  

It first appears that the inequality is not adapted for problems of the form \eqref{rig-eq: general energy} as in $|Eu|(\Omega)$ the jump height is involved and in \eqref{rig-eq: general energy} we only have control over the size of the crack. However, the main strategy of our approach is to prove that one may indeed find bounds on the jump heights after a suitable modification of the jump set and the displacement field whose total energy \eqref{rig-eq: general energy} almost coincides with the original energy. In particular, \eqref{rig-eq: help-korn} is then applicable and it can be shown that the distance of the modification from an infinitesimal rigid motion can already be controlled in terms of  the elastic energy.

The goal of this work is to prove the following Korn-Poincar\'e-type inequality for SBD$^2$ functions.  For $\mu>0$ let $Q_\mu=(-\mu,\mu)^2$ and by ${\rm diam}(R)$ denote the diameter of a rectangle $R \subset Q_\mu$.

\begin{theorem}\label{rig-th: kornpoin**}
Let $\theta> 0$. Then  there is a constant $C=C(\theta)>0$  such that for all $\eps>0$ and all  $u \in SBD^2(Q_\mu,\R^2)$ the following holds: We get  paraxial rectangles $R_1,\ldots,R_n$ with 
\begin{align}\label{eq:erst}
\sum\nolimits^n_{j=1} {\rm diam}(R_j) \le (1 + \theta)\big( {\cal H}^1(J_u) + \eps^{-1}\Vert e(u) \Vert^2_{L^2(Q_\mu)}\big)
\end{align}
 such that for $E :=\bigcup^n_{j=1} R_j$ and the square $\tilde{Q} = (-\tilde{\mu},\tilde{\mu})^2$ with $\tilde{\mu} = \max\lbrace \mu - 3\sum_j{\rm diam}(R_j), 0\rbrace$ we have $|E| \le  \sum_j({\rm diam}(R_j))^2$ and
\begin{align}\label{eq:erst2}
\Vert u - (A\,\cdot + b) \Vert^2_{L^2(\tilde{Q} \setminus E)}\le C \mu^{2} \big(\Vert e(u) \Vert^2_{L^2( Q_\mu)} +  \eps {\cal H}^1(J_u)\big)
\end{align}
for some $A \in \R^{2 \times 2}_{\rm skew}$ and $b \in \R^2$.
\end{theorem}

The specific choice $\eps = \Vert e(u) \Vert^2_{L^2(Q_\mu)} ({\cal H}^1(J_u))^{-1}$ shows that the distance from an infinitesimal rigid motion may indeed be controlled exclusively by the elastic part of the energy. As a direct consequence we obtain a modification $\tilde{u} \in SBD^2(Q_\mu,\R^2)$ with $\tilde{u} = u$ on $Q_\mu \setminus E$ and $\tilde{u}(x) = A\,x +b$ for $x \in E$ such that the estimate in Theorem \ref{rig-th: kornpoin**} still holds.   To the best of our knowledge our Korn-Poincar\'e inequality differs from other inequalities of this type available in the SBV-setting (see \cite{CarrieroLeaciTomarelli:2011, DeGiorgiCarrieroLeaci:1989}) as it is not based on a truncation of the function (which is forbidden in the SBD framework), but on a modification of the displacement field on an exceptional set $E=\bigcup_j R_j$. This set is associated to the parts of $Q_\mu$ being detached from the bulk part of $Q_\mu$ by $J_u$. In contrast to the recently established estimate in \cite{Chambolle-Conti-Francfort:2014},  Theorem \ref{rig-th: kornpoin**} provides an exceptional set with a rather simple geometry. In particular,  we have control over ${\cal H}^1(\partial E)$ which allows to apply compactness results for GSBD functions (see \cite{DalMaso:13}).

In revising the present contribution let us also mention that after submission of the original version of this article, similar inequalities \cite{Conti-Iurlano:15.2, Friedrich:15-3}  in a planar setting occurred also providing a bound on the perimeter of an exceptional set. In contrast to these results our modification technique  allows to obtain a very fine estimate  on the sum of the diameter of the rectangles if in \eqref{eq:erst} we choose  $\eps = \theta^{-1} \Vert e(u) \Vert^2_{L^2(Q_\mu)} ({\cal H}^1(J_u))^{-1}$ for given small $\theta$. In fact, also in related questions concerning the approximation  of SBD functions (see e.g. \cite{Chambolle:2004, Iurlano:13}) it turned out that the derivation of sharp estimates for the surface energy is a particularly challenging problem. The bound in \eqref{eq:erst}  will be  fundamental in the rigidity estimates in \cite{Friedrich-Schmidt:15}, where based on Theorem \ref{rig-th: kornpoin**} we construct a modification whose crack length essentially coincides with the surface energy of the original displacement field.

Apart from \eqref{eq:erst2} we also provide a trace estimate for a modification $\bar{u}$ of $u$ on the boundary $\partial E$ (see Section \ref{rig-sec: subsub,  modalg2} below), i.e. we do not only control the length of $\partial E$, but also the jump height of $\bar{u}$ along $\partial E$. This byproduct of our main result is (1) again essential in \cite{Friedrich-Schmidt:15} to estimate the fracture energy and (2) a basic ingredient for the derivation of a Korn-type inequality in SBD \cite{Friedrich:15-3}.

We comment that the original motivation for establishing Theorem \ref{rig-th: kornpoin**} was the investigation of quantitative geometric rigidity results in  SBD  and the derivation of linearized Griffith energies from nonlinear models (see \cite{ Friedrich:15-2, Friedrich-Schmidt:15}) generalizing the results in nonlinear elasticity theory (cf. \cite{DalMasoNegriPercivale:02, FrieseckeJamesMueller:02}) to the framework of brittle materials. Nevertheless, we believe that this inequality is of independent interest and may contribute to solve related problems in the future, especially concerning fracture models in the realm of linearized elasticity which are related to problems in SBV where Poincar\'e inequalities   (see \cite{DeGiorgiCarrieroLeaci:1989}) have proved to be useful. Moreover, although similarly as other results in this context (see \cite{CarrieroLeaciTomarelli:2011, Chambolle-Conti-Francfort:2014, DeGiorgiCarrieroLeaci:1989}) Theorem \ref{rig-th: kornpoin**} provides an estimate only for functions with small jump sets,  the modification technique presented in this article is adapted for problems involving functions with  possibly large jump sets (see \cite{Friedrich:15-4,  Friedrich-Schmidt:15,  Schmidt:15}) where arguments based on slicing are not expedient.

The derivation of our main result is very involved as apart from the fact that the body may be disconnected by the jump set one has to face the problems that e.g. (1) the body might be still connected but only in a small region where the elastic energy is possibly large, (2) the crack geometry might become extremely complex including highly oscillating crack paths, (3) there might be infinite crack patterns occurring on different scales.  The common difficulty  of all these phenomena is the possible high irregularity of the jump set whence there are no uniform bounds for the constant in \eqref{rig-eq: help-korn}. Therefore, our proof is based on a modification algorithm which analyzes the problem iteratively on regions of mesoscopic size gradually becoming larger. In each step we have to carefully balance the elastic and surface contributions as well as the crack geometry and decide whether to establish an estimate for the function using \eqref{rig-eq: help-korn} or to alter the jump set.

Although we in principle believe that the main strategy and a lot of techniques can be employed to treat the problem in arbitrary space dimension, we only prove the result in a planar setting in order to avoid even more technical difficulties concerning the topological structure of the crack geometries occurring in $d \ge 3$.

As the proof of Theorem \ref{rig-th: kornpoin**} is very technical and long, we give a short overview and highlight the principal proof strategies for   convenience of the reader at the end of this introduction. The rest of the paper is then organized as follows. In Section \ref{rig-sec: 1} we present the main modification results and show that Theorem \ref{rig-th: kornpoin**} may be derived herefrom. Section \ref{rig-sec: sub, neigh} is devoted to a fine investigation of   the crack geometry in the neighborhood of a jump component which is the starting point for the derivation of a trace estimate. In Section \ref{rig-sec: sub, proof kornpoin} we then present a modification algorithm  which iteratively modifies the jump set such that the estimates on the jump height may be applied. Section \ref{rig-sec: proofs} and Section \ref{rig-sec: sub, trace bc} contain the proofs of the technical results. Finally, in Appendix \ref{rig-sec: appe} we provide some auxiliary estimates including an analysis for the constant  of the Korn-Poincar\'e inequality in BD (see \eqref{rig-eq: help-korn}) for varying domains and a trace theorem in SBD which allows to control the $L^2$-norm of the functions on the boundary. Here we also give a list of frequently used notation.

\subsection*{Overview of the proof}

We give an outline of the principal proof ideas.

\smallskip

\textbf{(A)} First, by approximation of $u \in SBD^2(Q_\mu)$ we can assume that the jump set can be covered by  a finite number of rectangles (see Theorem \ref{th: sbd dense} and the proof of Theorem \ref{rig-th: kornpoin**XXXX}). The boundary of these sets, called \emph{boundary components} in the following, will be altered during an iterative procedure and we measure their length by a convex combination of the Hausdorff-measure ${\cal H}^1$ and the `diameter' defined by
$$
|\Gamma|_* := h_* {\cal H}^1(\Gamma) + (1-h_*) |\Gamma|_\infty \ \ \ \ \text{with} \ \ \ |\Gamma|_\infty  = \sqrt{|\pi_1 \Gamma|^2 + |\pi_2 \Gamma|^2}
$$
for $0 <h_* < 1$ small, where $\Gamma$ denotes the boundary component and $\pi_1,\pi_2$ the orthogonal projections onto the coordinate axes.

\smallskip

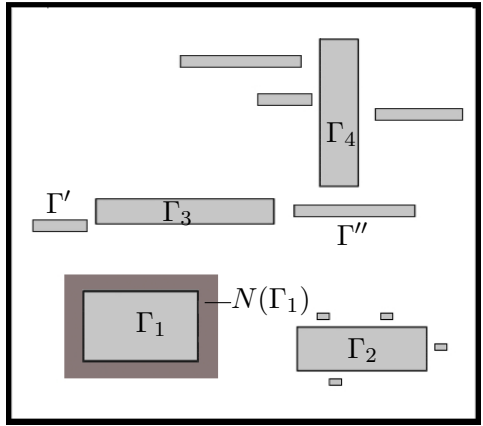
\begin{figure}[!htb]
\begin{minipage}[hbt]{7cm}

\textbf{(B)}  Beginning with small components the goal is to derive a trace estimate on each $\Gamma$ such that replacing the function $u$ by an appropriate infinitesimal rigid motion in the interior of the rectangle we get $[u] \sim \sqrt{|\Gamma|_\infty\eps}$ on   $\Gamma$. Here the scaling of $[u]$ reflects the intuition that the crack opening of small cracks is generically small. Then the estimate \eqref{eq:erst2} in  Theorem  \ref{rig-th: kornpoin**} can be established by \eqref{rig-eq: help-korn}. The fundamental ingredient in the proof  is a modification algorithm introduced in Section \ref{rig-sec: subsub,  modalg} which iteratively changes the jump set and determines the trace at boundary components once specific conditions in a neighborhood are fulfilled.

\end{minipage}
\hfill
\begin{minipage}[hbt]{7cm}
\begin{figure}[H]
\centering
\vspace{-1cm}
\begin{overpic}[width=0.9\linewidth,clip]{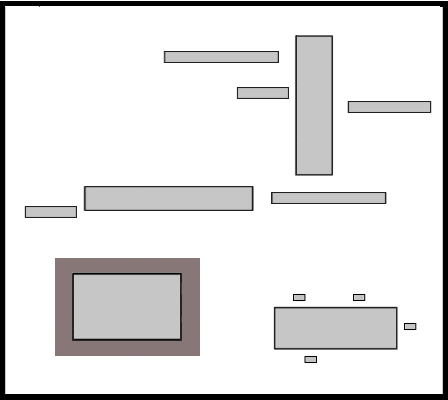}
\put(49,35){\small{$\Gamma_1$}}
\put(59,77){\small{$\Gamma_3$}}
\put(15,80){\small{$\Gamma'$}}
\put(125,68){\small{$\Gamma''$}}
\put(121,107){\small{$\Gamma_4$}}
\put(129,25){\small{$\Gamma_2$}}
\put(85,44){\small{$N(\Gamma_1)$}}
\put(75,46) {\line(1,0){10}}
\end{overpic}
\caption{\small Different boundary components whose interior is illustrated in light grey. In dark grey we depicted $N(\Gamma_1)$.} \label{1A}
\end{figure}

\end{minipage}
\end{figure}

\smallskip

\textbf{(C)} To derive a trace estimate for $\Gamma$ we consider a rectangular neighborhood $N=N(\Gamma)$ (see $\Gamma_1$ in Figure \ref{1A} and Section \ref{rig-sec: subsub,  rec-neigh}). If $N$ does not contain other components (see $\Gamma_1$) and the elastic energy $\Vert e(u) \Vert^2_{L^2(N)}$ is bounded by  $\sim |\Gamma|_\infty \eps$, the property essentially follows by a scaling argument in  \eqref{rig-eq: help-korn} and a corresponding trace theorem. Likewise, an estimate can be derived if $N$ contains only small components $(\Gamma_l)_l$, for which $[u] \sim \sqrt{|\Gamma_l|_\infty\eps}$ on $\Gamma_l$ has already been established and one has $\sum_l |\Gamma_l|_\infty \le C|\Gamma|_\infty$ (see $\Gamma_2$).     

         We will even see that the behavior of $u$ in $N$ is still sufficiently rigid if there are (at most) two large components $\Gamma'$, $\Gamma''$ intersecting $N$  and the elastic and surface energy in the two areas close to $\Gamma, \Gamma'$ and $\Gamma$, $\Gamma''$ are   small enough (see $\Gamma_3$). We refer to the beginning of Section \ref{rig-sec: sub, trace bc} for a more detailed overview of the proof of the trace estimate and to Section \ref{rig-sec: subsub, neigh,trace} where the statement is formulated and the required bounds on the energy in $N$ are listed.

\smallskip

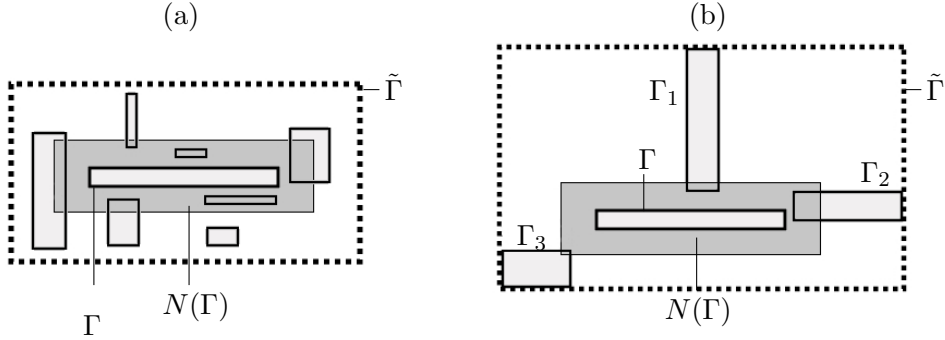
\begin{figure}[H]
\centering
\begin{overpic}[width=0.9\linewidth,clip]{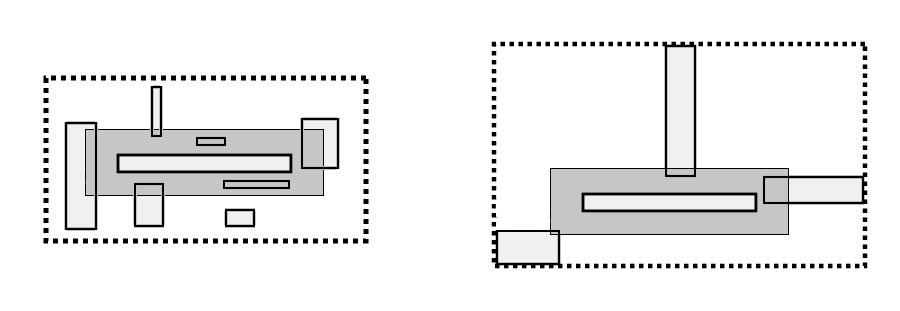}
\put(46,1){\small{$\Gamma$}}
\put(76,9){\small{$N(\Gamma)$}}
\put(76,119){\small{(a)}}
\put(275,119){\small{(b)}}
\put(160,90){\small{$\tilde{\Gamma}$}}
\put(152,94) {\line(1,0){6}}
\put(365,90){\small{$\tilde{\Gamma}$}}
\put(357,94) {\line(1,0){6}}
\put(86,20) {\line(0,1){30}}
\put(50,20) {\line(0,1){37}}

\put(278,18) {\line(0,1){20}}
\put(266,7){\small{$N(\Gamma)$}}
\put(256,64){\small{$\Gamma$}}
\put(260,90){\small{$\Gamma_1$}}
\put(340,59){\small{$\Gamma_2$}}
\put(210,35){\small{$\Gamma_3$}}
\put(258,48) {\line(0,1){14}}
\end{overpic}

\caption{\small (a) The components $\Gamma$ and $(\Gamma_l)_l$ are replaced by a larger dotted rectangle $\tilde{\Gamma}$ which in general strictly contains $N(\Gamma)$. (b) Exploiting the properties of $|\cdot|_\infty$ we find $|\tilde{\Gamma}|_{\infty} < |\Gamma|_{\infty} + \sum_{k=1,2,3} |\Gamma_k|_{\infty}$  and then also $|\tilde{\Gamma}|_{*} < |\Gamma|_{*} + \sum_{k=1,2,3} |\Gamma_k|_{*}$ for $h_*$ small.}  \label{1B}
\end{figure}

\textbf{(D)} In general, however, there is no control over the elastic energy   and the components intersecting $N$ (see $\Gamma_4$). We will show that in this case the jump set can be  changed without substantial  increase of $\sum_j |\Gamma_j|_*$, which will be crucial for the derivation of \eqref{eq:erst}. We essentially consider three situations (see case (a)-(c) in the proof of Theorem \ref{rig-th: derive prop} in Section \ref{rig-sec: subsub,  proofimain}):

\textbf{(D1)}  If  the size of the jump components $(\Gamma_l)_l$ in $N$ is much larger than $|\Gamma|_\infty$, then it is favorable to replace the component by a larger rectangle $\tilde{\Gamma}$ since $|\tilde{\Gamma}|_* \le |\Gamma|_* + \sum_l |\Gamma_l|_*$ (see Figure \ref{1B}\,(a)) and the sets $\Gamma$, $(\Gamma_l)_l$ will not be `used' anymore in the following iteration steps. A similar construction may be applied if the elastic energy    exceeds $\eps |\Gamma|_\infty$. (Here we see why in \eqref{eq:erst} also $\Vert e(u)\Vert^2_{L^2(Q_\mu)}$ occurs.)

\textbf{(D2)} Moreover, we show by exploiting the properties of $|\cdot|_\infty$ that it is never convenient if there are three  components $\Gamma_1,\Gamma_2,\Gamma_3 \cap N \neq \emptyset$ (or more)  whose size is comparable to $|\Gamma|_\infty$. As illustrated in Figure \ref{1B}\,(b), the strict convexity of the euclidian norm $|(\pi_1 \Gamma,\pi_2 \Gamma)|= \sqrt{|\pi_1 \Gamma|^2 + |\pi_2 \Gamma|^2}$ implies that the contribution of $|\Gamma|_*$, $(|\Gamma_k|_*)_{k=1}^3$ can be `used' to generate a new rectangular component $\tilde{\Gamma}$ without increase of the surface energy measured in terms of $|\cdot|_*$. Likewise, the position of the components with respect to $\Gamma$ can be investigated (see Figure \ref{E} below).

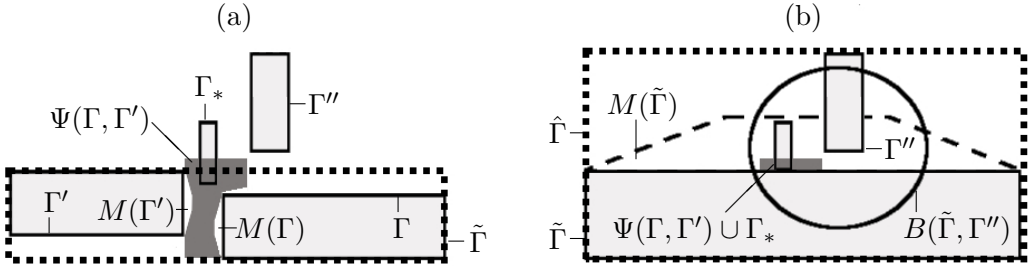
\begin{figure}[H]
\centering
\begin{overpic}[width=1.0\linewidth,clip]{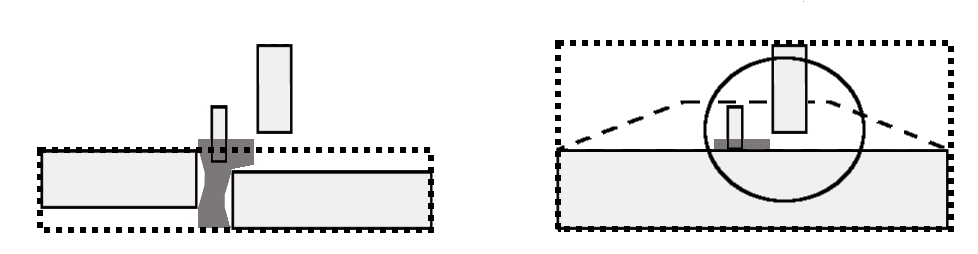}
\put(95,99){\small{(a)}}
\put(310,99){\small{(b)}}

\put(191,13){\small{$\tilde{\Gamma}$}}
\put(183,17) {\line(1,0){6}}

\put(130,65){\small{$\Gamma''$}}
\put(123,69) {\line(1,0){6}}

\put(87,72){\small{$\Gamma_*$}}
\put(91,63) {\line(0,1){6}}

\put(86,47) {\line(-1,1){15}}
\put(32,60){\small{$\Psi(\Gamma,\Gamma')$}}

\put(221,14){\small{$\tilde{\Gamma}$}}
\put(229,18) {\line(1,0){6}}

\put(221,54){\small{$\hat{\Gamma}$}}
\put(229,58) {\line(1,0){6}}

\put(245,18){\small{$\Psi(\Gamma,\Gamma') \cup \Gamma_*$}}
\put(288,28) {\line(1,1){18}}

\put(355,18){\small{$B(\tilde{\Gamma},\Gamma'')$}}
\put(360,28) {\line(0,1){8}}

\put(102,18){\small{$M(\Gamma)$}}
\put(96,22) {\line(1,0){7}}

\put(50,25){\small{$M(\Gamma')$}}
\put(80,29) {\line(1,0){4.5}}

\put(162,18){\small{$\Gamma$}}
\put(164,28) {\line(0,1){6}}

\put(30,28){\small{$\Gamma'$}}
\put(32,20) {\line(0,1){6}}

\put(348,48){\small{$\Gamma''$}}
\put(341,52) {\line(1,0){6}}

\put(243,65){\small{$M(\tilde{\Gamma})$}}
\put(254,48) {\line(0,1){15}}

\end{overpic}

\caption{\small (a) The newly generated rectangle $\tilde{\Gamma}$ including $\Gamma,\Gamma'$ does not completely contain $\Psi(\Gamma,\Gamma') \cup \Gamma_*$. Note that $\Psi(\Gamma,\Gamma')$ (the set in dark grey) does not intersect $M(\Gamma) \cup M(\Gamma')$.  (b) $\Psi(\Gamma,\Gamma') \cup \Gamma_*$  is contained in the neighborhood $M(\tilde{\Gamma})$ whose boundary has been illustrated in dashed lines. (For convenience only the relevant of the four connected component of  $M(\tilde{\Gamma})$ is depicted.) Moreover, $B(\tilde{\Gamma},\Gamma'')$ as well as $\Psi(\Gamma,\Gamma') \cup \Gamma_*$ are contained in the rectangle $\hat{\Gamma}$ which arises from the combination of $\tilde{\Gamma}$ and $\Gamma''$. Consequently, even if possibly $\Psi(\Gamma,\Gamma') \cap \Psi(\tilde{\Gamma},\Gamma'')\neq \emptyset$, the contribution of $\Psi(\Gamma,\Gamma') \cup \Gamma_*$ will not be involved in subsequent steps and is used in at most two applications of (D3).}  \label{1C}
\end{figure}

\textbf{(D3)} Finally, if the region $\Psi=\Psi(\Gamma,\Gamma')$ close to $\Gamma,\Gamma'$ contains too much jump or elastic energy, the components can be combined to a new component $\tilde{\Gamma}$ as depicted in Figure \ref{1C}\,(a). In contrast to (D1),(D2), however, we see that $\Psi$ in general is not included in the rectangle $\tilde{\Gamma}$ and therefore the energy contained therein is potentially `used' again  in a subsequent iteration step.  To keep track of this fact, we assign a weight $0<\omega(\Gamma_l)<1$ to all $\Gamma_l$ intersecting $\Psi$ indicating that the component has already been involved in a modification of this kind. In the following steps we then may not use the full contribution $|\Gamma_l|_*$, but only 
$$|\Gamma_l|_\omega:=\omega(\Gamma_l)(1-h_*)|\Gamma_l|_\infty + h_*|\Gamma_l|_{\cal H}.$$

\textbf{(E)} To establish an estimate of the form \eqref{eq:erst}, it is crucial that each component is only involved in a bounded number of applications of (D3). Although in general $\Psi(\Gamma,\Gamma')$ is  not contained in the new component $\tilde{\Gamma}$, we will see (cf. Figure \ref{1C}\,(b) and also Lemma \ref{rig-lemma: Asets}, Figure \ref{H} below) that it is always included in a neighborhood $M(\tilde{\Gamma})$  which differs from $N(\tilde{\Gamma})$ by removing the parts at the corners of the rectangle. (For a detailed illustration of the neighborhood we refer to Figure \ref{C} in Section \ref{rig-sec: subsub,  dod-neigh}.) There are in principle two different situations (see Lemma \ref{rig-lemma: Asets}):

\textbf{(E1)} If the mutual position of $\Gamma,\Gamma'$ is as illustrated in Figure \ref{1C}\,(a), then the set $\Psi(\Gamma,\Gamma')$, whose size and shape is essential for the derivation of the trace estimate, can be chosen such that $(M(\Gamma) \cup M(\Gamma')) \cap \Psi(\Gamma,\Gamma') = \emptyset$. Consequently, other sets $(\Psi_l)_l$ in $M(\Gamma) \cup M(\Gamma')$ given by previous steps have empty intersection with $\Psi(\Gamma,\Gamma')$  and therefore their contribution is not used again.

\textbf{(E2)} In the other possible situation depicted in Figure \ref{1C}\,(b) for components $\tilde{\Gamma},\Gamma''$ we observe that in general $\Psi(\tilde{\Gamma},\Gamma'') \cap (M(\tilde{\Gamma}) \cup M(\Gamma'')) \neq \emptyset$. However, in this case we find a ball $B(\tilde{\Gamma},\Gamma'') \supset \Psi(\tilde{\Gamma},\Gamma'')$ which is contained in the newly generated rectangle $\hat{\Gamma}$. This then implies that each component is involved in at most two applications of (D3).

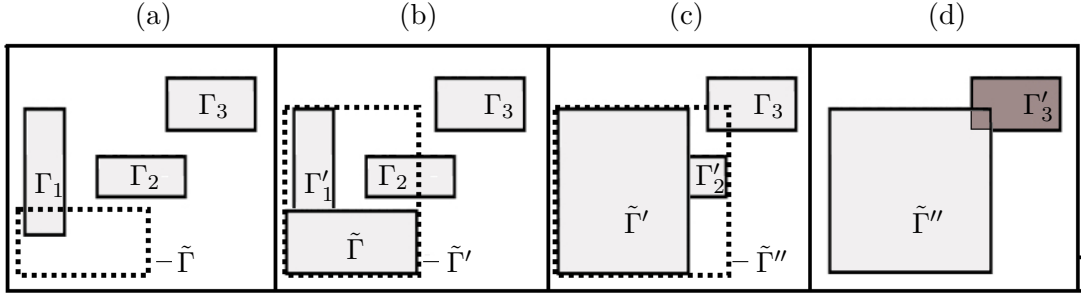
\begin{figure}[H]
\centering
\begin{overpic}[width=0.99\linewidth,clip]{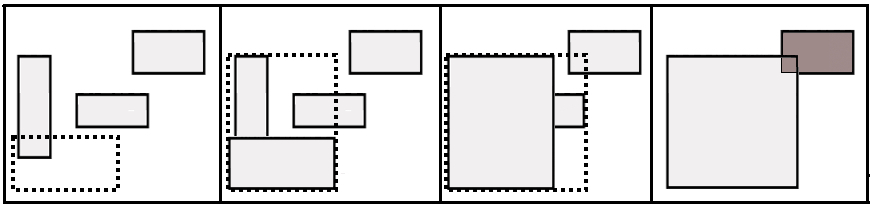}
\put(50,104){\small{(a)}}
\put(150,104){\small{(b)}}
\put(252,104){\small{(c)}}
\put(350,104){\small{(d)}}
\put(66,10){\small{$\tilde{\Gamma}$}}
\put(58,14) {\line(1,0){6}}
\put(285,10){\small{$\tilde{\Gamma}''$}}
\put(168,10){\small{$\tilde{\Gamma}'$}}
\put(159,14) {\line(1,0){6}}
\put(276,14) {\line(1,0){6}}

\put(130,15){\small{$\tilde{\Gamma}$}}
\put(235,25){\small{$\tilde{\Gamma}'$}}

\put(344,25){\small{$\tilde{\Gamma}''$}}

\put(12,40){\small{$\Gamma_1$}}
\put(114,40){\small{$\Gamma_1'$}} 

\put(48,43){\small{$\Gamma_2$}}
\put(142,43){\small{$\Gamma_2$}}
\put(262,43){\small{$\Gamma_2'$}}

\put(74,70){\small{$\Gamma_3$}}
\put(182,70){\small{$\Gamma_3$}}
\put(284,70){\small{$\Gamma_3$}}
\put(386,70){\small{$\Gamma_3'$}}

\end{overpic}

\caption{\small (a) Consider a newly generated rectangle $\tilde{\Gamma}$ and $\Gamma_1,\Gamma_2,\Gamma_3$ with $\omega(\Gamma_1)= 1$, $\omega(\Gamma_2), \omega(\Gamma_3)<1$. (b) $\tilde{\Gamma}$ and $\Gamma_1'$ are combined to $\tilde{\Gamma}'$. (c) Due to \eqref{eq: wei} we get $\omega(\Gamma_2') = 1$ and $\tilde{\Gamma}'$,  $\Gamma_2'$ are combined to $\tilde{\Gamma}''$. (d) As $\omega(\Gamma_3')<1$, we do not further combine. In dark grey the rectangle $\Gamma_3$ corresponding to $\Gamma_3'$ is depicted.}  \label{1D}
\end{figure}

\textbf{(F)} After the derivation  of a trace estimate in (C) or the modification of the jump set in (D) we pass to the next iteration step. We observe that possibly not all boundary components of the new configuration are pairwise disjoint and rectangular (cf. Figure \ref{1D}\,(a),(b)).

\textbf{(F1)} Therefore, similarly as in (D2) above we apply another modification combining different components with nonempty intersection (see Figure \ref{1D}\,(b),(c) and Lemma \ref{rig-lemma: derive prop2}). By iterative application of this argument we obtain a configuration where all components with weight $1$ are rectangular and pairwise disjoint (cf. Figure \ref{1A}). 

\textbf{(F2) }For $\Gamma_l$ with $\omega(\Gamma_l)<1$, however, a similar reasoning is not possible since we may only use the contribution $|\Gamma_l|_\omega$, where $|\Gamma_l|_\omega< |\Gamma_l|_*$. In this case we do not combine the sets, but consider the part ${\Gamma}'_l$ not covered by the newly generated rectangle as a new component (see  Figure \ref{1D}\,(d)). Since for $\Gamma_l$ a trace estimate has already been established in a previous step, it is indeed not necessary that ${\Gamma}'_l$ is rectangular. In this context, however, it is crucial to ensure that
\begin{align}\label{eq: wei}
\omega(\Gamma_l')|\Gamma_l'|_\infty   \ge c|\Gamma_l|_\infty
\end{align}
since $[u] \sim  \sqrt{|\Gamma_l|_\infty\eps}$ has been derived in terms of the original diameter $|\Gamma_l|_\infty$. To this end, we (1) always assign to ${\Gamma}'_l$ its original rectangle  $\Gamma_l$ (see \eqref{rig-eq: W prop}(i),(ii)) and (2) carefully adjust the weight $\omega({\Gamma}'_l)$ such that \eqref{eq: wei} still holds after further modifications (see Section \ref{rig-sec: sub, modification}).

\vspace{0.5cm}
 
\section{Basic notions and main modification theorem}\label{rig-sec: 1}

In this section we first introduce basic notions including the definition  of \emph{boundary components} and a corresponding  measure for the length of a component being a convex combination of the Hausdorff-measure ${\cal H}^1$ and the `diameter'. Afterwards, we present our main modification results and show how Theorem \ref{rig-th: kornpoin**} can be derived herefrom.

\subsection{Basic notions}\label{rig-sec: sub, prep}

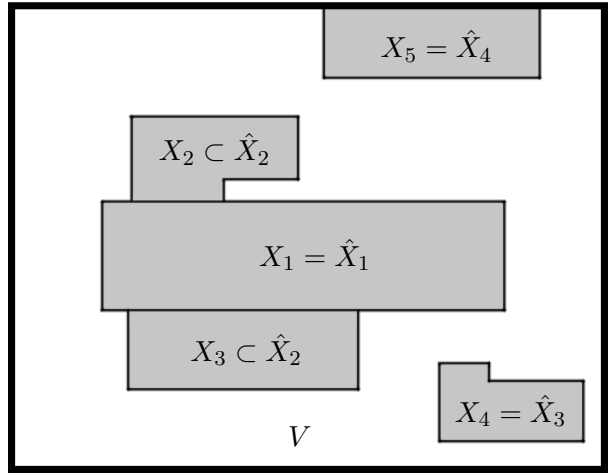
\begin{figure}[H]
\centering
\begin{overpic}[width=0.55\linewidth,clip]{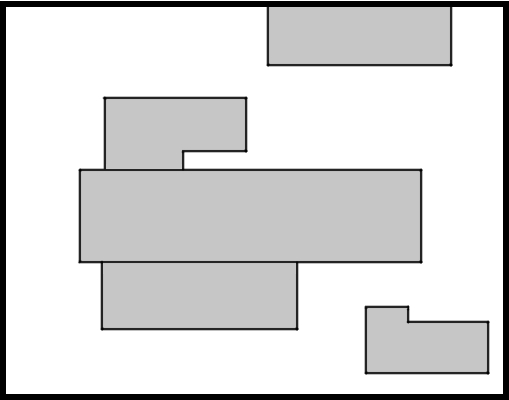}
\put(95,78){\small{$X_1= \hat{X}_1$}}
\put(141,157){\small{$X_5= \hat{X}_4$}}
\put(169,19){\small{$X_4= \hat{X}_3$}}
\put(69,42){\small{$X_3 \subset \hat{X}_2$}}
\put(57,118){\small{$X_2 \subset \hat{X}_2$}}
\put(106,10){\small{$V$}}
\end{overpic}
\caption{\small The square $Q_\mu$ with a subset $V$. The set $V$ has two representations $V_1 = Q_\mu \setminus \bigcup^5_{i=1} X_i$ and $V_2 = Q_\mu \setminus \bigcup^4_{i=1} \hat{X}_i$, where $\hat{X}_2 = X_2 \cup X_3$, which are regarded as different elements of ${\cal V}^s$. The corresponding set $U \in {\cal U}^s$ arises from $V$ by subtracting the black boundary lines $\bigcup^5_{i=1} \partial X_i$.} \label{D}
\end{figure}

\textbf{Sets and boundary components.} For $s>0$ we partition $\R^2$ up to a set of measure zero into squares $Q^s(p) = p + s(-1,1)^2$ for $p \in I^s := s(1, 1) + 2s\Z^2$ and define 
\begin{align}\label{rig-eq: calV-s def}
{\cal U}^s := \Big\{ U \subset \R^2: \exists I \subset I^s: U = \Big(\bigcup\nolimits_{p \in I} \overline{Q^s(p)} \Big)^\circ \Big\}.
\end{align}
Let $Q_\mu = (-\mu,\mu)^2$. We will concern ourselves with subsets  $V  \subset Q_\mu$ of the form
\begin{align}\label{rig-eq: calV-s def2}
{\cal V}^s := \lbrace V \subset Q_\mu: V = Q_\mu  \setminus \, \bigcup\nolimits^m_{i=1} X_i, \ \ X_i \in {\cal U}^s, \ X_i  \text{ pairwise disjoint} \rbrace
\end{align}
for $s >0$. Note that each set in $V \in {\cal V}^s$ coincides with a set $U \in {\cal U}^s$ up to subtracting a set of zero Lebesgue measure, i.e. $U \subset V$, ${\cal L}^2(V \setminus U) = 0$. The essential difference of $V$ and the corresponding $U$ concerns the connected components of the complements $Q_\mu \setminus V$ and $Q_\mu \setminus U$.  Observe that one may have $Q_\mu  \setminus \, \bigcup\nolimits^{m}_{i=1} X_i = Q_\mu  \setminus \, \bigcup\nolimits^{\hat{m}}_{i=1} \hat{X}_i$ with $(X_1,\ldots,X_{m}) \neq (\hat{X}_1,\ldots,\hat{X}_{\hat{m}})$, e.g. by combination of different sets (see Figure \ref{D}). In such a case we will regard $V_1 = Q_\mu  \setminus \, \bigcup\nolimits^{m}_{i=1} X_i$ and $V_2 = Q_\mu  \setminus \, \bigcup\nolimits^{\hat{m}}_{i=1} \hat{X}_i$ as different elements of ${\cal V}^s$. For this and the following sections we will always tacitly assume that all considered sets are elements of ${\cal V}^s$ for some small, fixed $s>0$.

\vspace{0.3cm}

Let  $W \in {\cal V}^s$ and arrange the components $X_1, \ldots, X_m$ of the complement such that $\partial X_i \subset Q_\mu$ for $1 \le i \le n$ and $\partial X_i \cap  \partial Q_\mu \neq \emptyset$ otherwise. Define $\Gamma_i(W) = \partial X_i$ for $i=1,\ldots,n$. In the following we will often refer to these sets as \emph{boundary components}  and say the component $\Gamma_i(W)$ is connected if $X_i$ is connected. Note that $\bigcup\nolimits^n_{i=1} \Gamma_i(W)$ might not cover $\partial W \cap Q_\mu$ completely if $n < m$. We frequently drop the subscript and write $\Gamma(W)$ or just $\Gamma$ if no confusion arises. We remark that for technical reasons, in particular concerning the application in \cite{Friedrich-Schmidt:15}, we here only consider the components not intersecting the boundary.  As a further  preparation, we define $H(W)\supset W \in {\cal V}^s$ as the `variant of $W$ without holes, i.e.
\begin{align}\label{rig-eq: no holes}
H(W) = W \cup \bigcup\nolimits^n_{j=1} X_j.
\end{align} 
\textbf{Surface measures.} In addition to the Hausdorff-measure $\vert \Gamma \vert_{\cal H} = {\cal H}^1(\Gamma)$ (we will use both notations) we define the `diameter' of a boundary component by $\vert \Gamma\vert_\infty:= \sqrt{|\pi_1 \Gamma|^2 + |\pi_2 \Gamma|^2}$, where $\pi_1$, $\pi_2$ denote the orthogonal projections onto the coordinate axes.

Note that by definition of ${\cal V}^s$ (in contrast to the definition of ${\cal U}^s$) two components in $(\Gamma_i)_i$ might not be disjoint. Therefore, we choose an (arbitrary) order $(\Gamma_i)^n_{i=1} = (\Gamma_i(W))^n_{i=1}$ of the boundary components of $W$, introduce 
\begin{align}\label{rig-eq: Xdef}
\Theta_i = \Theta_i(W) = \Gamma_i \setminus \bigcup\nolimits_{j<i} \Gamma_j
\end{align}
 for $i=1,\ldots,n$ and observe that the components $(\Theta_i)_i$ are pairwise disjoint. Again we will often drop the subscript if we consider a fixed boundary component. We now introduce a convex combination of $\vert \cdot \vert_\infty$ and $\vert \cdot \vert_{\cal H}$. For  $0 < h_*<1$  to be specified below we set
\begin{align}\label{rig-eq: h*}
\vert \Theta_i \vert_* =  h_* \vert \Theta_i \vert_{\cal H} + (1-h_*) \vert \Gamma_i \vert_\infty.
\end{align}
For sets $W \in {\cal V}^s$ we then define 
\begin{align*}
\Vert W \Vert_Z =  \sum\nolimits^n_{j=1} \vert \Theta_j(W) \vert_Z 
\end{align*}
for $Z={\cal H}, \infty,*$. Note that $\Vert W\Vert_\infty, \Vert W\Vert_{\cal H}$ and thus also $\Vert W\Vert_*$ are independent of the specific order which we have chosen in \eqref{rig-eq: Xdef}. Indeed, for $\Vert W\Vert_\infty$ this is clear, for $\Vert W\Vert_{\cal H}$ it follows from the fact that $\Vert W\Vert_{\cal H} = {\cal H}^1(\bigcup\nolimits^n_{i=1} \Gamma_i(W))$.  

We collect some elementary properties of $\vert \cdot \vert_*$ which follow directly from the definition.

\begin{lemma}\label{rig-lemma: inftyX}
Let $s>0$ and  $W \in {\cal V}^s$. Let $\Gamma=\Gamma(W)$ be a  boundary component with $\Gamma = \partial X$ and let $\Theta \subset \Gamma$ be the corresponding set defined in \eqref{rig-eq: Xdef}.   Moreover,  let  $V \in {\cal U}^s$ be a rectangle with $\overline{V} \cap \overline{X} \neq \emptyset$. Then
\begin{itemize} 
\item[(i)] $\vert \Theta \vert_* = \vert \Gamma \vert_* \Leftrightarrow \vert \Theta \vert_{\cal H} = \vert \Gamma \vert_{\cal H}$, \ \ \ $\vert \Gamma \vert_\infty \le \tfrac{1}{2}\vert \Gamma \vert_{\cal H}$
\item[(ii)] $|\partial (X \setminus \overline{V})|_\infty \le \vert \Gamma \vert_\infty$ and $|\Theta \setminus \overline{V}|_{\cal H} \le \vert \Theta \vert_{\cal H}$,
\item[(iii)] $ |\partial (V \cup X)|_\infty \le |\partial V|_\infty + |\Gamma|_\infty $ and $ |\partial (V \cup X)|_{\cal H} \le |\partial V|_{\cal H} + |\Gamma \setminus \partial V|_{\cal H} $
\item[(iv)] $ \frac{1}{\sqrt{2}}\vert \Gamma \vert_{\cal H} \le 2\vert \Gamma \vert_\infty \le \vert \Gamma \vert_{\cal H}$ if $\Gamma = \partial R$ for a rectangle $R \in {\cal U}^s$. 
\end{itemize}
\end{lemma}

\textbf{Weights.} As alluded to in (D3), it will be convenient for the modification scheme to introduce weights associated to the boundary components. Given a set $W \in {\cal V}^s$ we say $\omega$ is a corresponding weight if $\omega_{\rm min} \le \omega(\Gamma_i) \le 1$ for all $\Gamma_i = \Gamma_i(W)$, where $\frac{1}{2} \le \omega_{\rm min} <1$ to be specified below.  We define  
\begin{align}\label{rig-eq: h*, omega}
|\Theta_i|_{\omega} =  (1-h_*)\omega(\Gamma_i)|\Gamma_i|_\infty + h_* |\Theta_i|_{\cal H}, \ \ \ \ \Vert W \Vert_{\omega} := \sum\nolimits_{j=1}^n |\Theta_j|_{\omega}. 
\end{align}
For $\lambda \ge 0$ and fixed small $\upsilon >0$ to be specified below let ${\cal W}^s_\lambda \subset {\cal V}^s$ be the subset  consisting of the sets $W \in {\cal V}^s$ with a corresponding weight $\omega$  and an ordering of the boundary components $(\Gamma_i)^n_{i=1}$ such that the following properties are satisfied:
\begin{alignat}{2}\label{rig-eq: W prop}
(i) \ \ &  \Theta_i \subset \partial R_i, \,   \Gamma_i \subset \overline{R_i} \  \text{ for a rectangle } \overline{R_i} \subset Q_\mu \ \ \,
&&\forall \ \Gamma_i: \ \omega(\Gamma_i)< 1, \notag\\
(ii)  \ \  & \vert \partial R_i \vert_\infty \le \omega^{-1}_{\rm min} \, \omega(\Gamma_i) \vert \Gamma_i\vert_\infty && \forall \ \Gamma_i: \  \omega(\Gamma_i)< 1,\notag\\
(iii)  \ \ &    R_i \setminus X_j \text{ is connected for all } j=1,\ldots,n  \ \ \ && \forall \ \Gamma_i: \  \omega(\Gamma_i)< 1,\\
(iv)  \ \  & \omega(\Gamma_i)= 1   &&  \forall \ \Gamma_i: \ \vert\Gamma_i\vert_\infty  \ge  19\upsilon\lambda,\notag \\
 (v)  \ \ & \Gamma_i = \Theta_i = \partial R_i \text{ for a rectangle } R_i \ \ \ \ &&  \forall \ \Gamma_i: \ \omega(\Gamma_i)= 1.\notag 
\end{alignat}
 Observe that (iv),(v) imply that boundary components larger than $19\upsilon\lambda$ are always rectangular and pairwise disjoint (cf. (F1)). In particular, (i)-(iii) are trivially satisfied for ${\cal W}^s_0$  and all sets in  ${\cal W}^s_0$  have rectangular boundary components. Due to   (i) we see that the choice of a specific ordering for the  components is essential.  The conditions (i)-(ii) reflect the motivation given in (F2) and Figure \ref{1D}.  In the modification scheme presented in Section \ref{rig-sec: subsub,  modalg}  we will start with a set in ${\cal W}^s_0$ and see that by our modifications the properties in \eqref{rig-eq: W prop} remain true (see Lemma \ref {rig-lemma: derive prop2} below). In this context, $\lambda$ will stand for the size of the component being analyzed in the actual iteration step.

The diameter of $\Gamma_i$ and the corresponding rectangle $R_i$ are  comparable. More precisely, we have 
\begin{align}\label{rig-eq: new2}
(i) \ \ |\Gamma_i|_\infty \le |\partial R_i|_\infty \le 2|\Gamma_i|_\infty, \ \ \  \ \ (ii) \ \ |\Theta_i|_{\cal H} \le 4\sqrt{2}|\Gamma_i|_\infty,
\end{align}
where (i) follows from \eqref{rig-eq: W prop}(i),(ii),(v)  as well as $\omega_{\rm min} \ge \frac{1}{2}$ and (ii) is a consequence of  \eqref{rig-eq: W prop}(i), \eqref{rig-eq: new2}(i) and Lemma \ref{rig-lemma: inftyX}(iv).

\textbf{Modification.} Consider a set $W = Q_\mu \setminus \bigcup^{m}_{i=1} X_i \in {\cal W}^s_\lambda$, $\lambda \ge 0$, and a rectangle $V \in {\cal U}^s$ with $|\partial V|_\infty \ge \lambda$ and  $\overline{V} \subset Q_\mu$.  We define a \emph{modification} of $W$ by
\begin{align}\label{rig-eq: new1}
\tilde{W} = Q_\mu \setminus \bigcup\nolimits^{m}_{i=0} \tilde{X}_i,
\end{align}
where $\tilde{X}_i = X_i \setminus \overline{V}$ for $i=1, \ldots, m$ and $\tilde{X}_{0} = V$. We observe that $\tilde{W}  = (W \setminus V) \cup \partial V$ (as a subset of $\R^2$). Therefore, for shorthand we will write $\tilde{W}  = (W \setminus V) \cup \partial V$ to indicate the element of ${\cal V}^s$ which is given by \eqref{rig-eq: new1}. For the new set $\tilde{W}$  there is in general also a new corresponding weight. As its definition will first be needed in Section \ref{rig-sec: sub, proof kornpoin}, we omit it at this stage and refer to Section \ref{rig-sec: sub, modification} for details.

\textbf{Parameters.}  In this section we have already introduced the (small) parameters $h_*,\omega_{\rm min}, \upsilon$. In the following sections we will additionally consider $ q,r$. The relation between the parameters will be crucial for our analysis.  To avoid confusion, we state at this point that the parameters can be chosen in the order $h_*, q, \omega_{\rm min},  r, \upsilon$. More explicitly, we have the following dependencies for suitable $c_1= c_1(h_*) < c_2= c_2(h_*)$:
\begin{align}\label{eq: parameters}
\begin{split}
(i)& \ \ q=q(h_*),  \ \ \omega_{\rm min} = \omega_{\rm min}(h_*), \ \  r=r(h_*,q),  \ \ \upsilon = \upsilon(h_*,q,\omega_{\rm min},r),\\
(ii)& \ \ c_1 (1-\omega_{\rm min})^3 \le \upsilon \le c_2 (1-\omega_{\rm min})^3. 
\end{split}
\end{align}

\subsection{Main modification results}\label{rig-sec: subsub,  modalg2}


We  now state our main modification results from which we eventually will deduce   Theorem \ref{rig-th: kornpoin**}. For basic properties of  SBD functions we refer to Appendix \ref{rig-sec: sub, bd}. Recall that $|\cdot|_\infty$ denotes the diameter of a component and $Q_\mu = (-\mu,\mu)^2$ for $\mu>0$.

\begin{theorem}\label{rig-th: kornpoin**XXXX}
Let $\theta> 0$. Then  there is a constant $C=C(\theta)$  and a universal constant $c>0$ such that for all $\eps>0$, $\delta >0$ and all  $u \in SBD^2(Q_\mu,\R^2)$ the following holds: There are paraxial rectangles $R_1,\ldots,R_n$ with 
\begin{align}\label{eq:ers3}
\sum\nolimits^n_{j=1} |\partial R_j|_\infty \le (1 + c \theta)\big( {\cal H}^1(J_u) + \eps^{-1}\Vert e(u) \Vert^2_{L^2(Q_\mu)}\big)
\end{align}
and a modification $\bar{u} \in SBD^2(\tilde{Q},\R^2)$ with $J_{\bar{u}} \subset \bigcup^n_{j=1} \partial R_j$ and 
\begin{align}\label{eq:erst5}
\Vert \bar{u}- u\Vert_{L^2(\tilde{Q} \setminus E)} \le \delta, \ \  \ \ \Vert e(\bar{u})\Vert_{L^2(\tilde{Q})} \le  \Vert e(u)\Vert_{L^2(Q_\mu)} + \delta,
\end{align}
 where $E :=\bigcup^n_{j=1} R_j$ and $\tilde{Q} := (-\tilde{\mu},\tilde{\mu})^2$ with $\tilde{\mu} = \max\lbrace \mu - 3\sum_j |\partial R_j|_\infty, 0\rbrace$, such that for all measurable sets $D \subset \tilde{Q}$ we have
\begin{align}\label{eq:erst4}
(|E \bar{u}|(D))^2  \le c |D| \Vert e(\bar{u}) \Vert^2_{L^2( D )} +  C\eps {\cal H}^1(D \cap J_{\bar{u}} ) \sum\nolimits_{R_j \in {\cal R}(D)}\vert \partial R_j \vert^2_\infty,
\end{align}
where ${\cal R}(D):= \lbrace R_j: D \cap \overline{R_j} \neq \emptyset \rbrace$. Moreover, if $u \in L^\infty(Q_\mu)$, we can choose the modification such that $\Vert \bar{u} \Vert_{L^\infty(\tilde{Q} \setminus E)} \le \Vert u \Vert_{L^\infty(Q_\mu)}$.
\end{theorem}

The above theorem, particularly \eqref{eq:ers3} and \eqref{eq:erst4}, are not only the key ingredients for the proof of   Theorem \ref{rig-th: kornpoin**}, but also fundamental for the analysis of modifications in \cite{Friedrich-Schmidt:15}.  In this context, it is  important that we are free to choose $\eps$ and that $\theta$ may be selected arbitrarily small as hereby one can construct modifications whose crack lengths essentially coincide with ${\cal H}^1(J_u)$.  Moreover, the above estimates are  the starting point for the derivation of a Korn-type inequality  \cite{Friedrich:15-3}.

The main ingredient for the proof of Theorem \ref{rig-th: kornpoin**XXXX} will be the following technical modification result.  As a preparation we define   for $U \in {\cal W}^s_\lambda$ an interior square  $Q_{U} = (-\mu_U,\mu_U)^2$ with
\begin{align}\label{eq: Uqmu}
 \mu_U = \max\big\{\mu - 3\sum\nolimits_{j=1}^n |\partial R_j|_\infty - \sum\nolimits_{j=n+1}^m|\partial X_j(U)|_\infty,0\big\},
\end{align}
where $(R_j)_{j=1}^n$  are the rectangles corresponding to $(\Gamma_j(U))_{j=1}^n$ given by \eqref{rig-eq: W prop}(i) or \eqref{rig-eq: W prop}(v), respectively, and $(X_j(U))_{j=n+1}^m$ denote the components at the boundary (see before \eqref{rig-eq: no holes}). Recall also \eqref{rig-eq: no holes}.

\begin{theorem}\label{rig-th: derive prop}
Let $h_* >0$ sufficiently small and $\sigma>0$. Then there are constants $C_1=C_1(\sigma,h_*) \ge 1$   and $0 < C_2=C_2(\sigma,h_*) <1$    such that for all $\eps>0$ the following holds: For all $W \in {\cal V}^{s}$ with connected boundary components and $u \in H^1(W)$ there is a set $U \in {\cal W}^{C_2s}_{\lambda}$ for some $\lambda \ge 0$  with $|U \setminus W| = 0$, $H(U) \supset H(W)$ and a modification and extension $\bar{u}$ in SBD defined by
\begin{align}\label{rig-eq: extend def_new}
\bar{u}(x) = \begin{cases} A_j\,x  +b_j & x \in X_j \ \ \ \text{for all } X_j \text{ with } X_j \cap Q_U \neq \emptyset, \\
u(x) & \text{else,} \end{cases}
\end{align}
with $A_j \in \R^{2 \times 2}_{\rm skew}$, $b_j \in \R^2$  such that for all $\Gamma_j(U)  = \partial X_j$  with $X_j \cap Q_U \neq \emptyset$
\begin{align}\label{rig-eq: D1} 
\int\nolimits_{\Theta_j(U)} |[\bar{u}](x)|^2 \,d{\cal H}^1(x) \le C_1  \eps  \vert \Gamma_j(U) \vert^2_\infty.
\end{align}
Moreover, one has   $|W \setminus U| \le \Vert U \Vert^2_\infty$ and
\begin{align}\label{rig-eq: energy balance3}
\eps \Vert U \Vert_* + \Vert e(u) \Vert^2_{L^2(U)} \le (1+ \sigma)(\eps \Vert W \Vert_* + \Vert e(u) \Vert^2_{L^2(W)} ).
\end{align}
 \end{theorem}

\begin{rem}\label{rig-rem: z}
 {\normalfont
 (i) In the proof of Theorem \ref{rig-th: derive prop}    we will see that the constants $C_i=C_i(\sigma,h_*)$ have polynomial growth in $\sigma$: We find $z \in \N$ large enough such that $C_1(\sigma,h_*) \le C(h_*) \sigma^{-z}$  and $C_2(\sigma,h_*) \ge C(h_*) \sigma^{z}$.
 
 (ii) By \eqref{rig-eq: extend def_new} we particularly get that the domain of $\bar{u}$ contains  $Q_U$.
 
 (iii) More precisely,  a trace estimate \eqref{rig-eq: D1}  (with extension defined as in \eqref{rig-eq: extend def_new}) can be established for all boundary components $\Gamma_j(U)$ such that $R_j \subset H(W)$ and $\dist(\partial R_j, \partial H(W)) \ge \bar{C}C_2 |\partial R_j|_\infty$ for some $\bar{C}= \bar{C}(h_*)$ large enough with $\bar{C}C_2 \le h_*$, where  $R_j$ is the corresponding rectangle given by \eqref{rig-eq: W prop}(i) or \eqref{rig-eq: W prop}(v), respectively, and $H(W)$ as introduced in  \eqref{rig-eq: no holes}.
 
 (iv) The proof shows $|\partial R_j|_\infty \le (1+2\max\lbrace h_*,\sigma\rbrace)|\Theta_j(U)|_*$ for all $\Gamma_j(U)$.
 }
 \end{rem}

 The next sections are devoted to the proof of Theorem \ref{rig-th: derive prop} which relies on an iterative modification procedure. We show that in each step we either (1) find a boundary component satisfying certain conditions in a neighborhood (see Section \ref{rig-sec: sub, neigh}) such that a trace estimate of the form \eqref{rig-eq: D1} can be derived (see Section \ref{rig-sec: sub, trace bc}) or (2) we can modify the components in such a way that a bound of the form \eqref{rig-eq: energy balance3} also holds for the new configuration (see Section \ref{rig-sec: sub, proof kornpoin}).  
 
Before we start with the proof of Theorem \ref{rig-th: derive prop} we show that hereby we indeed obtain our main result.

\subsection{Proof of the Korn-Poincar\'e inequality}\label{rig-sec: subsub,  proofi}


In this section we prove our Korn-Poincar\'e-type inequality. We split the proof into  three steps and begin with a corollary of Theorem \ref{rig-th: derive prop}.

\begin{corollary}\label{rig-cor: kornpoin}
Let be given the situation of Theorem \ref{rig-th: derive prop}. Then there is   a universal constant $C>0$ such that for all measurable sets $D \subset Q_U$ we have
 
\begin{align*} 
(|E\bar{u}|(D))^2  \le C |D| \Vert e(u) \Vert^2_{L^2(U\cap D)} + CC_1\eps {\cal H}^1(D\cap J_{\bar{u}} ) \sum\nolimits_{R_j \in {\cal R}(D)}\vert \partial R_j \vert^2_\infty,
\end{align*}
where $C_1$ is the constant in Theorem \ref{rig-th: derive prop} and ${\cal R}(D):= \lbrace R_j: D \cap \overline{R_j} \neq \emptyset \rbrace$.
\end{corollary}

\Proof Let $D \subset Q_U$ be given. First, by H\"older's inequality we have
\begin{align}\label{rig-eq: split up}
 (|E\bar{u}|(D))^2  \le  C|D| \Vert e(\bar{u}) \Vert^2_{L^2(D)}  + C{\cal H}^1(D \cap J_{\bar{u}} ) \int_{D \cap J_{\bar{u}}} |[\bar{u}]|^2\,d{\cal H}^1.
\end{align}
By \eqref{rig-eq: extend def_new} we find $\Vert e(\bar{u}) \Vert_{L^2(D)}  = \Vert e(u) \Vert^2_{L^2(U \cap D)}$. Moreover, we have  $J_{\bar{u}} \cap D = \bigcup_j \Theta_j(U) \cap D$ by \eqref{rig-eq: extend def_new} and thus by \eqref{rig-eq: D1}   
\begin{align*}
\int_{D \cap J_{\bar{u}}} |[\bar{u}]|^2\,d{\cal H}^1 &\le  \sum\nolimits_{\Theta_j(U)\cap D \neq \emptyset}  \int_{\Theta_j(U)} |[\bar{u}]|^2\,d{\cal H}^1  \le C_1  \eps  \sum\nolimits_{\Theta_j(U)\cap D \neq \emptyset} \vert \Gamma_j(U) \vert^2_\infty \\&  \le  C_1 \eps \sum\nolimits_{R_j \in {\cal R}(D)}\vert \partial R_j \vert^2_\infty,
\end{align*}
where in the last step we used \eqref{rig-eq: new2}(i).  \eop

In \cite{Friedrich-Schmidt:15} we will also use the following properties. 
\begin{rem}\label{rem:rigi}
{\normalfont
(i) Here and in \eqref{eq:erst4} we can replace $|E\bar{u}|(D)$ by $\Vert e(\bar{u}) \Vert_{L^1(D)} + \int_{D \cap J_{\bar{u}}} |[\bar{u}]|\,d{\cal H}^1$.

(ii) By \eqref{rig-eq: extend def_new} we have ${\cal H}^1(D\cap J_{\bar{u}}) \le {\cal H}^1(D\cap \partial U)$. By $|\partial R_j|_\infty \le 2\sqrt{2}\mu$ and \eqref{rig-eq: new2} we get $\sum\nolimits_{j}\vert \partial R_j \vert^2_\infty \le 4\sqrt{2}\mu \sum_j |\Gamma_j|_\infty \le 4\sqrt{2}\mu \sum_j |\Gamma_j|_{\cal H} \le 8\sqrt{2}\mu {\cal H}^1(Q_\mu\cap \partial U)$.

}
\end{rem}

We are now in the position to prove the modification result.

\noindent {\em Proof of Theorem \ref{rig-th: kornpoin**XXXX}.} 
Let $\theta >0$ and $u \in SBD^2(Q_\mu, \R^2)$   be given. Clearly, is suffices to show the assertion for $0 < \delta \le \delta_0$ for $\delta_0$ small depending on $u$ and $\theta$. Fix $0 < \delta \le \delta_0$. Due to the fact that  $u \in L^2(Q_\mu)$ we can apply Theorem \ref{th: sbd dense} to find $\tilde{u} \in SBD^2(Q_\mu, \R^2)$ satisfying \eqref{rig-eq: V3.20} for $\delta$. As the jump set of $\tilde{u}$ consists of a finite number of closed, connected pieces of $C^1$-curves, we can cover $J_{\tilde{u}}$ up to a negligible set with finitely many, pairwise disjoint paraxial rectangles $(Q_j)^m_{j=1} \subset {\cal U}^s$ for some $s>0$ small enough such that
$$\sum\nolimits_{j=1}^m |\partial Q_j|_\infty \le {\cal H}^1(J_{\tilde{u}}) + \delta \le {\cal H}^1(J_{u}) + 2\delta \le {\cal H}^1(J_{u}) + 2\delta_0.$$
We first assume that $\overline{Q_j} \subset  Q_\mu$ for all $j=1,\ldots,m$ and indicate the adaption for the general case at the end of the proof. We  then get $W := Q_\mu \setminus \bigcup^m_{j=1} Q_j \in {\cal V}^s$ with connected boundary components and $H(W) = Q_\mu$. Choosing $h_* \le \frac{\theta}{2}$ we find for $\delta_0$ sufficiently small
\begin{align}\label{rig-eq: energy balance3XXXX}
 \Vert W \Vert_* \le (1+ch_*) \sum\nolimits^m_{j=1} |\partial Q_j|_\infty \le  (1+c\theta){\cal H}^1(J_u).
 \end{align}
We now apply Theorem \ref{rig-th: derive prop} on $\tilde{u}|_W \in H^1(W)$ to obtain a modification $\bar{u}$ and a set $U \in {\cal W}^{C_2s}_\lambda$, $\lambda \ge 0$, with $H(U) = Q_\mu$ as $H(U) \supset H(W)$. Let $(R_j)^n_{j=1}$ be the rectangles corresponding to $(\Gamma_j(U))_{j=1}^n$ as given by \eqref{rig-eq: W prop}(i),(v), respectively. Then choosing $\sigma \le   \frac{\theta}{2}$     we derive by  Remark \ref{rig-rem: z}(iv), \eqref{rig-eq: energy balance3} and \eqref{rig-eq: energy balance3XXXX} 
\begin{align*}
\sum^n_{j=1} |\partial R_j|_\infty &\le (1+\theta)\sum^n_{j=1} |\Theta_j|_* = (1+\theta)\Vert U \Vert_* \le (1 + \theta)^2\big( \Vert W\Vert_* + \tfrac{1}{\eps}\Vert e(\tilde{u}) \Vert^2_{L^2(W)}\big) \\ &\le (1 + c \theta)\big( {\cal H}^1(J_u) + \eps^{-1}\Vert e(u) \Vert^2_{L^2(Q_\mu)}\big),
\end{align*}
where in the last step we used \eqref{rig-eq: V3.20} for $\delta_0$ small enough.  Letting $E :=\bigcup^n_{j=1} R_j$ we see that \eqref{eq:erst5} follows from \eqref{rig-eq: W prop}(i), \eqref{rig-eq: extend def_new} and \eqref{rig-eq: V3.20}.  Likewise, $\Vert \bar{u} \Vert_{L^\infty(\tilde{Q} \setminus E)} \le \Vert u \Vert_{L^\infty(Q_\mu)}$ holds if $u \in L^\infty(Q_\mu)$. As $H(U) = Q_\mu$, \eqref{eq: Uqmu} implies $\tilde{Q} = Q_U$  and then \eqref{eq:erst4} is a consequence of Corollary \ref{rig-cor: kornpoin}.  Finally, $J_{\bar{u}} \subset \bigcup_j \partial R_j$  follows from \eqref{rig-eq: extend def_new} and \eqref{rig-eq: W prop}(i). 

It remains to treat the case when $H(W) \neq Q_\mu$. As $W$ is the union of squares, we can extend $\tilde{u}|_{W} \in H^1(W)$ from $W$ to $W \cup (\R^2 \setminus Q_{\mu})$, still denoted by $\tilde{u}$, and can choose a square $Q'_{\mu} \supset \overline{Q_\mu}$ such that $\Vert e(\tilde{u}) \Vert_{L^2(Q_\mu' \setminus Q_\mu)} \le \Vert \nabla \tilde{u} \Vert_{L^2(Q_\mu' \setminus Q_\mu)}\le\delta$. Then we may proceed as above with $W' := W \cup (Q_\mu' \setminus Q_{\mu})$ and $Q_\mu'$ in place of $W$ and $Q_\mu$. \eop

We can now finally give the proof of Theorem \ref{rig-th: kornpoin**}.

\noindent {\em Proof of Theorem \ref{rig-th: kornpoin**}.} 
We apply Theorem \ref{rig-th: kornpoin**XXXX} for $\frac{\theta}{c}$ in place of $\theta$ with $c$ as in Theorem \ref{rig-th: kornpoin**XXXX}. Then \eqref{eq:erst} directly follows and to see \eqref{eq:erst2}, we use \eqref{eq:erst4} for $D=\tilde{Q}$, \eqref{eq:erst5} and  Theorem \ref{rig-theo: korn} to get  $A \in \R^{2 \times 2}_{\rm skew}$, $b \in \R^2$ such that 
\begin{align*}
\Vert u - (A\, \cdot + b) \Vert^2_{L^2(\tilde{Q} \setminus E)} &\le C\delta^2 +  C\Vert \bar{u} - (A\, \cdot + b) \Vert^2_{L^2(\tilde{Q} \setminus E)} \\
& \le C\delta^2 + C \mu^2\Vert e(\bar{u}) \Vert^2_{L^2( \tilde{Q})} +  C\eps {\cal H}^1(\tilde{Q} \cap J_{\bar{u}} ) \sum\nolimits_{j=1}^n \vert \partial R_j \vert^2_\infty.
\end{align*} 
We can assume that  $\sum_j|\partial R_j|_\infty \le C\mu$ as otherwise the assertion trivially holds. This together with  ${\cal H}^1(\tilde{Q} \cap  J_{\bar{u}} )  \le \sum_{j=1}^n|\partial R_j|_{\cal H} \le C\mu$ and \eqref{eq:erst5} yields
\begin{align*}
\Vert u - (A\, \cdot + b) \Vert^2_{L^2(\tilde{Q} \setminus E)} & \le C\delta^2 + C \mu^2\Vert e(u) \Vert^2_{L^2(Q_\mu)} +  C\mu^2\eps \sum\nolimits_{j=1}^n|\partial R_j|_\infty.
\end{align*} 
We now conclude by \eqref{eq:ers3} and the fact that $\delta$ was arbitrary. \eop 
 
 \begin{rem}
 {\normalfont Using the density result \cite{Iurlano:13} instead of Theorem \ref{th: sbd dense} in the proof of Theorem \ref{rig-th: kornpoin**XXXX} we see that in Theorem \ref{rig-th: kornpoin**} and Theorem \ref{rig-th: kornpoin**XXXX} we can replace $SBD^2(Q_\mu,\R^2)$ by $GSBD^2(Q_\mu) \cap L^2(Q_\mu)$ (see \cite{DalMaso:13} for the definition of this space). }
 \end{rem}

\section{Neighborhoods of boundary components and trace estimate}\label{rig-sec: sub, neigh}

The goal of this section is to formulate conditions such that one can derive a trace estimate of the form \eqref{rig-eq: D1} for a boundary component. For the whole section we consider some $W \in {\cal W}^s_\lambda$ for $\lambda \ge 0$ and a component  $\Gamma = \Gamma(W)$ with $\omega(\Gamma)=1$  and $|\Gamma|_\infty \ge \lambda$. This implies that $\Gamma$ has rectangular shape  by \eqref{rig-eq: W prop}(v). As motivated in  (C), we will first concern ourselves with rectangular neighborhoods and investigate their properties (see Section \ref{rig-sec: subsub,  rec-neigh}). In (E) we have seen that in addition we have to consider neighborhoods where parts at the corners have been removed. These sets, which  in general have dodecagonal shape, will be discussed in Section \ref{rig-sec: subsub,  dod-neigh}.  Afterwards, in Section \ref{rig-sec: subsub, neigh,trace} we can formulate the trace estimate (see Theorem \ref{rig-theorem: D}). 

The main condition for both Theorem \ref{rig-theorem: D} and the analysis of the neighborhoods will be the following minimality condition for $\Vert \cdot \Vert_\omega$: We require
\begin{align}\label{rig-eq: property4}
\begin{split}
\Vert \tilde{W}\Vert_\omega   \ge \Vert W \Vert_\omega  \   \  \ \text{ for  all rectangles } V \in {\cal U}^s \text{  with }  \Gamma \subset \overline{V} \subset Q_\mu,
\end{split}
\end{align}
where $\tilde{W} = (W \setminus V) \cup \partial V$ as  defined in \eqref{rig-eq: new1}.  We recall that as discussed below \eqref{rig-eq: new1} for the configuration $\tilde{W}$ there are in  general also new weights which will be introduced below in Section \ref{rig-sec: sub, modification}. If \eqref{rig-eq: property4} is violated for some $\tilde{W}$, we will see that it is convenient to replace $W$ by $\tilde{W}$. (Recall the motivation in (D1),(D2) and see also case (a) in the proof of Theorem \ref{rig-th: derive prop}.)


\subsection{Rectangular neighborhood}\label{rig-sec: subsub,  rec-neigh}

This section is devoted to the definition and properties of rectangular neighborhoods of $\Gamma$.  Without restriction let $\Gamma = \partial X$ with $X = (-l_1,l_1) \times (-l_2,l_2)$ for $0 <l_2 \le l_1$ and $l_1,l_2 \in s \N$. For $t \in s\N$  with $t$ small with respect to  $l_1$ we set
\begin{align*}
& N^t(\Gamma) :=  (-t - l_1 , l_1 + t) \times  (-t - l_2 , l_2 + t) \setminus \overline{X},   \\
& N^t_{j,\pm}(\Gamma) := N^t(\Gamma) \cap \lbrace  \pm x_j \ge l_j \rbrace \ \ \text{ for } j=1,2.
\end{align*}
(In the following we will use  $\pm$ for shorthand if something holds for sets with index $+$ and $-$.) We drop $\Gamma$ in the brackets if no confusion arises.

\begin{figure}[H]
\centering
\begin{overpic}[width=0.65\linewidth,clip]{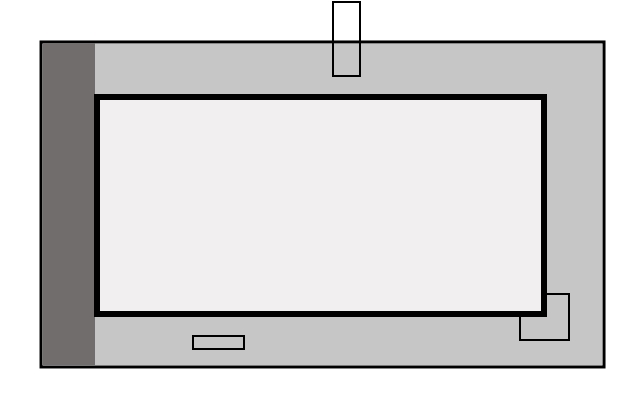}
\put(135,75){\small{$X$}}
\put(49,75){\small{$\Gamma$}}
\put(55,134){\small{$N^t(\Gamma)$}}
\put(-31,84){\small{$N^t_{1,-}(\Gamma)$}}
\put(10,88) {\line(1,0){10}}
\put(156,134){\small{$\Gamma_{j_1}$}}
\put(200,23){\small{$\Theta_{j_2}$}}
\put(106,20){\small{$\Gamma_{j_3}$}}
\put(40,78) {\line(1,0){6}}
\end{overpic}
\caption{\small Neigborhood $N^t(\Gamma)$ with other small boundary components. The part $N^t_{1,-}(\Gamma)$ is colored in dark grey.} \label{A}
\end{figure}
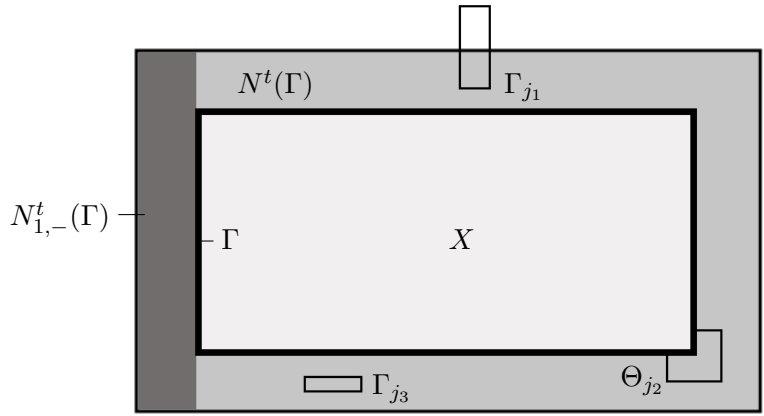

For shorthand we define $\bar{\tau} = \upsilon \vert \Gamma \vert_\infty$ with  $\upsilon > 0  $ as in  \eqref{eq: parameters}  and we will assume that (possibly by passing to a smaller $s$)
\begin{align}\label{rig-eq: tau def}
 \bar{\tau} = \upsilon \vert \Gamma \vert_\infty \in s\N \ \ \ \text{ and }  \ \ \ \bar{\tau}  \gg s.
\end{align}
This ensures that all the neighborhoods we consider below can be chosen as elements of ${\cal U}^s$.

\textbf{Properties of rectangular neighborhoods.} We now present some properties of $N^t(\Gamma)$ induced by condition \eqref{rig-eq: property4}. In particular, we will see that the neighborhood contains at most two other `large' boundary components. For the proofs of the following lemmas we refer to Section \ref{rig-sec: proofs2}. Let $(\Gamma_j)_j = (\Gamma_j(W))_j$ be the boundary components of $W$ and recall \eqref{rig-eq: no holes}.  We will always add a subscript to avoid a mix up with $\Gamma$. The following lemma shows that one can control the size of other components in the neighborhood.

\begin{lemma}\label{rig-lemma: crack bound}
Let $\lambda\ge 0$  and $0< h_*<1$ small. Let  $W \in {\cal W}^s_\lambda$ and let $\Gamma$ be a boundary component with $\omega(\Gamma)=1$ and $|\Gamma|_\infty \ge \lambda$. Assume that \eqref{rig-eq: property4} holds. Then for all $t>0$ such that $\overline{N^t} \subset H(W)$ one has  
\begin{align*}
\vert \partial W \cap N^{t} \vert_{\cal H} \le \tfrac{8t}{h_*}.
\end{align*}
\end{lemma}

The proof is essentially based on the observation that \eqref{rig-eq: property4} is violated for $V= N^{t}(\Gamma) \cup \overline{X}$ if the property does not hold. Note that $h_*$ appears in the denominator due to the definition of $|\cdot|_*$ in \eqref{rig-eq: h*}.   We now formulate the main lemma about the number and position of large components   intersecting $N^t$.

\begin{lemma}\label{rig-lemma: two components}
 For $h_*$, $\upsilon$ and $1 - \omega_{\rm min}$  as in \eqref{eq: parameters} small enough  the following holds: For all $\lambda\ge 0$, $W \in {\cal W}^s_\lambda$ and boundary components $\Gamma$ satisfying $\omega(\Gamma)=1$, $|\Gamma|_\infty \ge \lambda$ and \eqref{rig-eq: property4} we have for all  $\bar{\tau}\le t \le 22\bar{\tau}$ that there are at most two boundary components $\Gamma_1$ and $\Gamma_2$ with $\vert \Gamma_i \vert_\infty \ge 19 t$ having nonempty intersection with $N^{t}$. \\
 If $\Gamma_1$, $\Gamma_2$ exist, $\Gamma_1\cup\Gamma_2$ intersects both $N^{t}_{1,+}$ and $N^{t}_{1,-}$ or both $N^{t}_{2,+}$ and $N^{t}_{2,-}$. Additionally, if $l_2 \le \frac{l_1}{2}$ then $\Gamma_1\cup \Gamma_2$ intersects both $N^{t}_{1,+}$, $N^{ t}_{1,-}$ and $|\pi_1 \Gamma_k| \ge \frac{1}{2}\vert \pi_2 \Gamma_k \vert $ for $k=1,2$.
 \end{lemma}

\begin{figure}[H]
\centering
\begin{overpic}[width=1.0\linewidth,clip]{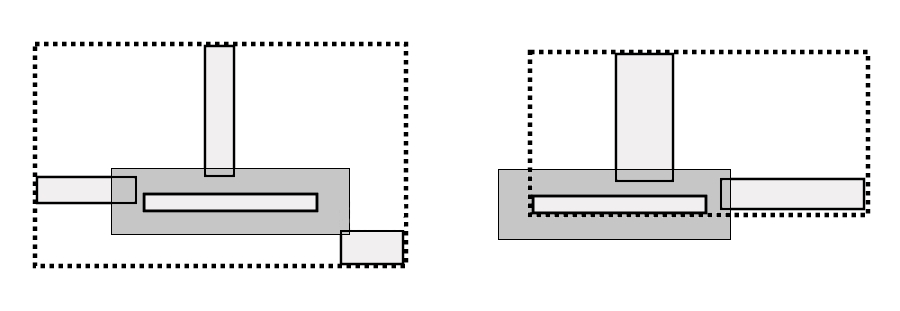}
\put(136,55){\small{$\Gamma$}}
\put(96,5){\small{$N^t(\Gamma)$}}
\put(170,40){\small{$\Gamma_1$}}
\put(30,38){\small{$\Gamma_2$}}
\put(110,90){\small{$\Gamma_3$}}
\put(194,95){\small{$V$}}
\put(108,16) {\line(0,1){25}}
\put(187,98) {\line(1,0){6}}

\put(278,21) {\line(0,1){20}}
\put(266,10){\small{$N^t(\Gamma)$}}
\put(236,98) {\line(1,0){6}}
\put(228,95){\small{$V$}}
\put(256,54){\small{$\Gamma$}}
\put(270,90){\small{$\Gamma^1$}}
\put(360,63){\small{$\Gamma^2$}}

\end{overpic}

\caption{\small The left picture shows three boundary components $\Gamma_1, \Gamma_2,  \Gamma_3$ intersecting $N^t(\Gamma)$. In the proof we argue that such a configuration violates \eqref{rig-eq: property4} for $\tilde{W} = (W \setminus V) \cup \partial V$, where $V$ is  dotted rectangle. Indeed, one might have $|\partial V|_{\cal H} > |\Gamma|_{\cal H} +\sum_{k=1,2,3} |\Gamma_k|_{\cal H}$, but we can show that one always has $|\partial V|_{\infty} < |\Gamma|_{\infty} + \sum_{k=1,2,3} |\Gamma_k|_{\infty}$ exploiting the properties of $|\cdot|_\infty$. Hereby we obtain  $|\partial V|_{*} < |\Gamma|_{*} + \sum_{k=1,2,3} |\Gamma_k|_{*}$ for $h_*$ sufficiently small. Likewise, we can control the position of the at most two large components $\Gamma^1,\Gamma^2$ in $N^t(\Gamma)$: A configuration depicted on the right, where $\Gamma^1,\Gamma^2$ do not intersect opposite parts of $N^t(\Gamma)$, violates \eqref{rig-eq: property4} for the dotted rectangle $V$. }  \label{E}
\end{figure}
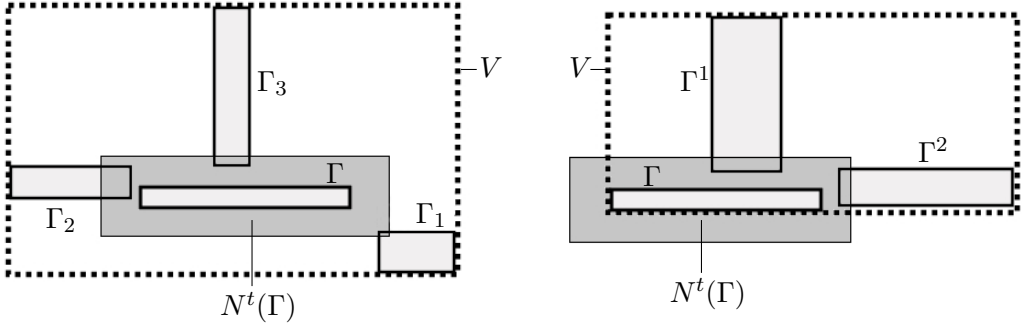

 We remark that the additional statement $|\pi_1 \Gamma_k| \ge \frac{1}{2}\vert \pi_2 \Gamma_k \vert $ also holds if only one $\Gamma_k$ exists.  The properties can be established by exploiting elementary geometric arguments and essential ideas of the procedure are exemplarily illustrated in Figure \ref{E}.

We close this first part by introducing a specific size of the neighborhood which will be needed in the following.  By Lemma \ref{rig-lemma: two components} for $t = \bar{\tau}$ we find (at most) two $\Gamma_i$, $i=1,2$, with $|\Gamma_1|_\infty, |\Gamma_2|_\infty \ge 19\bar{\tau}$ intersecting $N^{\bar{\tau}}$. We   can choose
\begin{align}\label{rig-eq: tau-bartau}
\tfrac{1}{800}\bar{\tau} \le \tau \le \tfrac{1}{2}\bar{\tau}
\end{align}
in such a way that  the neighborhood $N^\tau= N^\tau(\Gamma)$ satisfies  $\Gamma_i \cap N^{\tau/20} \neq \emptyset$ or  $\Gamma_i \cap N^{\tau} = \emptyset$ for $i=1,2$: If $\Gamma_1, \Gamma_2 \cap N^{\bar{\tau}/800} \neq \emptyset$ choose $\tau =  \frac{\bar{\tau}}{2}$, if $\Gamma_1, \Gamma_2 \cap N^{\bar{\tau}/800} = \emptyset$ choose $\tau = \frac{\bar{\tau}}{800}$, otherwise choose either $\tau=\frac{\bar{\tau}}{2}$ or $\tau = \frac{1}{40}\bar{\tau}.$ 

\textbf{Covering and  projections.} We now discuss a generalization of Lemma \ref{rig-lemma: two components}. We show that there are at most two small exceptional sets $K_1, K_2 \subset N^\tau$, where $|K_1 \cap \partial W|_{\cal H}$, $|K_2 \cap \partial W|_{\cal H}$ is `large' in a sense to be specified below. In particular, we will see that the at most two large components $\Gamma_1,\Gamma_2$ intersecting $N^\tau$ satisfy  $(\Gamma_1 \cup \Gamma_2) \cap N^\tau \subset K_1 \cup K_2$. 

\begin{figure}[H]
\centering
\begin{overpic}[width=0.7\linewidth,clip]{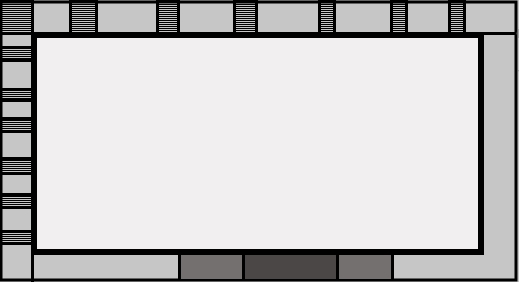}
\put(136,75){\small{$X$}}
\put(28,75){\small{$\Gamma$}}
\put(41,134){\small{$\underbrace{\ \ \ \ \ \ \ \ \ \ \ \ \ \ \ \   }$}}
\put(129,134){\scriptsize{$\underbrace{ \ \  \    }$}}
\put(251,158){\small{$\overbrace{\ \ \ \ \ \ \ \ \ \   }$}}
\put(101,18){\small{$\overbrace{\ \ \ \ \ \ \ \ \ \  \ \  \ \ \ \ \ \ \ \ \ \ \  \ \ \ \ \ \ \ \ \ }$}}
\put(135,-1){\small{$\underbrace{\ \  \ \ \ \ \ \ \  \  \ \ \ \ \ }$}}
\put(64,115){\small{$Y^2_{2,+}$}}
\put(112,115){\small{$Y^3_{2,+} \cap Y^4_{2,+}$}}
\put(261,171){\small{$Y^7_{2,+}$}}
\put(151,31){\small{$\bar{Y}^t_i$}}
\put(160,-17){\small{$Y^t_i$}}
\put(20,78) {\line(1,0){6}}
\end{overpic}
\vspace{0.7cm}
\caption{\small On the upper left side of the neigborhood $N^t(\Gamma)$ one can see elements of the partition ${\cal C}^t(\Gamma)$ (which are not necessarily of the same size). The sets where two elements overlap are striped. In the lower part an element $Y^t_i$ and the corresponding enlarged set $\bar{Y}^t_i$ are highlighted.} \label{B}
\end{figure}
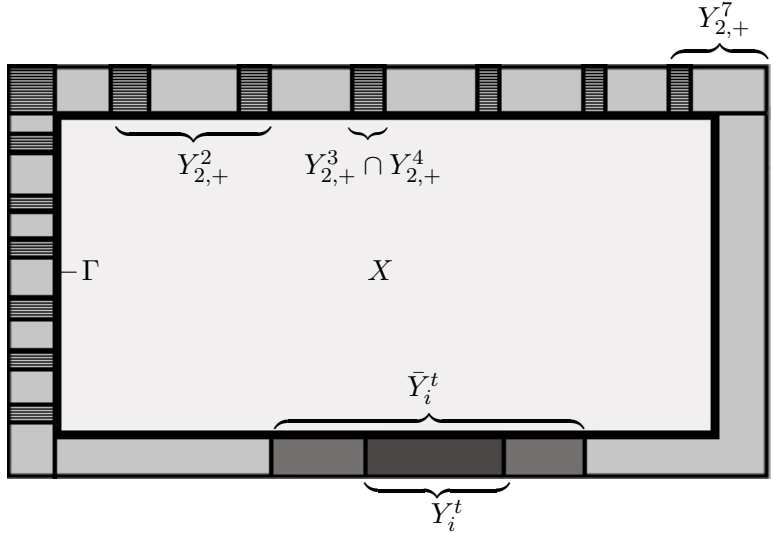

More precisely, we will construct a covering of the neighborhood (see Figure \ref{B}) and show that on  the elements not intersecting $K_1$, $K_2$ the   projection $\Vert \cdot \Vert_\pi$ (see \eqref{rig-eq: VertP}) can be controlled. As the construction provided below will in principle first be needed in Section \ref{rig-sec: sub, trace bc} for a slicing argument in the proof of Theorem \ref{rig-theorem: D}, the reader is invited to skip the details of the remainder of this section. For the following sections  we essentially only need the existence of the exceptional sets $K_1$, $K_2$ (see Lemma \ref{rig-lemma: slicing bound2}).

 \bigskip 
 
 We now proceed with the definition of the covering.  We cover $N^t_{2,\pm}$  up to a set of measure $0$ with disjoint translates of a `quasi square' $(0,\tilde{t}) \times (0,t)$, $\frac{\tilde{t}}{t} \approx 1$. If $l_2 \ge \frac{t}{2}$ we cover $N^t_{1,\pm} \setminus (N^t_{2,-} \cup N^t_{2,+})$ with translates of the rectangle $(0,t) \times (0,a)$ with $\frac{1}{2}t \le a \le t$. By $E^t_{\pm,\pm}$ we denote the four squares in the corners whose boundaries contain the points $(\pm l_1, \pm l_2)$, respectively. For $l_2 < \frac{t}{2}$ we cover each $N^t_{1,\pm}$ by itself, i.e. by a translate of the rectangle $(0,t) \times (0,2t + 2l_2)$. For convenience we will often refer to these sets as `squares' in the following. We number the squares by $Q^t_0, Q^t_1, \ldots, Q^t_n=Q^t_0$ such that $\overline{Q^t_j} \cap \overline{Q^t_{j+1}} \neq \emptyset$ for $j=0,\ldots,n-1$ and let $J^t = \lbrace Q^t_1, \ldots, Q^t_n \rbrace $.

We now introduce a coarser covering of $N^t$:  Let ${\cal Y}^t$ be the union of connected sets $Y$ having the form $Y = \big(\bigcup^k_{i=j} \overline{Q^t_i}\big)^\circ$ for $Q^t_i \in J^t$. Cover each set $N^t_{2,\pm}$ with seven sets $Y^j_{2,\pm}$ such that

\begin{align}\label{rig-eq: covering def}
|Y^j_{2,\pm}| \ge \bar{C}t\vert \Gamma\vert_\infty, \ \ \ \ \ \tfrac{1}{8}\bar{C} t\vert \Gamma\vert_\infty \le |Y^j_{2,\pm}\cap Y^{j+1}_{2,\pm}| \le \tfrac{1}{4}\bar{C}  t\vert \Gamma\vert_\infty
\end{align}
for a constant $\bar{C}>0$. If $l_2 \ge \frac{l_1}{ 2}$ we proceed likewise for $N^t_{1,\pm}$ passing possibly to  a smaller constant $\bar{C}$. If $l_2 < \frac{l_1}{2}$ we cover $N^t_{1,\pm}$ by itself. Denote the covering by ${\cal C}^t = {\cal C}^t(\Gamma) = \lbrace Y^t_0, \ldots Y^t_{m-1} \rbrace$ and order the sets in a way that $Y^t_i \cap Y^t_{i+1} \neq \emptyset$ for all $i=0,\ldots,m-1$, where by convention $Y^t_{i}=Y^t_{i \, {\rm mod}m}$.  In particular, \eqref{rig-eq: covering def} implies $Y_{i}^t \cap Y_{j}^t = \emptyset$ for $|i-j| \ge 2$. For later purposes, for a set $Y^t_i \in {\cal C}^t$ we also introduce the enlarged set
\begin{align}\label{eq: label large}
\bar{Y}^t_i = \overline{\bigcup\nolimits_{|l|\le 1} Y^t_{i+l}}.
\end{align}

\vspace{0.4cm}

Let $(R_j)_j$ be the rectangles corresponding to the components $(\Gamma_j)_j$ given by \eqref{rig-eq: W prop}(i),(v). For $Y \subset N^t$ we set ${\cal R}(Y) = \lbrace R_j: \overline{R_j} \cap Y \neq \emptyset \rbrace$ and define
\begin{align}\label{rig-eq: VertP}
\begin{split}
\vert \partial R_j\vert_\pi =   
 \min \lbrace |\partial R_j |_\infty, \ t -  \max\nolimits_{i=1,2}\dist(\pi_i R_j, \pi_i \Gamma)  \rbrace  
\end{split}
\end{align} 
for all $R_j \in {\cal R}(N^t)$. Observe that $|\partial R_j|_\pi \le |\partial R_j|_\infty$ and  that the definition of $|\cdot|_\pi$ depends on $t$. For a set $Y \subset N^t$ we then let $\Vert Y \Vert_\pi = \sum_{R_j  \in {\cal R}(Y)}  \vert \partial R_j \vert_\pi$. The projection $\Vert \cdot \Vert_\pi$ is one essential object we will need to  apply a slicing argument in the investigation of the jump heights in Section \ref{rig-sec: sub, trace bc}.

Recall the specific choice of $\tau$ in \eqref{rig-eq: tau-bartau}. 
We close this section with a lemma providing bounds on $\Vert \cdot \Vert_\pi$.

\begin{lemma}\label{rig-lemma: slicing bound2}
There are universal $c>0$ small  and $C>0$ large  such that the following holds for $h_*$, $\upsilon$ and  $1 - \omega_{\rm min}$ as in \eqref{eq: parameters} small enough: For all $\lambda \ge 0$, $W \in {\cal W}^s_\lambda$ and boundary components $\Gamma$ satisfying $\omega(\Gamma)=1$, $|\Gamma|_\infty \ge \lambda$, \eqref{rig-eq: property4} and $\overline{N^{\bar{\tau}}}\subset H(W)$  we obtain  sets $K_1, K_2 \in {\cal Y}^\tau$ with $|K_j| \le C\frac{\tau^2}{h_*}$, $j=1,2$, and $\dist(K_1,K_2) \ge c|\Gamma|_\infty$    such that
\begin{itemize}
\item[(i)] The covering $\lbrace \hat{Y}_1, \ldots, \hat{Y}_k \rbrace$ of $N^\tau\setminus (K_1 \cup K_2)$ consisting of the connected components of $\lbrace Y^\tau \setminus (K_1 \cup K_2): Y^\tau \in {\cal C}^\tau\rbrace$, satisfies $\Vert \hat{Y}_i \Vert_\pi \le \frac{19}{20}\tau$ for all $i=1, \ldots, k$. 
\item[(ii)] $\Gamma_i \cap N^\tau \subset K_1 \cup K_2$ for the (at most two) components $\Gamma_i$ with $\vert \Gamma_i \vert_\infty \ge 19\bar{\tau}$.
\end{itemize} 
\end{lemma}

The basic proof idea is very similar to the one in Lemma \ref{rig-lemma: two components}  (cf. also Figure \ref{E}): One shows that there are at most two sets $Y^1,Y^2 \in {\cal C}^\tau$, where $\Vert \cdot \Vert_\pi$ exceeds $\frac{19}{20}\tau$. For later purposes, we remark that the proof shows that if $K_i \cap (\Gamma_1 \cup \Gamma_2) \neq \emptyset$ for $i=1,2$, then $K_i$ is contained in one of the sets $N^\tau_{j,\pm}$, $j=1,2$.

\subsection{Dodecagonal neighborhood}\label{rig-sec: subsub,  dod-neigh}

We now introduce neighborhoods of $\Gamma$ which in general have dodecagonal shape  and differ from $N^t(\Gamma)$ near the corners of $\Gamma$. These neighborhoods will be essential in the modification algorithm below (see Section \ref{rig-sec: subsub,  modalg}) as we have to treat the modification near the corners of a boundary component with special care.  For further motivation we also refer to (D3) and (E).

\textbf{Definition of neighborhoods.} For  $t>0$ we define
\begin{align}\label{rig-eq: MMhatdef2}
\hat{M}^t(\Gamma) =  \bigcup\nolimits_{i=1,2} \lbrace x \in N^t(\Gamma): \   |x_i + l_i| \ge q h_*^{-1} t, |x_i - l_i| \ge q h_*^{-1} t  \rbrace  
\end{align}
for $q \ge 1$ sufficiently large to be specified below (cf. \eqref{eq: parameters}).    Moreover, for $\tilde{l} = l_1 + \min\lbrace t,q^{-1} h_* l_2\rbrace$ let
\begin{align}\label{rig-eq: MMhatdef}
\begin{split}
M^t(\Gamma) :=  \text{co}(\hat{M}^t(\Gamma) \cup \Gamma  \cup (\tilde{l},0) \cup (-\tilde{l},0)) \cap N^t(\Gamma),
\end{split}
\end{align}
where $\text{co}(\cdot)$ denotes the convex hull of a set.  Observe that $M^t(\Gamma) \supset \hat{M}^t(\Gamma)$ and that $M^t(\Gamma)$, $\hat{M}^t(\Gamma)$ may differ by some triangles.  Moreover, the shape of $M^t(\Gamma)$ is dodecagonal for $l_2 > q h_*^{-1} t$ and decagonal otherwise, cf. Figure \ref{F}.  For later reference we also define
\begin{align}\label{rig-eq: new7}
M^t_k(\Gamma) = M^t(\Gamma) \cap (N^t_{k,+}(\Gamma) \cup N^t_{k,-}(\Gamma)) \ \ \ \text{ for } k=1,2.
\end{align}
Recall the definition of $\bar{\tau}$ in \eqref{rig-eq: tau def} and the choice $\tau$ in \eqref{rig-eq: tau-bartau}. Let $K_1,K_2 \in {\cal Y}^\tau$ be the sets constructed in Lemma  \ref{rig-lemma: slicing bound2}. Let $\Gamma_m = \Gamma_m(W)$ be another boundary component satisfying $\Gamma_m \cap K \neq \emptyset$ for some $K \in \lbrace K_1,K_2 \rbrace$ and $\vert \Gamma_m\vert_\infty \ge \frac{q^2 \bar{\tau}}{h_*}$ with $q$ given in  \eqref{rig-eq: MMhatdef2}. For $q$ large enough we have $|\Gamma_m|_\infty \ge 19\bar{\tau}$ and thus $\omega(\Gamma_m)=1$ by \eqref{rig-eq: W prop}(iv). Moreover, \eqref{rig-eq: tau-bartau} implies that $\Gamma_m$  is one of the (at most) two rectangular boundary components given by Lemma \ref{rig-lemma: two components}. Recalling the remark below Lemma \ref{rig-lemma: slicing bound2},  $K$ is contained in one of the sets $N^\tau_{j,\pm}$, $j=1,2$. Let $X_m \in {\cal U}^s$ be the corresponding component of $Q_\mu \setminus W$.

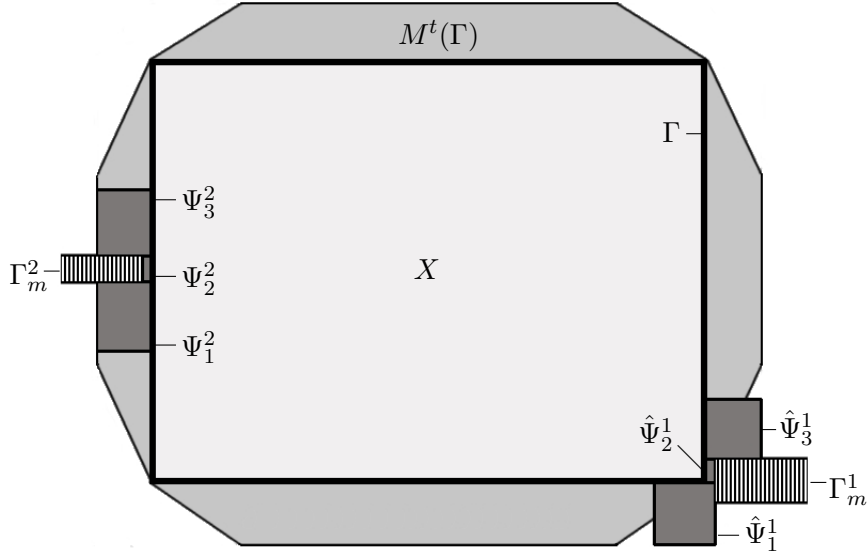
\begin{figure}[H]
\centering
\begin{overpic}[width=0.7\linewidth,clip]{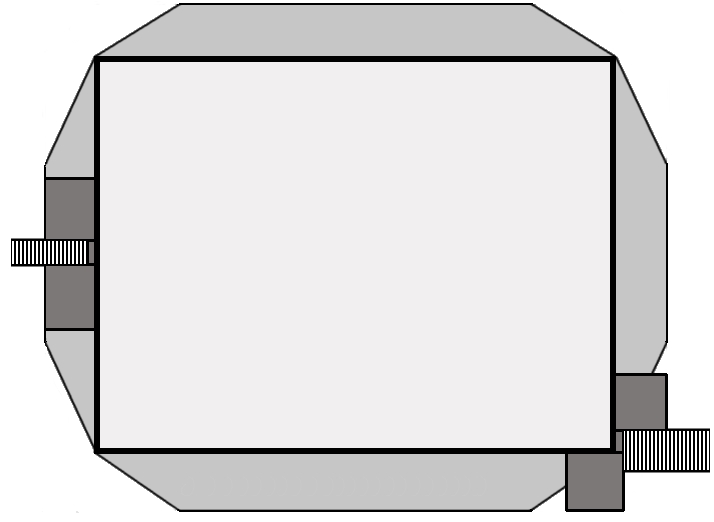}
\put(138,105){\small{$X$}}
\put(132,194){\small{$M^t(\Gamma)$}}
\put(232,157){\small{$\Gamma$}}
\put(241,160){\line(1,0){6}}

\put(-2,108){\line(1,0){6}}
\put(-15,103){\small{$\Gamma_m^2$}}

\put(288,28){\line(1,0){6}}
\put(295,22){\small{$\Gamma_m^1$}}

\put(269,48){\line(1,0){6}}
\put(276,45){\small{$\hat{\Psi}^1_3$}}

\put(238,42){\line(1,-1){10}}
\put(223,43){\small{$\hat{\Psi}^1_2$}}

\put(253,8){\line(1,0){8}}
\put(263,5){\small{$\hat{\Psi}^1_1$}}

\put(40,80){\line(1,0){8}}
\put(50,76){\small{$\Psi^2_1$}}

\put(40,106){\line(1,0){8}}
\put(50,102){\small{$\Psi^2_2$}}

\put(40,135){\line(1,0){8}}
\put(50,131){\small{$\Psi^2_3$}}

\end{overpic}
\caption{\small Neigborhood $M^t(\Gamma)$ with two other boundary components $\Gamma_m^1$, $\Gamma_m^2$ (the interiors $X^1_m$, $X^2_m$ are striped) and corresponding  sets $\hat{\Psi}^1$ and $\Psi^2$.} 
\label{C}
\end{figure}

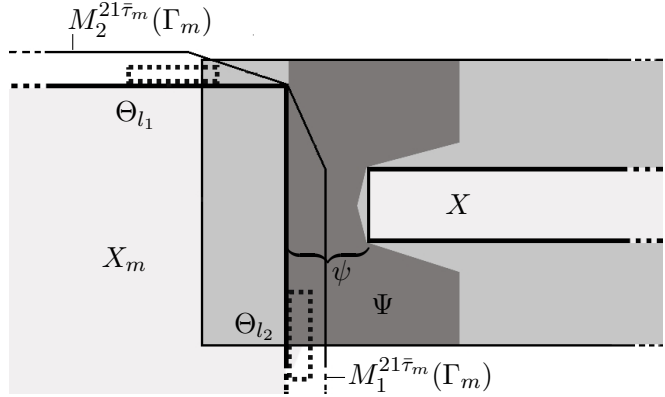
\begin{figure}[H]
\centering
\begin{overpic}[width=0.6\linewidth,clip]{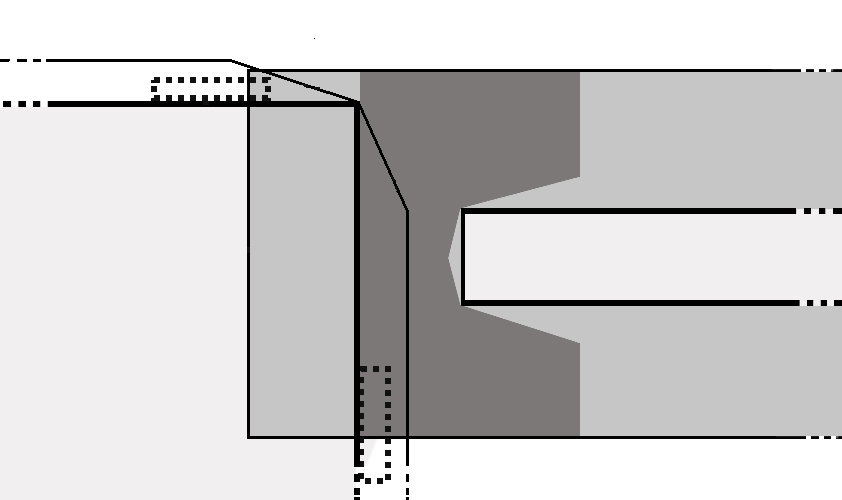}
\put(165,68){\small{$X$}}
\put(35,48){\small{$X_m$}}
\put(22,138){\small{$M^{21\bar{\tau}_m}_2(\Gamma_m)$}}
\put(128,3.5){\small{$M^{21\bar{\tau}_m}_1(\Gamma_m)$}}
\put(40,103){\small{$\Theta_{l_1}$}}
\put(85,23){\small{$\Theta_{l_2}$}}
\put(24,130){\line(0,1){6}}
\put(120,6){\line(1,0){7}}

\put(105,56){\small{$\underbrace{\ \ \  \ \, \ \ \ \ }$}}
\put(122,43){\small{$\psi$}}

\put(137,30){\small{$\Psi$}}

\end{overpic}
\caption{\small Sketch of $\Psi$ (grey) in the case  (II)(ii), where only parts of the boundary components $\Gamma$,$\Gamma_m$ are depicted. In particular $M^{21\bar{\tau}_m}(\Gamma_m) \cap \Psi \neq \emptyset$ and $M^{21\bar{\tau}}(\Gamma) \cap \Psi = \emptyset$. Also note that $M^{21\bar{\tau}_m}(\Gamma_m)$ is dodecagonal, whereas $M^{21\bar{\tau}}(\Gamma)$ is decagonal. Moreover, for later reference (see proof of Lemma \ref{rig-lemma: Asets})  we have also drawn two boundary components $\Theta_{l_1}, \Theta_{l_2}  \subset M^{21\bar{\tau}_m}(\Gamma_m)$ in dashed lines.} 
\label{G}
\end{figure}

 \textbf{Definition of sets $\Psi$ near $\Gamma_m$.} We now define sets near $\Gamma$ and $\Gamma_m$ for which we will formulate a specific condition in Section \ref{rig-sec: subsub, neigh,trace} necessary for the derivation of the trace estimate. The properties of these sets will  be exploited in Lemma \ref{lemma: rig-eq: calB prop}, Lemma \ref{rig-lemma: Asets} and Section \ref{rig-sec: sub, trace bc}  below. For the formulation of the trace estimate (Section \ref{rig-sec: subsub, neigh,trace}) and the modification algorithm (Section \ref{rig-sec: subsub,  modalg}) we will essentially only need that such sets $\Psi$ and the corresponding `distance between $\Gamma$ and $\Gamma_m$', denoted by $\psi$, exist, where there are in principle three different situations depicted in Figure \ref{C} and Figure \ref{G}. Consequently, the reader may skip the details below on first reading.

\bigskip

 For shorthand we write $M = M^\tau(\Gamma)$, $\hat{M} = \hat{M}^\tau(\Gamma)$ and $N = N^\tau(\Gamma)$, $N_{j,\pm} = N^\tau_{j,\pm}(\Gamma)$ in the following. We  treat two different cases depending on whether $K$ is near  a corner of $\Gamma$ or not:

(I) Assume $K \cap \hat{M} \neq \emptyset$. As  $K$ is contained in one of the sets $N_{j,\pm}$, $j=1,2$, we assume e.g. $K \subset N_{1,-}$. As $|K| \le C\frac{\tau^2}{h_*}$ by Lemma \ref{rig-lemma: slicing bound2}, we find $|\pi_2\Gamma_m| \le C\frac{\tau}{h_*}$ and thus recalling \eqref{rig-eq: tau-bartau} and $\vert \Gamma_m\vert_\infty \ge \frac{q^2 \bar{\tau}}{h_*}$   we get that $\vert \Gamma_m\vert_\infty \ge \frac{q^2 \bar{\tau}}{2 h_*}$ for $q$ sufficiently large. Consequently,  for $q$ sufficiently large we have $|\pi_1\Gamma_m|$ is large with respect to $\bar{\tau}$ which implies
\begin{align}\label{rig-eq: -5tau}
\Gamma_m \cap \lbrace -l_1 - 21\bar{\tau}\rbrace \times \R \neq \emptyset.
\end{align} 
Let $Q_1, Q_2 \in J^\tau$ be the neighboring squares of $K$, i.e. $Q_i \cap K = \emptyset$ and $\partial K \cap \partial Q_i \neq \emptyset$ for $i=1,2$. Let $\Psi = (\overline{Q_1 \cup K \cup Q_2} \setminus X_{ m})^\circ$  and observe that $\Psi \subset N_{1,-}$ as  $K \cap \hat{M} \neq \emptyset$. By \eqref{rig-eq: -5tau} the set $\Psi = \Psi_1 \cup \Psi_2 \cup \Psi_3$ decomposes into three rectangles, where (up to translation and sets of measure zero) $\Psi_1 = (0,\tau) \times (0,\tau +a_1)$, $\Psi_2 = (0,\psi) \times (0,\hat{\psi})$  and  $\Psi_3 = (0,\tau) \times (0,\tau+a_3)$ for $ -\frac{1}{2}\tau \le a_1,a_3 \le \tau$.  Furthermore, let 
$$\Phi = \lbrace x \in  Q_\mu: \dist(x,\Psi) \le 20\bar{\tau}\rbrace.$$
Before we go on with case (II) we state two observations. We say that two sets are $C$-Lipschitz equivalent if they are related through a bi-Lipschitzian homeomorphism with Lipschitz constants of both the homeomorphism itself and its inverse bounded by $C$.

\begin{lemma}\label{rig-lemma: two observations1}
Let $\Gamma, \Gamma_m$ with $\omega(\Gamma_m) =\omega(\Gamma) = 1$, $\overline{N^{\bar{\tau}}}\subset H(W)$ and $|\Gamma|_\infty \ge \lambda$, $\vert \Gamma_m \vert_\infty \ge  q^2 \frac{ \bar{\tau}}{h_*}$ be given. In the situation of (I) the following holds:
\begin{itemize}
\item[(i)] Let $V \in {\cal U}^s$ be the smallest rectangle containing $X$ and $X_m$. Then $\Phi \subset V$.
\item[(ii)] $\hat{\psi} \le  \frac{16\psi}{h_*} $. In particular, there is a suitable set $\Psi_2 \subset \Psi_2^* \subset \Psi$ such that each set $\Psi_1, \Psi^*_2, \Psi_3$ is  $C(h_*)$-Lipschitz equivalent to a square.
\end{itemize}
\end{lemma}

\Proof (i) As $K \cap \hat{M} \neq \emptyset$, we get that $\Phi \subset N^*:= N^{21\bar{\tau}}_{1,-} \setminus (N^{21\bar{\tau}}_{2,+} \cup N^{21\bar{\tau}}_{2,-})$ if we again choose $q$ large enough. By \eqref{rig-eq: -5tau} we have $\partial N^* \cap \Gamma_m \neq \emptyset$.  Therefore, the smallest rectangle $V$ containing $\Gamma$ and $\Gamma_m$ satisfies $N^* \subset \overline{V}$ which gives the assertion.

(ii) By Lemma \ref{rig-lemma: crack bound} we obtain $\hat{\psi} \le |\Gamma_m \cap N^{2\psi}|_{\cal H} \le \frac{16\psi}{h_*}.$ If also $ \hat{\psi} \ge h_*\psi$, we set $\Psi_2^* = \Psi_2$, otherwise we choose some $\Psi_2^* \supset \Psi_2$ with $|\pi_2 \Psi_2^*| =  \psi$. \eop

The situation of property (i) is  described in (E2) (see Figure \ref{1C}\,(b)).

(II) Assume now $K \cap \hat{M} = \emptyset$. (i)  We first treat the case $l_2 \ge  C'\frac{\tau}{h_*}$  for $C'$ large enough and  similarly as in (I) suppose without restriction that $K \subset N_{1,-}$. Again let $Q_1, Q_2$ be the neighboring squares of $K$ and set $\hat{\Psi} = (\overline{Q_1 \cup K \cup Q_2} \setminus X_m)^\circ$. If $Q_j \subset N_{1,-}$ for $j=1,2$, the set $\hat{\Psi}$ decomposes as before in (I). 

Otherwise,  we may assume that e.g. $Q_1 \subset N_{2,-}\setminus N_{1,-}$. Observe that then $K,Q_2 \subset N_{1,-}$ as $|K| \le  C\frac{\tau^2}{h_*}$ and $l_2 \ge  C'\frac{\tau}{h_*}$. As indicated in Figure \ref{C}, the set $\hat{\Psi}$ contains  three rectangles $\hat{\Psi}_1, \hat{\Psi}_2, \hat{\Psi}_3$, where (up to translation and sets of measure zero) $\hat{\Psi}_1 = (0,\tau+\psi) \times (0,\tau)$, $\hat{\Psi}_2 = (0,\psi) \times (0,\hat{\psi})$ and  $\hat{\Psi}_3 = (0,\tau) \times (0,\tau+a_3)$ for $0 \le a_3 \le \tau$. Note that $\hat{\psi}= 0$ is possible and that an argumentation as in Lemma \ref{rig-lemma: two observations1} yields $\hat{\psi} \le \frac{16\psi}{h_*}$. Now let 

$$\Psi_j = \hat{\Psi}_j \setminus (M^{21\bar{\tau}}(\Gamma) \cup M^{21\bar{\tau}_m}(\Gamma_m)), \  j=1,2,3, \ \ \ \ \Psi = \Big(\bigcup\nolimits^3_{j=1} \overline{\Psi_j} \Big)^\circ, $$
where $\bar{\tau}_m = \upsilon|\Gamma_m|_\infty$. Furthermore, let $\Phi = \lbrace x \in  Q_\mu: \dist(x,\Psi) \le 20\bar{\tau}\rbrace$.  The fact that we subtract $M^{21\bar{\tau}}(\Gamma) \cup M^{21\bar{\tau}_m}(\Gamma_m)$ from $\hat{\Psi}_j$ has been motivated in (E1) (cf. Figure \ref{1C}\,(a)).
\vspace{.4cm}

(ii) We finally treat the case that $l_2$ is small with respect to $l_1$   (i.e. $l_2 \le C'\frac{\tau}{h_*}$) which particularly implies that $M^{21\bar{\tau}}(\Gamma)$ is decagonal. Suppose without restriction that $K \subset N_{1,-}$. If $K\cap N_{2,+} = \emptyset$ or $K \cap N_{2,-} = \emptyset$, we may proceed as before in (II)(i). Otherwise, the set $\hat{\Psi} \supset \hat{\Psi}_1 \cup \hat{\Psi}_2 \cup \hat{\Psi}_3$ contains three rectangles, where (up to translation and sets of measure zero) $\hat{\Psi}_1 = (0,\tau+\psi) \times (0,\tau)$, $\hat{\Psi}_2 = (0,\psi) \times (0,2l_2)$  and  $\hat{\Psi}_3 = (0,\tau + \psi) \times (0,\tau)$   (cf. Figure \ref{G}). The same argumentation as in Lemma \ref{rig-lemma: two observations1}(ii) yields $2l_2 \le \frac{16\psi}{h_*}$. We let $\Psi_j = \hat{\Psi}_j \setminus M^{21\bar{\tau}}(\Gamma) $ for $j=1,2,3$.  Observe that in contrast to case (II)(i) we only subtract the set $M^{21\bar{\tau}}(\Gamma)$. We again let $\Phi = \lbrace x \in  Q_\mu: \dist(x,\Psi) \le 20\bar{\tau}\rbrace$. We now have the following properties.

\begin{lemma}\label{rig-lemma: two observations2}
Let $\Gamma, \Gamma_m$ with $\omega(\Gamma_m) =\omega(\Gamma) = 1$, $\overline{N^{\bar{\tau}}}\subset H(W)$ and $|\Gamma|_\infty \ge \lambda$, $\vert \Gamma_m \vert_\infty \ge  q^2 \frac{ \bar{\tau}}{h_*}$ be given. In the situation of (II) the following holds:
\begin{itemize}
\item[(i)] Let $V \in {\cal U}^s$ be the smallest rectangle containing  $X$ and $X_m$. Then we have $\Phi\cap  \lbrace x : x_{1} \ge -l_1 - \psi \rbrace \cap M^{21\bar{\tau}_m}(\Gamma_m) \subset V$.
\item[(ii)] In the cases (II)(i),(ii) we have $\hat{\psi} \le \frac{16\psi}{h_*}$ and  $2l_2 \le \frac{16\psi}{h_*}$, respectively.  Moreover, there is a suitable set $\Psi_2 \subset \Psi_2^* \subset \Psi$ such that each set $\Psi_1, \Psi^*_2, \Psi_3$ is  $C(h_*)$-Lipschitz equivalent to a square.
\end{itemize}
\end{lemma}

\Proof (i) It suffices to note that $\lbrace x: x_1 \ge -l_1 -\psi\rbrace \cap M^{21\bar{\tau}_m}(\Gamma_m) \subset [-l_1-\psi,\infty ) \times \pi_2 \Gamma_m$ and $\pi_1 \Phi \subset (-\infty, l_1]$ (cf. Figure \ref{G}). 

(ii) The bounds on $\hat{\psi}$ and $l_2$ were already discussed above. As in the proof of Lemma \ref{rig-lemma: two observations1}(ii) we can choose $\hat{\Psi}^*_2 \supset \hat{\Psi}_2$ such that $\hat{\Psi}^*_2 $ is $C(h_*)$-Lipschitz equivalent to a square. Let $\Psi^*_2 = \hat{\Psi}^*_2 \setminus  (M^{21\bar{\tau}}(\Gamma) \cup M^{21\bar{\tau}_m}(\Gamma_m))$ or $\Psi^*_2 = \hat{\Psi}^*_2 \setminus  M^{21\bar{\tau}}(\Gamma)$, respectively, depending on the cases (II)(i) and (II)(ii). For $q$ sufficiently large in  \eqref{rig-eq: MMhatdef2} it is elementary to see that $\Psi_1, \Psi^*_2, \Psi_3$  are $C(h_*)$-Lipschitz equivalent to a square.  \eop 

We remark that property (i) is a `variant' of the situation described in (E2) which we did not discuss explicitly in the introduction. The basic idea, however, is the same in the sense that the part of $M^{21\bar{\tau}_m}(\Gamma_m)$ intersecting $\Psi$ is contained in the larger rectangle $V$ (see Figure \ref{G}).

\subsection{Conditions for boundary components and trace estimate}\label{rig-sec: subsub, neigh,trace}


After the introduction of the neighborhoods we are now in the position to formulate the trace estimate. As before we consider $W_i \in {\cal W}^s_\lambda$ for $\lambda\ge 0$  with the corresponding weight $\omega$ and a specific ordering of the boundary components $(\Gamma_j(W_i))^n_{j=1}$ (cf. \eqref{rig-eq: W prop}).  (In Section \ref{rig-sec: sub, proof kornpoin} the subscript $i$ will indicate the iteration step and we already use the notation here for later reference.) Consider  $\Gamma = \Gamma(W_i)$ with  $|\Gamma|_\infty \ge \lambda$ and recall that $\Gamma = \Theta$ is rectangular by \eqref{rig-eq: W prop}(v). 

Suppose there is another set $W \in {\cal V}^s$ with $|W_i \setminus W|=0$ and assume that  $u \in H^1(W)$ is given. ($W$ will be the starting point for our modification.) As in the following we will consider the elastic energy on various subsets, it is convenient to set for shorthand  $\alpha(U) = \Vert e(u)\Vert^2_{L^2(U)}$  for $U \subset W$.

In \eqref{rig-eq: property4} we have already introduced a first condition whereby the number and position of other components in a neighborhood can be controlled. We now present additional conditions concerning (1) the elastic energy in a neighborhood, (2) the behavior in the small sets $\Psi$ constructed in Section \ref{rig-sec: subsub,  dod-neigh}   and (3) the property that a trace estimate has already been established for smaller components (cf. (C)).  
Consider $N^{2\hat{\tau}} =  N^{2\hat{\tau}}(\Gamma)$ with
\begin{align}\label{rig-eq: hat tau def}
\hat{\tau} = q^2\bar{\tau}h_*^{-1} =  q^2\upsilon h_*^{-1} |\Gamma|_\infty \ \ \ \ \text{ with }  \ \ \ \  \hat{\tau} \le c h_* |\Gamma|_\infty
\end{align}
for some $c>0$ small with $q$ from  \eqref{rig-eq: MMhatdef2} and $\bar{\tau}$ as defined in \eqref{rig-eq: tau def}. The latter inequality holds in view of  \eqref{eq: parameters} for $\upsilon$ small enough with respect to $h_*,q$.  Recall that $\hat{\tau}$ is the least length of boundary components considered in Section \ref{rig-sec: subsub,  dod-neigh}.   For $\eps >0$ and for $D =D(h_*)$ sufficiently large to be specified later (see case (b),(c) in the proof of Theorem \ref{rig-th: derive prop}) we require 
\begin{align}\label{rig-eq: property,b}
\alpha(N^{2\hat{\tau}} \cap W_i) + \eps\vert \partial W_i\cap N^{2\hat{\tau}}\vert_{\cal H} \le D\eps \hat{\tau}. 
\end{align}
Moreover, let $\Psi^j$ and $\psi^j$, $j=1,2$, be defined as in Section \ref{rig-sec: subsub,  dod-neigh} (I),(II) (cf. Figure \ref{C} and Figure \ref{G}) corresponding to the sets  $K_j$, $j=1,2$, provided by Lemma  \ref{rig-lemma: slicing bound2}. We introduce the condition
\begin{align}\label{rig-eq: property,a}
\alpha(\Psi^j \cap W_i) + \eps \vert\partial W_i \cap \Psi^j\vert_{\cal H} \le D  (1-\omega_{\rm min})^{-1} \eps \psi^j
\end{align}
for $j=1,2$, where  $D=D(h_*)$ as in \eqref{rig-eq: property,b}. 
For   $\eta \ge 0$ we let ${\cal T}_{\eta}(W_i, W)$ be the set of $\Gamma_l(W_i)$ satisfying $\vert \Gamma_l(W_i) \vert_\infty \le \eta$ and (recall \eqref{rig-eq: no holes})
\begin{align}\label{eq: etatau2}
\overline{N^{2\hat{\tau}_l}(\partial R_l)} \subset  H(W),
\end{align}    
where  $\hat{\tau}_l = q^2 \upsilon|\Gamma_l|_\infty h_*^{-1}$ (cf. \eqref{rig-eq: hat tau def}) and $R_l$ is the corresponding rectangle given in \eqref{rig-eq: W prop}(i) or \eqref{rig-eq: W prop}(v), respectively. We assume that for all $\Gamma_l(W_i) \in {\cal T}_{\hat{\tau}}(W_i, W)$ there are $A_l \in \R^{2 \times 2}_{\rm skew}$, $b_l \in \R^2$  such that for the modification $\bar{u}$ in SBD defined by
\begin{align}\label{rig-eq: extend def}
\begin{split}
\bar{u}(x) = \begin{cases}  A_l\, x + b_l & x \in X_l \ \text{ for } \ \Gamma_l(W_i) \in {\cal T}_{\hat{\tau}}(W_i, W)\\ u(x) & \text{else},  \end{cases}
\end{split}
\end{align}
we have for all $\Gamma_l(W_i) \in {\cal T}_{\hat{\tau}}(W_i, W)$  the trace estimate 
\begin{align}\label{rig-eq: property,c}
\begin{split}
\int\nolimits_{\Theta_l(W_i) } |\bar{u}(x) - (A_l\,x + b_l)|^2 \, d{\cal H}^1(x) & = \int\nolimits_{\Theta_l(W_i) } |[\bar{u}](x)|^2 \, d{\cal H}^1(x) \\ & \le C_*\frac{\eps} {\upsilon^4} \vert \Gamma_l(W_i) \vert^2_\infty
\end{split}
\end{align}
for some $C_*  =C_*(h_*) > 0$ sufficiently large.  (The left hand side has to be understood as the trace of $\bar{u}|_{W_i }$ on $\Theta_l(W_i)$.) We now state that under the above conditions also $\Gamma = \Gamma(W_i)$ with $\vert \Gamma \vert_\infty  \ge \lambda$ satisfies an estimate similar to \eqref{rig-eq: property,c}.

\begin{theorem}\label{rig-theorem: D}
Let $h_*, \omega_{\rm min},q > 0$ as in \eqref{eq: parameters}. Then there is a constant $\hat{C} = \hat{C}(h_*) >0$ independent of $\omega_{\rm min}$ such that for $\upsilon$ sufficiently small (depending on $h_*$, $\omega_{\rm min}$ and $q$)  the following holds: For all  $\eps>0$, $\lambda>0$, for all  $W_i \in {\cal W}^s_\lambda$, $W \in {\cal V}^s$ with $|W_i \setminus W|=0$,  for all $ u \in H^1(W)$ and boundary components $\Gamma=\Gamma(W_i)$ with $|\Gamma|_\infty  \ge \lambda$ such that \eqref{rig-eq: property4}, \eqref{rig-eq: property,b}, \eqref{rig-eq: property,a}, $\overline{N^{2\hat{\tau}}(\Gamma)} \subset  H(W)$ hold and \eqref{rig-eq: property,c} is satisfied for ${\cal T}_{\hat{\tau}}(W_i, W)$ one has (in the sense of traces)
\begin{align}\label{rig-eq: D0} 
\int\nolimits_{\Gamma} |\bar{u}(x) - (A\,x + b)|^2 \,d{\cal H}^1(x) \le \Big(\hat{C} + \frac{C_*}{2} \Big)   \frac{\eps}{\upsilon^4}  \vert \Gamma\vert^2_\infty
\end{align}
for suitable $A \in \R^{2 \times 2}_{\rm skew}$, $b \in \R^2$.
\end{theorem}

As the proof of this assertion is very technical and involves several steps,  we  postpone it to Section \ref{rig-sec: sub, trace bc}. In the next section we proceed with the formulation of the modification algorithm and the proof of Theorem \ref{rig-th: derive prop}.

\section{Proof of the modification result}\label{rig-sec: sub, proof kornpoin}

This section is devoted to the proof of Theorem \ref{rig-th: derive prop} which relies on an iterative modification procedure. We we will first give a precise meaning to the modification of sets and the adjustment of the corresponding weights (Section \ref{rig-sec: sub, modification}). Afterwards, we present the modification algorithm in Section \ref{rig-sec: subsub,  modalg} and give the proof in Section \ref{rig-sec: subsub,  proofimain}.

\subsection{Modification of sets}\label{rig-sec: sub, modification}

As motivated in (D),(F) a basic idea in our construction is the replacement of certain boundary components by (larger) rectangles. As before we consider $W = Q_\mu \setminus \bigcup^{m}_{i=1} X_i\in {\cal W}^s_\lambda$ for $\lambda\ge 0$  with the corresponding weight $\omega$ and a specific ordering of the boundary components $(\Gamma_j(W))^n_{j=1}$ (cf. Section \ref{rig-sec: sub, prep}).

Recall  from \eqref{rig-eq: new1} that for a rectangle $V \in {\cal U}^s$ with $|\partial V|_\infty \ge \lambda$,  $\overline{V} \subset Q_\mu$ we introduced the  modification $\tilde{W}  = (W \setminus V) \cup \partial V$, or more precisely $\tilde{W} = Q_\mu \setminus \bigcup\nolimits^{m}_{i=0} \tilde{X}_i$, where $\tilde{X}_i = X_i \setminus \overline{V}$ for $i=1, \ldots, m$ and $\tilde{X}_{0} = V$.

\smallskip

\textbf{Components and weights.} We will now give a precise meaning to the boundary components of $\tilde{W}$, their ordering and the corresponding weight. First let $\Gamma_0(\tilde{W}) = \partial V$ (it is convenient to start with index $0$) and  for $j \ge 1$ we have by construction $\Gamma_j(\tilde{W}) = \partial (X_j \setminus \overline{V})$. Observe that some boundary components may be empty and therefore reordering the indices we let $(\Gamma_j(\tilde{W}))^{\tilde{n}}_{j=1}$ for $\tilde{n} \le n$ be the nonempty boundary components. Clearly, for each   $\Gamma_j(\tilde{W})$,  $j\ge 1$, there is exactly one corresponding $\partial X_{i_j} = \Gamma_{i_j}(W)$ such that $\Gamma_j(\tilde{W}) = \partial(X_{i_j} \setminus \overline{V})$.  (This mapping is injective.) We order the components of $\tilde{W}$ such that $1 \le j_1 < j_2$ if and only if $i_{j_1} < i_{j_2}$, i.e. we preserve the ordering of $W$.

We now define the corresponding subsets as in \eqref{rig-eq: Xdef} and obtain  $\Theta_0(\tilde{W}) =\partial V$ as well as \begin{align}\label{rig-eq: new1XXXX}
\Theta_j(\tilde{W}) = \Theta_{i_j}(W) \setminus \overline{V}
\end{align}
for $j \ge 1$. Moreover, we choose the same corresponding rectangles as given for $W$ by \eqref{rig-eq: W prop}(i), i.e. for $\Gamma_j(\tilde{W})$ with $\omega(\Gamma_j(\tilde{W}))< 1$ we define $R_j(\tilde{W}) = R_{i_j}(W)$.

From now on for  notational convenience we may assume that $i_j = j$ for all $j \ge 1$. Recall that for our analysis it is also crucial to adjust the weights (cf. (F2)). We define the following `new' weights: Set $\omega(\Gamma_0(\tilde{W})) = 1$ and for $j \ge 1$

\begin{align}\label{rig-eq: weights}
\omega(\Gamma_j(\tilde{W})) = \begin{cases} 1 & \text{if } \omega(\Gamma_{j}(W)) = 1, \\  \min\Big\{ \frac{\vert \Gamma_{j}(W) \vert_\infty}{\vert \Gamma_j(\tilde{W}) \vert_\infty}\omega(\Gamma_{j}(W))  ,1\Big\} & \text{else.} \end{cases}
\end{align}

\textbf{Basic properties.} We note that 
\begin{align}\label{rig-eq: weights2}
\omega(\Gamma_j(\tilde{W})) \ge \omega(\Gamma_j(W)) \ \ \ \text{and} \ \ \ \omega(\Gamma_j(\tilde{W}))|\Gamma_j(\tilde{W})|_\infty \le \omega(\Gamma_j(W))|\Gamma_j(W)|_\infty 
\end{align}
for all $j\ge 1$. To see this, it suffices to observe $|\Gamma_j(\tilde{W})|_\infty \le|\Gamma_j(W)|_\infty $ which follows  from $\Gamma_j(\tilde{W}) = \partial(X_j \setminus \overline{V})$ and Lemma   \ref{rig-lemma: inftyX}(ii). We now state a basic estimate that we will frequently use. Recall  \eqref{rig-eq: h*, omega}.

\begin{lemma}\label{rig-lemma: derive prop3}
Let $\lambda \ge 0$ and $W \in {\cal W}^{s}_\lambda$. Let $\tilde{W} = (W \setminus V) \cup \partial V$ for a rectangle $V \in {\cal U}^s$ with $\overline{V} \subset Q_\mu$ and let ${\cal L}_V \subset (\Gamma_j(W))_{j=1}^n$ be the components satisfying $\Gamma_j(W) \subset \overline{V}$. (i) We have
\begin{align*}
\Vert \tilde{W}\Vert_\omega   &\le |\partial V|_* + \Vert W \Vert_\omega - (1-h_*)\sum_{\Gamma_j(W) \in {\cal L}_V}\omega(\Gamma_j(W))\vert \Gamma_j(W)\vert_\infty- h_*  \vert \partial_{\cal V} W \cap \overline{V}\vert_{\cal H},
\end{align*}
where we adjusted the weights as in \eqref{rig-eq: weights} and for shorthand $\partial_{\cal V} W = \bigcup_{j=1}^n \Theta_j(W)$.  

(ii) In particular, we get
\begin{align*}
\Vert \tilde{W}\Vert_\omega   &\le |\partial V|_* + \Vert W \Vert_\omega - \sum\nolimits_{\Gamma_l(W) \in {\cal L}_V}\vert \Theta_j(W)\vert_\omega.
\end{align*}
\end{lemma}

\Proof  We first observe  
\begin{align*}
\Vert \tilde{W}\Vert_\omega & =   |\partial V|_* + \sum_{j=1}^n \vert \Theta_j(\tilde{W})\vert_\omega = |\partial V|_* + \Vert W\Vert_\omega + \sum_{j=1}^n (\vert \Theta_j(\tilde{W})\vert_\omega-\vert \Theta_j(W)\vert_\omega),
\end{align*}
where $\Theta_j(\tilde{W}) = \Theta_j(W) \setminus \overline{V}$. For $\Gamma_j(W) \in {\cal L}_V$ we have $\Theta_j(\tilde{W}) = \emptyset$ and thus
\begin{align*}
\vert \Theta_j(\tilde{W})\vert_\omega-\vert \Theta_j(W)\vert_\omega  
  = -(1-h_*)\omega(\Gamma_j(W))|\Gamma_j(W)|_\infty - h_*|\Theta_j(W) \cap \overline{V}|_{\cal H}\notag
\end{align*}
by \eqref{rig-eq: h*, omega}. By \eqref{rig-eq: weights2} we find $\vert \Theta_j(\tilde{W})\vert_\omega-\vert \Theta_j(W)\vert_\omega \le h_* (\vert \Theta_j(\tilde{W})\vert_{\cal H}-\vert \Theta_j(W)\vert_{\cal H}) \le 0$ for $\Gamma_j(W) \notin {\cal L}_V$. Now (i) follows from  $\Theta_j(\tilde{W}) = \Theta_j(W) \setminus \overline{V}$ and the fact that $(\Theta_j(W))_j$ are pairwise disjoint (see \eqref{rig-eq: Xdef}). Moreover, as $\Theta_j(\tilde{W}) = \emptyset$ for $\Gamma_j(W) \in {\cal L}_V$ and $\vert \Theta_j(\tilde{W})\vert_\omega-\vert \Theta_j(W)\vert_\omega \le 0$ for $\Gamma_j(W) \notin {\cal L}_V$,  (ii) holds. \eop

\textbf{Modifications.} We observe that $\tilde{W}$ might not be an element of ${\cal W}^s_\lambda$. One can show, however, that $\tilde{W}$ can be modified to a set in ${\cal W}^s_\lambda$.

\begin{lemma}\label{rig-lemma: derive prop2}
Let $\lambda \ge 0$ and $W \in {\cal W}^{s}_\lambda$. Let $\tilde{W} = (W \setminus V) \cup \partial V$ for a rectangle $V \in {\cal U}^s$ with $|\partial V|_\infty \ge \lambda$  and $\overline{V} \subset Q_\mu$. Then there is another rectangle $V' \in {\cal U}^s$ with  $\overline{V} \subset \overline{V'} \subset Q_\mu$ such that $U := (W \setminus V') \cup \partial V' \in {\cal W}^{s}_\lambda$,   $H(U) = H(W) \cup \overline{V'}$ and  
\begin{align}\label{rig-eq: energy balance2}
\Vert U\Vert_\omega  \le \Vert \tilde{W} \Vert_\omega,
\end{align}
where for both sets $\tilde{W}$, $U$ we adjusted the weights as in \eqref{rig-eq: weights}. 
\end{lemma}

The proof idea is as follows: If $\tilde{W} \notin {\cal W}^{s}_\lambda$, we observe that one of the conditions \eqref{rig-eq: W prop}(iii),(v) is violated and so we particularly find a boundary component $\Gamma_j(\tilde{W})$ having nonempty intersection with $V$. We then modify $\tilde{W}$ by cutting out the smallest rectangle containing $V$ and $\Gamma_j(\tilde{W})$ (see Figure \ref{1D}). This procedure possibly has to be repeated iteratively. To see   \eqref{rig-eq: energy balance2} we then  exploit the properties stated in Lemma \ref{rig-lemma: inftyX}. 

\begin{rem}\label{rem: modi}
{\normalfont
In the proof of Lemma \ref{rig-lemma: derive prop2} we see that each $\Gamma(W)$ with $\omega(\Gamma(W)) = 1$ (and so particularly each $\Gamma(W)$ with $|\Gamma(W)|_\infty \ge \lambda$) is either left unchanged, i.e.  $\Gamma(W) = \Gamma(U)$ for a component $\Gamma(U)$ of $U$ with $\Gamma(U) \cap \overline{V'} = \emptyset$ or $\Gamma(W) \subset \overline{V'}$.

}
\end{rem}

As a  consequence of the above result, we derive that sets in ${\cal V}^s$ (cf. \eqref{rig-eq: calV-s def2}) can be modified such that the boundary components have rectangular form.

\begin{corollary}\label{rig-cor: derive prop2}
Let $W \in {\cal V}^s$ with connected boundary components. Then there is a subset $U \subset W$ with $H(U) \supset H(W)$ such that  $|W \setminus U| \le \Vert U\Vert_\infty^2$ and all boundary components of $U$ are rectangular and pairwise disjoint. Moreover, we have
\begin{align*}
\Vert U\Vert_*  \le \Vert W \Vert_*.
\end{align*}
\end{corollary}

 In particular, if we introduce  a weight $\omega$ corresponding to $U$ by $\omega(\Gamma_j(U))  = 1$ for all $j$ and define an (arbitrary) ordering of the boundary components we obtain $U \in {\cal W}^{s}_0$. The proofs of these modification results will be given in Section \ref{rig-sec: proofs1}.

\subsection{Modification algorithm}\label{rig-sec: subsub,  modalg}


In this section we present the modification algorithm needed for the proof of Theorem \ref{rig-th: derive prop}.  Let be given $W \in {\cal V}^s$ with connected boundary components and  $u \in H^1(W)$. For shorthand we again write $\alpha(U) = \Vert e(u)\Vert^2_{L^2(U)}$  for $U \subset W$.

As a   preparation we introduce a further notion needed in the formulation of the modification algorithm. For  $W_i \in {\cal W}^s_\lambda$, $\lambda \ge 0$, and a component $\Gamma_l(W_i)$ we define 
\begin{align}\label{rig-eq: S shorthand}
S(\Gamma_l(W_i)) = \bigcup\nolimits_{\Gamma_k(W_i) \in {\cal S}(\Gamma_l(W_i))}\Gamma_k(W_i)  
\end{align}
with ${\cal S}(\Gamma_l(W_i)) = \lbrace \Gamma_k(W_i): \omega(\Gamma_k(W_i)) = 1 , |\Gamma_k(W_i)|_\infty > |\Gamma_l(W_i)|_\infty  \rbrace$.  Before we start, we specify the constant $C_2$ in Theorem \ref{rig-th: derive prop} to be 
\begin{align}\label{rig-eq: C2 choice}
C_2 = C'\upsilon
\end{align}
for a universal constant $C'$ small enough and consider $W \in {\cal V}^{s}$ as an element of ${\cal V}^{C_2s}$ such that \eqref{rig-eq: tau def} (with $C_2 s$ in place of $s$) is satisfied for all boundary components $\Gamma_l(W)$.  From now on we will always tacitly assume that all involved sets lie in ${\cal V}^{C_2 s}$ and write $ {\cal W}_\lambda$ instead of ${\cal W}_\lambda^{C_2 s}$.  In the proof below we will show that $C_2$ is in fact a constant only depending on $h_*$ and $\sigma$ as claimed in Theorem \ref{rig-th: derive prop}.

We set $W_0 =  \hat{W}$, where $\hat{W}$ is the modification given by Corollary \ref{rig-cor: derive prop2}.   Choosing an (arbitrary) ordering of the boundary components and setting $\omega(\Gamma_j(W_0)) = 1$ for all $j$ we obtain $W_0 \in {\cal W}_0$. Moreover, we let $\lambda_0 = 0$, $B^0_0 = \emptyset$.  Assume that $\lambda_i \ge \lambda_0 $ and  $W_i \in {\cal W}_{\lambda_i}$ with   $|W_i \setminus  W_0| = 0$ and $H(W_i) \supset H(W)$ (recall \eqref{rig-eq: no holes}) are given. Suppose there are sets     $\lbrace B^i_j: j=0,\ldots, i\rbrace$ which  will represent the sets where we already `used' the `energy lying in the set' to modify $W$  (see (D),(E)).    

Suppose that in the iteration step $i$ the following conditions are satisfied: 
\begin{align}\label{rig-eq: energy balance4}
\eps\Vert W_i \Vert_{\omega} + \alpha(W_i) \le \eps \Vert W \Vert_{*} + \alpha(W) +  h_*(1-\omega_{\rm min}) \sum\nolimits^i_{j=0} \alpha(B^i_j)
\end{align}
as well as  
\begin{align}\label{rig-eq: Asets}
\begin{split}
(i) & \ \ \text{Each $x \in Q_\mu$ lies in at most two different  $B^i_{j_1}$, $B^i_{j_2}$ and} \\ & \,  \ \ \text{each $x \in W_i$ lies in at most one $B^i_j$,  $j,j_1,j_2 \in \lbrace 0,\ldots,i\rbrace$,}\\
(ii) & \ \ \text{Either } \Theta_l(W_i) \subset B^i_j \text{ for some } 0 \le j \le i \text{ or }\\ &  \ \ \Gamma_l(W_i) \in {\cal G}_i := \big\{ \Gamma_l: \Gamma_l\cap \bigcup\nolimits^i_{j=0} B^i_j = \emptyset, \ \omega(\Gamma_l) = 1 \big\} \  \text{ for all } \Gamma_l(W_i),\\
(iii)   &\ \  \text{Each } B^i_j \text{ with } B^i_j \cap W_i\neq \emptyset   \text{ satisfies } B^i_j \cap W_i \subset M^{\eta^i_l}_k(\Gamma_l(W_i)) \setminus S(\Gamma_l(W_i)),   \\
& \ \ \  \ \ \  \ \ \  \ \ \  \ \ \ \text{ for some } \Gamma_l(W_i) \in {\cal G}_i \text{ and } k \in \lbrace 1,2\rbrace, \ j=0,\ldots,i. 
\end{split}
\end{align}
Here $\eta^i_l := 21\upsilon \min\lbrace |\Gamma_l(W_i)|_\infty, \lambda_i \rbrace = \min \lbrace 21\bar{\tau}_l, 21\upsilon \lambda_i\rbrace$, the neighborhood  $M_k$ was defined in \eqref{rig-eq: new7}  and $S(\Gamma_l(W_i))$ as in \eqref{rig-eq: S shorthand}.  Moreover, recalling \eqref{rig-eq: hat tau def} and the definition before \eqref{eq: etatau2}   we suppose  
\begin{align}\label{rig-eq: property,d}
\begin{split}
&\alpha( N^{\hat{\tau}_l}(\Gamma_l(W_i)) \cap W_i) + \eps | N^{\hat{\tau}_l}(\Gamma_l(W_i)) \cap  \big(\partial W_i \setminus S(\Gamma_l(W_i)) \big)|_{\cal H} \le D \eps\hat{\tau}_l \\
& \ \ \ \ \ \ \ \ \ \ \ \ \ \ \ \ \ \ \ \ \ \ \ \ \ \ \ \ \ \ \ \ \ \ \ \ \ \ \ \ \ \ \ \ \ \ \text{ for all } \Gamma_l(W_i) \in {\cal G}_i \cap {\cal T}_{\lambda_i}(W_i, W),
\end{split}
\end{align}
where $D$   as defined in \eqref{rig-eq: property,b}.  Furthermore,  we assume that there is an extension $\bar{u}_i$ in SBD defined as in   \eqref{rig-eq: extend def} with suitable $A_l \in \R^{2 \times 2}_{\rm skew}$ and $b_l \in \R^2$ for all $\Gamma_l(W_i) \in {\cal T}_{\lambda_i}(W_i, W)$ such that all $\Gamma_l(W_i) \in {\cal T}_{\lambda_i}(W_i, W)$ satisfy
\begin{align}\label{rig-eq: D2}
\int_{\Theta_l(W_i)} |\bar{u}_i - (A_l\, x + b_l)|^2 \,d{\cal H}^1 \le \hat{C}\sum^i_{n=0}   \Big( \frac{2}{3} \Big)^n \ {\frac{\omega(\Gamma_l(W_i))^2}{\hat{\omega}_i(\Gamma_l(W_i))^2}}\frac{\eps}{\upsilon^4}  \vert \Gamma_l(W_i) \vert^2_\infty, 
\end{align}
 where  $\hat{\omega}_i(\Gamma_l(W_i)) := 1 - \frac{1-\omega_{\rm min}}{2}\, \#\lbrace j=0,\ldots,i:\Theta_l(W_i) \subset B^i_j \rbrace $ and  $\hat{C}$ is the constant from \eqref{rig-eq: D0}. Finally, we assume 
\begin{align}\label{rig-eq: new5}
\begin{split}
(i)   &\ \,   \omega(\Gamma_l(W_i))\ge \hat{\omega}_i(\Gamma_l(W_i)),  \ \ \ \ \ \  \ \ \ \ \ \ \ \ \ \ \ \ \ \ \ \ \ \ \ \ \ \ \ \ \ \ \ \ \forall \ \Gamma_l,\\
 (ii)   &\ \,  |\partial R_l(W_i)|_\infty \le \omega(\Gamma_l(W_i)) (\hat{\omega}_i(\Gamma_l(W_i)))^{-1} \vert \Gamma_l(W_i) \vert_\infty, \, \forall \ \Gamma_l: \ \omega(\Gamma_l)< 1.
\end{split} 
\end{align}
Recall that $\partial W_i \cap Q_\mu \subset W_i$ by definition (see \eqref{rig-eq: calV-s def2}).  This particularly implies that each $\Theta_l(W_i)$ is contained in at most one set $B_j^i$ (see \eqref{rig-eq: Asets}(i),(ii)). Consequently, we have $\hat{\omega}_i(\Gamma_l(W_i)) \ge \omega_{\rm min}$ and see that for components satisfying \eqref{rig-eq: D2} also \eqref{rig-eq: property,c} holds if we replace $C_*$ by $\hat{C} \omega_{\rm min}^{-2}\sum\nolimits^i_{n=0}  (2/3)^n$.

Condition \eqref{rig-eq: energy balance4} yields a bound on the length of the boundary components ultimately leading to \eqref{rig-eq: energy balance3}. In  \eqref{rig-eq: Asets} the properties of the sets `used' for the modification are described. In particular, \eqref{rig-eq: Asets}(i) is discussed in (E2) and \eqref{rig-eq: Asets}(iii) reflects the observation in (E). Moreover,  \eqref{rig-eq: property,d} and \eqref{rig-eq: D2} state that a trace estimate already holds for smaller components and that the energy in a neighborhood can be controlled (cf. (C)). The second condition in \eqref{rig-eq: new5} is a refinement of \eqref{rig-eq: W prop}(ii).

In the proof of Theorem  \ref{rig-th: derive prop}  it will be crucial that the above stated conditions are preserved under modification. To state this property, it is convenient to assume that there is an additional set $B^{i}_{i+1} \subset Q_\mu$ and to define
\begin{align}\label{eq:omega*}
\hat{\omega}^*_i(\Gamma_l(W_i)) := 1 - \tfrac{1}{2}(1-\omega_{\rm min})\, \#\lbrace j=0,\ldots,i+1:\Theta_l(W_i) \subset B^i_j \rbrace. 
\end{align}
Clearly, if $B^i_{i+1} = \emptyset$, then $\hat{\omega}^*_i(\Gamma_l(W_i))= \hat{\omega}_i(\Gamma_l(W_i))$ for all $\Gamma_l(W_i)$.

\begin{lemma}\label{rig-lemma: modification}
Let $\eps >0$ and $\lambda \ge 0$. Let $W_i \in {\cal W}_\lambda$, $\bar{u}_i$ and  $\lbrace B^i_j: j=0,\ldots, i\rbrace$ be given such that  \eqref{rig-eq: energy balance4}-\eqref{rig-eq: property,d} hold for $W_i$ and  $\lambda$ in place of $\lambda_i$. Moreover, suppose that \eqref{rig-eq: D2}-\eqref{rig-eq: new5} are satisfied for $\bar{u}_i$ and $\hat{\omega}^*_i$ in place of $\hat{\omega}_i$ for an additional set $B^{i}_{i+1} \subset Q_\mu$. For a rectangle $ \overline{V} \subset Q_\mu$ with $\vert \partial V \vert_\infty >  \lambda$, let $\tilde{W}_i = (W_i \setminus V) \cup \partial V$ and assume that  (recall \eqref{rig-eq: new1}, \eqref{rig-eq: weights})
$$\eps\Vert \tilde{W}_i \Vert_\omega + \alpha(\tilde{W}_i)\le \eps\Vert W_i\Vert_\omega + \alpha(W_i).$$
Let $W_{i+1}  = (W_i\setminus V') \cup \partial V'  \in {\cal W}_{\lambda}$ be the set given by Lemma \ref{rig-lemma: derive prop2}. Define $B^{i+1}_j = B^i_j \setminus \partial V'$ for $j=0,\ldots,i+1$  and $\bar{u}_{i+1} = \bar{u}_i$. Then the following holds:
\begin{itemize}
\item[(i)]  \eqref{rig-eq: energy balance4} and \eqref{rig-eq: property,d}-\eqref{rig-eq: new5} are satisfied for $W_{i+1}$, $\bar{u}_{i+1}$, $\hat{\omega}_{i+1}$ and $\lambda$ in place of $\lambda_{i+1}$.
\item[(ii)] If $B^{i}_{i+1} = \emptyset$, then \eqref{rig-eq: Asets} holds for $W_{i+1}$  and $\lambda$ in place of $\lambda_{i+1}$.
\end{itemize}
\end{lemma}

We will first give the proof of Theorem \ref{rig-th: derive prop} to see how the above algorithm is applied and postpone the proof of Lemma \ref{rig-lemma: modification} to the end of Section \ref{rig-sec: subsub,  proofimain}. Moreover, the corresponding property (ii) in the case $B^{i}_{i+1} \neq \emptyset$ will be stated below in Lemma \ref{rig-lemma: Asets} since it relies on the construction of $B^{i}_{i+1}$ in \eqref{eq: BII}.

\subsection{Proof of Theorem \ref{rig-th: derive prop}}\label{rig-sec: subsub,  proofimain}


We are now in a position to prove Theorem \ref{rig-th: derive prop}.

\noindent {\em Proof of Theorem \ref{rig-th: derive prop}.}
Using Corollary \ref{rig-cor: derive prop2} we first see that \eqref{rig-eq: energy balance4} holds for $W_0 =  \hat{W}$ and \eqref{rig-eq: Asets}-\eqref{rig-eq: new5}  are trivially satisfied for $\lambda_0 = 0$, $B^0_0 = \emptyset$. Assume that $\lambda_i$, $W_i \in {\cal W}_{\lambda_i}$,   $\lbrace B^i_j: j=0,\ldots, i\rbrace$ and  $\bar{u}_i$ have already been constructed and that \eqref{rig-eq: energy balance4}-\eqref{rig-eq: new5} hold. 

If now all $\Gamma_l(W_i)$ with $\overline{N^{2\hat{\tau}_l}(\partial R_l(W_i))} \subset H(W)$  satisfy $|\Gamma_l(W_i)|_\infty \le \lambda_i$ we stop  and set $U=W_i$.    We  will see that in this case \eqref{rig-eq: D1} holds for $C_1 = \hat{C}\sum^\infty_{n=0}(2/3)^n \omega_{\rm min}^{-2} \upsilon^{-4}$ by \eqref{rig-eq: D2}. Otherwise, there is some smallest $\Gamma=\Gamma(W_i)$ with respect to $\vert \cdot \vert_\infty$ satisfying $\vert \Gamma \vert_\infty >  \lambda_i$ and $\overline{N^{2\hat{\tau}}(\Gamma)} \subset H(W)$  with $\hat{\tau} = q^2\upsilon h_*^{-1} |\Gamma|_\infty$ as defined in \eqref{rig-eq: hat tau def}.  To simplify the exposition, we will suppose that the choice of $\Gamma$ is unique. At the end of the proof we briefly indicate the necessary changes if there are several components of the same size.
  
 Choose $\omega_{\rm min} \ge \sqrt{\tfrac{3}{4}}$.  We observe that ${\cal T}_{\hat{\tau}}(W_i, W) \subset {\cal T}_{\lambda_i}(W_i, W)$. Indeed, for $\hat{\tau} \le \lambda_i$ it is obvious  and for $\hat{\tau} > \lambda_i$ it follows from  the choice of $\Gamma$ with respect to $|\cdot|_\infty$.  Thus, by \eqref{rig-eq: D2} we get that \eqref{rig-eq: property,c}  is satisfied  replacing $C_*$ by $  \frac{4}{3}\hat{C}\sum\nolimits^i_{n=0}  (2/3)^n$.   
 
 \smallskip
 \textbf{Trace estimate.}  If $\Gamma$ additionally fulfills \eqref{rig-eq: property4}, \eqref{rig-eq: property,b} and \eqref{rig-eq: property,a}, we may apply Theorem \ref{rig-theorem: D}. Therefore,  recalling $\omega(\Gamma) = 1$ we get that for suitable $A \in \R^{2 \times 2}_{\rm skew}$, $b\in \R^2$
\begin{align*}
\int_\Gamma |\bar{u}_i(x) - (Ax+b)|^2 \, dx &\le \Big(\hat{C} +  \frac{1}{2}  \cdot \frac{4}{3}\hat{C} \sum\nolimits^i_{n=0}   \Big( \frac{2}{3} \Big)^n \Big) \frac{\eps}{\upsilon^4}|\Gamma|^2_\infty \\ &\le  \hat{C}\sum\nolimits^{i+1}_{n=0}   \Big( \frac{2}{3} \Big)^n\, \frac{\omega(\Gamma)^2}{\hat{\omega}_{i+1}(\Gamma)^2} \frac{\eps}{\upsilon^4}|\Gamma|^2_\infty.
\end{align*}
Thus, \eqref{rig-eq: D2} holds, as desired. We define $\bar{u}_{i+1}(x)= Ax + b$ for $x \in X$ and $\bar{u}_{i+1} = \bar{u}_{i}$ else, where $\partial X = \Gamma$.  Moreover, we set $W_{i+1} = W_i$,  $\lambda_{i+1} = \vert \Gamma \vert_\infty$, $B^{i+1}_{j} = B^i_j$ for $j=0,\ldots,i$ and $B^{i+1}_{i+1} = \emptyset$. Now \eqref{rig-eq: energy balance4}-\eqref{rig-eq: new5} still hold due to choice of $\Gamma$ with respect to $|\cdot|_\infty$. In particular,  for \eqref{rig-eq: property,d}, \eqref{rig-eq: D2} we note that ${\cal T}_{\lambda_{i+1}}(W_{i+1}, W) = {\cal T}_{\lambda_{i}}(W_{i}, W) \cup \lbrace \Gamma \rbrace$ and in \eqref{rig-eq: property,d}  we additionally take \eqref{rig-eq: property,b} for the component $\Gamma$ into account.  Moreover, \eqref{rig-eq: Asets}(iii) still holds as $\lambda_{i+1} \ge \lambda_i$.  As also \eqref{rig-eq: W prop} is satisfied for $\lambda_{i+1}$, we get $W_{i+1} \in {\cal W}_{\lambda_{i+1}}$. We continue with the next iteration step.

\smallskip

Otherwise (a) \eqref{rig-eq: property4}, (b) \eqref{rig-eq: property,b} or (c) \eqref{rig-eq: property,a} is violated. We refer to (D1)-(D3) where the basic situation of each case is described.

\smallskip
\textbf{Case (a).} We find some $\overline{V}  \supset\Gamma$ with $\partial V \neq \Gamma$  and $\overline{V} \subset Q_\mu$ such that setting $\tilde{W}_{i} = (W_i \setminus V) \cup \partial V$, we get $\Vert \tilde{W}_{i} \Vert_\omega \le \Vert W_i \Vert_\omega$ and therefore
\begin{align*}
\eps \Vert \tilde{W}_{i} \Vert_\omega + \alpha(\tilde{W}_{i}) & \le \eps \Vert W_i \Vert_\omega + \alpha(W_i).
\end{align*}
Here we adjusted the weights as in \eqref{rig-eq: weights}. Define $B^i_{i+1} = \emptyset$ and observe that $\hat{\omega}^*_i = \hat{\omega}_i$ with $\hat{\omega}^*_i$ as defined in \eqref{eq:omega*}. Then Lemma \ref{rig-lemma: modification}  applied for $\lambda = \lambda_i$  yields a set $W_{i+1} \in {\cal W}_{\lambda_{i}}$ with $|W_{i+1} \setminus \tilde{W}_i|=0$  and $H(W_{i+1}) \supset H(W_i) \supset H(W)$, corresponding $(B^{i+1}_j)^{i+1}_{j=0}$ and $\bar{u}_{i+1} = \bar{u}_i$ such that  \eqref{rig-eq: energy balance4}-\eqref{rig-eq: new5} hold for $W_{i+1}$, $\bar{u}_{i+1}$, $\hat{\omega}_{i+1}$ and $\lambda_{i}$. Let $\lambda_{i+1} = \vert \Gamma \vert_\infty$ and observe  that  $W_{i+1} \in {\cal W}_{\lambda_{i+1}}$  satisfies \eqref{rig-eq: energy balance4}-\eqref{rig-eq: new5}  also for $\lambda_{i+1}$.   Indeed, \eqref{rig-eq: property,d} and \eqref{rig-eq: D2}    follow from the choice of $\Gamma$ with respect to $|\cdot|_\infty$, Remark \ref{rem: modi} and the fact that $|\partial V|_\infty > \lambda_{i+1}$.  For the other properties we may argue as before.  We now continue with the next iteration step.

\smallskip

\textbf{Case (b).} First, suppose $\vert \partial W_i\cap N^{2\hat{\tau}}\vert_{\cal H} > \frac{16\hat{\tau}}{h_*}$. In view of  $\overline{N^{2\hat{\tau}}} \subset H(W) \subset H(W_i)$ this contradicts the assertion of Lemma \ref{rig-lemma: crack bound}   and shows that also \eqref{rig-eq: property4} is violated. Thus, we can proceed as in case (a). Otherwise, choosing  $D=D(h_*)  \ge \frac{32}{h_*}$, the fact that \eqref{rig-eq: property,b}  does not hold yields $\alpha(N^{2\hat{\tau}} \cap W_i) > \frac{D}{2}\eps \hat{\tau}$.

 We set $\tilde{W}_{i} = (W_i \setminus V) \cup \partial V$, where $V$ is the smallest rectangle containing $N^{2\hat{\tau}}$. We observe that $\overline{V} \subset Q_\mu$  and a short computation yields $|\partial V|_* - |\Gamma|_*\le 8 \cdot 2\hat{\tau}$. Thus,  using Lemma \ref{rig-lemma: derive prop3}(ii) and noting that $\Gamma \in {\cal L}_V$ we derive  
\begin{align*}
\eps \Vert \tilde{W}_{i} \Vert_\omega + \alpha(\tilde{W}_{i}) & \le \eps \Vert W_i \Vert_\omega +  \alpha(W_i) + 16\eps  \hat{\tau} - \alpha(V \cap W_i)  \le \eps \Vert W_i \Vert_\omega   +  \alpha(W_i).
\end{align*}
We now may proceed exactly as in case (a) and then continue with the next iteration step.

\smallskip
\textbf{Case (c).} Let $\Psi_i = \Psi^j$ and $\psi_i = \psi^j$, where $\Psi^j$ is  a set such that \eqref{rig-eq: property,a} is violated.  As discussed above in Section \ref{rig-sec: subsub,  dod-neigh}, we find a boundary component $\Gamma_m=\Gamma_m(W_i)$ with $|\Gamma_m|_\infty \ge \hat{\tau}$, $\omega(\Gamma_m) = 1$. Moreover, there is a rectangle $T \subset Q_\mu$ with $\vert \partial T \vert_{\cal H} \le 4 \psi_i$ and $\overline{T} \cap \Gamma \neq \emptyset$, $\overline{T} \cap\Gamma_m \neq \emptyset$ (cf. Figure \ref{H} below). 

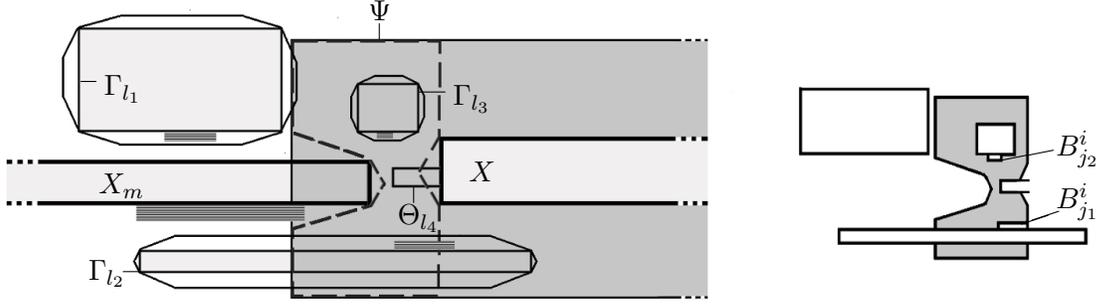
\begin{figure}[H]
\centering
\begin{overpic}[width=1.0\linewidth,clip]{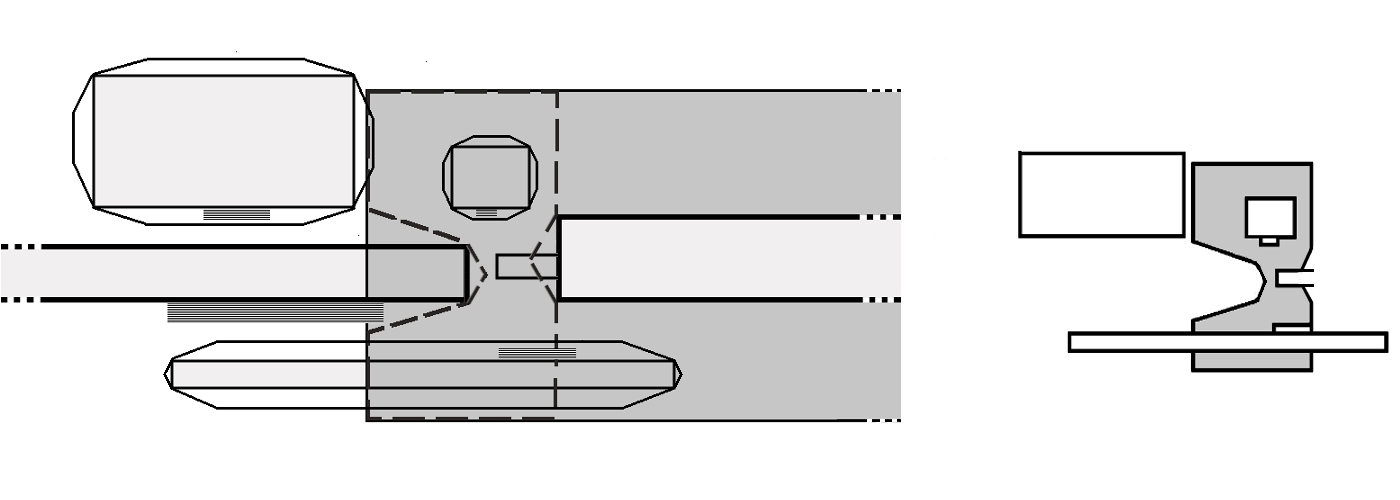}
\put(175,67){\small{$X$}}
\put(35,62.5){\small{$X_m$}}

\put(37,100){\small{$\Gamma_{l_1}$}}
\put(28,105){\line(1,0){6}}
\put(169,96){\small{$\Gamma_{l_3}$}}
\put(156,101){\line(1,0){10}}
\put(44,33){\line(1,0){6}}
\put(31,30){\small{$\Gamma_{l_2}$}}

\put(154,59){\line(0,1){6}}
\put(148,51){\small{$\Theta_{l_4}$}}

\put(397,57){\small{$B^i_{j_1}$}}
\put(397,77){\small{$B^i_{j_2}$}}
\put(384,51){\line(3,1){14}}
\put(376,75){\line(5,1){20}}
\put(141,121){\line(0,1){6}}
\put(137,128){\small{$\Psi$}}

\end{overpic}
\caption{\small On the left side $\Psi$ (the set surrounded by the dashed grey line) and parts of $\Gamma$, $\Gamma_m$ are sketched. Observe that $M^{21\bar{\tau}_m}(\Gamma_m) \cap \Psi = M^{21\bar{\tau}}(\Gamma) \cap \Psi = \emptyset$. Moreover, the picture includes several boundary components with corresponding dodecagonal or  decagonal neighborhoods as well as four striped sets $B^i_{j_1},\ldots,B^i_{j_4}$. On the right hand side the resulting $B^i_{i+1}$ is drawn, where $\partial B^i_{i+1}$ is black and the interior $(B^i_{i+1})^\circ$ is grey. Observe that in general only parts of $\partial B^i_{i+1}$ are contained in $B^i_{i+1}$.} \label{F}
\end{figure}

Let ${\cal A}_{i} \subset (\Gamma_l(W_i))_l \setminus \lbrace \Gamma, \Gamma_m \rbrace$ be the boundary components with $ \Theta_l(W_i) \cap \Psi_i \neq \emptyset$ or,  if $\Gamma_l(W_i) = \Theta_l(W_i) \in {\cal G}_i$, with $M^{\eta_l^i}(\Gamma_l(W_i)) \cap \Psi_i \neq \emptyset$.  (Recall the definition of ${\cal G}_i$ in \eqref{rig-eq: Asets}(ii) and see \eqref{rig-eq: MMhatdef} for the neighborhood.) We now define an additional set $B^{i}_{i+1}$, where we will `use the energy' to modify $W_i$.  Let
\begin{align}\label{eq: BII}
B^{i}_{i+1} = \Big((\Psi_i \cap W_i)\cup \bigcup\nolimits_{\Gamma_l(W_i) \in {\cal A}_{i}}  \Theta_l(W_i)\Big) \setminus  \bigcup\nolimits_{B^{i}_j \in {\cal B}_i} B_j^i,
\end{align}
where ${\cal B}_i:= \lbrace B_j^i: B_j^i \cap W_i \subset M^{\eta_l^i}(\Gamma_l(W_i)) \ \text{ for some } \Gamma_l(W_i) \in {\cal A}_{i} \cap {\cal G}_i\rbrace.$ In the definition of $B^{i}_{i+1}$ it is essential to subtract the set on the right hand side such that we will be able to ensure \eqref{rig-eq: Asets}(i).

Note that by \eqref{rig-eq: Asets} we have  for all $\Gamma_l(W_i) \in {\cal A}_{i}$  either $\Theta_l(W_i) \subset B^{i}_{i+1}$ or $\Theta_l(W_i) \cap B^{i}_{i+1} = \emptyset$ depending on whether $\Theta_l(W_i) \cap \bigcup\nolimits_{B^i_j \in {\cal B}_i} B^i_j = \emptyset$ or $\Theta_l(W_i) \subset B^i_j \in {\cal B}_i$.  Moreover, the components $\Gamma_l(W_i) \notin {\cal A}_{i}$  clearly satisfy $\Theta_l(W_i) \cap B^{i}_{i+1} = \emptyset$ by the definition of ${\cal A}_{i}$. Denote by $\tilde{\cal A}_{i} \subset {\cal A}_{i}$ the boundary components completely contained in $B^{i}_{i+1}$ and observe that $ {\cal G}_i \cap {\cal A}_i \subset \tilde{\cal A}_i $ (see \eqref{rig-eq: Asets}(ii)).  This definition   implies
\begin{align}\label{eq: BIII}
\bigcup\nolimits_{\Gamma_l(W_i) \in \tilde{\cal A}_i} \Theta_l(W_i) =  \partial W_i \cap B^i_{i+1}.
\end{align}
One can show that $|\Gamma_l(W_i)|_\infty < 19\bar{\tau}$ and $\dist(\Gamma_l(W_i),\Psi_i ) \le \bar{\tau}$ for all  $\Gamma_l(W_i) \in {\cal A}_i$ which we postpone to Lemma \ref{lemma: rig-eq: calB prop} below for convenience. This together with \eqref{rig-eq: tau-bartau}, \eqref{rig-eq: hat tau def} and the fact that $\Psi_i \subset N^\tau(\Gamma) \subset N^{\bar{\tau}}(\Gamma)$ implies $\overline{N^{2\hat{\tau}_l}(\partial R_l)} \subset N^{2\hat{\tau}}(\Gamma)  \subset  H(W)$ for all $\Gamma_l(W_i) \in {\cal A}_i$  for $q$ large enough, where $R_l$ are as usual the corresponding rectangles. Then, as $|\Gamma_l(W_i)|_\infty < |\Gamma|_\infty$, by the choice of $\Gamma$ with respect to $|\cdot|_\infty$ we obtain $|\Gamma_l(W_i)|_\infty \le \lambda_i$ for all $\Gamma_l(W_i) \in {\cal A}_i$ and thus  ${\cal A}_i \subset {\cal T}_{\lambda_i}(W_i, W)$.

Consequently, by \eqref{rig-eq: property,d} we have   $\gamma(M^{21\bar{\tau}_l}(\Gamma_l(W_i)) \setminus S(\Gamma_l(W_i))) \le D \eps\hat{\tau}_l$ for all $\Gamma_l(W_i) \in {\cal A}_i \cap {\cal G}_i$, where
\begin{align*}
\gamma(A) := \alpha(A \cap W_i \cap \Psi_i) + \eps |A \cap  (\partial W_i \cap \Psi_i)|_{\cal H}
\end{align*}
for $A \subset \R^2$. Observe  that $\eta^i_l = 21 \bar{\tau}_l$ as $|\Gamma_l(W_i)|_\infty \le \lambda_i$ for all $\Gamma_l(W_i) \in {\cal A}_i$.  Using the definition of ${\cal B}_i$ below \eqref{eq: BII} as well as \eqref{rig-eq: hat tau def}, \eqref{rig-eq: Asets}(iii)   and recalling $D=D(h_*)$ we find for $\upsilon$ small enough (with respect to $h_*$ and $q$, cf. \eqref{eq: parameters})
\begin{align*}
\sum\nolimits_{B^i_j \in {\cal B}_i} \gamma(B_j^i ) &\le \sum\nolimits_{\Gamma_l(W_i) \in {\cal A}_i \cap {\cal G}_i} \gamma(M^{21\bar{\tau}_l}(\Gamma_l(W_i)) \setminus S(\Gamma_l(W_i))) \\ 
&\le \sum\nolimits_{\Gamma_l(W_i) \in \tilde{\cal A}_i \cap {\cal G}_i}  \eps|\Gamma_l(W_i)|_\infty \le  \eps\vert\partial W_i \cap B^{i}_{i+1}\vert_{\cal H}.
\end{align*}
In the first step we used that $B^i_{j_1} \cap B^i_{j_2} \cap W_i= \emptyset$ for $j_1 \neq j_2$ by \eqref{rig-eq: Asets}(i) and $\partial W_i \cap Q_\mu \subset W_i$. The last step follows from \eqref{eq: BIII}.   The fact that \eqref{rig-eq: property,a} is violated  and the definition of $B^i_{i+1}$ in \eqref{eq: BII} then imply  
\begin{align}\label{rig-eq: new 4} 
\begin{split}
D (1-\omega_{\rm min})^{-1}\eps \psi_i &< \alpha(\Psi_i \cap W_i) + \eps\vert\partial W_i \cap \Psi_i\vert_{\cal H} \\
& \le\alpha(\Psi_i \cap W_i) + \eps\vert\partial W_i \cap \Psi_i\vert_{\cal H}  - \sum\nolimits_{B^i_j \in {\cal B}_i} \gamma(B_j^i) \\ & \ \ \  + \eps\vert\partial W_i \cap B^{i}_{i+1}\vert_{\cal H}\\
& \le \alpha(B^{i}_{i+1}) + 2\eps\vert\partial W_i \cap B^{i}_{i+1}\vert_{\cal H}. 
\end{split}
\end{align}
We adjust the weights for components in $\tilde{\cal A}_{i}$: Let $W^*_i = W_i$ and $\omega(\Gamma_l(W^*_i)) = \omega(\Gamma_l(W_i)) - \frac{1 -\omega_{\rm min}}{2}$ for $\Gamma_l(W^*_i) = \Gamma_l(W_i) \in \tilde{\cal A}_{i}$ and $\omega(\Gamma_l(W^*_i)) = \omega(\Gamma_l(W_i))$ otherwise. (The set as a subset of $\R^2$ is left unchanged, we have only changed the weights of the boundary components.) By \eqref{rig-eq: h*, omega}, \eqref{rig-eq: new2}(ii) and \eqref{eq: BIII} we derive
\begin{align}\label{rig-eq: weights3}
\begin{split}
\Vert W^*_i \Vert_\omega &= \Vert W_i \Vert_\omega  - \tfrac{1}{2}(1-h_*)(1-\omega_{\rm min})\sum\nolimits_{\Gamma_l(W_i) \in \tilde{\cal A}_{i}} |\Gamma_l(W_i)|_\infty\\
& \le  \Vert W_i \Vert_\omega -  h_*(1-\omega_{\rm min}) \vert \partial W_i \cap B^{i}_{i+1}\vert_{\cal H}
\end{split}
\end{align}
for $h_*$ small enough. We briefly note that \eqref{rig-eq: D2}, \eqref{rig-eq: new5} are still satisfied for $W^*_i$ if we replace $\hat{\omega}_i$ by  $\hat{\omega}^*_{i}$ as given in  \eqref{eq:omega*}.  Indeed, as $\omega(\Gamma_l(W_i))\ge \hat{\omega}_i(\Gamma_l(W_i))$  by \eqref{rig-eq: new5} we find 
$$\frac{\omega(\Gamma_l(W^*_i))}{\hat{\omega}^*_{i}(\Gamma_l(W^*_i))} = \frac{\omega(\Gamma_l(W_i))- (1-\omega_{\rm min})/2}{\hat{\omega}_{i}(\Gamma_l(W_i))- (1-\omega_{\rm min})/2} \ge  \frac{\omega(\Gamma_l(W_i))}{\hat{\omega}_{i}(\Gamma_l(W_i))}  \ge 1$$
 for  $\Gamma_l(W_i) \in \tilde{\cal A}_i$. We set $\tilde{W}_{i} = (W^*_i \setminus V) \cup \partial V$, where  $\overline{V} \subset Q_\mu$ is the smallest rectangle containing $\Gamma$, $\Gamma_m$ and $T$. As usual we define $\omega(\partial V) = 1$ and adjust the other weights as in \eqref{rig-eq: weights}. Observe that $|\partial V|_* \le |\Gamma|_* + |\Gamma_m|_* + |\partial T|_{\cal H}$. We then derive by \eqref{rig-eq: new 4}, \eqref{rig-eq: weights3} and Lemma \ref{rig-lemma: derive prop3}(ii) in view of $\lbrace \Gamma, \Gamma_m \rbrace \subset {\cal L}_V$  
\begin{align*}
\Vert \tilde{W}_{i}\Vert_\omega &\le \Vert W^*_i \Vert_\omega   + \vert\partial T\vert_{\cal H} \le \Vert W_i \Vert_\omega  - h_*  (1-\omega_{\rm min}) \vert \partial W_i \cap B^{i}_{i+1}\vert_{\cal H} + \vert\partial T\vert_{\cal H} \\
& \le \Vert W_i \Vert_\omega + h_*(1-\omega_{\rm min}) \tfrac{1}{2\eps} \alpha(B^{i}_{i+1}) - \tfrac{1}{2} h_*D\psi_i +\vert\partial T\vert_{\cal H}.
\end{align*}
Recall $\vert \partial T \vert_{\cal H} \le 4 \psi_i$. Now choosing $D \ge \frac{8}{h_*}$  we conclude
$$\eps \Vert \tilde{W}_{i}\Vert_\omega + \alpha(\tilde{W}_{i})  \le \eps \Vert W_i\Vert_\omega + \alpha(W_i) + h_*(1-\omega_{\rm min}) \alpha(B^{i}_{i+1}).$$
Let $\lambda_{i+1} = \vert \Gamma \vert_\infty$. We observe that \eqref{rig-eq: energy balance4}-\eqref{rig-eq: new5} hold for $W_i^*$ for all $\lambda \in [\lambda_i,\lambda_{i+1})$, where we particularly use ${\cal T}_{\lambda}(W_i^*, W) ={\cal T}_{\lambda_i}(W_i^*, W)$ due to the choice of $\Gamma$ with respect to $|\cdot|_\infty$. Choosing  $\lambda' \in [\lambda_i,\lambda_{i+1})$ large enough we get that $W_i^* \in {\cal W}_{\lambda'}$. In fact, \eqref{rig-eq: W prop}(iv) follows from the definition of the weights and Lemma \ref{lemma: rig-eq: calB prop} below in view of $\bar{\tau} = \upsilon \lambda_{i+1}$.  Moreover, \eqref{rig-eq: W prop}(ii) is a consequence of \eqref{rig-eq: new5}(ii). 

Recall that \eqref{rig-eq: D2} and  \eqref{rig-eq: new5} hold for $W_i^*$ and $\hat{\omega}^*_{i}$.  By Lemma \ref{rig-lemma: modification} for $\lambda = \lambda'$ we find a set $W_{i+1} =  (W_i^* \setminus  V') \cup \partial V' \in {\cal W}_{\lambda'}$  for a rectangle $\overline{V} \subset  \overline{V'} \subset Q_{\mu}$ with $|W_{i+1} \setminus \tilde{W}_i|=0$ and $H(W_{i+1}) \supset H(W_i) \supset H(W)$. Moreover, with $\bar{u}_{i+1} = \bar{u}_i$ and the sets $B^{i+1}_j = B^{i}_j \setminus \partial V'$ for $j=0,\ldots,i+1$ we find that \eqref{rig-eq: energy balance4} and \eqref{rig-eq: property,d}-\eqref{rig-eq: new5} hold for $W_{i+1}$, $\bar{u}_{i+1}$, $\hat{\omega}_{i+1}$ and $\lambda'$. Observe that \eqref{rig-eq: Asets}   does not follow from Lemma \ref{rig-lemma: modification} as in general $B^{i}_{i+1} \neq \emptyset.$  We postpone the proof of \eqref{rig-eq: Asets} to Lemma \ref{rig-lemma: Asets} below. Then arguing as before in case (a) we   get  that  $W_{i+1} \in {\cal W}_{\lambda_{i+1}}$  satisfies \eqref{rig-eq: energy balance4}-\eqref{rig-eq: new5}  also for $\lambda_{i+1}$. We continue with the next iteration step.

\smallskip
\textbf{Conclusion.} In each iteration step either the number of components satisfying \eqref{rig-eq: D2} increases or the volume of $W_i$ decreases by at least $(2C_2s)^2$. Consequently, after a finite number of steps, denoted by $i^*$, we find a set $U=W_{i^*} \in {\cal W}^{C_2s}_{\lambda_U}$, $\lambda_U \ge 0$, with $H(U) \supset H(W)$ satisfying \eqref{rig-eq: D2} for all boundary components $\Gamma_l(U)$ with $\overline{N^{2\hat{\tau}_l}(\partial R_l(U))} \subset H(W)$.  Let $\bar{u} = \bar{u}_{i^*}$. Then \eqref{rig-eq: D1} holds for a constant $C_1 = C(h_*) \upsilon^{-4}$   for all such boundary components  as $\hat{C} = \hat{C}(h_*)$, $\sum_n  (2/3)^n < \infty$ as well as $\hat{\omega}(\Gamma_l(U)) \ge \omega_{\rm min }$ for all $\Gamma_l(U)$ by \eqref{rig-eq: Asets}(i). (Below we will show that indeed $C_1=C_1(h_*,\sigma)$.) 

In view of \eqref{eq: parameters}, \eqref{rig-eq: hat tau def} and \eqref{rig-eq: C2 choice} we see that Remark \ref{rig-rem: z}(iii) holds for some $\bar{C}= \bar{C}(h_*)$ large enough with $\bar{C}C_2 \le h_*$. In particular, this also implies that  \eqref{rig-eq: D1} holds for $\Gamma_l(U)$ with $X_l \cap Q_{U} \neq \emptyset$, where $Q_{U} = (-\mu_U,\mu_U)^2$, $\mu_U =  \max\big\{\mu - 3\sum\nolimits_{l=1}^n |\partial R_l|_\infty - \sum\nolimits_{l=n+1}^m|\partial X_l(U)|_\infty,0\big\}$ (see \eqref{eq: Uqmu}). Indeed, $\dist(\partial Q_{\mu_U},\partial H(W)) \ge 2\sum\nolimits_{l=1}^n |\partial R_l|_\infty$ and then the claim follows from Remark \ref{rig-rem: z}(iii) for $h_*$ small.

We now show \eqref{rig-eq: energy balance3}. By \eqref{rig-eq: Asets}(i) we find $\sum^{i^*}_{j=0} \alpha(B^{i^*}_j) \le 2 \alpha(W)$  and  by \eqref{rig-eq: Asets}(i), \eqref{rig-eq: new5} we get ${\omega}(\Gamma_l(U)) \ge \omega_{\rm min }$ for all $\Gamma_l(U)$. Possibly passing to a larger $\omega_{\rm min}$, we may assume that $\sigma':=2(1-\omega_{\rm min}) \le \sigma$  with $\sigma$ as in the assertion. Then using \eqref{rig-eq: energy balance4}  and $h_* \le 1$ we conclude 
$$\eps\Vert U \Vert_* + \alpha(U) \le (1-\tfrac{1}{2}\sigma')^{-1}\eps\Vert W \Vert_* + (1+\sigma')\alpha(W)\le (1+ \sigma)(\eps\Vert W \Vert_* + \alpha(W)).$$
In view of \eqref{rig-eq: h*},\eqref{rig-eq: W prop}(ii) we get $|\partial R_l|_\infty \le \omega_{\rm min}^{-1}(1-h_*)^{-1}|\Theta_l(U)|_*$  and by the choice of $\sigma$ we also get Remark \ref{rig-rem: z}(iv) for $h_*$ small enough. To  see that $C_1$ and $C_2$ indeed only depend on $h_*$ and $\sigma$, we recall that $C_1 = C_1(h_*)\upsilon^{-4}$ and  $C_2 = C'\upsilon$ (see \eqref{rig-eq: C2 choice}). Then the claim follows from the fact that $\omega_{\rm min} = \omega_{\rm min}(h_*,\sigma)$ and $\upsilon =  \upsilon(h_*,\omega_{\rm min})$ by \eqref{eq: parameters}. Likewise,  Remark \ref{rig-rem: z}(i) also follows from \eqref{eq: parameters} and the choice of $\omega_{\rm min}$ with respect to $\sigma$.   Finally, to find  $|W\setminus U| \le  \Vert U \Vert_\infty^2$, we observe $W\setminus U \subset \bigcup_{l=1}^n X_l(U)$ and compute 
  \begin{align}\label{eq: before}
  |W \setminus U| \le \sum\nolimits_l |X_l(U)| \le  \sum\nolimits_l |\Gamma_l(U)|^2_\infty \le \Vert U \Vert^2_\infty.
  \end{align}

\textbf{Adaptions for the general case.}  It remains to indicate the necessary changes if in some iteration step $i$ the choice of $\Gamma$ is not unique.  If there are several components $\Gamma_1,\ldots,\Gamma_m$ with $\lambda_{i+1} := |\Gamma_j|_\infty >\lambda_i$ for $j=1,\ldots,m$, we choose an order such that $\Gamma_1,\ldots,\Gamma_{m'}$, $m' \le m$, are the components satisfying \eqref{rig-eq: property4}, \eqref{rig-eq: property,b}, \eqref{rig-eq: property,a}. We now apply Theorem \ref{rig-theorem: D} successively  on each $\Gamma_j$, $j=1,\ldots,m'$, and replace ${\cal T}_{\lambda_{i+1}}(W_{i+1},W)$ in \eqref{rig-eq: property,d}, \eqref{rig-eq: D2} by ${\cal T}^j_{\lambda_{i+1}}(W_{i+1},W) := {\cal T}_{\lambda_{i}}(W_{i},W) \cup \bigcup^j_{k=1}\lbrace \Gamma_k \rbrace$. For each $\Gamma_j$, $j=m'+1,\ldots,m$, we proceed as in one of the cases (a)-(c) and let ${\cal T}^j_{\lambda_{i+1}}(W_{i+1},W) := {\cal T}_{\lambda_{i}}(W_{i},W) \cup \bigcup^{m'}_{k=1}\lbrace \Gamma_k \rbrace$ in \eqref{rig-eq: property,d}, \eqref{rig-eq: D2}. 
\eop
 
 \bigskip
 
 Recall that we have first proven  Theorem \ref{rig-th: derive prop} and have postponed the proof of some lemmas. We now   concern ourselves with Lemma   \ref{rig-lemma: modification}. 
 
 \noindent {\em Proof of Lemma \ref{rig-lemma: modification}.} Let $\tilde{W} = (W_i \setminus V) \cup \partial V$ for some rectangle $V $ with $\vert \partial V \vert_\infty > \lambda$ and choose $\overline{V} \subset \overline{V'} \subset Q_\mu$ such that $W_{i+1} := (W_i \setminus V') \cup \partial V' \in {\cal W}_\lambda$ as in  Lemma \ref{rig-lemma: derive prop2}.   Due to the definition of boundary components in Section \ref{rig-sec: sub, modification}  (see \eqref{rig-eq: new1XXXX}) we observe that for each $\Gamma_l(W_{i+1}) \neq \partial V'$ there is a (unique) corresponding $\Gamma_l(W_i)$ (for notational convenience we use the same index) such that $\Theta_l (W_{i+1}) =  \Theta_l(W_i) \setminus \overline{V'}$. In particular, as $B^{i+1}_j = B^i_j \setminus \partial V'$ for $j=0,\ldots,i+1$, this implies
  \begin{align}\label{eq:forlater}
 \Theta_l(W_i) \subset B^i_j \ \ \ \Rightarrow \ \ \ \Theta_l(W_{i+1}) \subset B^{i+1}_j.
 \end{align}  
We begin with the proof of (i). First,  by assumption and \eqref{rig-eq: energy balance2} we have 
 $\eps\Vert W_{i+1} \Vert_\omega + \alpha(W_{i+1})\le \eps\Vert W_i\Vert_\omega + \alpha(W_i),$
 where we adjust the weights  as described in \eqref{rig-eq: weights}. As $|B^i_j \setminus B^{i+1}_j| = 0$ for all $j=0,\ldots,i$, \eqref{rig-eq: energy balance4}  is trivially satisfied.

We now confirm \eqref{rig-eq: new5}(i), i.e. $\omega(\Gamma_l(W_{i+1}))\ge \hat{\omega}_{i+1}(\Gamma_l(W_{i+1}))$ for all $\Gamma_l(W_{i+1})$. First, we  observe that $\hat{\omega}_{i+1}(\partial V') = \omega(\partial V') = 1$. Moreover,  recalling \eqref{eq:omega*} we get  in view of \eqref{eq:forlater} 
\begin{align}\label{eq:forlater2}
\hat{\omega}_{i+1}(\Gamma_l(W_{i+1})) \le \hat{\omega}^*_{i}(\Gamma_l(W_{i}))
\end{align} 
for all $\Gamma_l(W_{i+1}) \neq \partial V'$,  where $\Gamma_l(W_{i})$ is the unique corresponding component.  This together with the fact that $\omega(\Gamma_l(W_{i+1})) \ge \omega(\Gamma_l(W_{i}))$ (see \eqref{rig-eq: weights2}) yields the desired property.  
We now show that \eqref{rig-eq: D2} remains true. As a preparation we find
\begin{align}\label{rig-eq: new41}
\begin{split}
(i)& \  \omega(\Gamma_l(W_{i+1})) \vert \Gamma_l(W_{i+1})) \vert_\infty = \omega(\Gamma_l(W_i)) \vert \Gamma_l(W_i) \vert_\infty, \  \text{ for all } \  \Gamma_l(W_{i+1}) \neq \partial V',\\
(ii)&  \ |\Gamma_l(W_{i})|_\infty \le \lambda  \ \ \ \text{ for all } \ \ \ \Gamma_l(W_{i+1}) \in {\cal T}_{\lambda}(W_{i+1}, W),
\end{split}
\end{align}
where $\Gamma_{l}(W_i)$ is the unique corresponding component.   If $\Gamma_l(W_{i+1}) \cap \partial V' = \emptyset$, then $\Gamma_l(W_{i+1}) = \Gamma_l(W_i)$ and therefore $\omega(\Gamma_l(W_{i+1})) =\omega(\Gamma_l(W_i))$ by  \eqref{rig-eq: weights} and in addition, if $\Gamma_l(W_{i+1}) \in {\cal T}_{\lambda}(W_{i+1}, W)$, we have $|\Gamma_l(W_i)|_\infty = |\Gamma_l(W_{i+1})|_\infty \le \lambda$.  Otherwise, as $W_{i+1} \in {\cal W}_{\lambda}$, we deduce  $\omega(\Gamma_l(W_{i+1})) < 1$ by \eqref{rig-eq: W prop}(v) and then (i) also follows from \eqref{rig-eq: weights}. Note that $|\Gamma_{l}(W_{i+1})|_\infty \le 19\upsilon \lambda$ by \eqref{rig-eq: W prop}(iv)  and thus \eqref{rig-eq: new41}(i)   implies  $|\Gamma_{l}(W_{i})|_\infty \le \lambda$ for $\upsilon$ small enough.

Moreover, we observe that also \eqref{eq: etatau2} holds since in the modification procedure described in Section \ref{rig-sec: sub, modification} we never change the rectangles $\partial R_l$ given by \eqref{rig-eq: W prop}(i)(v). Consequently, we get $\Gamma_{l}(W_i) \in {\cal T}_\lambda(W_i, W)$ if $\Gamma_{l}(W_{i+1}) \in {\cal T}_\lambda(W_{i+1}, W)$. This together with the fact that $\bar{u}_{i+1} = \bar{u}_i$  as well as  \eqref{eq:forlater2} and \eqref{rig-eq: new41}(i) yields \eqref{rig-eq: D2}. By the same argument as in \eqref{eq:forlater2}, \eqref{rig-eq: new41}(i) and the fact that the rectangles $\partial R_l$ are never changed, property \eqref{rig-eq: new5}(ii) is satisfied.    

To conclude the proof of (i), it remains to show \eqref{rig-eq: property,d}. For a given $\Gamma_l(W_{i+1}) \in {\cal T}_{\lambda}(W_{i+1}, W) \cap {\cal G}_{i+1}$ we deduce $\Gamma_l(W_{i+1}) \cap \partial V' = \emptyset$  by \eqref{rig-eq: W prop}(v) and thus $\Gamma_l(W_{i+1}) = \Gamma_l(W_i) \in {\cal T}_{\lambda}(W_{i}, W) \cap {\cal G}_{i}$. The assertion now follows from the $i$-th iteration step of \eqref{rig-eq: property,d}.  In fact, for the left part it suffices to recall $|W_{i+1} \setminus  W_i|=0$.  For the right part we note $\partial V' \subset S(\Gamma_l(W_{i+1}))$ as $|\partial V'|_\infty > \lambda$ and $\omega(\partial V') =1$. Moreover, Remark \ref{rem: modi} implies  
$(S(\Gamma_l(W_{i})) \setminus S(\Gamma_l(W_{i+1}))) \cap  \partial W_{i+1} = \emptyset$. Consequently, as   $\partial W_{i+1} \setminus \partial W_i \subset \partial V'$, an elementary computation shows $\partial W_{i+1} \setminus S(\Gamma_l(W_{i+1})) \subset \partial W_{i} \setminus S(\Gamma_l(W_{i}))$.

\smallskip

We now show (ii). Clearly, \eqref{rig-eq: Asets}(i) still holds as  $B^{i+1}_j \subset B^i_j$ for all $j=0,\ldots,i$ and $B^{i+1}_{i+1} = \emptyset$.  To see \eqref{rig-eq: Asets}(ii), we first observe   $\partial V' \in {\cal G}_{i+1}$ by definition of $(B^{i+1}_j)_{j=1}^{i+1}$ and it thus suffices to consider $\Gamma_l (W_{i+1}) \neq \partial V'$.   If the corresponding $\Gamma_l(W_i)$ lies in ${\cal G}_{i}$, we immediately get $\Gamma_l(W_{i+1}) \in {\cal G}_{i+1}$ by \eqref{rig-eq: weights2} and the fact that the sets $(B^i_j)^{i}_{j=1}$ do not become larger. Consequently, we can assume that $\Theta_l(W_i) \subset B^i_j$ for some $j$. Then \eqref{eq:forlater} implies \eqref{rig-eq: Asets}(ii).

 Due to  Remark \ref{rem: modi}  for all $\Gamma_l(W_i) \in {\cal G}_i$ we find a $\Gamma_j(W_{i+1}) \in {\cal G}_{i+1}$ such that  $\Gamma_l(W_i) \subset \overline{X_j(W_{i+1})}$, where $\partial X_j(W_{i+1}) = \Gamma_j(W_{i+1})$. In fact, one can choose either $\Gamma_l(W_i)$ itself or $\partial V'$. (Note that both are elements of ${\cal G}_{i+1}$.)  Therefore, $M^{\eta^i_l}_k(\Gamma_l(W_i)) \subset M^{\eta^{i+1}_j}_k(\Gamma_j(W_{i+1}))$  for $k=1,2$. Moreover, in both cases we have $( S(\Gamma_j(W_{i+1}))\setminus S(\Gamma_l(W_{i}))) \cap (W_{i+1} \setminus \partial V') = \emptyset$ by Remark \ref{rem: modi}. This together with $W_{i+1} \setminus W_i \subset \partial V'$ and $B^{i+1}_j = B^i_j \setminus \partial V'$ yields \eqref{rig-eq: Asets}(iii).  \eop

 We now prove a property for sets in ${\cal A}_i$ needed in the proof of the main result. For the following lemma recall the construction of the sets $\Psi_i$ in Section \ref{rig-sec: subsub,  dod-neigh}.

\begin{lemma}\label{lemma: rig-eq: calB prop}
Let be given the situation of case (c) in the proof of Theorem \ref{rig-th: derive prop}. Then we have 
\begin{align}\label{eq: rig-eq: calB prop}
|\Gamma_l(W_i)|_\infty  < 19\bar{\tau} \ \text{ and } \ \dist(\Gamma_l(W_i),\Psi_i ) \le \bar{\tau}\ \ \  \text{ for all } \  \Gamma_l(W_i) \in {\cal A}_i.
\end{align}
\end{lemma}

\Proof Let $\Gamma_m$ be the component considered in case (c) (cf. before \eqref{eq: BII}). Recall that  $\Gamma_m \cap \overline{\Psi}_i \neq \emptyset$ and $|\Gamma_m|_\infty \ge \hat{\tau} \ge 19 \cdot 22\bar{\tau}$ for $q$ large enough (see \eqref{rig-eq: tau def} and \eqref{rig-eq: hat tau def}). For sets $\Gamma_l=\Gamma_l(W_i) \in {\cal A}_i$ intersecting $\Psi_i \subset N^{\bar{\tau}}(\Gamma)$ the assertion is clear by construction of $\Psi_i$ and Lemma \ref{rig-lemma: two components}.  (We can assume that property \eqref{rig-eq: property4} holds and Lemma \ref{rig-lemma: two components} is applicable as otherwise we would have applied case (a).) 

Now assume $\Gamma_l \cap  \Psi_i = \emptyset$ but $M^{\eta_l^i}(\Gamma_l) \cap \Psi_i \neq \emptyset$ for $\Gamma_l \in {\cal A}_i \cap  {\cal G}_i$, which implies $\dist(\Gamma_l,\Psi_i) < \eta^l_i \le 21\upsilon\lambda_i \le 21\bar{\tau}$ since $|\Gamma|_\infty > \lambda_i$. In particular, this yields $\Gamma_l \cap N^{22\bar{\tau}}(\Gamma) \neq \emptyset$.  In view of  $|\Gamma_m|_\infty \ge \hat{\tau} \ge 19 \cdot 22\bar{\tau}$ we apply Lemma \ref{rig-lemma: two components}  for $t = 22\bar{\tau}$  and derive that $|\Gamma_l|_\infty \le 19 \cdot 22 \bar{\tau}$. 
Repeating the above arguments we obtain $\dist(\Gamma_l,\Psi_i) \le \eta^l_i \le 21\bar{\tau}_l \le \upsilon  \cdot 21\cdot 19\cdot 22\bar{\tau} < \frac{\bar{\tau}}{2}$ for $\upsilon$ small enough. This gives the second part of \eqref{eq: rig-eq: calB prop}. Moreover, we have $\Gamma_l \cap N^{\bar{\tau}}(\Gamma) \neq \emptyset$ as  $\dist(\Gamma_l,N^\tau(\Gamma)) < \frac{\bar{\tau}}{2} $ and $\tau \le \frac{\bar{\tau}}{2}$ by \eqref{rig-eq: tau-bartau}. Then the first part of \eqref{eq: rig-eq: calB prop} follows again from Lemma \ref{rig-lemma: two components}. \eop

 \begin{figure}[H]
\centering
\begin{overpic}[width=.9\linewidth,clip]{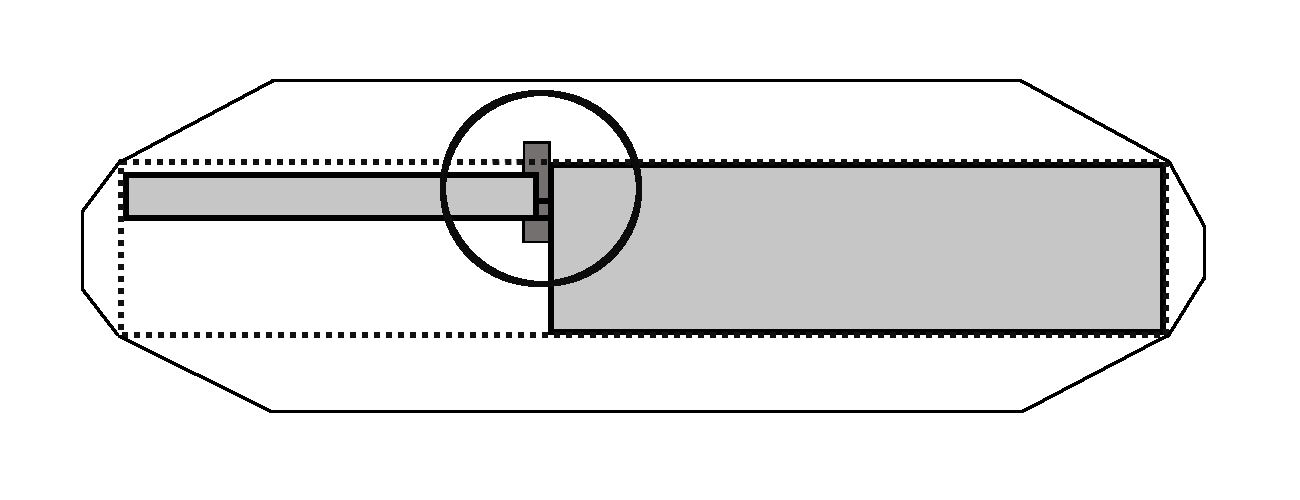}
\put(225,79){\small{$\Gamma$}}
\put(228,89){\line(0,1){6}}

\put(75,64){\small{$\Gamma_m$}}
\put(78,74){\line(0,1){6}}

\put(100,55){\small{$\partial V$}}
\put(105,47){\line(0,1){6}}

\put(138,98){\small{$\Psi$}}
\put(149,101){\line(-1,0){3.5}}

\put(180,107){\small{$B$}}
\put(173,111){\line(1,0){6}}

\put(231,120){\line(0,1){6}}
\put(227,128){\small{$M^{21\upsilon \lambda'}_2(\partial V)$}}

\put(328,32){\small{$M^{21\upsilon \lambda' }_1(\partial V)$}}
\put(343,43){\line(0,1){18}}

\put(164,78){\small{$\partial T$}}
\put(158,82){\line(1,0){6}}

\end{overpic}
\caption{\small Sketch of the components $\Gamma$, $\Gamma_m$, $\partial T$ and in dotted lines the corresponding rectangle $V$ (which in this example coincides with $V'$). The ball $B$ is chosen large enough such that $\Phi \subset B$. (The proportions were adapted for illustration purposes.)} \label{H}
\end{figure}
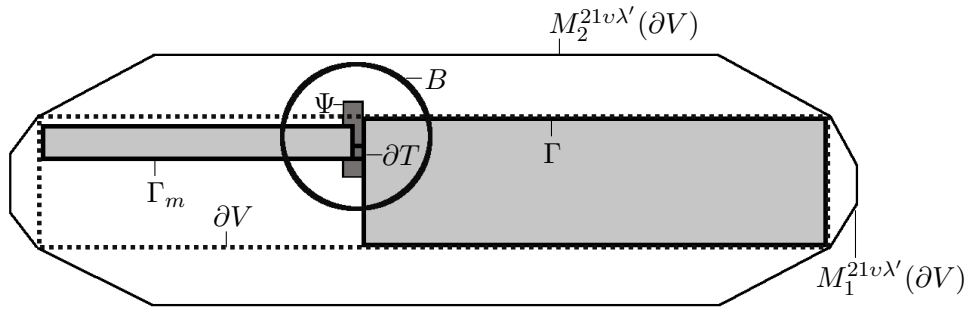

It remains to formulate and prove the analog of Lemma \ref{rig-lemma: modification}(ii)   in case (c). Recall the construction of the sets in case (c) of the proof of Theorem \ref{rig-th: derive prop}.   We refer to  (E) and Figure \ref{1C} where the basic idea has been described.

\begin{lemma}\label{rig-lemma: Asets}
If in the $i$-th iteration step of the above modification procedure case (c) is applied,  then for   $\lambda' \in [\lambda_i,\lambda_{i+1})$ sufficiently large we get that  \eqref{rig-eq: Asets} holds for $W_{i+1}$  and $\lambda'$ in place of $\lambda_{i+1}$.
\end{lemma}

\Proof  Note that by Lemma \ref{rig-lemma: modification} we have $W_{i+1} = (W^*_i \setminus V') \cup \partial V'$ for a rectangle $V'$ which contains $\Gamma$, $\Gamma_m$ and $T$.  Moreover, recall that $W_i = W_i^*$ only differ by the definition of the weights.

First of all, to see \eqref{rig-eq: Asets}(ii) it suffices to show that either $\Theta_l(W^*_i) \subset B^{i}_{j}$ for some $0 \le j \le i+1$ or $\Theta_l(W^*_i) \cap \bigcup^{i+1}_{j=0} B_j^i = \emptyset$ and $\omega(\Gamma_l(W^*_i)) = 1$. In fact, we can then follow exactly the argumentation in the proof of Lemma \ref{rig-lemma: modification}, particularly using \eqref{eq:forlater}, to obtain the desired property also for the sets $B^{i+1}_j = B^i_j \setminus \partial V'$, $j=0,\ldots,i+1$.  

Recall that $\Theta_l(W^*_i) \subset B^{i}_{i+1}$ or $\Theta_l(W^*_i) \cap B^{i}_{i+1} = \emptyset$ for all $\Gamma_l(W^*_i)$  (see before \eqref{eq: BIII}). Thus, if $\Theta_l(W^*_i) \not\subset B^{i}_{j}$ for some $0 \le j \le i+1$ we find $\Theta_l(W^*_i) \cap \bigcup^{i+1}_{j=0} B_j^i = \emptyset$ by \eqref{rig-eq: Asets}(ii) for iteration step $i$. This particularly implies $\Gamma_l(W^*_i) \notin \tilde{\cal A}_i$ as $\Theta_l(W^*_i) \cap B^{i}_{i+1} = \emptyset$. Again by \eqref{rig-eq: Asets}(ii) and the construction of the weights in \eqref{rig-eq: weights3} we then get  $\omega(\Gamma_l(W^*_i))  =1$, as desired.

We concern ourselves with \eqref{rig-eq: Asets}(i).  First, the assertion is clear for $x \in W_{i+1} \setminus W_i$ as $W_{i+1} \setminus W_i \subset \partial V'$  and $\partial V' \cap \bigcup_{j=0}^{i+1} B^{i+1}_j = \emptyset$.  For $x \notin W_{i+1} \setminus W_i$ it is enough to show the property for $(B^i_j)^{i+1}_{j=0}$ since $B^{i+1}_j \subset B^i_j$ for $j=0,\ldots,i+1$. As   $\bigcup^n_{l=1} \Theta_l(W_i) \subset W_i$ and thus $B^i_{i+1} \subset W_i$ (see \eqref{eq: BII}), it is elementary to see that it suffices to confirm $B^{i}_{i+1} \cap \bigcup\nolimits^i_{j=0} B^{i}_j \subset W_i \setminus W_{i+1}$.  Recall that $\Gamma,\Gamma_m \notin {\cal A}_i$. We now show that 
\begin{align}\label{eq: BIIIIII}
B^{i}_{i+1} \cap \bigcup\nolimits^i_{j=0} B^{i}_j \subset B^{i}_{i+1}  \cap \big( f(M^{21\bar{\tau}}(\Gamma)) \cup f(M^{21\bar{\tau}_m}(\Gamma_m) )\big),
\end{align}
where $f(A) = A$ if $A \cap \Psi_i \neq \emptyset$ and $f(A) = \emptyset$ else for $A \subset \R^2$.   (The possible different cases can be seen in Figure \ref{C}, \ref{G}, \ref{F}.) To see this, let ${\cal A}^* \subset \lbrace \Gamma,\Gamma_m\rbrace$ such that the boundary component is contained in ${\cal A}^*$ if the corresponding neighborhood intersects $\Psi_i$.  Fix $j$. Observe that if $B_j^i \cap B_{i+1}^i \neq \emptyset$, then  by \eqref{rig-eq: Asets}(iii) (for $W_i$) and the definition of $B^{i}_{i+1}$ in \eqref{eq: BII}  we get  
\begin{align}\label{eq: BIIIIIIX}
B_j^i \cap B_{i+1}^i \subset B_j^i \cap W_i \subset M^{\eta^i_k}(\Gamma_k(W_i)) \  \ \ \text{for some} \ \ \   \Gamma_k(W_i) \in {\cal G}_i \setminus {\cal A}_i.
\end{align}
On the other hand, by \eqref{rig-eq: Asets}(ii) and the definition of ${\cal A}_i$ we derive that each $\Gamma_l(W_i) \in {\cal A}_i$ with $\Theta_l(W_i) \subset B_j^i$ satisfies $\Gamma_l(W_i) \notin {\cal G}_i$ and thus $\Theta_l(W_i) \cap \Psi_i \neq \emptyset$. Therefore, $\Theta_l(W_i) \not\subset M^{\eta^i_k}(\Gamma_k(W_i))$ for all $\Gamma_k(W_i) \in {\cal G}_i \setminus ({\cal A}_i \cup {\cal A}^* )$. Likewise, we get $\Psi_i \not\subset M^{\eta^i_k}(\Gamma_k(W_i))$ for all  $\Gamma_k(W_i) \in {\cal G}_i \setminus ({\cal A}_i \cup {\cal A}^* )$ and thus    $(B_j^i \cap B_{i+1}^i) \not\subset M^{\eta^i_k}(\Gamma_k(W_i))$  for all  $\Gamma_k(W_i) \in {\cal G}_i \setminus ({\cal A}_i \cup {\cal A}^* )$.This implies $B_j^i \cap B_{i+1}^i \subset M^{\eta^i_k}(\Gamma_k(W_i))$ for some $\Gamma_k(W_i) \in {\cal A}^*$ by \eqref{eq: BIIIIIIX} and shows \eqref{eq: BIIIIII} as $j$ was arbitrary.

Setting $\Phi:= \lbrace x \in Q_\mu: \dist(x,\Psi_i) \le 20\bar{\tau} \rbrace$ and recalling \eqref{eq: BII}, \eqref{eq: rig-eq: calB prop} we then find by \eqref{eq: BIIIIII}
$$B^{i}_{i+1} \cap \bigcup\nolimits^i_{j=0} B^{i}_j \subset \Phi \cap \big( f(M^{21\bar{\tau}}(\Gamma)) \cup f(M^{21\bar{\tau}_m}(\Gamma_m) )\big).$$
 We are now in the position to confirm $B^{i}_{i+1} \cap \bigcup\nolimits^i_{j=0} B^{i}_j \subset W_i \setminus W_{i+1}$ and   differ the cases (I) and (II) as considered in Lemma \ref{rig-lemma: two observations1}, \ref{rig-lemma: two observations2}. In case (I) we get $\Phi \cap W_{i+1} = \emptyset$  as by Lemma \ref{rig-lemma: two observations1}(i) the rectangle $V'$  satisfies $\Phi \subset V'$. In (II)(i) the assertion follows as $ M^{21\bar{\tau}}(\Gamma), M^{21\bar{\tau}_m}(\Gamma_m) \cap \Psi_i = \emptyset$ . Finally, in (II)(ii)  it suffices to derive
\begin{align}\label{rig-eq: new6}
B^{i}_{i+1} \cap \bigcup\nolimits^i_{j=0} B^{i}_j \subset \Phi \cap\lbrace x:  x_1 \ge -l_1 -\psi\rbrace \cap M^{21\bar{\tau}_m}(\Gamma_m),
\end{align}
 where without restriction we treat the case $\Gamma_m \cap N^\tau(\Gamma) \subset N^\tau_{1,-}(\Gamma)$. Then Lemma \ref{rig-lemma: two observations2}(i) gives $\Phi \cap\lbrace x: x_1 \ge -l_1 -\psi\rbrace \cap M^{21\bar{\tau}_m}(\Gamma_m) \subset V'$ which finishes the proof of \eqref{rig-eq: Asets}(i).
To see \eqref{rig-eq: new6}, first note that $f(M^{21\bar{\tau}}(\Gamma)) = \emptyset$ and $\Psi_i \subset \lbrace x:  x_1 \ge -l_1 -\psi\rbrace$.  Consequently, recalling \eqref{rig-eq: Asets}(ii),(iii), if the assertion was wrong, there would be  some $\Gamma_l(W_i) \in  {\cal A}_i \setminus {\cal G}_i$ which satisfies $\Theta_l(W_i) \subset M^{21\bar{\tau}_m}(\Gamma_m)$ and $\Theta_l(W_i) \cap  \lbrace x:  x_1 < -l_1 -\psi\rbrace \neq \emptyset$. Again by \eqref{rig-eq: Asets}(iii) we then get $\Theta_l(W_i) \subset M^{21\bar{\tau}_m}_2(\Gamma_m)$ (see Figure \ref{G}) and therefore $\Theta_l(W_i) \cap \Psi_i = \emptyset$. This implies $\Gamma_l(W_i) \notin {\cal A}_i$ and yields a contradiction.

Finally, we now show \eqref{rig-eq: Asets}(iii).  It suffices to consider  $B^{i+1}_{i+1}$ as for the other sets the property follows by an argumentation as in the proof of Lemma \ref{rig-lemma: modification}.  Without restriction we set $V' = (-v_1, v_1) \times (-v_2,v_2)$. We first observe  that $ \Phi \setminus  \overline{V'}   \subset N^{21\upsilon\lambda'}(\partial V')$ for $\lambda' \in [\lambda_i,\lambda_{i+1})$ large enough.  This is a consequence of the definition of $\Phi$ and the fact that $\bar{\tau}= \upsilon|\Gamma|_\infty = \upsilon\lambda_{i+1}$. We may assume that $|\pi_1 \Gamma_m| \ge \frac{1}{2}|\pi_2\Gamma_m|$.  In fact,  if $l_2 \le \frac{l_1}{2}$ this follows from Lemma \ref{rig-lemma: two components}  and the remark below Lemma \ref{rig-lemma: two components},  if $l_2 \ge \frac{l_1}{2}$ then $l_1$, $l_2$ are comparable and the assumption holds possibly after a rotation of the components by $\frac{\pi}{2}$.  (As before, $l_1$, $l_2$ denote the sidelengths of the rectangle $\Gamma$.) As $|\Gamma_m|_\infty \ge \hat{\tau}$, we   obtain $|\pi_1 \Gamma_m| \ge \frac{1}{\sqrt{5}} \hat{\tau}= \frac{q^2\bar{\tau}}{\sqrt{5}h_*}$. Recall that $|\pi_1 \Gamma_m \cap \pi_1\Gamma| \le 8\frac{\bar{\tau}}{h_*}$ by Lemma \ref{rig-lemma: crack bound} (for $t= \bar{\tau}$). We now find for all  $x \in \Phi \setminus V' $  with $x_1,x_2 \ge 0$
$$ \frac{v_1 - x_1}{x_2 - v_2} \ge \frac{\min\lbrace |\pi_1\Gamma|, |\pi_1 \Gamma_m| \rbrace -  8\bar{\tau}h_*^{-1} }{21\bar{\tau}} \ge \frac{q^2\bar{\tau} (\sqrt{5}h_*)^{-1} -  8\bar{\tau}h_*^{-1} }{21\bar{\tau}} \ge q h_*^{-1} $$
for  $q$ sufficiently large and  may proceed likewise for $\pm x_1, \pm x_2 \ge 0$. Thus, recalling \eqref{rig-eq: MMhatdef2}  and \eqref{rig-eq: MMhatdef}, we obtain $B^{i}_{i+1}  \cap W_{i+1} \subset \Phi  \setminus  \overline{V'} \subset M^{21\upsilon\lambda'}_{2}(\partial V')$  (cf. Figure \ref{H}).  Moreover, $B^i_{i+1} \cap S(\partial V') = \emptyset$ due to \eqref{eq: rig-eq: calB prop} with $S(\partial V')$ as defined in \eqref{rig-eq: S shorthand} (as a component of $W_{i+1}$). Consequently, as  $\upsilon |\partial V'|_\infty \ge \upsilon\lambda'$, $\omega(\partial V')= 1$ and $\partial V' \cap \bigcup^{i+1}_{j=0} B_j^{i+1} = 0$,  we finally obtain \eqref{rig-eq: Asets}(iii).  \eop 
 
 We close this section with a remark useful for the analysis in \cite{Friedrich-Schmidt:15}.

 \begin{rem}\label{rig_rem: connect}
 {\normalfont

(i) During the modification process in Theorem \ref{rig-th: derive prop} the components  $X_{n+1}(W), \ldots, X_{m}(W)$  at the boundary of $Q_\mu$ might be changed and the corresponding components of $U$ are given by $X_j(U) = X_j(W) \setminus \overline{H(U)}$ for $j=n+1,\ldots,m$. In particular, by Lemma \ref{rig-lemma: inftyX}(ii) we observe $|\partial X_j(U)|_* \le |\partial X_j(W)|_*$.
 
(ii) In view of  the modification process described in Section \ref{rig-sec: sub, modification}, in particular Lemma \ref{rig-lemma: derive prop2}, we find that each rectangle $R_j(U)$ given by \eqref{rig-eq: W prop}(i),(v) is contained in $\overline{H(U) \setminus U}$.
  
(iii) In general, the components of the set $U$ might not be connected as they can  be separated by other components. However, by combining boundary components with nonempty intersection we can pass to a configuration $U' \in {\cal V}^{C_2s}$ with connected boundary components satisfying $U' \subset U$, $|U \setminus U'| = 0$, $H(U) = H(U')$ and $\Vert U' \Vert_* \le \Vert U \Vert_*$, where the last property follows from Lemma \ref{rig-lemma: inftyX}. Then by Corollary \ref{rig-cor: derive prop2} we obtain a set $U'' \subset U$ with $\Vert U'' \Vert_* \le \Vert U \Vert_*$ and $|U\setminus U''| \le  \Vert U'' \Vert_\infty^2 \le C\mu \Vert U'' \Vert_\infty$ such that all components of $U''$ are pairwise disjoint and rectangular.  As in (ii) we find that for each $\Gamma_j(U)$ the corresponding rectangle $R_j(U)$ is contained in a boundary component of $U''$.

 }
 \end{rem}

\section{Proofs of lemmas}\label{rig-sec: proofs}

In this section we prove the lemmas stated in Section  \ref{rig-sec: subsub,  rec-neigh} and Section  \ref{rig-sec: sub, modification}.

\subsection{Modifications of sets}\label{rig-sec: proofs1}

Before we give the proof of Lemma \ref{rig-lemma: derive prop2} and Corollary \ref{rig-cor: derive prop2}, we collect some further properties of   $\vert \cdot \vert_*$ (see \eqref{rig-eq: h*}).

\begin{lemma}\label{rig-lemma: infty}
Let $W\in {\cal V}^s$ and  let $\Gamma=\Gamma(W)$ be a  boundary component with $\Gamma = \partial X$. Moreover,  let  $V \in {\cal U}^s$ be a rectangle with $\overline{V} \cap \overline{X} \neq \emptyset$.  Then for $h_*$ sufficiently small and $c>0$ small  (independent of $W$, $\Gamma$ and $V$) the following holds:
\begin{itemize} 
\item[(i)] If $\Gamma$ is connected, then $|\Gamma|_* \ge |\partial R(\Gamma)|_*$, where $R(\Gamma)$ denotes the smallest (closed) rectangle such that $\Gamma \subset R(\Gamma)$.
\item[(ii)] If $X$ is a rectangle, $|\Gamma|_\infty \le c|\partial V|_\infty$ and $\Gamma \setminus \overline{V}$ is not connected, then we have $|\partial (V \cup X)|_* \le |\partial V|_*  + \frac{1}{2}(1-h_*)|\Gamma|_\infty  $.
\end{itemize}
\end{lemma}

\Proof To see (i) it suffices to observe that for $\Gamma$ connected one has  $|\Gamma|_{\cal H} \ge |\partial R(\Gamma)|_{\cal H}$ and $|\Gamma|_{\infty} = |\partial R(\Gamma)|_{\infty}$.  Consider (ii). After translation and rotation we can assume $V = (-a,a) \times (0,b)$ and $\pi_1 X = (-d_1,d_2)$ with $d_1,d_2\ge a$  as well as $\pi_2 X \subset (0,b)$, where $\pi_1 X, \pi_2 X$ denote the orthogonal projections onto the coordinate axes. Letting $d = d_1 + d_2$ we derive $|\partial (V \cup X)|_\infty = \sqrt{d^2 + b^2} \le b + \frac{d^2}{2b} \le b + \frac{d}{4} \le |\partial V|_\infty + \frac{1}{4}|\Gamma|_\infty$. Here we used that $2d \le b$ for $c$ small enough.  As  $|\partial (V \cup X)|_{\cal H} \le |\partial V|_{\cal H} + |\Gamma|_{\cal H} \le |\partial V|_{\cal H} + 2\sqrt{2}|\Gamma|_{\infty}$ by Lemma \ref{rig-lemma: inftyX}(iii),(iv) the claim now follows from \eqref{rig-eq: h*}   for $h_*$ small enough. \eop

\smallskip

We now give the proof of Lemma  \ref{rig-lemma: derive prop2}. Recall \eqref{rig-eq: new1} and the definition of the components and the weights in \eqref{rig-eq: new1XXXX} and \eqref{rig-eq: weights}, respectively.

\noindent {\em Proof of Lemma \ref{rig-lemma: derive prop2}.} Without restriction we can assume $\overline{V}\cap W \neq \emptyset$ as otherwise there is nothing to show.   Let $\tilde{W} = (W \setminus V) \cup \partial V$. Due to the definition of boundary components in Section \ref{rig-sec: sub, modification}  (see \eqref{rig-eq: new1XXXX}) we observe that for each $\Gamma_j(\tilde{W}) \neq \partial V$ there is a (unique) corresponding $\Gamma_j(W)$ (for notational convenience we use the same index) such that $\Theta_j (\tilde{W}) =  \Theta_j(W) \setminus \overline{V}$. Then $\tilde{W}$ clearly satisfies \eqref{rig-eq: W prop}(i),(iv), where \eqref{rig-eq: W prop}(iv)  follows from \eqref{rig-eq: weights} and for \eqref{rig-eq: W prop}(i) we use the fact that $R_j(\tilde{W}) = R_j(W)$, i.e we take the same rectangles as for the boundary components of $W$ (cf. below \eqref{rig-eq: new1XXXX}). To see \eqref{rig-eq: W prop}(ii), it suffices to note that for a given $\Gamma_j(\tilde{W})$ with $\omega(\Gamma_j(\tilde{W}))<1$, \eqref{rig-eq: weights} implies $\omega(\Gamma_j(\tilde{W})) \vert \Gamma_j(\tilde{W}) \vert_\infty = \omega(\Gamma_j(W)) \vert \Gamma_j(W) \vert_\infty$. 

Possibly  \eqref{rig-eq: W prop}(iii) or \eqref{rig-eq: W prop}(v) are violated:  Let ${\cal F}_0^1$ be the set of components  $\Gamma_{j}(\tilde{W})$ with $\omega(\Gamma_{j}(\tilde{W})) =1$  such that $\Gamma_{j}(\tilde{W})$ is not rectangular or $\Gamma_{j}(\tilde{W}) \neq \Theta_{j}(\tilde{W})$  and let ${\cal F}_0^2$ be the set of components $\Gamma_{j}(\tilde{W})$ such that for the corresponding rectangles $R_{j}$  given by \eqref{rig-eq: W prop}(i) one has that ${R}_{j} \setminus X$ is disconnected for a suitable component $X$ of $\tilde{W}$.  As   \eqref{rig-eq: W prop}(iii),(v) are satisfied for $W$, we have $\Gamma_j(\tilde{W}) \cap \partial V \neq \emptyset$ for $\Gamma_j(\tilde{W}) \in {\cal F}^1_0$ and $R_j \setminus V$ is disconnected for $\Gamma_j(\tilde{W}) \in {\cal F}^2_0$.   The goal is to modify $\tilde{W}$ iteratively to obtain a set satisfying  \eqref{rig-eq: W prop}(iii) and \eqref{rig-eq: W prop}(v). 

Set $W_0 = \tilde{W}$  and $V_0 = V$. Assume $W_i = (\tilde{W} \setminus V_i) \cup \partial V_i$ has been constructed, where $V_i \in {\cal U}^s$ is a rectangle with  $\overline{V} \subset \overline{V_i} \subset Q_\mu$, and assume that the boundary  component  $\partial V_i = \Gamma_0(W_i)$ satisfies $\omega(\Gamma_0(W_i))=1$. Moreover, suppose  that  \eqref{rig-eq: W prop}(i),(ii),(iv) hold for $W_i$  and that
there is a subset of components ${\cal G}_i \subset (\Gamma_j(\tilde{W}))_j$ such that 
\begin{align}\label{rig-eq: energy balanceXXXX}
(i)& \ \ \Gamma_j(\tilde{W}) \subset \overline{V_i} \ \ \text{for all} \ \ \Gamma_j(\tilde{W}) \in {\cal G}_i,\\
(ii)& \ \ (1-h_*){\hspace{-0.2cm}}\sum_{\Gamma_j(\tilde{W}) \in {\cal G}_i} \omega(\Gamma_j(\tilde{W}))|\Gamma_j(\tilde{W})|_\infty + h_*|\partial_{\cal V} \tilde{W} \cap (\overline{V_i} \setminus \overline{V})|_{\cal H} + |\partial V|_*\ge |\partial V_i|_*,\notag
\end{align}
where for shorthand $\partial_{\cal V} \tilde{W} = \bigcup_j \Theta_j(\tilde{W})$. As before by ${\cal F}^1_i$ we denote the set of components $\Gamma_j(W_i)$ with $\omega(\Gamma_j(W_i)) = 1$  which are not rectangular or satisfy $\Gamma_j(W_i) \neq \Theta_j(W_i)$ and ${\cal F}^2_i$ denotes the set of boundary components  for which the corresponding rectangle is disconnected. We suppose that $\Gamma_j(W_i) \cap  \partial V_i \neq \emptyset$  for $\Gamma_j(W_i) \in {\cal F}^1_i$ and $R_j \setminus V_i$ is disconnected for  $\Gamma_j(W_i) \in {\cal F}^2_i$.

\smallskip

Observe that $|\partial V_i|_\infty \ge \lambda$  as $|\partial V|_\infty \ge \lambda$.  If now $W_i \in {\cal W}^s_\lambda$ (i.e ${\cal F}^1_i = {\cal F}^2_i = \emptyset$), we stop and set $U=W_i$. Otherwise, we choose $\Gamma_k(W_i) \in {\cal F}_i$  and let  $X_k$ be the component of $Q_\mu \setminus W_i$ corresponding to $\Gamma_k(W_i)$.  Note that by \eqref{rig-eq: new1XXXX} there is a component $\Gamma_k(\tilde{W})$ such that $\Gamma_k(W_i)= \Gamma_k(\tilde{W}) \setminus \overline{V_i}$ (we again use the same index for convenience). If $\Gamma_k(W_i) \in {\cal F}^1_i$, we let $V_{i+1} \in {\cal U}^s$ be the smallest rectangle whose closure contains $V_i$ and $\Gamma_k(W_i)$. By Lemma  \ref{rig-lemma: infty}(i) applied on $\Gamma = \partial (V_{i} \cup X_k)$  we get $|\partial V_{i+1}|_* \le |\partial (V_{i} \cup X_k)|_*$. Moreover, \eqref{rig-eq: h*} and Lemma  \ref{rig-lemma: inftyX}(iii)   yield
\begin{align*}
 |\partial (V_{i} \cup X_k)|_* &\le  (1-h_*)(|\partial V_i|_\infty + |\Gamma_k(W_i)|_\infty)  + h_*(|\partial V_i|_{\cal H} + |\Gamma_k(W_i) \setminus \partial  V_i|_{\cal H}) \\
&=  |\partial V_i|_* + (1-h_*) \omega( \Gamma_k(W_i))| \Gamma_k(W_i)|_\infty + h_* | \Gamma_k(W_i)  \setminus \partial V_i|_{\cal H}, 
\end{align*}
where we used  $ \omega( \Gamma_k(W_i)) = 1$. In view of \eqref{rig-eq: new1} we have $\Gamma_k(W_i) \setminus  \partial V_i = \Gamma_k(\tilde{W})  \cap (\overline{V_{i+1}} \setminus \overline{V_i})$. Then  by \eqref{rig-eq: weights2}  we get 
\begin{align}\label{eq:step1}
|\partial V_{i+1}|_* \le |\partial V_i|_* + (1-h_*)\omega( \Gamma_k(\tilde{W}))|\Gamma_k(\tilde{W})|_\infty +  h_*|\partial_{\cal V} \tilde{W} \cap (\overline{V_{i+1}} \setminus \overline{V_i})|_{\cal H}.
\end{align}
If $\Gamma_k(W_i) \in {\cal F}_i^2 \setminus {\cal F}_i^1$ we let $V_{i+1} \in {\cal U}^s$ be the smallest rectangle whose closure contains $V_i$ and $\partial R_k$, where $R_k$ is the rectangle given by \eqref{rig-eq: W prop}(i) associated to both sets $ \Gamma_k(\tilde{W})$ and $ \Gamma_k(W_i)$. As \eqref{rig-eq: W prop}(ii),(iv) hold for $W_i$, we get  $|\partial R_k|_\infty \le 2| \Gamma_k(W_i)|_\infty \le 38\upsilon |\partial V_i|_\infty$ using $\omega_{\rm min} \ge \frac{1}{2}$ and $|\partial V_i|_\infty \ge \lambda$. Consequently, due to the fact that $R_k \setminus V_i$ is disconnected by assumption   Lemma \ref{rig-lemma: infty}(i),(ii) then yields for $\upsilon$ sufficiently small by \eqref{rig-eq: W prop}(ii) and $\omega_{\rm min} \ge \frac{1}{2}$
\begin{align}\label{eq:step2}
\begin{split}
|\partial V_{i+1}|_* &\le    |\partial (V_{i} \cup R_k)|_* \le |\partial V_{i}|_* + \tfrac{1}{2}(1-h_*)|\partial R_k|_\infty \\&\le |\partial V_{i}|_* + (1-h_*)\omega(\Gamma_k(W_i))|\Gamma_k(W_i)|_\infty   \\&  \le |\partial V_{i}|_* + (1-h_*)\omega(\Gamma_k(\tilde{W}))|\Gamma_k(\tilde{W})|_\infty,
\end{split} 
\end{align}
where in the last step we used \eqref{rig-eq: weights2}. We now define the set $W_{i+1} = (\tilde{W} \setminus V_{i+1}) \cup \partial V_{i+1}$ (recall \eqref{rig-eq: new1}) and adjust the  weights of the boundary components of $W_{i+1}$ as in \eqref{rig-eq: weights}. Moreover, we let ${\cal G}_{i+1} = {\cal G}_i \cup \lbrace \Gamma_k(\tilde{W}) \rbrace$.  In view of   $\overline{V_i} \subset \overline{V_{i+1}}$, $\Gamma_k(\tilde{W}) \subset \overline{V_{i+1}}$ and \eqref{eq:step1}-\eqref{eq:step2} we  observe that \eqref{rig-eq: energy balanceXXXX} holds.

Repeating the arguments at the beginning of the proof for $W_{i+1}$ in place of $\tilde{W}$ we see that  $W_{i+1}$ satisfies  \eqref{rig-eq: W prop}(i),(ii),(iv). Moreover, there are sets ${\cal F}^1_{i+1}$ and ${\cal F}^2_{i+1}$ defined as below \eqref{rig-eq: energy balanceXXXX} such that $\Gamma_j(W_{i+1}) \cap \partial V_{i+1} \neq \emptyset$ for $\Gamma_j(W_{i+1}) \in {\cal F}^1_{i+1}$ and $R_j \setminus V_{i+1}$ is disconnected for  $\Gamma_j(W_{i+1}) \in {\cal F}^2_{i+1}$. We continue with the next iteration step.

After a finite number of steps $i^*$ we find a rectangle $V_{i^*}$ with  $\overline{V} \subset \overline{V_{i^*}} \subset Q_\mu$ and a set $W_{i^*} = (\tilde{W} \setminus V_{i^*}) \cup \partial V_{i^*} \in {\cal W}^s_\lambda$ as in each step the number of boundary components decreases. Letting ${\cal G}' = {\cal G}_{i^*} \cup \lbrace \partial V \rbrace$ we derive by  \eqref{rig-eq: h*, omega}, \eqref{rig-eq: energy balanceXXXX}(ii) and the fact that $\omega(\partial V) = 1$
\begin{align*}
|\partial V_{i^*}|_* &\le (1-h_*)\sum_{\Gamma_j(\tilde{W}) \in {\cal G}'} \omega(\Gamma_j(\tilde{W}))|\Gamma_j(\tilde{W})|_\infty + h_*|\partial V|_{\cal H} + h_*|\partial_{\cal V} \tilde{W} \cap (\overline{V_{i^*}} \setminus \overline{V})|_{\cal H}  \\ &
\le (1-h_*)\sum\nolimits_{\Gamma_j(\tilde{W}) \in {\cal G}'} \omega(\Gamma_j(\tilde{W}))|\Gamma_j(\tilde{W})|_\infty  +   h_*|\partial_{\cal V} \tilde{W} \cap \overline{V_{i^*}}|_{\cal H}.
\end{align*} 
We apply   Lemma \ref{rig-lemma: derive prop3}(i) for $W_{i^*} = (\tilde{W} \setminus V_{i^*}) \cup \partial V_{i^*}$. Recalling that ${\cal G}' \subset {\cal L}_{V_{i^*}}$    by  \eqref{rig-eq: energy balanceXXXX}(i) we get $
\Vert W_{i^*}\Vert_\omega  \le  \Vert \tilde{W}\Vert_\omega$. Define $V' = V_{i^*}$ and  $U = W_{i^*}$ as an element of ${\cal W}_\lambda^s$, i.e. with the same weights and the same order for the components.  Then we clearly have \eqref{rig-eq: energy balance2} and $H(U) = H(W) \cup \overline{V'}$. 

 It remains to show that $U$ coincides with $(W \setminus V') \cup \partial V'$ in the sense of \eqref{rig-eq: new1}, i.e. $X_j(U) =X_j(W) \setminus \overline{V'}$ for all $j$, and that the weights are adjusted as in \eqref{rig-eq: weights}. First, we see that $\omega(\Gamma_0(U)) = 1$ with $\Gamma_0(U) = \partial V'$ as required. Fix another component $\Gamma_j(U) =\partial X_j(U)$ and let $\Gamma_j(W) = \partial X_j(W)$, $\Gamma_j(\tilde{W}) = \partial X_j(\tilde{W})$ be the corresponding components. We observe
$$ X_j(\tilde{W}) = X_j(W)  \setminus \overline{V}, \ \ \ \ X_j(U) =  X_j(\tilde{W}) \setminus \overline{V'}.$$
This gives  $X_j(U) =X_j(W) \setminus \overline{V'}$ and establishes \eqref{rig-eq: new1}. To see   \eqref{rig-eq: weights}, we use that for the sets $\tilde{W} = (W \setminus V) \cup \partial V$, $U = (\tilde{W} \setminus V') \cup \partial V'$ the weights have been adjusted as in \eqref{rig-eq: weights} and argue as follows:

If $\Gamma_j(U) \cap \partial V' = \emptyset$, the component has not been changed during the modification process and thus in view of  \eqref{rig-eq: weights} we have  $\omega(\Gamma_j(U)) = \omega(\Gamma_j(\tilde{W}))$ and $\omega(\Gamma_j(\tilde{W})) = \omega(\Gamma_j(W))$. This gives   $\omega(\Gamma_j(U)) = \omega(\Gamma_j(W))$ and confirms \eqref{rig-eq: weights}.

Otherwise, if  $\Gamma_j(U) \cap \partial V' \neq  \emptyset$,  \eqref{rig-eq: W prop}(v) implies that $\omega(\Gamma_j(U)) < 1$ since components with weight $1$ are pairwise disjoint. Then again by \eqref{rig-eq: weights}  we have $\omega(\Gamma_j(U))|\Gamma_j(U)|_\infty = \omega(\Gamma_j(\tilde{W}))|\Gamma_j(\tilde{W})|_\infty$. Due to \eqref{rig-eq: weights2} we have $\omega(\Gamma_j(\tilde{W})) \le \omega(\Gamma_j(U)) < 1$  and thus using \eqref{rig-eq: weights} once more we get $\omega(\Gamma_j(\tilde{W}))|\Gamma_j(\tilde{W})|_\infty = \omega(\Gamma_j(W))|\Gamma_j(W)|_\infty$. This yields  $\omega(\Gamma_j(U))|\Gamma_j(U)|_\infty = \omega(\Gamma_j({W}))|\Gamma_j({W})|_\infty$ and confirms \eqref{rig-eq: weights} in this case, as desired. \eop

\smallskip

We observe that Remark \ref{rem: modi} follows from the construction above since for each $\Gamma(W)$ with $\omega(\Gamma(W)) = 1$ and $\Gamma(W) \cap \overline{V'} \neq \emptyset$, $\Gamma(W) \not\subset \overline{V'}$ the corresponding component $\Gamma(U)$ (satisfying $\omega(\Gamma(U))=1$) would have nonempty intersection with the component $\partial V'$ which contradicts \eqref{rig-eq: W prop}(v).  An adaption of the above arguments leads to the proof of Corollary \ref{rig-cor: derive prop2}.

\noindent {\em Proof of Corollary \ref{rig-cor: derive prop2}.}   Set $W_0 = W$ and assume $W_i \subset W$ with connected boundary components  has been constructed with $\Vert W_i \Vert_* \le \Vert W \Vert_*$. If $W_i \in {\cal W}^s_0$, we stop, otherwise we find a component $\Gamma = \Gamma(W_i) \in (\Gamma_j(W_i))_{j=1}^n$ which is not rectangular  or satisfies $\Gamma \neq \Theta$.  Let ${\cal G}_i$ be the components $\Gamma_j(W_i)$ having nonempty intersection with  $\Gamma$ and set ${\cal G}' = \lbrace \Gamma\rbrace \cup {\cal G}_i$. Let $W'_{i+1} = (W_i \setminus V) \cup \partial V$, where $\overline{V} \subset Q_\mu$ is the smallest closed rectangle  which contains the components in ${\cal G}'$.  By \eqref{rig-eq: h*} and arguing as in  Lemma \ref{rig-lemma: inftyX}(iii), Lemma \ref{rig-lemma: infty}(i) we get 
\begin{align*}
|\partial V|_* &\le (1-h_*) \sum\nolimits_{\Gamma_j(W_i) \in {\cal G}'} |\Gamma_j(W_i)|_\infty + h_*\big|\bigcup\nolimits_{\Gamma_j(W_i) \in {\cal G}'} \Gamma_j(W_i)\big|\\
&\le (1-h_*) \sum\nolimits_{\Gamma_j(W_i) \in {\cal G}'} |\Gamma_j(W_i)|_\infty + h_*\big|\bigcup\nolimits_{j=1}^n \Theta_j(W_i) \cap \overline{V}\big|.
\end{align*}
Noting that all sets in ${\cal G}'$ are contained in $\overline{V}$ and repeating the calculation in  Lemma \ref{rig-lemma: derive prop3}(i)  (the assumption $W_i \in {\cal W}^s_\lambda$ is not needed as no weights are introduced) we then get $\Vert W'_{i+1} \Vert_* \le \Vert W_i \Vert_* \le \Vert W \Vert_*$. By combining all boundary components of $W'_{i+1}$ having nonempty intersection with $\partial V$ we obtain a set $W_{i+1}$ with connected boundary components with $W_{i+1} \subset W_{i+1}'$, $|W'_{i+1} \setminus W_{i+1}|=0$ and $\Vert W_{i+1} \Vert_* \le \Vert W'_{i+1} \Vert_*$ arguing as in Lemma \ref{rig-lemma: inftyX}. We now continue with iteration step $i+1$ and note that we find the desired set $U$ after a finite number of iterations. We observe that $H(U) \supset H(W)$ and as $W\setminus U \subset \bigcup_{l=1}^n X_l(U)$  we conclude $|W \setminus U| \le \sum_{l=1}^n |X_l(U)| \le  \sum_{l=1}^n |\Gamma_l(U)|_\infty^2 \le \Vert U\Vert_\infty^2$.   \eop

\subsection{Neighborhoods}\label{rig-sec: proofs2}

We now prove the properties of rectangular neighborhoods announced in Section \ref{rig-sec: subsub,  rec-neigh}. Let $W \in {\cal W}^s_\lambda$ for $\lambda \ge 0$ be given and a component  $\Gamma  = \partial X = \Gamma(W)$ with $\omega(\Gamma)=1$,   $|\Gamma|_\infty \ge \lambda$  and $X= (-l_1,l_1) \times (-l_2,l_2)$. Recall the definition of the neighborhood before \eqref{rig-eq: tau def}.  Let $(\Gamma_j)_j = (\Gamma_j(W))_j$ be the boundary components of $W$ with the corresponding subsets $(\Theta_j)_j$ and rectangles $(R_j)_j$ as given in \eqref{rig-eq: Xdef} and \eqref{rig-eq: W prop}, respectively. We will always add a subscript to avoid a mix up with $\Gamma$. We begin with the proof of Lemma \ref{rig-lemma: crack bound}.

\noindent {\em Proof of Lemma \ref{rig-lemma: crack bound}.}  Let $t>0$ be given and define   $V = (-l_1 - t , l_1 + t) \times  (-l_2 - t , l_2 + t) \in {\cal U}^s$ with $\overline{V} \subset Q_\mu$.  Let $\tilde{W} = (W \setminus V) \cup \partial V$ and adjust the weights as in \eqref{rig-eq: weights}. It is not hard to see that  $|\partial V|_* \le |\Gamma|_* + 8t$. We apply Lemma \ref{rig-lemma: derive prop3}(i) and derive in view of $\Gamma  \in {\cal L}_V$  and $\omega(\Gamma)=1$
\begin{align*}
\begin{split}
\Vert \tilde{W}\Vert_\omega & \le  \Vert W \Vert_\omega +   |\partial V|_* - (1-h_*)|\Gamma|_\infty - h_*  \big|\bigcup\nolimits_j \Theta_j(W) \cap \overline{V}\big|_{\cal H} \\ 
&  \le  \Vert W \Vert_\omega +   |\partial V|_* - (1-h_*)|\Gamma|_\infty - h_*|\Gamma|_{\cal H}   - h_*|\partial W \cap N^t|_{\cal H} \\&\le  \Vert W \Vert_\omega +   8t   - h_* |\partial W \cap N^t|_{\cal H},
\end{split}
\end{align*} 
 where in the second step we used that $N^t \cap \Gamma = \emptyset$ and $\partial W \cap N^t = \bigcup\nolimits_j \Theta_j(W) \cap N^t$ since $\overline{N^t} \subset H(W)$. Since $\Vert W \Vert_\omega  \le \Vert \tilde{W} \Vert_\omega$ by condition \eqref{rig-eq: property4}, we find $\vert \partial W \cap N^{t} \vert_{\cal H} \le \frac{8t}{h_*}$. \eop

Before we proceed with the proof of  Lemma \ref{rig-lemma: two components}, we state an adaption of Lemma \ref{rig-lemma: derive prop3}.

\begin{lemma}\label{rig-lemma: slicing bound4}
Let $c'>0$. Then there is a constant $C=C(c')>0$ such that for $h_*$, $1 - \omega_{\rm min}$ as in \eqref{eq: parameters} small enough depending on $c'$ the following holds: For all $\lambda\ge 0$,  $W \in {\cal W}^s_\lambda$ and boundary components $\Gamma$ satisfying  $\omega(\Gamma)=1$, $|\Gamma|_\infty \ge \lambda$ and  \eqref{rig-eq: property4} we get that for all $t \ge c' \bar{\tau}$ and for all subsets ${\cal R} \subset \lbrace R_j: R_j \cap N^t \neq \emptyset \rbrace$
\begin{align*}
\Vert \tilde{W} \Vert_{\omega}   &\le \Vert W \Vert_{\omega} + (1-h_*)(|\partial V|_\infty  - |\Gamma|_\infty - \sum\nolimits_{R_j \in {\cal R}} |\partial R_j|_\infty ) + Cth_*,
\end{align*}
where $\tilde{W} = (W \setminus V) \cup \partial V$ with $V$ being the smallest rectangle  whose closure contains $\Gamma$ and the rectangles ${\cal R}$. The weights have been adjusted as in \eqref{rig-eq: weights}. 
\end{lemma}

\Proof  First, we control the difference of $|\partial R_j|_* $ and $|\Theta_j|_\omega$ for $R_j \in {\cal R}$. Let ${\cal R}' = \lbrace R_j \in {\cal R}: |\partial R_j|_\infty \le 38\bar{\tau} \rbrace$.  By \eqref{rig-eq: W prop}(ii),(iv),(v),  \eqref{rig-eq: new2}(i) and the fact that $\bar{\tau} \ge \upsilon \lambda$ we have  $\vert \Theta_j\vert_\omega = \vert \partial R_j\vert_{*}$ if $ R_j \notin {\cal R}'$ and for $R_j \in {\cal R}'$
\begin{align}\label{eq: inter*}
\begin{split}
\vert \Theta_j\vert_\omega &= (1-h_*)\omega(\Gamma_j) \vert \Gamma_j\vert_\infty + h_*|\Theta_j|_{\cal H} \ge    (1-h_*)\omega_{\rm min} \vert \partial R_j\vert_{\infty} \\&\ge  (1-3h_*) \omega_{\rm min}|\partial R_j|_* \ge  (1-4h_*)  |\partial R_j|_*.
\end{split}
\end{align}
Here the third step holds by Lemma \ref{rig-lemma: inftyX}(iv) for $h_*$ small enough and in the last step we have chosen $1 - \omega_{\rm min}$ small with respect to $h_*$ (cf. \eqref{eq: parameters}). Let $V' \in {\cal U}^s$ be the smallest rectangle whose closure contains $\Gamma$ and the rectangles ${\cal R}'$. By \eqref{rig-eq: W prop}(i) we get $\overline{V'} \subset Q_\mu$ and a short computation yields $|\partial V'|_* \le |\Gamma|_* + 8(t + 38\bar{\tau}) \le |\Gamma|_* + Ct$ for $C=C(c')$ large enough. Defining $W' = (W \setminus V') \cup\partial V'$ we then get by \eqref{rig-eq: property4} and Lemma \ref{rig-lemma: derive prop3}(ii)  
\begin{align}\label{eq: inter}
0 \le \Vert W' \Vert_\omega - \Vert W \Vert_\omega \le Ct - \sum\nolimits_{R_j \in {\cal R}'} |\Theta_j|_\omega,
\end{align}
where we used that $\lbrace \Gamma \rbrace \cup \lbrace \Gamma_j: R_j \in {\cal R}' \rbrace \subset {\cal L}_{V'}$. For notational convenience we define $\Vert W \Vert_{\omega,{\cal R}}  = \Vert W \Vert_\omega + \sum_{R_j \in {\cal R}} (|\partial R_j|_* - |\Theta_j|_\omega) $.  We obtain by \eqref{eq: inter*} and  \eqref{eq: inter}
\begin{align}\label{rig-eq: ess-three-compYX}
\begin{split}
\Vert W \Vert_{\omega,{\cal R}} &\le \Vert W \Vert_{\omega}  + \sum\nolimits_{R_j \in {\cal R}'} (  (1-4h_*)^{-1} - 1)|\Theta_j|_\omega \\ 
&\le \Vert W \Vert_{\omega} + Ct( (1-4h_*)^{-1} - 1)\le \Vert W \Vert_{\omega} + Cth_*
\end{split}
\end{align}
for $h_*$ small enough for a constant $C=C(c')$.  Now let $V \in {\cal U}^s$ be the smallest rectangle whose closure contains $\Gamma$ and ${\cal R}$. We obtain $\overline{V} \subset Q_\mu$ by \eqref{rig-eq: W prop}(i) and   $|\partial V|_{\cal H} \le |\Gamma|_{\cal H} + \sum_{R_j \in {\cal R}} |\partial R_j|_{\cal H}  + 8t$.  Then applying again Lemma \ref{rig-lemma: derive prop3}(ii) and noting that $\lbrace \Gamma \rbrace \cup \lbrace \Gamma_j: R_j \in {\cal R}\rbrace \subset {\cal L}_{V}$   we get 
$$\Vert \tilde{W} \Vert_\omega \le \Vert W \Vert_\omega + |\partial V|_* - |\Gamma|_* -  \sum\nolimits_{R_j \in {\cal R}} |\Theta_j|_\omega.$$
Using $\Vert W \Vert_{\omega,{\cal R}} = \Vert W \Vert_\omega + \sum_{R_j \in {\cal R}} (|\partial R_j|_* - |\Theta_j|_\omega)$  we then find by \eqref{rig-eq: h*}, \eqref{rig-eq: ess-three-compYX}
\begin{align*}
\Vert \tilde{W} \Vert_{\omega} - \Vert W \Vert_{\omega}  &\le \Vert \tilde{W} \Vert_{\omega} - \Vert W \Vert_{\omega,{\cal R}} + Cth_* \\
& \le |\partial V|_* - |\Gamma|_*  - \sum\nolimits_{R_j \in {\cal R}} |\partial R_j|_*  +C t h_* \notag \\& \le (1-h_*)(|\partial V|_\infty  - |\Gamma|_\infty - \sum\nolimits_{R_j \in {\cal R}} |\partial R_j|_\infty )  + 8t h_* + Cth_*.\notag
\end{align*}
 \eop

As a preparation for the proof of Lemma \ref{rig-lemma: two components} and Lemma \ref{rig-lemma: slicing bound2} we formulate a lemma which shows that $\Vert \cdot \Vert_\pi$ can be controlled in a suitable way. Recall the construction of the covering ${\cal C}^t$ in \eqref{rig-eq: covering def} consisting of sets in ${\cal Y}^t$ and the definition of  $\Vert \cdot \Vert_\pi$ in \eqref{rig-eq: VertP}. Let $\frac{1}{800} \bar{\tau} \le t \le 20\cdot 22\bar{\tau}$ be given with $\bar{\tau}$ as in \eqref{rig-eq: tau def} and assume  $\overline{N^{t}}\subset H(W)$. We note that the  construction in  \eqref{rig-eq: VertP}  implies  for $\upsilon$ sufficiently small (with respect to $h_*$)
\begin{align}\label{rig-eq: covering def2}
R_j \cap Y^t_i \neq \emptyset \ \ \ \Rightarrow \ \ \ R_j \cap Y^t_{i+l} = \emptyset \ \ \ \text{for } |l| \ge 3
\end{align}
for all rectangles $R_j$ given by \eqref{rig-eq: W prop} and $i=0,\ldots,m-1$. To see this, we first observe that 
\begin{align}\label{rig-eq: maxcracklength}
|\partial R_j \cap N^t|_{\cal H} \le Ct h_*^{-1}  \le  C \upsilon h_*^{-1} |\Gamma|_\infty
\end{align}
 for $C>0$ large enough. Indeed, if $|\Gamma_j|_\infty < 19\bar{\tau}$, we  obtain $|\partial R_j|_\infty < 38\bar{\tau}$ by \eqref{rig-eq: new2}(i) and thus  $|\partial R_j|_{\cal H} \le 2\sqrt{2} |\partial R_j|_\infty  \le Ct$. Otherwise,  recalling $|\Gamma|_\infty \ge \lambda$, by \eqref{rig-eq: W prop}(iv),(v) we have $\partial R_j = \Gamma_j$,  $\omega(\Gamma_j)=1$  and thus  employing  Lemma \ref{rig-lemma: crack bound} we get $|\partial R_j \cap N^t|_{\cal H} \le 8t h^{-1}_*$.

Let $\bar{C}$ as in \eqref{rig-eq: covering def}. If now $\dist(Y^t_i,Y^t_{i+l}) \ge \bar{C} |\Gamma|_\infty$ for some $|l|\ge 3$, \eqref{rig-eq: covering def2} follows as by \eqref{rig-eq: maxcracklength} we get that $|\partial R_j \cap N^t|_{\cal H}$ is small with respect to $\bar{C}|\Gamma|_\infty$ for $\upsilon$ small enough  (depending on $h_*, \bar{C}$). 

On the other hand, suppose $\dist(Y^t_i,Y^t_{i+l}) \le \bar{C} |\Gamma|_\infty$. This is only possible  in the case $l_2 \le \frac{l_1}{2}$ if (up to  interchanging $+$ and $-$) $Y^t_i \subset N^t_{2,+} \setminus (N^t_{1,-} \cup N^t_{1,+})$, $Y^t_{i+l} \subset N^t_{2,-}$  and $\dist(Y^t_{i},N^t_{1,\pm})\ge c|\Gamma|_\infty$ or $\dist(Y^t_{i+l},N^t_{1,\pm})\ge c|\Gamma|_\infty$ for a universal $c>0$.  Now assume that \eqref{rig-eq: covering def2} was wrong. Then by \eqref{rig-eq: covering def} and the fact that $|\partial R_j \cap N^t|_{\cal H}$ is small with respect to $\bar{C}|\Gamma|_\infty$   (for $\upsilon$ small) this would imply  $R_j \cap N^t_{2,\pm} \neq \emptyset$ and $R_j \cap (N^t_{1,-} \cup N^t_{1,+}) = \emptyset$. But then we would get that $R_j \setminus X$  is not connected, where as before $X$ is the component with $\partial X = \Gamma$. This, however, contradicts \eqref{rig-eq: W prop}(iii).

After this preparation we now prove the following lemma yielding a bound on $\Vert \cdot \Vert_\pi$. Recall the definition of the enlarged sets in \eqref{eq: label large}.
\begin{lemma}\label{rig-lemma: slicing bound}
For $h_*$, $\upsilon$ and  $1 - \omega_{\rm min}$ as in \eqref{eq: parameters} small enough the following holds: For all $\lambda\ge 0$, $W \in {\cal W}^s_\lambda$ and boundary components $\Gamma$ satisfying $\omega(\Gamma)=1$, $|\Gamma|_\infty \ge \lambda$ and  \eqref{rig-eq: property4}  we get  for all $\frac{1}{800} \bar{\tau}\le t \le 20\cdot 22\bar{\tau}$ with  $\overline{N^{t}}\subset H(W)$ that there are two sets  $Y^1, Y^2 \in {\cal C}^{t}$  such that $ \Vert Y^{ t} \Vert_\pi \le \frac{19}{20}t$ for all $Y^{t} \in {\cal C}^{ t}$ with $Y^{t} \cap (\bar{Y}^1 \cup \bar{Y}^2) = \emptyset$. \\
Additionally,  if $\Vert Y^1 \Vert_\pi, \Vert Y^2 \Vert_\pi \ge \frac{19}{20}t$, then  $\bar{Y}^1\cup\bar{Y}^2$ intersects  both $N^{t}_{1,+}$ and $N^{t}_{1,-}$ or both $N^{t}_{2,+}$ and $N^{t}_{2,-}$. If $l_2 \le \frac{l_1}{2}$, then $\bar{Y}^1\cup\bar{Y}^2$ intersects $N^{t}_{1,+}$ and $N^{t}_{1,-}$.
\end{lemma}

 We briefly remark that by similar arguments the additional statement can also be proved without the extra assumption $\Vert Y^1 \Vert_\pi, \Vert Y^2 \Vert_\pi \ge \frac{19}{20}t$. We omit the proof of this fact here as we will not need it in the following.

\Proof  For convenience we drop the superscript $t$ in the following proof.   We proceed in two steps: 

In a) we first show that it is not possible that there are three sets $Y^1,Y^2,Y^3  \in {\cal C}$ such that $\bar{Y}^k \cap Y^l = \emptyset$ if $k \neq l$ and $\Vert Y^k\Vert_\pi > \frac{19}{20} t$ for $k,l=1,2,3$. Provided that a) is proven we can then select the two desired  sets $Y^1,Y^2$ as follows:

(1) If $\Vert Y\Vert_\pi \le \frac{19}{20} t$ for all $Y \in {\cal C}$, we can choose arbitrary sets $Y^1,Y^2$ satisfying the additional condition. Otherwise, we can assume that there is some $Y^*$ with $\Vert Y^*\Vert_\pi > \frac{19}{20} t$. 

(2) If $\Vert Y\Vert_\pi \le \frac{19}{20} t$ for all $Y \in {\cal C}^t$ with $Y \cap \bar{Y}^* = \emptyset$, we set $Y^1 = Y^*$ and choose $Y^2$ arbitrarily such that the additional condition holds. 

(3) Otherwise, we set $Y^1 = Y^*$ and choose $Y^2$ with $\Vert Y^2\Vert_\pi > \frac{19}{20} t$ and $Y^2 \cap \bar{Y}^* = \emptyset$. Now a) indeed shows that $\Vert Y\Vert_\pi \le \frac{19}{20} t$ for all $Y \in {\cal C}^t$ with $Y \cap (\bar{Y}^1 \cup \bar{Y}^2) = \emptyset$.

In step b) we concern ourselves with the additional assertions on the position of $\bar{Y}^1  \cup \bar{Y}^2$  in case (3). 

\bigskip

a) Suppose that there are three sets $Y^1,Y^2,Y^3  \in {\cal C}$ such that $\bar{Y}^k \cap Y^l = \emptyset$ if $k \neq l$ and $\Vert Y^k\Vert_\pi > \frac{19}{20} t$ for $k,l=1,2,3$. First note that the assumption implies that if e.g. $Y^1 = Y_i$, then $Y^2,Y^3 \notin \lbrace Y_{i-2},\ldots,Y_{i+2}\rbrace$. Let $V$ be the smallest rectangle whose closure contains  $\Gamma$ and the sets ${\cal R} := \bigcup^3_{k=1} {\cal R}(Y^k)$  with ${\cal R}(Y^k)$ as defined before \eqref{rig-eq: VertP}. Define $\tilde{W} = (W \setminus V) \cup \partial V$ (recall \eqref{rig-eq: new1}). We will show that 
\begin{align}\label{rig-eq: ess-three-comp}
|\partial V|_\infty   \le |\Gamma|_\infty + \sum\nolimits_{R_j \in {\cal R}} |\partial R_j|_\infty  - \tfrac{1}{50}t.
\end{align}
 Then applying Lemma \ref{rig-lemma: slicing bound4}  we get 
 $\Vert \tilde{W} \Vert_{\omega} - \Vert W \Vert_{\omega} \le -(1-h_*)\frac{1}{50}t + Ch_* t< 0$
 for $h_*$ small enough. This then gives a contradiction to \eqref{rig-eq: property4} and  concludes the proof of a).

We now proceed to show \eqref{rig-eq: ess-three-comp}.  (We refer to Figure \ref{E} where the essential idea is illustrated.) Assume $V= (-a_{1,-}-l_1 ,l_1 + a_{1,+}) \times (-a_{2,-} - l_2,l_2 + a_{2,+})$ and select (not necessarily pairwise different) $R_{k,\pm} \in {\cal R}$ such that $\pm (l_k + a_{k,\pm}) \in \pi_k  \partial R_{k,\pm}$ for $k=1,2$. (If $a_{k,\pm} = 0$ then $R_{k,\pm} = \emptyset$.) We find by \eqref{rig-eq: VertP}
\begin{align}\label{rig-eq: pi est}
|\pi_k R_{k,\pm}|  \ge a_{k,\pm} - \dist(\pi_k R_{k,\pm}, \pi_k \Gamma) \ge a_{k,\pm} - t + |\partial R_{k,\pm}|_\pi.
\end{align}
We suppose for the moment that  $R_{k,+} \neq R_{k,-}$ for $k=1,2$.  (In particular, this implies that three rectangles never coincide.) At the end of the proof we will briefly indicate how the following arguments can be adapted to the general case. We first assume that two rectangles coincide, e.g. $R = R_{1,-} = R_{2,-}$.  By \eqref{rig-eq: pi est} and an elementary computation we obtain 
$$\sqrt{a_{1,-}^2 + a_{2,-}^2} \le  |\partial R|_\infty + \sqrt{2}(t - |\partial R|_\pi) \le  |\partial R|_\infty + \sqrt{2}t - |\partial R|_\pi.$$ 
Otherwise, if e.g. $R_{1,-} \neq R_{2,-}$, again applying \eqref{rig-eq: pi est} we get 
\begin{align*}
\sqrt{a_{1,-}^2 + a_{2,-}^2} &\le \sqrt{(t + (|\partial R_{1,-}|_\infty - |\partial R_{1,-}|_\pi))^2 + (t + (|\partial R_{2,-}|_\infty - |\partial R_{2,-}|_\pi))^2} \\
&\le   \sqrt{2}t + |\partial R_{1,-}|_\infty - |\partial R_{1,-}|_\pi + |\partial R_{2,-}|_\infty - |\partial R_{2,-}|_\pi.
\end{align*}
Consequently, we obtain 
\begin{align}\label{rig-eq: Fest}
F \le 2\sqrt{2}t + \sum\nolimits_{k,\pm} (|\partial R_{k,\pm}|_\infty  - |\partial R_{k,\pm}|_\pi), 
\end{align}
where each rectangle is only counted once in the sum and 
 $$F = \sqrt{a_{1,-}^2 + a_{2,-}^2} + \sqrt{a_{1,+}^2 + a_{2,+}^2} \ \ \ \text{or} \ \ \ F = \sqrt{a_{1,-}^2 + a_{2,+}^2} + \sqrt{a_{1,+}^2 + a_{2,-}^2}.$$
Moreover, note that by assumption on the sets $Y^1,Y^2,Y^3$ and \eqref{rig-eq: covering def2} each rectangle $R_j$ intersects at most one of the three sets $Y^k$. Therefore, as $|\cdot|_\infty \ge |\cdot|_\pi$ we obtain
\begin{align}\label{eq-rig: new3}
\begin{split}
|\partial V|_\infty & \le |\Gamma|_\infty +F \le |\Gamma|_\infty + 2\sqrt{2}t + \sum\nolimits_{R_j \in {\cal R}} (|\partial R_j|_\infty - |\partial R_j|_\pi) \\
& = |\Gamma|_\infty + \sum\nolimits_{R_j \in {\cal R}} |\partial R_j|_\infty + 2\sqrt{2}t - \sum\nolimits^3_{k=1} \Vert Y^k\Vert_\pi.
\end{split}
\end{align}
As $\sum\nolimits_{k=1}^3 \Vert Y^k\Vert_\pi\ge 3\cdot\frac{19}{20}t$, this gives \eqref{rig-eq: ess-three-comp}.

\bigskip

b)  Suppose that $Y^1, Y^2 \in {\cal C}$  with $\Vert Y^1 \Vert_\pi, \Vert Y^2 \Vert_\pi \ge \frac{19}{20}t$  have been chosen according to  case (3) above. We show that $\bar{Y}^1 \cup \bar{Y}^2$ intersect both $N_{1,+}$ and $N_{1,-}$ or both $N_{2,+}$ and $N_{2,-}$ ($N_{1,\pm}$ if $l_2 \le \frac{l_1}{ 2}$). Let $V$ be the smallest rectangle whose closure contains $\Gamma$ and the sets ${\cal R} := \bigcup^2_{k=1} {\cal R}(Y^k)$. Set $\tilde{W} = (W \setminus V) \cup \partial V$. Observe that by assumption and \eqref{rig-eq: covering def2} each boundary component $R_j$ intersects at most one of the two sets $Y^k$.

(i) First we assume $\bar{Y}_1 \cup \bar{Y}_2$ intersects at most two adjacent parts of the neighborhood, e.g. $(\bar{Y}_1\cup \bar{Y}_2)\cap (N_{1,+} \cup N_{2,+}) =\emptyset$. This implies $V = (-a_1 - l_1,l_1) \times (-a_2 - l_2,l_2)$. Selecting (not necessarily different) $R_1,R_2$ such that $-l_k - a_k \in \pi_k  \partial R_k$ for $k=1,2 $ and proceeding as in \eqref{rig-eq: Fest} we obtain
$$\sqrt{a_1^2 + a_2^2} \le \sqrt{2}t  + \sum\nolimits_{k=1,2} (|\partial R_k|_\infty - |\partial R_k|_\pi) $$ 
and therefore 
\begin{align}\label{rig-eq: (i)}
\begin{split}
|\partial V|_\infty &\le |\Gamma|_\infty + \sum\nolimits_{R_j \in {\cal R}} |\partial R_j|_\infty + \sqrt{2}t - \sum\nolimits_{k = 1,2} \Vert Y^k\Vert_\pi  
\end{split}
\end{align}
Applying Lemma \ref{rig-lemma: slicing bound4} and using  $\sum\nolimits_k \Vert Y^k\Vert_\pi  \ge 2\cdot\frac{19}{20}t$ we again find $\Vert \tilde{W} \Vert_\omega -  \Vert W \Vert_{\omega}  <0$ for $h_*$ small enough which contradicts \eqref{rig-eq: property4}.

(ii) We finally show the additional statement that $\bar{Y}_1\cup \bar{Y}_2$ intersects $N_{1,\pm}$ in the case $l_2 \le \frac{l_1}{2}$. Assume without restriction that $(\bar{Y}_1 \cup \bar{Y}_2) \cap N_{1,+}(\Gamma) = \emptyset$. Then we have $V = (-a_1 - l_1,l_1) \times (-a_- - l_2,l_2 + a_+)$ and select rectangles $R_1, R_-,R_+$ as before.  (For the moment we assume $R_{k,+} \neq R_{k,-}$ for $k=1,2$.) Similarly as in a), we find
$$\sqrt{a_1^2 +a_-^2} \le \sqrt{2}t + \sum\nolimits_{k=1,-} (|\partial R_k |_\infty - |\partial R_k |_\pi), \ \ \ a_+ \le t + |\partial R_+ |_\infty - |\partial R_+ |_\pi.$$
Then
\begin{align*}
|\partial V|_\infty = \max\nolimits_{x \in [0,1]} \big( \sqrt{1-x^2}|\pi_1\Gamma| + x|\pi_2\Gamma| + \sqrt{1-x^2}a_1 +  x(a_+ + a_-)\big).
\end{align*}
We define $f(x) = \frac{2}{\sqrt{5}}(\sqrt{1-x^2} + \frac{1}{2}x)$ for $x \ge \frac{1}{\sqrt{5}}$ and $f(x) =1$ else. As $l_2 \le \frac{l_1}{2}$ an elementary argument yields $\sqrt{1-x^2}|\pi_1\Gamma| + x|\pi_2\Gamma| \le f(x)|\Gamma|_\infty$. Using $|\Gamma|_\infty \ge  (20\cdot 22\upsilon)^{-1} t$ we then obtain
\begin{align}\label{rig-eq: c split}
\begin{split}
|\partial V|_\infty & \le \max\nolimits_{x \in [0,1]} \big( f(x)|\Gamma|_\infty +  \sqrt{2}t + xt + \sum\nolimits_{k=1,\pm} (|\partial R_k |_\infty - |\partial R_k |_\pi) \big) \\
& \le|\Gamma|_\infty  + t\max\nolimits_{x \in [0,1]} r(x) + \sum\nolimits_{R_j \in {\cal R}} |\partial R_j|_\infty  -\sum\nolimits_{k=1,2} \Vert Y^k\Vert_\pi,
\end{split}
\end{align}
where  $r(x) = (f(x)-1) (20\cdot 22\upsilon)^{-1} + \sqrt{2}  +x$.  A computation yields $r(x) \le \frac{19}{10}-\frac{1}{100}$ for $x \le \sqrt{\frac{9}{40}}$ as $f\le 1$. Otherwise we have $\max_{[\sqrt{\frac{9}{40}},1]} ( f(x) - 1) < 0 $ and thus for $\upsilon$ sufficiently small  we also obtain $r(x) \le \frac{19}{10}-\frac{1}{100}$ for $x \in [ \sqrt{\frac{9}{40}},1]$.  Consequently, due to $\sum\nolimits_k \Vert Y^k\Vert_\pi \ge 2\cdot\frac{19}{20}t$ we have 
\begin{align}\label{rig-eq: c split2}
\begin{split}
|\partial V|_\infty &\le    |\Gamma|_\infty + \sum\nolimits_{R_j \in {\cal R}} |\partial R_j|_\infty - \tfrac{1}{100}t.
\end{split}
\end{align}
By Lemma \ref{rig-lemma: slicing bound4} we now again derive a contradiction to \eqref{rig-eq: property4} for $h_*$ small enough.

To finish the proof we briefly indicate how to proceed if e.g.  $R_{2,-} = R_{2,+}$. This may happen in the cases a) and b)(ii) above if $l_2 \ll l_1$. In this case we reduce the problem to the above treated situation by applying a translation argument: We replace $R_j$ by $R_j' := R_j - a_{2,+}\e_2$ for all $R_j \in {\cal R}$ as well as $V$ by $V' := V - a_{2,+}\e_2$.  Then we may set $R_{2,+} = \emptyset$ and can repeat the arguments above to derive  \eqref{rig-eq: ess-three-comp}, \eqref{rig-eq: (i)} and \eqref{rig-eq: c split2}, respectively, for $V'$ and $R'_{k,\pm}$. But then  \eqref{rig-eq: ess-three-comp}, \eqref{rig-eq: (i)} and \eqref{rig-eq: c split2}  also hold for the original sets $V$ and $R_j \in {\cal R}$ as $|V'|_* = |V|$ and $|R'_{k,\pm}|_* = |R_{k,\pm}|_*$. Consequently, we may then proceed as before and use Lemma \ref{rig-lemma: slicing bound4} to derive a contradiction to \eqref{rig-eq: property4}. \eop

We now can derive Lemma \ref{rig-lemma: two components} as a special case of the previous lemma.

\noindent {\em Proof of Lemma \ref{rig-lemma: two components}.}  Let  $\bar{\tau}\le t \le 22\bar{\tau}$ be given and assume that there are three components $\Gamma_k$, $k=1,2,3$, intersecting $N^{t}$ with $\vert \Gamma_k \vert_\infty \ge 19 t$.  By \eqref{rig-eq: W prop}(iv),(v), \eqref{rig-eq: tau def} and $t \ge \bar{\tau}$ we see that $\Gamma_k$ are rectangular with $\omega(\Gamma_k) =1$ for $k=1,2,3$. Set $\bar{t} =20 t$ and  recalling \eqref{rig-eq: VertP} we observe that $| \Gamma_k |_\pi \ge 19 t = \frac{19}{20}\bar{t}$.  (Here $|\cdot|_\pi$ has to be understood with respect to the neighborhood $N^{\bar{t} }$.) In view of $\bar{t} \le 20 \cdot 22\bar{\tau}$   we now may follow the lines of the proof of Lemma \ref{rig-lemma: slicing bound} for the neighborhood $N^{\bar{t}}$ with the essential difference that we replace the set of rectangles ${\cal R} = \bigcup^3_{k=1} {\cal R}(Y^k)$ (see beginning of step a))  by ${\cal R} = \lbrace \Gamma_1 \rbrace \cup \lbrace \Gamma_2 \rbrace \cup \lbrace \Gamma_3 \rbrace$ and in \eqref{eq-rig: new3} we replace $\sum^3_{k=1} \Vert Y^k\Vert_\pi$ by $\sum^3_{k=1} | \Gamma_k|_\pi$.  (Observe that in this case the assumption $\overline{N^{\bar{t}}} \subset H(W)$ can be dropped as it was only needed for \eqref{rig-eq: covering def2}.) Noting that $\sum^3_{k=1} | \Gamma_k|_\pi \ge 3\cdot\frac{19}{20}\bar{t}$ we again obtain a contradiction to \eqref{rig-eq: property4} and thus there are at most two large components $\Gamma_k$, $k=1,2$, in $N^{t}$. Likewise, we can proceed to determine the possible position of the two sets.

It remains to show that $|\pi_1 \Gamma_k| < \frac{1}{2}\vert \pi_2 \Gamma_k \vert$ leads to a contradiction if $l_2 \le \frac{l_1}{2}$. Let $V$ be the smallest rectangle whose closure contains  $\Gamma$, $\Gamma_k$ and derive similarly as in \eqref{rig-eq: c split}
\begin{align*}
|\partial V|_\infty &\le \max\nolimits_{x \in [0,1]} \big( \sqrt{1-x^2}|\pi_1\Gamma| + x|\pi_2\Gamma| + \sqrt{1-x^2}(t + |\pi_1\Gamma_k|) +  x(t + |\pi_2\Gamma_k|)\big)\\
&\le  \max\nolimits_{x \in [0,1]} \big(f(x)|\Gamma|_\infty + \sqrt{2}t + \sqrt{1-x^2}|\pi_1\Gamma_k| +  x|\pi_2\Gamma_k|)\big)
\end{align*} 
with $f$ as defined before \eqref{rig-eq: c split}, where we used $\sqrt{1-x^2}|\pi_1\Gamma| + x|\pi_2\Gamma| \le f(x)|\Gamma|_\infty$ due to the fact that $|\pi_2 \Gamma| \le \frac{1}{2}|\pi_1 \Gamma|$. Likewise, using the assumption $|\pi_1 \Gamma_k| < \frac{1}{2} |\pi_2 \Gamma_k|$ we find $\sqrt{1-x^2}|\pi_1\Gamma_k| + x|\pi_2\Gamma_k| \le f(\sqrt{1-x^2})|\Gamma_k|_\infty$ and thus obtain
\begin{align*}
|\partial V|_\infty \le |\Gamma|_\infty + |\Gamma_k|_\infty + \sqrt{2}t + \max_{x \in [0,1]} \big((f(x)-1) \upsilon^{-1}\tfrac{t}{22}  + (f(\sqrt{1-x^2})-1) 19 t\big),
\end{align*} 
where we used $|\Gamma|_\infty = \upsilon^{-1}\bar{\tau} \ge \upsilon^{-1}\frac{t}{22}$ and $|\Gamma_k|_\infty \ge 19 t$. Again separating the cases $x \ge \sqrt{\frac{9}{40}}$,  where $\max_{x\in [\sqrt{\frac{9}{40}},1]} ( f(x) - 1) < 0 $, and $x \le \sqrt{\frac{9}{40}}$,  where we find $(f(\sqrt{1-x^2})-1) 19t \le  -3 t$, we obtain for $\upsilon$ small enough $|\partial V|_\infty \le |\Gamma|_\infty + |\Gamma_k|_\infty -  t$. As $|\partial V|_{\cal H} \le |\Gamma|_{\cal H} + 4 t + \vert \Gamma_k \vert_{\cal H}$, we derive $\Vert \tilde{W} \Vert_\omega - \Vert W \Vert_\omega <0$ for $h_*$ small enough  by Lemma \ref{rig-lemma: derive prop3}(ii) since $\lbrace \Gamma, \Gamma_k\rbrace \subset {\cal L}_V$, where $\tilde{W} = (W \setminus V) \cup \partial V.$ This gives a contradiction to \eqref{rig-eq: property4} and finishes the proof. \eop

We close this section with the proof of Lemma \ref{rig-lemma: slicing bound2}. Recall the specific choice of $\tau$ in \eqref{rig-eq: tau-bartau} and set $N = N^\tau(\Gamma)$, $J = J^\tau$, ${\cal Y} = {\cal Y}^\tau$,  ${\cal C} = {\cal C}^\tau$ for shorthand.

\noindent {\em Proof of Lemma \ref{rig-lemma: slicing bound2}.}  As $\frac{1}{800}\bar{\tau} \le \tau \le \frac{1}{2}\bar{\tau}$, we can apply Lemma \ref{rig-lemma: slicing bound} to obtain that there are two sets $Y^1,Y^2 \in {\cal C}$ with $\bar{Y}^1 \cap \bar{Y}^2 = \emptyset$ such that  $\Vert Y \Vert_\pi \le \frac{19}{20}\tau$ for $Y \in {\cal  C}$ with $Y \cap (\bar{Y}^1 \cup \bar{Y}^2) = \emptyset$. We only construct the set $K_1$. Choose $Y_i = Y^1$ and set $S_l = Y_{i+l} \in {\cal C}$ for $|l|\le 3$. In particular, we have $S_l \cap S_0 \neq \emptyset$ for $l=-1,1$  and $\Vert S_l \Vert_\pi \le \frac{19}{20} \tau$ for $l=-3,3$. Set $S = \bigcup^2_{l=-2} S_l$. 

Let ${\cal R}(S)$ be defined as before \eqref{rig-eq: VertP}.  Arguing as in \eqref{rig-eq: (i)} or  \eqref{rig-eq: c split2} for $t=\tau$, respectively, depending on whether $S$ is contained in at most two adjacent parts of the neighborhood or $S$ intersects three parts of the neighborhood (possible for $l_2 \le \frac{l_1}{2})$, we derive $|\partial V|_\infty \le |\Gamma|_\infty + \sum\nolimits_{R_j \in {\cal R}(S)} |\partial R_j|_\infty + (\frac{19}{10}- \frac{1}{100})\tau - \Vert S\Vert_\pi$, where $V$ is the smallest rectangle whose closure contains  $\Gamma$ and ${\cal R}(S)$.  (Note that in the above calculation we possibly have to repeat the translation argument indicated at the end of the proof of Lemma \ref{rig-lemma: slicing bound2}.)  Thus,  by Lemma \ref{rig-lemma: slicing bound4} for $t=\tau$ and   \eqref{rig-eq: property4} 
\begin{align}\label{rig-eq: Ses}
\begin{split}
0  \le \Vert \tilde{W} \Vert_{\omega} - \Vert W \Vert_{\omega} &\le (1- h_*)\big((\tfrac{19}{10}- \tfrac{1}{100})\tau - \Vert S\Vert_\pi\big) + Ch_* \tau \\
&\le (1- h_*) \tfrac{19}{10}  \tau   - (1-h_*) \Vert S\Vert_\pi 
\end{split}
\end{align}
for $h_*$ small enough, where $\tilde{W} = (W \setminus V) \cup \partial V$.   We now construct the set $K_1$ and the corresponding  (at most) two connected components $T_1$, $T_2$ of $S\setminus K_1$ by distinction of the two following cases: 

a) If there is some $R_j$ with $| \partial R_j |_\pi \ge \frac{19}{20} \tau$ we choose $K_1 \in {\cal Y}$ as the smallest set such that $R_j \cap N \subset K_1$. Then the  (at most) two connected components  $T_1, T_2$ of $S \setminus K_1$ satisfy $\Vert T_i \Vert_\pi \le \frac{19}{20}\tau$ by \eqref{rig-eq: Ses}. Using \eqref{rig-eq: maxcracklength} we derive that $|K_1| \le C\frac{\tau^2}{h_*}$, as desired. 

b) Otherwise, we choose $K_1$ as follows. Assume $S = (\bigcup^{n'}_{i=1} \overline{Q_i})^\circ$ for $Q_i \in J$ and let $k \in \lbrace 0,\ldots,n'\rbrace$ be the index (if existent) such that $\Vert (\bigcup^k_{i=1} \overline{Q_i})^\circ\Vert_\pi\le \frac{19}{20}\tau$ and $\Vert (\bigcup^{k+1}_{i=1} \overline{Q_i})^\circ\Vert_\pi > \frac{19}{20}\tau$. Now define $T_1 = (\bigcup^k_{i=1} \overline{Q_i})^\circ$ and choose $K_1 = (\bigcup^l_{i=k+1} \overline{Q_i})^\circ$ for $l$ large enough such that $|K_1| \ge \bar{c}\frac{\tau^2}{h_*}$. Finally, let $T_2 = S \setminus (\overline{T_1 \cup {K_1}})$ and observe that for $\bar{c}$ large enough also $\Vert  T_2\Vert_\pi\le \frac{19}{20}\tau$ by \eqref{rig-eq: Ses}  and \eqref{rig-eq: maxcracklength}  since each rectangle can intersect at most one of the sets $T_1$, $T_2$.

Let $S^1_l$, $S^2_l$ be the connected components of $S_l \setminus K_1$  for $l=-2,-1,0,1,2$. Both cases a),b) above imply $\Vert S^i_l \setminus K_1 \Vert_\pi \le \max_{k=1,2} \Vert  T_k\Vert_\pi \le \frac{19}{20}\tau$ for $l=-2,\ldots,2$, $i=1,2$, which gives assertion (i). Assertion (ii) follows from the construction of the  set $K_1$ and definition \eqref{rig-eq: tau-bartau}.  Indeed, if $\Gamma_i \cap N  \neq \emptyset$, then $\Gamma_i \cap N^{\tau/20}$ and thus recalling \eqref{rig-eq: VertP} we find $\vert \Gamma_i \vert_\pi \ge \frac{19}{20}\tau$ and then $\Gamma_i \cap N \subset \bar{Y}^1 \cup \bar{Y}^2$.  This also shows, as remarked at the end of Section \ref{rig-sec: subsub,  rec-neigh}, that in this case $K_1$ is contained in one of the sets $N_{j,\pm}$, $j=1,2$.

Finally, in the case  $\Vert Y^1\Vert_\pi, \Vert Y^2\Vert_\pi \ge \frac{19}{20}\tau$ the statement $\dist(K_1,K_2) \ge c|\Gamma|_\infty$ follows directly from the fact that the set $\bar{Y}^1 \cup \bar{Y}^2$  intersects both $N_{1,+}$ and $N_{1,-}$ or both $N_{2,+}$ and $N_{2,-}$ ($N_{1,\pm}$ if $l_2 \le \frac{l_1}{ 2}$).  Otherwise, we can set $K_1 = \emptyset$ or $K_2 = \emptyset$ and the property trivially holds.\eop

\section{Trace estimates for boundary components}\label{rig-sec: sub, trace bc}

 This section is entirely devoted to the proof of Theorem \ref{rig-theorem: D} and we start with an outline of the proof.  We first give some preliminary estimates including an approximation of $u$ by a piecewise infinitesimal rigid motion. Here we also discuss the passage from an estimate in the neighborhood to a trace estimate. Afterwards, we proceed in several steps.

In Step 1   we assume that in the neighborhood $N$ of $\Gamma$ only small cracks $\Gamma_l$ are present. This indeed induces that $|Eu|(N)$ is sufficiently small as on each $\Gamma_l$ we have by assumption $[u] \sim \sqrt{|\Gamma_l|_\infty \eps}$.  In general, the idea is to construct thin long paths in $N$ which avoid cracks being to large. We then first measure the distance of the function from an infinitesimal rigid motion only on this path and may apply this result to estimate the distance in the whole set $N$ afterwards. A major drawback of such a technique is that the constant in \eqref{rig-eq: help-korn} crucially depends on the domain and blows up for sets getting arbitrarily thin. Consequently, in this context we have to carry out a careful  quantitative analysis how the constant in \eqref{rig-eq: help-korn} depends on the shape of the domain (see Section \ref{rig-sec: sub, pre}).

In Step 2 we suppose that we have a bound on the projection  $\Vert \cdot \Vert_\pi$ (recall definition \eqref{rig-eq: VertP}). We  observe that  the paths in general cannot be selected in such  a way that they only intersect sufficiently small cracks. Nevertheless, it can be shown that boundary components being too large for a direct application of the above ideas occupy only a comparably small region. In this region we then do not use the Korn-Poincar\'e inequality in \eqref{rig-eq: help-korn}, but  circumvent the estimate of the surface energy by a slicing technique and the fact that by the bound on $\Vert \cdot \Vert_\pi$  we   find  small stripes in the neighborhood, which do not intersect the jump set at all. The assertion then follows as this exceptional set can then be taken arbitrarily small by an iterative application of the slicing method.   
  
In Step 3 we repeat the ideas of Step 2 near the corners of $\Gamma$. Finally, in Step 4 we present the general proof taking into account the possible existence of sets $\Psi^1$, $\Psi^2$ discussed in Section  \ref{rig-sec: subsub,  dod-neigh}. We remark that the result crucially depends on a suitable L$^2$- trace theorem for SBD functions (see Lemma \ref{rig-sec: sub, trace}) which can be established in our framework due to the sufficiently regular jump set. Moreover, it is essential that there are at most two large cracks in a neighborhood. Already with three or four cracks the configurations might be significantly less rigid.

\subsection{Preliminary estimates}\label{rig-sec: subsub,  pre-est}

We remark that in this section we will frequently use the estimates derived in Appendix \ref{rig-sec: sub, pre}  and \ref{rig-sec: sub, trace}.  Moreover, recall the construction of the neighborhoods in Section \ref{rig-sec: sub, neigh}. We start with a discussion about the parameters.  

Assume $h_*, q, \omega_{\rm min} >0$ as in \eqref{eq: parameters} have been chosen in the previous sections and $\upsilon >0$ is small enough with respect to $h_*, q$. (Note that in the proofs up to now the choice of $\upsilon$ was independent of $\omega_{\rm min}$.) Below we will introduce $r =r(h_*,q) \in (0,1)$ and we then possibly pass   to a smaller $\upsilon>0$  such that $r (1-\omega_{\rm min})^3\ge \upsilon$. This indeed implies $\upsilon=\upsilon(h_*,q,\omega_{\rm min},r)$ and also shows \eqref{eq: parameters}(ii) as $r=r(h_*)$ since $q=q(h_*)$.  

As in the assumption of Theorem \ref{rig-theorem: D} we  let $\eps>0$, $\lambda > 0$, $W \in {\cal V}^s$, $W_i \in {\cal W}^s_\lambda$  and let $u \in H^1(W)$. Let $\Gamma = \Gamma(W)$ with $\vert \Gamma \vert_\infty \ge \lambda$ and the corresponding neighborhoods $N^\tau = N^\tau(\Gamma)$ and $N^{2\hat{\tau}}= N^{2\hat{\tau}}(\Gamma)$ be given. Suppose that \eqref{rig-eq: property4},  \eqref{rig-eq: property,b}, \eqref{rig-eq: property,a} and \eqref{rig-eq: property,c} for ${\cal T}_{\hat{\tau}}(W_i,W)$ hold. Assume that $\overline{N^{2\hat{\tau}}} \subset H(W)$.

In addition, we define the neighborhood $\tilde{N} = N(\Gamma) \setminus (X^1 \cup X^2)$, where  $\partial X^1 = \Gamma^1$, $\partial X^2 = \Gamma^2$ are the boundary components satisfying $|\Gamma^i|_\infty \ge \hat{\tau} = q^2 h_*^{-1} \upsilon \vert \Gamma \vert_\infty$ and $\Gamma_i \cap N^\tau \neq \emptyset$, see  Lemma  \ref{rig-lemma: two components}  (note that $X_1, X_2 = \emptyset$ is possible).  In the following we will write for shorthand  $N$, $N^{2\hat{\tau}}$ and $\tilde{N}$ if no confusion arises. By Remark \ref{rig-rem: 1}(i) it is not restrictive to assume that $J = J^\tau(\Gamma)$ as defined  before equation  \eqref{rig-eq: covering def}  consists of  (almost) squares.

Note that the inclusion $W_i \cap N \subset \tilde{N}$ may be strict due to boundary components $\Gamma_l$ with  $\vert \Gamma_l \vert_\infty  < \hat{\tau}$ having nonempty intersection with $N$.  Observe that by \eqref{rig-eq: W prop}(iv),(v) and \eqref{rig-eq: new2}(i) for $q$ large we have $\vert \partial R_l \vert_\infty  < \hat{\tau}$  for these  $\Gamma_l$, where $R_l$ is the corresponding rectangle given by \eqref{rig-eq: W prop}(i),(v). Then recalling \eqref{rig-eq: hat tau def}, for $q$ sufficiently large and $h_*$ small we get 
\begin{align}\label{rig-eq: new9}
\overline{N^{2\hat{\tau}_l} (\partial R_l)} \subset N^{2\hat{\tau}}(\Gamma) \subset H(W)
\end{align}
and thus $\Gamma_l \in  {\cal T}_{\hat{\tau}}(W_i,W)$. Consequently, we can extend $u$ as an SBD function from $N \cap W_i$ to  $\tilde{N}$ as defined in \eqref{rig-eq: extend def}.    As before, we will write for shorthand $\alpha(U) = \Vert e(\bar{u})\Vert^2_{L^2(U)}$ as well as 
\begin{align}\label{rig-eq: shorthand}
{\cal E}(U) =  \int_{U} |e(\bar{u})| + \int_{J_{\bar{u}} \cap U}|[\bar{u}]|\,d{\cal H}^1, \ \  \ \   \hat{\alpha}(U) = \int_{U} |e(\bar{u})|
\end{align}
for $U \subset \tilde{N}$. Note that $\alpha(U) = \alpha(U \cap W_i)$.

We now begin with some preliminary estimates.   First assume $\tilde{N} = N$. We apply Theorem \ref{rig-theo: korn} on each $Q \in J$  recalling that the constant is invariant under rescaling of the domain: This yields functions $\bar{A}: N \to \R^{2 \times 2}_{\rm skew}$, $\bar{b}: N \to \R^2$ being constant on each $Q \in J$ such that by \eqref{rig-eq: property,b} and \eqref{rig-eq: property,c}  we obtain
\begin{align}\label{rig-eq: D1.5}
\begin{split}
\int_N & |\bar{u}(x)- (\bar{A}\,x +  \bar{b})|^2 \, dx  \le C  \sum\nolimits_{Q \in J} ({\cal E}(Q))^2 \le C  ({\cal E}(N))^2 \\
& \le C  \upsilon\vert \Gamma \vert^2_\infty  \alpha(N) + C\big(\sum\nolimits_l \vert \Theta_l \cap N\vert_{\cal H}^{1/2}  \Vert [\bar{u}]\Vert_{L^2(\Theta_l)}  \big)^2  \\
& \le C   \upsilon^2\vert \Gamma \vert^3_\infty  \eps  + CC_*  \upsilon^{-4}  \eps \vert \partial W_i \cap N\vert_{\cal H}\big(\sum\nolimits_l \vert \Gamma_l \vert_{\infty} \big)^2  \\
& \le C   \upsilon^2\vert \Gamma \vert^3_\infty  \eps   + CC_*  \upsilon^{-4}  \eps \vert \partial W_i \cap  N^{2\hat{\tau}}\vert^3_{\cal H}  \le  C(1+C_*)  \upsilon^{-1} \vert \Gamma \vert^3_\infty  \eps
\end{split}
\end{align}
for some  $C =C(h_*, q)$. Here and in the following the sum always runs over the boundary components having nonempty intersection with $N$. In the third step we  employed H\"older's inequality and  $|N| \le C\upsilon |\Gamma|_\infty^2$. In the penultimate step we used $\sum_l |\Gamma_l|_{\infty} \le \sum_l | \Gamma_l|_{\cal H}\le 2|\partial W_i \cap N^{2\hat{\tau}}|_{\cal H}$  by Lemma \ref{rig-lemma: inftyX}(i) and \eqref{rig-eq: new9}. The constants used in this section may as usual vary from line to line depending on $h_*$ and $q$, but are always independent of the parameters $r$,  $\omega_{\rm min}$ and $\upsilon$.

In the general case, recall the definition of $\Psi^1$ and $\Psi^2$ in Section  \ref{rig-sec: subsub,  dod-neigh}. By Lemma \ref{rig-lemma: two observations1}(ii), Lemma \ref{rig-lemma: two observations2}(ii) and Remark \ref{rig-rem: 1}(i) it is not restrictive to assume that $\Psi^{i}_1$, $\Psi^{i,*}_2$, $\Psi^{i}_3$ are squares for $i=1,2$.  For simplicity we write $\Psi^{i}_2$ instead of $\Psi^{i,*}_2$ in the following. Similarly as in the previous estimate we obtain functions $\bar{A}: N \to \R^{2 \times 2}_{\rm skew}$, $\bar{b}: N \to \R^2$ being constant on each $Q \in J$ with $Q \cap (\Psi^1 \cup \Psi^2) = \emptyset$ and $\Psi^i_j$, $i=1,2$, $j=1,2,3$, such that
\begin{align}\label{rig-eq: D1.2}
\begin{split}
(i)& \ \ \int_{\tilde{N} \setminus (\Psi^1 \cup \Psi^2)}  |\bar{u}(x)- (\bar{A}\,x +  \bar{b})|^2 \, dx  \le C(1+C_*) \upsilon^{-1} \vert \Gamma \vert^3_\infty  \eps ,\\
(ii) & \ \ \int_{\Psi^i_j}  |\bar{u}(x)- (\bar{A}\,x +  \bar{b})|^2 \, dx  \le C(1+C_* )  \upsilon^{-3}\vert \Gamma \vert^2_\infty   \psi^i \eps. 
\end{split}
\end{align}
To see (ii), we apply Theorem \ref{rig-theo: korn} on the sets $\Psi^{i}_1$, $\Psi^{i}_2$, $\Psi^{i}_3$ and follow the lines of the previous estimate to obtain that the left hand side is bounded by $C  \upsilon\vert \Gamma \vert^2_\infty  \alpha(\Psi^i_j) + CC_*  \upsilon^{-4}  \eps \vert \partial W_i \cap \Psi^i_j\vert_{\cal H} \vert \partial W_i \cap  N^{2\hat{\tau}}\vert^2_{\cal H}$. We then use $\alpha(\Psi^i_j) \le  D\eps (1- \omega_{\rm min})^{-1} \psi^i \le D\eps \upsilon^{-1} \psi^i$ and $|\partial W_i \cap \Psi^i|_{\cal H}\le C  \upsilon^{-1}\psi^i$ by  \eqref{rig-eq: property,b}, \eqref{rig-eq: property,a}. 

The goal will be to replace the functions $\bar{A}$, $\bar{b}$ in \eqref{rig-eq: D1.2} by constants $A \in \R^{2 \times 2}_{\rm skew}$ and $b \in \R^2$ such that 
\begin{align}\label{rig-eq: D1.2.2}
\begin{split}
(i)& \ \ \int_{\tilde{N} \setminus (\Psi^1_2 \cup \Psi^2_2)}  |\bar{u}(x)- (A\,x +  b)|^2 \, dx  \le C(1+rC_*)  {\upsilon^{-3}}\vert \Gamma \vert^3_\infty  \eps, \\
(ii) & \ \ \int_{\Psi^i_2}  |\bar{u}(x)- (A\,x +  b)|^2 \, dx   \le C(1+ rC_* )  \upsilon^{-4}\vert \Gamma \vert^2_\infty    \psi^i \eps, 
\end{split}
\end{align}
for $i=1,2$ and for $r  (1- \omega_{\rm min})^3\ge \upsilon$. Then the trace theorem applied  on each square  (if $J$ or $\Psi^j$, $j=1,2$, consist also of rectangles, they can be covered by possibly overlapping squares) implies the assertion as follows: 

To satisfy the assumptions of Lemma \ref{rig-lemma: trace}, the jump set has to be the union of rectangle boundaries. Therefore, we extend $J_{\bar{u}} \cap \tilde{N}$ to $\tilde{J}_{\bar{u}} = \bigcup_l \partial R_l \cap \tilde{N}$  by $[\bar{u}](x)=0$ for $x \in \tilde{J}_{\bar{u}} \setminus J_{\bar{u}}$, where $R_l$ are the corresponding rectangles given in \eqref{rig-eq: W prop}(i),(v).  We observe that by    Lemma  \ref{rig-lemma: inftyX}(i),(iv) and \eqref{rig-eq: new2}(i) we get for a universal $C>0$
\begin{align}\label{rig-eq: new8}
\sum\nolimits_l |\partial R_l|_{\cal H} \le C\sum\nolimits_l |\Gamma_l|_{\cal H}, \ \ \ \sum\nolimits_l |\partial R_l|^{-1}_{\cal H} |\Gamma_l|_\infty^2 \le C\sum\nolimits_l |\Gamma_l|_{\cal H}.
\end{align}
Note that every boundary component $\Gamma_l$ intersects   at most $C(h_*,q)$ different squares  since $|\Gamma_l|_\infty < \hat{\tau}$ (see before \eqref{rig-eq: new9}).    Then by  Lemma \ref{rig-lemma: trace},  applied on each square in $J$, either for $\mu \sim \upsilon |\Gamma|_\infty$ or $\mu \sim \psi^i$,  taking the sum over all squares we obtain  by \eqref{rig-eq: D1.2.2} and \eqref{rig-eq: property,c}   for $\upsilon$ small enough 
\begin{align*}
& \int_\Gamma  |\bar{u}(x)- (A\,x + b)|^2 \, d{\cal H}^1(x)  \\
&\le C\upsilon \vert \Gamma \vert_\infty \alpha(\tilde{N}) + CC_*  \sum\nolimits_l \vert \partial R_l\vert_{\cal H} \sum\nolimits_{l}  \vert \partial R_l\vert^{-1}_{\cal H}  \eps \upsilon^{-4}    \vert \Gamma_l \vert^2_{\infty}\\
& \ \, + C (\upsilon \vert \Gamma \vert_\infty)^{-1} \Vert \bar{u} - (A \, \cdot + b)\Vert^2_{L^2(\tilde{N}\setminus (\Psi^1_2 \cup \Psi^2_2)} + \sum\nolimits^2_{i=1} C (\psi^i)^{-1} \Vert \bar{u} - (A \, \cdot + b)\Vert^2_{L^2(\Psi^i_2)}\\
& \le  C(1+rC_*) \upsilon^{-2} \vert \Gamma \vert^2_\infty  \eps + C(1+ rC_* )  \upsilon^{-4}\vert \Gamma \vert^2_\infty     \eps \le C(1+ rC_* ) \eps \upsilon^{-4}\vert \Gamma \vert^2_\infty   ,
\end{align*}
where for the first two terms we proceeded similarly as in \eqref{rig-eq: D1.5}, also taking  \eqref{rig-eq: new8} into account.  Finally, choosing $\hat{C} = C  = C(h_*,q)$ and $r =  r(h_*,q)$ small enough (i.e. also $\upsilon$ small enough) such that $r\hat{C} \le \frac{1}{2}$ we get \eqref{rig-eq: D0},  as desired.  Consequently, it suffices to establish \eqref{rig-eq: D1.2.2}.

\subsection{Step 1: Small boundary components}\label{rig-sec: subsub,  small-bound}

We first treat the case that only small components $\Gamma_l$ lie in $N$. For $1 > r \ge \upsilon >0$  define $T =  \lfloor \text{log}_r(\upsilon) \rfloor$ and let  for all $Q \in J$
\begin{align}\label{eq:S}
{\cal S}_t(Q) = \lbrace \Gamma_l: \Gamma_l \cap Q \neq \emptyset, \ \upsilon^{4} r^{-2t}\vert \Gamma \vert_\infty  < \vert \Gamma_l \vert_\infty \le \upsilon^{4} r^{-2t-2}\vert \Gamma \vert_\infty\rbrace
\end{align}
for all $ t \in \N$ and  ${\cal S}_0(Q) = \lbrace \Gamma_l: \Gamma_l \cap Q \neq \emptyset, \vert \Gamma_l \vert_\infty \le \upsilon^{4} r^{-2}\vert \Gamma \vert_\infty\rbrace$.  Moreover, by ${\cal Y}'$ we denote the set of subsets of $N$ consisting of squares in $J$.  (In contrast to ${\cal Y}^\tau$ as defined before \eqref{rig-eq: covering def} the connectedness of the sets is not required.)

\begin{lemma}\label{rig-lemma: D1}
Theorem \ref{rig-theorem: D} holds under the additional assumption that there is some $\frac{T}{4} + 2 \le t \le  \frac{T}{2} - 1 $ such that $ \bigcup_{s > t} {\cal S}_s(Q) = \emptyset$ for all $Q \in J$ and $\sum_{Q} \#{\cal S}_t(Q) \le \upsilon^{-3} r^{2t + 3} $. 
\end{lemma}

\Proof We first observe that the assumption implies $\tilde{N} =N$. Let $\frac{T}{4} + 2 \le t \le \frac{T}{2}-1$  with the above properties be given and write $\hat{\upsilon} = \upsilon^{2} r^{ - \frac{2t+1}{2}}$ for shorthand. For later we note that
\begin{align}\label{rig-eq: D1.7}
 \upsilon^\frac{7}{4} r^{-\frac{3}{2}} \le \hat{\upsilon} \le \upsilon^{\frac{3}{2}} \sqrt{r}. 
\end{align}
We cover $N$ with squares $\hat{Q}(\xi) =  \hat{Q}^{\hat{\upsilon} \vert \Gamma \vert_\infty}(\xi)$ of length $2\hat{\upsilon} \vert \Gamma \vert_\infty$ and midpoint $\xi$.  (If the sets in $J=J(\Gamma)$ constructed in Section \ref{rig-sec: subsub,  rec-neigh} are not perfect squares, the sets $\hat{Q}(\xi)$ shall be chosen appropriately. The difference in the  possible shapes, however, does not affect the following estimates by Remark \ref{rig-rem: 1}(i).)  We will now consider a rectangular path, i.e. a path  $\xi = (\xi_0, \ldots \xi_n=\xi_0)$ of square midpoints intersecting all $Q \in J$ such that there are indices $i_1,i_2,i_3$ with $\xi_j - \xi_{j-1} = \pm 2\hat{\upsilon} \vert \Gamma \vert_\infty \e_1$ for all $0 \le j \le i_1$, $i_2 \le j \le i_3$ and $\xi_j - \xi_{j-1} = \pm 2\hat{\upsilon} \vert \Gamma \vert_\infty \e_2$ else. Observe that the number of squares in a path satisfies $n \le C\hat{\upsilon}^{-1}$ and that we can find $\sim \upsilon \hat{\upsilon}^{-1}$ disjoint rectangular paths in $N$. Consequently, by assumption and \eqref{rig-eq: property,b} we can find at least one rectangular path $P:= \overline{\bigcup_j \hat{Q}(\xi_j)}$
such that
\begin{align}\label{rig-eq: D1.4}
\alpha(P) \le C\hat{\upsilon} \eps \vert \Gamma \vert_\infty, \ \ \,  \sum_{\Gamma_l \in \hat{\cal S}(P)} \vert \Gamma_l\vert_\infty \le C\hat{\upsilon} \vert \Gamma \vert_\infty, \ \ \,   \# \hat{\cal S}^t(P) \le C\tfrac{\hat{\upsilon}}{\upsilon} \upsilon^{-3} r^{2t + 3} 
\end{align}
for some sufficiently large constant $C=C(h_*,q)$, where $\hat{\cal S}(P)=  \lbrace \Gamma_l: \Gamma_l \cap P \neq \emptyset \rbrace$ and $\hat{\cal S}^t(P) = \hat{\cal S}(P) \cap \bigcup_{Q} {\cal S}_t(Q)$. Here we used that each $\Gamma_l \in \bigcup_{Q \in J} \bigcup_{s\le t}{\cal S}_s(Q)$ intersects at most four adjacent squares $\hat{Q}$ because $| \partial R_l|_\infty \le C\hat{\upsilon}^2r^{-1}|\Gamma|_\infty  \ll \hat{\upsilon}|\Gamma|_\infty$ by  \eqref{rig-eq: new2}(i), \eqref{rig-eq: D1.7} and  $\Gamma_l \subset \overline{R_l}$ (see  \eqref{rig-eq: W prop}(i)). Observe that the above path can be chosen in the way that also $|P \cap Q| \ge C\frac{\hat{\upsilon}}{\upsilon}|Q|$ for $Q = E_{\pm,\pm}$, where $E_{\pm,\pm}$ denote the squares in the corners of $N$, i.e.  $(\pm l_1,\pm l_2) \cap  \overline{E_{\pm,\pm}} \neq \emptyset$.  (Recall the construction of the neighborhood in Section \ref{rig-sec: subsub,  rec-neigh}.) This implies $|P \cap Q| \ge C\frac{\hat{\upsilon}}{\upsilon}|Q|$ for all $Q \in J$. It is convenient to write the above estimate in the form 
\begin{align}\label{rig-eq: D1.6}
\upsilon^4 r^{-2t-1 \pm 1} = \hat{\upsilon}^2 r^{\pm 1} ,\ \ \ \# \hat{\cal S}^t(P) \upsilon^{4} r^{-2t-2} \le  C\hat{\upsilon} r.
\end{align}
We now apply Lemma  \ref{rig-lemma: weak rig}(ii) with $s= \hat{\upsilon} \vert \Gamma \vert_\infty$ and $U = P$ with $|P| \le C\hat{\upsilon}|\Gamma|_\infty^2$. For later purpose in Section \ref{rig-sec: subsub, proj-bound1}, we consider general subsets $Z \in {\cal Y}'$ consisting of squares in $J$. Recall that  we get $\Vert [\bar{u}] \Vert_{L^1(\Theta_l)} \le \sqrt{|\Theta_l|_{\cal H}} \Vert [\bar{u}] \Vert_{L^2(\Theta_l)} \le C\sqrt{C_* \eps\upsilon^{-4}}|\Gamma_l|^{3/2}_{\infty}$  by \eqref{rig-eq: new2}(ii), \eqref{rig-eq: property,c} and H\"older's inequality.  Arguing similarly as in \eqref{rig-eq: D1.5} we find $A \in \R^{2 \times 2}_{\rm skew}$ and $b \in \R^2$ such that by  \eqref{rig-eq: new9},  \eqref{rig-eq: D1.4} and \eqref{rig-eq: D1.6}      ($\dashint_{P \cap Z}$ stands for $\frac{1}{|P \cap Z|}\int_{P \cap Z}$)
\begin{align}\label{rig-eq: D1.1}
\begin{split}
 |P|\dashint_{P\cap Z}& |\bar{u}(x)- (A\,x +   b)|^2 \, dx  \le C  \hat{\upsilon}^{-3} ({\cal E}(P))^2 \\
& \le C   \hat{\upsilon}^{-3}\hat{\upsilon}\vert \Gamma \vert^2_\infty \alpha(P)   + CC_*  \hat{\upsilon}^{-3}\frac{\eps}{\upsilon^4} \big(\sum\nolimits_{\Gamma_l \in \hat{\cal S}(P)} \vert \Gamma_l \vert_\infty^{3/2} \big)^2 \\
& \le C   \hat{\upsilon}^{-1}\vert \Gamma \vert^3_\infty  \eps  + CC_*  \hat{\upsilon}^{-1}r^{-1} \vert \Gamma\vert_\infty\frac{\eps}{\upsilon^4} \big(\sum\nolimits_{\Gamma_l \in \hat{\cal S}^t(P)} \vert \Gamma_l \vert_\infty \big)^2 \\
& \ \ \ + CC_*  \hat{\upsilon}^{-1}r \vert \Gamma\vert_\infty\frac{\eps}{\upsilon^4} \big(\sum\nolimits_{\Gamma_l \in \hat{\cal S}(P) \setminus \hat{\cal S}^t(P)} \vert \Gamma_l \vert_\infty \big)^2 \\
& \le C\hat{\upsilon}^{-1}\vert \Gamma \vert^3_\infty  \eps + CC_* \hat{\upsilon} r  \frac{\eps}{\upsilon^4}\vert \Gamma\vert^3_\infty   + CC_* \hat{\upsilon} r \frac{\eps}{\upsilon^4} \vert \Gamma\vert^3_\infty.
\end{split}
\end{align}
Observing that $\hat{\upsilon}^{-1} \le \hat{\upsilon} \upsilon^{-4} $ by  definition of $\hat{\upsilon}$, we derive
\begin{align}\label{rig-eq: D1.3}
 |P|\dashint_{P\cap Z} |\bar{u}(x)- (A&\,x +   b)|^2 \, dx \le C(1 + C_* r) \, \hat{\upsilon}  \frac{\eps}{\upsilon^4} \vert \Gamma \vert^3_\infty =: F.
\end{align}
We now pass from an estimate on $P$ to an estimate on $N$. Since  $|P \cap Q| \ge C\frac{\hat{\upsilon}}{\upsilon}|Q|$ for all $Q \in J$ we find $\frac{|Z|}{|N|}\ge C \frac{|Z\cap P|}{|P|} \ge C\upsilon$. Therefore, by \eqref{rig-eq: D1.5}  (recall that $\tilde{N} = N$) and $\upsilon^2 \le \hat{\upsilon}  r$ (see \eqref{rig-eq: D1.7}) we also have 
$$\int_Z |\bar{u}(x)- (\bar{A}\,x +   \bar{b})|^2 \, dx \le |Z \cap P||P|^{-1}CF \le C|Z||N|^{-1}F.$$
 We apply the triangle inequality, Lemma \ref{lemma: diff1}(ii) and Remark \ref{rem: diff1}(i) on each $Q \subset Z$ with $D_1 = P \cap Q$, $D_2 = Q$ noting that $\bar{A},\bar{b}$ are constant on each square. As  $|P \cap Q| \ge C\frac{\hat{\upsilon}}{\upsilon}|Q|$ for all $Q \in J$  and ${\rm diam}(D_1) \ge C {\rm diam}(D_2)$ we obtain 
$$\Vert  (\bar{A}\cdot + \bar{b}) - (A\cdot + b)\Vert^2_{L^2(Q)}  \le C\frac{\upsilon}{\hat{\upsilon}}  \big( \Vert\bar{u}- (\bar{A}\cdot +   \bar{b})\Vert^2_{L^2(Q)}  + \Vert\bar{u}- (A\cdot+   b)\Vert^2_{L^2(P\cap Q)}\big)$$
and thus summing over all $Q \subset Z$ we derive  again by the triangle inequality
\begin{align}\label{rig-eq: P to N}
|N|\dashint_Z |\bar{u}(x)- (A\,x + b)|^2 \, dx \le C\hat{\upsilon}^{-1} \upsilon F = C(1 + C_* r) \,   \frac{\eps}{\upsilon^3} \vert \Gamma \vert^3_\infty.
\end{align}
Consequently, setting $Z = N$, \eqref{rig-eq: D1.2.2}(i) is established, as desired. \eop

\subsection{Step 2: Subset with small projection of components}\label{rig-sec: subsub, proj-bound1}

The next step will be the case that $\Vert \cdot \Vert_\pi$ is  not too large.  For that purpose, recall \eqref{rig-eq: VertP} and the definition of ${\cal Y} = {\cal Y}^\tau$ (see before \eqref{rig-eq: covering def}). Consider some $U \in {\cal Y}$. In this section we show that for all $Z \subset U$, $Z \in  {\cal Y}'$, one has 
\begin{align}\label{rig-eq: D2.3} 
|U|\dashint_{Z}  |\bar{u}(x)- (A_U\,x +  b_U)|^2 \, dx  \le C(1+rC_*)  \upsilon^{-3}\vert \Gamma \vert^3_\infty  \eps 
\end{align}
for $A_U \in \R^{2 \times 2}_{\rm skew}$, $b_U \in \R^2$. Recall that $E_{\pm,\pm}$ denote the squares at the corners of $\Gamma$ (see construction before  \eqref{rig-eq: covering def}).   We start with an auxiliary result.

\begin{lemma}\label{rig-lemma: D1,5}
Let $r \ge \upsilon >0$ as in \eqref{eq: parameters}  small enough and $C'>0$. Then there is a constant $C=C(C')>0$ independent of $r,\upsilon$ such that for all $U \in {\cal Y}$ with  $|U| \ge C' \upsilon |\Gamma|_\infty^2$,  $U \cap E_{\pm,\pm} = \emptyset$ and $ \Vert U\Vert_\pi \le \frac{19}{20} \tau$ there is a subset $U' \subset U$, 
 $U' \in {\cal Y}'$, with $|U\setminus U'|  \le Cr|U|$ such that \eqref{rig-eq: D2.3} holds for all $Z \subset U'$, $Z \in  {\cal Y}'$.
\end{lemma}

\Proof Let $U \in  {\cal Y}$ be given with $ \Vert U\Vert_\pi \le \frac{19}{20} \tau$ and assume without restriction $U \subset N_{2,+}\setminus (N_{1,-} \cup N_{1,+})$.  By the choice of $\tau$ in \eqref{rig-eq: tau-bartau}  and the definition of $\Vert \cdot \Vert_\pi$ in \eqref{rig-eq: VertP} we obtain that all  $\Gamma_l$ having nonempty intersection with $U$ satisfy $\vert \Gamma_l\vert_\infty < 19\bar{\tau}$. In particular, this implies $U \cap \tilde{N} = U$. Let $(\partial R_l)_l$ be the rectangles corresponding to $(\Gamma_l)_l$  as given by \eqref{rig-eq: W prop}(i),(v).  We first prove that there is a $\frac{T}{4} +  2 \le t\le \frac{T}{2} -1 $ such that $\sum_{Q \subset U} \#{\cal S}_t(Q) \le \upsilon^{-3}r^{2t+3}$ as in the assumption of Lemma \ref{rig-lemma: D1}. If the claim were false, we would have (assume without restriction that $T \in 4\N$)
\begin{align}\label{rig-eq: D2.6}
\vert \partial W_i \cap N^{ \tau + 38\bar{\tau}} \vert_{\cal H} & \ge C \sum\nolimits^{\frac{T}{2} -1 }_{t=\frac{T}{4}+ 2}\sum\nolimits_{Q \subset U} \#{\cal S}_t(Q) \upsilon^{4}r^{ - 2t} \vert \Gamma \vert_\infty \ge CT \upsilon r^3 \vert \Gamma \vert_\infty\notag\\
& \ge C\text{log}_r(\upsilon) r^3 \upsilon\vert \Gamma \vert_\infty,  
\end{align}
 where in the first step we used that $|\partial R_l|_\infty <  38\bar{\tau}$  by \eqref{rig-eq: new2}(i) which implies $\Gamma_l \subset \overline{R_l}\subset N^{\tau +38\bar{\tau}}$ and ensures that $\Gamma_l$ intersects only a uniformly bounded number of different squares $Q \subset U$ (independently of $r,\upsilon$).   As for $\upsilon$ small enough (with respect to $r = r(h_*,q)$, $q$ and $h_*$) this is larger than $D\hat{\tau}$, we get a contradiction to \eqref{rig-eq: property,b}.  As before, we define $\hat{\upsilon} = \upsilon^{2} r^{ - \frac{2t+1}{2}}$ for shorthand.

As in the previous proof we will select a path in $U$ with certain properties.  Recalling \eqref{rig-eq: VertP} it is not hard to see that $|\pi_2 (R_l\cap U)| \le \vert \partial R_l \vert_\pi$. As by assumption $\Vert U \Vert_\pi \le \frac{19}{20} \tau$, we find a set $S \subset (l_2, l_2 + \tau)$  being the union of intervals $2k's + (-s,s), k'\in \Z$, with $|S| \ge \frac{\tau}{20}$ such that the stripe $\hat{U} = U \cap (\R \times S)$ satisfies $\partial W_i \cap \hat{U} = \emptyset$. We cover $U$ by  $k$ horizontal paths ${\cal P} = (P_i)_{i=1}^k$, consisting  of $\hat{Q}(\xi) = \hat{Q}^{\hat{\upsilon}\vert \Gamma\vert_\infty}(\xi)$, i.e. $k = \lceil(2\hat{\upsilon}\vert \Gamma\vert_\infty)^{-1} \tau\rceil$ as $|\pi_2 U|=\tau$. We can find a subset $\hat{\cal P}_1 \subset {\cal P}$ with $\# \hat{\cal P}_1 \ge c_1 k$ for $c_1$  small enough such that
\begin{align}\label{rig-eq: PS}
\Gamma_l \cap P_i = \emptyset \ \ \  \text{ for all } |\Gamma_l|_\infty \ge \bar{C}^2\hat{\upsilon}|\Gamma|_\infty \ \text{ and } \ P_i \in \hat{\cal P}_1
\end{align}
and  $|S \cap \pi_2 \bigcup_{P_i \in \hat{\cal P}_1} P_i| \ge \frac{\tau}{21}$, if  $\bar{C}=\bar{C}(h_*,q)$ is chosen sufficiently large. Indeed, for $\lbrace \Gamma_l: |\pi_2 \Gamma_l|\ge  \bar{C} \hat{\upsilon}|\Gamma|_\infty \rbrace$ this follows by an elementary argument  in view of  $\Vert U \Vert_\pi \le \frac{19}{20} \tau$.  On the other hand,  by \eqref{rig-eq: new2}(i) we see that each component in 
$${\cal G} := \lbrace \Gamma_l: |\Gamma_l|_\infty \ge \bar{C}^2\hat{\upsilon}|\Gamma|_\infty, |\pi_2 \Gamma_l|\le \bar{C}\hat{\upsilon}|\Gamma|_\infty  \rbrace$$
intersects at most $\sim \frac{\bar{C} \hat{\upsilon} |\Gamma|_\infty}{\hat{\upsilon} |\Gamma|_\infty} = \bar{C}$ different $P_i \in {\cal P}$ and thus using \eqref{rig-eq: property,b} the components in ${\cal G}$ intersect at most  $\bar{C} \frac{D\hat{\tau}}{\bar{C}^2 \hat{\upsilon} |\Gamma|_\infty}  \sim \frac{D\hat{\tau}}{\hat{\upsilon}|\Gamma|_\infty\bar{C}} \ll \frac{\tau}{|\Gamma|_\infty\hat{\upsilon}}$ different $P_i \in {\cal P}$.

Moreover, we can then find a subset $\hat{\cal P}_2 \subset \hat{\cal P}_1$ with $\#  \hat{\cal P}_2 \ge c_2k$ for $c_2$ sufficiently small such that 
\begin{align}\label{rig-eq: D2.2}
|\pi_2(P_i \cap \hat{U})|\ge \tfrac{1}{22} 2\hat{\upsilon} \vert \Gamma\vert_\infty \ \ \ \text{ for all } \ P_i \in \hat{\cal P}_2.
\end{align}
 Recall that we have already found a  $\frac{T}{4} +  2 \le t\le \frac{T}{2} -1 $ such that $\sum_{Q \subset U} \#{\cal S}_t(Q) \le \upsilon^{-3}r^{2t+3}$. Using \eqref{rig-eq: PS} we can now choose a path $P  =\overline{\bigcup_j \hat{Q}(\xi_j)} \in \hat{\cal P}_2$ such that \eqref{rig-eq: D1.4} is satisfied possibly passing to a larger constant $C>0$ depending on $\bar{C}$.  (The essential difference to the argument developed in \eqref{rig-eq: D1.4} is the fact that every boundary component may intersect not only four squares but a number depending on $\bar{C}$.) Observe that  $n \le C \hat{\upsilon}^{-1}$ for a universal $C>0$,  where $n$ denotes the number of squares in the path $P$. Recall \eqref{eq:S}, $\hat{\cal S}(P)=  \lbrace \Gamma_l: \Gamma_l \cap P \neq \emptyset \rbrace$ and let
\begin{align}\label{rig-eq: D2.7}
\hat{\cal S}^t_>(P)=  \lbrace \Gamma_l: \Gamma_l \cap P \neq \emptyset, \Gamma_l \in \bigcup\nolimits_{Q \in J} \bigcup\nolimits_{s>t} {\cal S}_s(Q) \rbrace.
\end{align}
Moreover, define ${\cal K} = \lbrace \hat{Q} = \hat{Q}(\xi_j): \hat{Q} \cap \Gamma_l = \emptyset \ \text{ for all } \ \Gamma_l \in \hat{\cal S}^t_>(P) \rbrace$. By \eqref{rig-eq: D1.4} it is elementary to see that $\# \hat{\cal S}^t_>(P) \le C \hat{\upsilon} \vert \Gamma \vert_\infty (\upsilon^{4} r^{ - 2t-2}\vert \Gamma \vert_\infty)^{-1} = C\hat{\upsilon}^{-1} r$. Consequently, as by \eqref{rig-eq: PS} every $\Gamma_l \in \hat{\cal S}^t_>(P)$ intersects only a uniformly bounded number of adjacent sets, we find 
\begin{align}\label{rig-eq: D2.5}
n - \# {\cal K} \le C\# \hat{\cal S}^t_>(P) \le C\hat{\upsilon}^{-1} r.
\end{align}
Consider two squares $\hat{Q}(\xi_0),\hat{Q}(\xi_m) \in {\cal K}$ and the path $(\xi_0, \xi_1, \ldots, \xi_m)$. Define $D = \overline{\bigcup^m_{j=0} \hat{Q}(\xi_j)}$. Without restriction we assume $\hat{Q}_0 := \hat{Q}(\xi_0) = \mu(-1,1)^2$ and $\hat{Q}_m := \hat{Q}(\xi_m) = \mu((2m,0) + (-1,1)^2)$, where for shorthand we write $\mu = \hat{\upsilon} \vert \Gamma\vert_\infty$. We will now derive an estimate  such that Lemma \ref{rig-lemma: weak rig2}(i) is applicable. First of all, Theorem \ref{rig-theo: korn}, Theorem \ref{rig-th: tracsbv} and a scaling argument  show (recall \eqref{rig-eq: shorthand})
\begin{align}\label{rig-eq: D2.4}
\Vert \bar{u} - (A^i\, \cdot + b^i)\Vert_{L^2(\hat{Q}_i)} \le C{\cal E}(\hat{Q}_i), \ \ \ \Vert \bar{u} - (A^i\, \cdot + b^i)\Vert_{L^1( \partial \hat{Q}_i)} \le C{\cal E}(\hat{Q}_i)
\end{align}
for $A^i \in \R^{2 \times 2}_{\rm skew}$, $b^i \in \R^2$, $i=0,m$, and a constant independent of $\mu$. We claim that
\begin{align}\label{rig-eq: D2.1}
\begin{split}
&\mu^2  |a^0 - a^m| + \mu |b^0_1 - b^m_1|  \le C({\cal E}(\hat{Q}_0) + {\cal E}(\hat{Q}_m)) + C \hat{\alpha}(D),\\
&\mu  |b^0_2 - b^m_2| \le Cm({\cal E}(\hat{Q}_0) + {\cal E}(\hat{Q}_m)) + C m\hat{\alpha}(D),
\end{split}
\end{align}
where $b^i_j$ denotes the $j$-th component of $b^i$, $i=0,m$, and  $a^0,a^m$ are defined such that $A^i = \begin{pmatrix} 0 & a^i \\ -a^i & 0 \end{pmatrix}$. By \eqref{rig-eq: D2.2} we find two (measurable) sets $B_1,B_2 \subset \mu (-1,1)$ such that $|B_j| \ge \frac{\mu}{ 44}$, $\dist(B_1,B_2) \ge \frac{\mu}{ 22}$ and $E:= \mu  (-1,2m+1) \times B_1 \cup B_2 \subset W_i$.  We  apply a slicing argument in the first coordinate direction and obtain 
\begin{align}\label{rig-eq: slicing1}
\begin{split}
\int_{B_1 \cup B_2} |\bar{u}_1 &(\mu (2m-1),y) - \bar{u}_1(\mu,y) | \, dy \\ &\le \int_{B_1 \cup B_2} \Big| \int^{\mu(2m-1)}_{\mu }  \partial_1 \bar{u}_1(t,y)  \, dt\Big|\,dy  \le C \hat{\alpha}(E). 
\end{split}
\end{align}
This together with \eqref{rig-eq: D2.4} and the triangle inequality yields 
$$\Vert (a^0 - a^m)\,\cdot + (b^0_1 - b^m_1) \Vert_{L^1(B_1 \cup B_2)} \le C({\cal E}(\hat{Q}_0) + {\cal E}(\hat{Q}_m)) + C \hat{\alpha}(E).$$
Choose $f: B_1 \to \R$ such that $\id + f: B_1 \to B_2$ is piecewise constant and bijective. Thanks to $|B_1| \ge\frac{\mu}{44}$ and $f(y) \ge \frac{\mu}{22}$ for $y\in B_1$ we derive
\begin{align*}
\mu^2 &|a^0 -a^m| \le  C\Vert (a^0 - a^m)\,f(\cdot) \Vert_{L^1(B_1)}  \le   C \Vert (a^0 - a^m)\,\cdot + (b^0_1 - b^m_1) \Vert_{L^1(B_1)} \\
& \ \ + C \Vert (a^0 - a^m)\,(\cdot + f(\cdot))  + (b^0_1 - b^m_1) \Vert_{L^1(B_1)}  \le C({\cal E}(\hat{Q}_0) + {\cal E}(\hat{Q}_m)) + C\hat{\alpha}(E)
\end{align*}
and likewise $\mu |b^0_1 - b^m_1| \le C({\cal E}(\hat{Q}_0) + {\cal E}(\hat{Q}_m)) + C\hat{\alpha}(E)$. This gives the first bound in \eqref{rig-eq: D2.1} since $E \subset D$.  Analogously, we slice in $\zeta = \mu(2m-2,d)$ direction for $0<d <1$. By \eqref{rig-eq: D1.4} we find $|\pi_2 (\partial W_i \cap P_i)|\le C\mu$. Consequently, choosing $d$ small enough and recalling \eqref{rig-eq: D2.2}, we find a set $B_3 \subset \mu(-1,1-d)$ with $|B_3| \ge \frac{1}{ 24} 2\hat{\upsilon}|\Gamma|_\infty=\frac{\mu}{ 12}$ such that $ \lbrace \mu \rbrace \times B_3 + [0,1]\zeta \subset W_i$. Letting $\bar{\zeta} = \frac{\zeta}{|\zeta|}$ we get
$$ \int_{B_3} |\bar{u}(\mu,y)\cdot \bar{\zeta} - \bar{u}((\mu, y) + \zeta)\cdot \bar{\zeta} | \, dy  \le C \hat{\alpha}(E)$$ 
similarly to \eqref{rig-eq: slicing1} and thus, using \eqref{rig-eq: D2.4} and the fact that $A^m\,\bar{\zeta} \cdot \bar{\zeta} = 0$, we derive
 $$ \int_{B_3} | (A^0 - A^m) \,(\mu,y)^T \cdot \bar{\zeta} + (b^0 - b^m) \cdot \bar{\zeta} | \, dy  \le C({\cal E}(\hat{Q}_0) + {\cal E}(\hat{Q}_m)) + C \hat{\alpha}(E).$$
This together with the  first part of \eqref{rig-eq: D2.1} then leads to
$$\mu |b^0_2 - b^m_2 + (a^m - a^0)\mu|\ \tfrac{d}{|\zeta|} \le C \hat{\alpha}(E) + C({\cal E}(\hat{Q}_0) + {\cal E}(\hat{Q}_m)) $$ 
and implies the second part of \eqref{rig-eq: D2.1} as $E \subset D$. Summarizing, \eqref{rig-eq: D2.1} yields
\begin{align}\label{rig-eq: D2.8}
 \Vert b^0 - b^m + (A^0 - A^m)\,\cdot \Vert_{L^2(\hat{Q}_0 \cup \hat{Q}_m)}    \le C\hat{\upsilon}^{-1}({\cal E}(\hat{Q}_0)+ {\cal E}(\hat{Q}_m) + \hat{\alpha}(D)),
\end{align}
where we used that $m \le n \le C\hat{\upsilon}^{-1}$.  Define $\tilde{P}=  \overline{\bigcup_{\hat{Q} \in {\cal K}} \hat{Q}}$ and    apply Lemma \ref{rig-lemma: weak rig}(i) for $s = \hat{\upsilon}|\Gamma|_\infty$, $U =\tilde{P}$ with $|\tilde{P}| \le C \hat{\upsilon}|\Gamma|^2_\infty$ and $Z \subset U$, $Z \in  {\cal Y}'$ to derive by \eqref{rig-eq: D2.4} and \eqref{rig-eq: D2.8} 
\begin{align}\label{rig-eq: D1.1*}
 \Vert \bar{u} -  (A\, \cdot +  b)\Vert^2_{L^2(\tilde{P} \cap Z)} \le C|\tilde{P} \cap Z||\tilde{P}|^{-1} \hat{\upsilon}^{-3} (({\cal E}(\tilde{P}))^2 + (\hat{\alpha}(P))^2)
\end{align}
for suitable $A \in \R^{2 \times 2}_{\rm skew}$ and $b \in \R^2$. Note that the difference to the estimate in the proof of Lemma \ref{rig-lemma: D1} is that due to the above slicing argument it suffices to consider the elastic part of the energy in the connected components of  $P \setminus \tilde{P}$ (cf. \eqref{rig-eq: D1.1}). Recall the definition of ${\cal K}$ (cf. \eqref{rig-eq: D2.7}) and note that $\Gamma_l \cap \tilde{P} = \emptyset$ for all $\Gamma_l \in \hat{\cal S}^t_>$. Proceeding as in \eqref{rig-eq: D1.1} and \eqref{rig-eq: D1.3}, in particular using  \eqref{rig-eq: D1.4} and \eqref{rig-eq: D1.6}, we then obtain
\begin{align}\label{rig-eq: D2.14}
 \Vert \bar{u} -  (A\, \cdot +  b)\Vert^2_{L^2(\tilde{P} \cap Z)} \le C|\tilde{P} \cap Z||\tilde{P}|^{-1} (1 + C_* r) \, \hat{\upsilon}  \frac{\eps}{\upsilon^4} \vert \Gamma \vert^3_\infty.
\end{align}
 Now let $J' \subset J$ be the set of squares such that $|Q \cap \tilde{P}| \ge   C\upsilon \hat{\upsilon} \vert \Gamma \vert^2_\infty \ge C\frac{\hat{\upsilon}}{\upsilon} |Q|$ for all $Q \in J'$. Setting $U' = \big(\bigcup_{Q \in J'} \overline{Q}\big)^\circ  \in {\cal Y}'$ it is not hard to see that $|U\setminus U'|  \le Cr|U|$ for $r$ small enough  and a constant $C=C(C')$ as $|P \setminus \tilde{P}| \le Cr|P|$ by \eqref{rig-eq: D2.5}. Then applying  Lemma \ref{lemma: diff1}(ii),  \eqref{rig-eq: D1.5} and arguing as in \eqref{rig-eq: P to N} we derive  for $Z \subset U'$, $Z \in  {\cal Y}'$
$$|U'|\dashint_{Z} |\bar{u}(x)- (A\,x + b)|^2 \, dx \le C (1 + C_* r)  \frac{\eps}{\upsilon^3} \vert \Gamma \vert^3_\infty.$$
As $|U \setminus U'| \le Cr|U|$, this gives \eqref{rig-eq: D2.3}  for $r$ small enough, as desired. \eop

The next step will be to replace $U'$ by $U$ in Lemma \ref{rig-lemma: D1,5} and to drop the assumption $|U| \ge C' \upsilon |\Gamma|_\infty^2$. To this end, we will  apply the above arguments iteratively.

\begin{lemma}\label{rig-lemma: D2}
 Let $r \ge \upsilon >0$   as in \eqref{eq: parameters} small enough. Then for all  $U \in {\cal Y}$ with   $U \cap E_{\pm,\pm} = \emptyset$ and  $ \Vert U\Vert_\pi \le \frac{19}{20} \tau$ the estimate \eqref{rig-eq: D2.3} holds.
\end{lemma}

\Proof We first treat the case that $|U|\ge C'\upsilon|\Gamma|_\infty^2$ with $C'$ as in Lemma \ref{rig-lemma: D1,5} and indicate the adaptions for the general case at the end of the proof. Define $U_1  =U'$ and $J_1 = J'$ as given in Lemma \ref{rig-lemma: D1,5}  (see below \eqref{rig-eq: D2.14}).  Assume that $U_i =  \big(\bigcup_{Q \in J_i} \overline{Q}\big)^\circ$, $J_i \subset J$, with $U_{1} \subset \ldots \subset U_i$ is given such that  for $\bar{C} >0$ sufficiently large
\begin{align}\label{rig-eq: r-i}
|U\setminus U_i|\le \bar{C}r^{i}|U|
\end{align}
and for all $Z \subset U_i$, $Z \in  {\cal Y}'$, one has
\begin{align}\label{rig-eq: D2.10}
\int_{Z} |\bar{u}(x)- (A\,x + b)|^2 \, dx \le |Z||U|^{-1}  C\prod\nolimits^{i-1}_{j=0} \big(1+ \bar{C} r^{\frac{j}{8}}\big) \ G,
\end{align}
where $A \in \R^{2 \times 2}_{\rm skew}$, $b \in \R^2$ as in Lemma \ref{rig-lemma: D1,5} and $G:=(1 + C_* r)  \frac{\eps}{\upsilon^3} \vert \Gamma \vert^3_\infty$.

Observe that \eqref{rig-eq: r-i}, \eqref{rig-eq: D2.10} hold for $i=1$ by Lemma \ref{rig-lemma: D1,5}. We now pass from $i $ to $i+1$ and suppose $i \le T +  2$. First, it is not restrictive to assume that $|U \setminus U_i| \ge  \bar{C}r^{i+1}|  U|$  for $\bar{C}>0$ as above since otherwise we may set $U_{i+1} = U_i$. We cover $U \setminus U_i$ with pairwise disjoint, connected sets $N^1_{i}, \ldots,  N^m_{i} \in {\cal Y}$, such that 
\begin{align}\label{rig-eq: D2.12}
 \tfrac{1}{2} r^{\frac{i}{8}} |N^k_{i} | \le |N^k_{i} \setminus U_i| \le 2 r^{\frac{i}{8}}|N^k_{i} |
\end{align}
for all $k=1,\ldots,m$. This can be done in the following way: Let $V_0 = U = (\bigcup^n_{j=1} \overline{Q_j})^\circ$. First, to construct  an auxiliary set $\tilde{N}^1_{i}$ let $l_1$, $l_2$ be the smallest and largest index, respectively, such that $Q_l \subset U \setminus U_i$ and choose $l = l_1$ if $l_1< n-l_2$ and $l=l_2$ otherwise. Then add neighbors $Q_{l-1}, Q_{l+1} \subset U$, $Q_{l-2}, Q_{l+2} \subset U$, $\ldots$  until  $|\tilde{N}^1_{i} \setminus U_i| \le 2 r^{\frac{i}{8}}|\tilde{N}^1_{i} |$ holds. (I.e. the right inequality in \eqref{rig-eq: D2.12} is satisfied.) This is possible due to the fact that $|V_0 \setminus U_i| \le \bar{C}r^i|V_0| \le \frac{1}{2} r^{\frac{i}{8}}|V_0|$ by \eqref{rig-eq: r-i} for $r$ sufficiently small. Then note that also $ r^{\frac{i}{8}} |\tilde{N}^1_{i} | \le |\tilde{N}^1_{i} \setminus U_i|$ holds, in particular the left inequality in \eqref{rig-eq: D2.12} is fulfilled. We now define $V_1$ as the connected component of $V_0 \setminus \tilde{N}_i^1$ which is not completely contained in $U_i$. (If both are contained in $U_i$ we have finished.) We repeat the procedure on sets $V_j$ to define $\tilde{N}^j_{i}$, $1 \le j \le k$, satisfying \eqref{rig-eq: D2.12}, where $k$ is the smallest index such that $|V_k \setminus U_i| > \frac{1}{2}r^{\frac{i}{8}}|V_k|$. We now define $N_i^j = \tilde{N}_i^j$ for $j<k$ and $N_i^k := \tilde{N}_i^k \cup V_k$. 

It remains to show that also $N_i^k$ satisfies \eqref{rig-eq: D2.12}. Recall $|V_{k-1} \setminus U_i| \le \frac{1}{2}r^{\frac{i}{8}}|V_{k-1}|$. As due to the choice of $l$ and the fact that $ r^{\frac{i}{8}} |\tilde{N}^{k-1}_{i} | \le |\tilde{N}^{k-1}_{i} \setminus U_i|$ we have $|V_k| \ge \frac{1}{2}|V_{k-1}| - |\tilde{N}^{k-1}_i|$ and $|V_k \setminus U_i| = |V_{k-1} \setminus U_i| - |\tilde{N}^{k-1}_i \setminus U_i| \le \frac{1}{2}r^{\frac{i}{8}}|V_{k-1}| - r^{\frac{i}{8}}|\tilde{N}_i^{k-1}|$, we find $|V_{k} \setminus U_i| \le  r^{\frac{i}{8}}|V_{k}|$. This together with \eqref{rig-eq: D2.12} for $\tilde{N}_i^k$ implies the desired property for $N_i^k$.

Let $N_{i} = \bigcup^m_{k=1} N^k_{i}$. Similarly as in \eqref{rig-eq: D2.6} we find some $\frac{T}{4}  + 2 \le t \le \frac{T}{2}  - 1$ such that for $t_i = t +  \frac{9}{8} \cdot \frac{i}{2}$ we have $\sum_{Q \subset N_{i}} \#{\cal S}_{t_i}(Q) \le \upsilon^{-3}r^{2t_i + 3}$. Again set $\hat{\upsilon} = \upsilon^{2} r^{ - \frac{2t+1}{2}}$. Arguing as in \eqref{rig-eq: D2.2}  we can find a horizontal path $P_{i}$  consisting of $\hat{Q}(\xi_j) = \hat{Q}^{\hat{\upsilon} \vert \Gamma \vert_\infty}(\xi_j)$, $j=1,\ldots,n_{i}$,  and lying in $N_{i}$ such that \eqref{rig-eq: D1.4}, \eqref{rig-eq: PS} and \eqref{rig-eq: D2.2} are satisfied  replacing $t$ by $t_i$.   By \eqref{rig-eq: r-i} and \eqref{rig-eq: D2.12} we obtain 
\begin{align}\label{rig-eq: n length}
 \bar{C} r^{i+1}\hat{\upsilon}^{-1} \le n_{i} \le  C\bar{C}r^{ i-\frac{i}{8}}\hat{\upsilon}^{-1}.
\end{align}
Clearly, in general the path $P_{i}$ is not connected. Define $\hat{\cal S}^{t_i}_>(P_{i})$  and ${\cal K}_{i}$ similarly as in \eqref{rig-eq: D2.7}. By \eqref{rig-eq: D1.4} we get
$$\# \hat{\cal S}^{t_i}_>(P_{i})  \le C \hat{\upsilon} \vert \Gamma \vert_\infty (\upsilon^{4} r^{ - 2t_i-2}\vert \Gamma \vert_\infty)^{-1} \le C\hat{\upsilon}^{-1}  r^{ \frac{9}{8}i}r =  C\hat{\upsilon}^{-1}  r^{\frac{9}{8}i+1}.$$
Therefore, letting $\tilde{P}_{i} = \overline{\bigcup_{\hat{Q} \in {\cal K}_{i}} \hat{Q}}\subset P_{i}$ we find by \eqref{rig-eq: PS} (cf. \eqref{rig-eq: D2.5})
\begin{align}\label{rig-eq: D2.9}
n_{i} - \# {\cal K}_{i}\le C\hat{\upsilon}^{-1} r^{\frac{9}{8}i+1} \ \ \ \text{ and } \ \ \  |P_{i}| - |\tilde{P}_{i}|  \le Cr^{ \frac{9}{8}i+1}\hat{\upsilon} |\Gamma|^2_\infty.
\end{align}
We now repeat the slicing arguments above on each $N^k_{i}$ and obtain expressions similar to \eqref{rig-eq: D2.8}. As before, applying Lemma \ref{rig-lemma: weak rig}(i) we get  (cf. \eqref{rig-eq: D1.1*})
$$ |\tilde{P}_{i} \cap N^k_{i}|\dashint_{\tilde{P}_{i} \cap N^k_{i} \cap Z} |\bar{u} -  (A^k\,x +  b^k)|^2\, dx \le C(n^k)^{3} ({\cal E}(\tilde{P}_{i} \cap N^k_{i}) + \hat{\alpha}(P_{i} \cap N^k_{i}))^2 $$
for suitable $A^k \in \R^{2 \times 2}_{\rm skew}$, $b^k \in \R^2$ and $Z \in {\cal Y}'$, $Z \subset N^k_i$  (recall the definition in \eqref{rig-eq: shorthand}). Here $n^k$ denotes  the number of squares forming the path $P_{i} \cap N^k_{i}$, particularly $n_i = \sum^m_{k=1} n^k$.  We observe that by \eqref{rig-eq: D1.4} the estimate in \eqref{rig-eq: D1.6} can now be replaced by 
$$\upsilon^4 r^{-2t_i-1 \pm 1} = \hat{\upsilon}^2 r^{-{ \frac{9}{8}}i \pm 1}, \ \ \ \ \ \# \hat{\cal S}^{t_i}(P) \upsilon^{4} r^{-2t_i-2} \le { C}\hat{\upsilon} r.$$
Consequently, recalling  $n_{i} \le  C\bar{C}r^{ \frac{7}{8}i}\hat{\upsilon}^{-1}$ by \eqref{rig-eq: n length},  $t_i = t+ \frac{9}{8} \cdot \frac{i}{2}$ and  following the arguments in  \eqref{rig-eq: D1.1}, \eqref{rig-eq: D1.1*} and \eqref{rig-eq: D2.14} we obtain   for $Z \subset N_i$, $Z \in  {\cal Y}'$
\begin{align}\label{eq: XX}
\begin{split}
 \sum\nolimits_k  &  (n^k)^{-1}    |N^k_{i} \cap \tilde{P}_{i}|\dashint_{N^k_{i} \cap \tilde{P}_{i} \cap Z} |\bar{u}(x) - (A^k \, x + b^k)|^2 \, dx  \\ &\le  C  \sum\nolimits_k (n^k)^2({\cal E}(N^k_{i} \cap \tilde{P}_{i}) + \hat{\alpha}(N^k_{i} \cap P_{i}))^2 \le  C n^2_i ({\cal E}(\tilde{P}_{i}) + \hat{\alpha}( P_{i}))^2\\
&\le  C  \hat{\upsilon}   r^{ \frac{7}{4}i} \hat{\upsilon}^{-3} ({\cal E}(\tilde{P}_{i}) + \hat{\alpha}( P_{i}))^2 \le  C \hat{\upsilon} r^{ \frac{7}{4}i} \, r^{- \frac{9}{8}i}F \le C \hat{\upsilon} r^{ \frac{i}{2}}F,
\end{split}
\end{align}
where $F$ was defined in \eqref{rig-eq: D1.3}.  Observe that in the calculation the additional $r^{-\frac{9}{8}i}$ in front of $F$ occurs as in \eqref{rig-eq: D1.6} $\upsilon^4 r^{-2t-1 \pm 1}$ was replaced by $\upsilon^4 r^{-2t-1 \pm 1} r^{-\frac{9}{8}i}$. Define $J^k_{i} \subset J$ such that $|Q \cap (\tilde{P}_{i} \cap N^k_{i})| \ge  C\upsilon \hat{\upsilon} \vert \Gamma \vert^2_\infty \ge C \frac{\hat{\upsilon}}{\upsilon}|Q|$  for $Q \in J^k_i$ and set $\hat{N}^k_{i} = \bigcup_{Q \in J^k_{i}} Q$.  Assume $\hat{N}^k_{i} \cap Z \neq \emptyset$ which implies $|\hat{N}^k_{i} \cap Z| \ge \upsilon^2  |\Gamma|^2_\infty$. Observe $\hat{\upsilon} \ge \upsilon^{\frac{7}{4}}r^{- \frac{3}{2}} \ge \upsilon^2 r^{-\frac{i}{4}{ -1}}$ by \eqref{rig-eq: D1.7} and the fact that $i \le T+{2}$. As $|N^k_i| (n^k)^{-1} \le C\hat{\upsilon} \upsilon|\Gamma|^2_\infty$, we find by \eqref{rig-eq: D1.5}
\begin{align*}
\sum\nolimits_k  &(n^k)^{-1} |\hat{N}^k_{i}|\dashint_{\hat{N}^k_{i} \cap Z} |\bar{u}(x) - (\bar{A} \, x + \bar{b})|^2 \, dx \\
& \le C\hat{\upsilon}\upsilon^{-1} \int_{U} |\bar{u}(x) - (\bar{A} \, x + \bar{b})|^2 \, dx \le C\hat{\upsilon}\upsilon^{-1} \ \upsilon^{3} (\hat{\upsilon}  r)^{-1} F \le C  \hat{\upsilon} r^{\frac{i}{4}} F.
\end{align*}
Again arguing as in \eqref{rig-eq: P to N},  in particular applying Lemma \ref{lemma: diff1}(ii), we derive by the triangle inequality and \eqref{eq: XX}
\begin{align}\label{rig-eq: D2.11}
\sum\nolimits_k  (n^k)^{-1} |\hat{N}^k_{i}|  \dashint_{Z \cap\hat{N}^k_{i} } |\bar{u}(x) - (A^k \, x + b^k)|^2 \, dx \le  C \hat{\upsilon} r^{\frac{i}{4}} \ \hat{\upsilon}^{-1} \upsilon F  = C r^{\frac{i}{4}}  \upsilon F. 
\end{align}
We set $U^k_{i} = \hat{N}^k_{i}$ if  $|N^k_{i} \setminus \hat{N}^k_{i}| \le r^{\frac{i}{8}}|N^k_{i}|$ and $U^k_{i} = \emptyset$ else for all $k=1,\ldots,m$. We now estimate the difference between $A,b$ given in \eqref{rig-eq: D2.10} and $A^k, b^k$ for $k=1,\ldots,m$.
Consider $U^k_i$ such that  $U^k_i = \hat{N}^k_{i}$.  Then $|U_i^k| \ge(1-r^{\frac{i}{8}})|N^k_i|$ and thus by \eqref{rig-eq: D2.12} we have $$|U^k_{i} \cap U_i| \ge |N^k_{i} \cap U_i| - |N^k_{i} \setminus \hat{N}^k_{i}| \ge (1-Cr^{\frac{i}{8}})|N^k_i| \ge (1-Cr^{\frac{i}{8}})|U_i^k|$$
for $r$ sufficiently small and some $C>0$.  We are now in the position to apply Lemma \ref{rig-lemma: weak rig4}  for $B_2 = U_i^k$, $B_1 = U_i^k \cap U_i$ and $s= \frac{\tau}{2}$, $\delta = Cr^{\frac{i}{8}}$, where we observe that for $C$ large enough we have  ${\rm diam}(U^k_i) \ge \frac{s}{\delta}$ by \eqref{rig-eq: D2.12} and  ${\rm diam}(U^k_i \cap U_i) \ge (1-\delta) {\rm diam} (N^k_{i}) \ge (1-\delta){\rm diam} (U^k_{i})$ by the previous estimate and the fact that $N^k_i$ is connected.  Set $\bar{C}_i = C\prod\nolimits^{i-1}_{j=0} \big(1+ \bar{C} r^{\frac{j}{8}}\big)$ (cf. \eqref{rig-eq: D2.10}). Using \eqref{rig-eq: D2.10} and \eqref{rig-eq: D2.11}, in particular recalling that the sets $(\hat{N}_i^k)_k$ are pairwise disjoint, we find for $Z \subset U$, $Z \in  {\cal Y}'$
\begin{align*}
&\Vert \bar{u} - (A\, \cdot + b)\Vert^2_{L^2(U_i^k \cap U_i \cap Z)} \le |U_i^k \cap U_i \cap Z||U_i^k|^{-1} H^k_1, \\
 & \Vert \bar{u} - (A^k\, \cdot + b^k)\Vert^2_{L^2(U_i^k \cap Z)} \le |U_i^k \cap Z||U_i^k|^{-1} H^k_2,
\end{align*}
where $H_1^k = |U_i^k||U|^{-1} \bar{C}_iG$ and $H_2^k = C n^kr^{\frac{i}{4}}\upsilon F$. Therefore, Lemma \ref{rig-lemma: weak rig4} yields 
\begin{align}\label{rig-eq: D2.11.2}
\begin{split}
\Vert \bar{u} - (A \, \cdot + b)\Vert^2_{L^2(U^k_{i} \cap Z)} &\le  |U^k_{i} \cap Z||U^k_{i}|^{-1} (1+Cr^{\frac{i}{8}})|U_i^k||U|^{-1} \bar{C}_i G \\
& \ \ \ +  |U^k_{i} \cap Z||U^k_{i}|^{-1}Cr^{-\frac{i}{8}} n^kr^{\frac{i}{4}}\upsilon F.
\end{split}
\end{align}
For shorthand we write $U^* = \big(\bigcup^m_{k=1} \overline{U^k_{i}}\big)^\circ$ and define $U_{i+1} =  (\overline{U_i} \cup \overline{U^*})^\circ$.
 We recall $N_i = \bigcup_{k=1}^m N_i^k$ as constructed in \eqref{rig-eq: D2.12}. We claim
\begin{align}\label{rig-eq: r-i.2}
|N_i \setminus U^*| \le Cr^{i+1}|U|
\end{align}
and postpone the proof of this assertion to the end of the proof.  Then \eqref{rig-eq: n length} for  $\bar{C}$ sufficiently large implies $|U^*| \ge  |N_i| - Cr^{i+1}|U| \ge  cn_{i} \hat{\upsilon} |U| $.   As for $U^k_i \neq \emptyset$ we have $| N^k_{i}| \le  (1-r^{\frac{i}{8}})^{-1}|U^k_{i}|$, it is not hard to see that $\frac{|U^*|}{|U^k_{i}|}\le C \frac{n_{i}} {n^k}$  and thus $n^k|U_i^k|^{-1} \le C|U^*|^{-1} n_i \le C|U|^{-1} \hat{\upsilon}^{-1}$.  Let $Z \subset U_{i+1}$, $Z \in  {\cal Y}'$. Now by \eqref{rig-eq: D2.11.2},  the  fact that the sets $U^k_i$ are pairwise disjoint and $F \le C \hat{\upsilon} \upsilon ^{-1} G$ we derive
\begin{align*}
\sum\nolimits_{k=1}^m\Vert \bar{u} - (A \, \cdot + b)\Vert^2_{L^2(U^k_{i} \cap Z)} &\le  |U^* \cap Z||U|^{-1}\big( \bar{C}_i(1+Cr^{\frac{i}{8}}) G + C r^{\frac{i}{8}} G \big) \\ 
& \le |U^* \cap Z||U|^{-1} \bar{C}_{i+1} G.
\end{align*}
The last estimate follows for $\bar{C}$ sufficiently large.  By \eqref{rig-eq: D2.10} we now conclude for $Z \subset U_{i+1}$ 
\begin{align*}
\Vert \bar{u} - (A \, \cdot + b)\Vert^2_{L^2(Z)} &  = \Vert \bar{u} - (A \, \cdot + b)\Vert^2_{L^2(Z \setminus U^*)} + \sum\nolimits_k\Vert \bar{u}- (A \, \cdot + b)\Vert^2_{L^2(U^k_{i}\cap Z)}\\
& \le |Z \setminus U^*||U|^{-1} \bar{C}_iG +  |U^*\cap Z| |U|^{-1} \bar{C}_{i+1}G \\
& \le | Z||U|^{-1} \bar{C}_{i+1}G.
\end{align*}
 This yields \eqref{rig-eq: D2.10}. To see \eqref{rig-eq: r-i} for $i+1$, we apply  \eqref{rig-eq: r-i.2} to obtain  $|U \setminus U_{i+1}| \le |(U \setminus U_i)\setminus N_i| + |N_i \setminus U^*| \le 0 + \bar{C}r^{i+1}|U|$. Here we used that  $U \setminus U_i \subset N_i$. 

Finally, we choose $i_* \le T+ 2$ large enough such that $|U \setminus U_{i_*}| \le \bar{C}r \upsilon |U| \ll (\upsilon\vert\Gamma\vert_\infty)^2$ for $r$ sufficiently small which implies $U_{i_*} = U$. Consequently, thanks to \eqref{rig-eq: D2.10},  \eqref{rig-eq: D2.3}  holds.

 To finish the proof in the case $|U|\ge C'\upsilon|\Gamma|_\infty^2$, it remains to show \eqref{rig-eq: r-i.2}. First, by \eqref{rig-eq: D2.9} and the construction of $\hat{N}^k_{i} $ we have $|\bigcup_k N^k_{i} \setminus \bigcup_k \hat{N}^k_{i}| = \sum_k |N^k_{i} \setminus\hat{N}^k_{i}|  \le Cr^{\frac{9}{8}i+1}\upsilon|\Gamma|_\infty^2 \le Cr^{ \frac{9}{8}i+1}|U|$. Therefore, it suffices to prove 
\begin{align}\label{rig-eq: D2.13} 
\sum\nolimits_k |\hat{N}^k_{i} \setminus U^k_{i}| \le  r^{-\frac{i}{8}}\sum\nolimits_k |N^k_{i} \setminus \hat{N}^k_{i}| 
\end{align}
as then we conclude $|N_i \setminus U^*|\le \sum_k| N^k_{i} \setminus \hat{N}^k_{i}|  + \sum_k| \hat{N}^k_{i} \setminus U^k_i| \le Cr^{i+1}|U|$.

To see \eqref{rig-eq: D2.13}  we  observe that if $\hat{N}^k_i \neq U^k$, then $|\hat{N}^k_i| \le |N^k_i| < r^{-\frac{i}{8}} |N^k_i \setminus \hat{N}^k_i|$.  Consequently, we calculate $\sum\nolimits_k |\hat{N}^k_i \setminus U^k_i| = \sum_{k: U^k_i =\emptyset} |\hat{N}^k_i| \le r^{-\frac{i}{8}}  \sum_{k: U^k_i =\emptyset} |N^k_i \setminus \hat{N}^k_i|  \le r^{-\frac{i}{8}} \sum_{k} |N^k_i \setminus \hat{N}^k_i| $, as desired.

It remains to indicate the adaptions for the (general) case that $U \in {\cal Y}$ is of arbitrary size. We choose $0 \le i_0 \le T+2$ such that $Cr^{i_0+1} \upsilon |\Gamma|_\infty^2 < |U| \le Cr^{i_0} \upsilon |\Gamma|_\infty^2$ and begin the induction in  \eqref{rig-eq: r-i}, \eqref{rig-eq: D2.10} not for $i = 1$, but for $i = i_0$. For the first step $i = i_0$ we do not apply Lemma \ref{rig-lemma: D1,5}, but follow the lines of the above proof for one single set $N^1_{i_0} = U$.   \eop

We now drop the assumption that $U \in {\cal Y}$ does not intersect a corner of $\Gamma$.

\begin{corollary}\label{rig-cor: D2}
Let $r \ge \upsilon >0$   as in \eqref{eq: parameters} small enough. Then for all  $U \in {\cal Y}$, $U \subset N_{j,\pm}$, $j=1,2$,  with   $ \Vert U\Vert_\pi \le \frac{19}{20} \tau$ the estimate \eqref{rig-eq: D2.3} holds.
\end{corollary}

\Proof  Assume without restriction $E_{+,+} \subset U$ and define $U' = U \setminus E_{+,+}$. Using Lemma \ref{rig-lemma: D2}   we find
$$|U'|\dashint_{Z}  |\bar{u}(x)- (A\,x +  b)|^2 \, dx  \le CG $$
for $Z \subset U'$, $Z \in  {\cal Y}'$, where $G:=(1 + C_* r)  \frac{\eps}{\upsilon^3} \vert \Gamma \vert^3_\infty$. Let $Q \in J$, $Q \subset U'$ such that $\partial Q \cap E_{+,+} \neq \emptyset$. Setting $Z=Q$ in the above inequality and arguing as in \eqref{rig-eq: D1.5} we find $\hat{A} \in \R^{2 \times 2}_{\rm skew}$, $\hat{b} \in \R^2$ such that
$$\int_Q |\bar{u}(x)- (A\,x +  b)|^2 \, dx \le C \upsilon G, \ \ \  \int_{Q\cup E_{+,+}} |\bar{u}(x)- (\hat{A}\,x +  \hat{b})|^2 \, dx \le C \upsilon^2  r^{-1} G.$$
Applying  the triangle inequality and Lemma \ref{lemma: diff1}(ii) on $D_1 = Q$ and $D_2 = Q \cup E_{+,+}$ we find $\Vert \bar{u}- (A\,\cdot +  b)\Vert^2_{L^2(Q\cup E_{+,+})}  \le C \upsilon G$ as $\upsilon \le r$. Now it is not hard to see that \eqref{rig-eq: D2.3} is satisfied. \eop

\subsection{Step 3: Neighborhood with small projection of components}\label{rig-sec: subsub, proj-bound2}

 Recall the covering ${\cal C}$ of the neighborhood $N$ introduced in \eqref{rig-eq: covering def} (see also Figure \ref{B}).  We now treat the case that $\Vert U \Vert_\pi$ is small for all $U \in {\cal C}$. It is essential that adjacent elements of the covering overlap sufficiently. Therefore, we introduce another covering $\hat{\cal C}$ as follows. First assume $l_2 \ge \frac{l_1}{2}$. If some $U \in {\cal C}$ intersects only one of the four sets $N_{j,\pm}$, $j=1,2$, we let $U \in \hat{\cal C}$. Then eight sets $U^1_{\pm,\pm}, U^2_{\pm,\pm}$ remain where $U^i_{\pm,\pm} \cap E_{\pm,\pm} \neq \emptyset$ and $U_{\pm,\pm}^i \subset N_{i,-} \cup N_{i,+}$ for $i=1,2$. As before $E_{\pm,\pm}$ denote the sets at the corners of $\Gamma$.  Add the four sets $U^1_{\pm,\pm} \cup U^2_{\pm,\pm}$ to $\hat{\cal C}$. If $l_2 < \frac{l_1}{2}$ we proceed likewise with the only difference that  instead of $U^1_{\pm,\pm} \cup U^2_{\pm,\pm}$ we add the two sets $U^2_{k,+} \cup U^1_{k,+} \cup U^1_{k,-} \cup U^2_{k,-}$, $k=+,-$, to $\hat{\cal C}$. (Note that by definition of ${\cal C}$ we have $U^1_{\pm,+} = U^1_{\pm,-}  = N_{1,\pm}$ in this case.)

\begin{lemma}\label{rig-lemma: D3}
Theorem \ref{rig-theorem: D} holds under the additional assumption that $ \Vert U\Vert_\pi \le \frac{19}{20} \tau$ for all $U \in {\cal C}$. 
\end{lemma}

\Proof It suffices to show that for all $U \in \hat{\cal C}$ there are $A_U \in \R^{2 \times 2}_{\rm skew}$,  $b_U \in \R^2$ such that
\begin{align}\label{rig-eq: D2.3.2} 
|U|\dashint_{Z}  |\bar{u}(x)- (A_U\,x +  b_U)|^2 \, dx  \le C(1+rC_*)  \upsilon^{-3}\vert \Gamma \vert^3_\infty  \eps 
\end{align}
holds for all $Z \subset U$, $Z \in  {\cal Y}'$. Indeed, the desired result then follows from the construction of the covering $\hat{\cal C}$ and  Lemma \ref{lemma: diff1}:  Write $\hat{\cal C} = \lbrace U_1, \ldots, U_n \rbrace$ with $U_{i-1} \cap U_{i} \neq \emptyset$ for all $i=1,\ldots,n$, where $U_0 = U_n$.

As $|U_i \cap U_{i+1}|\ge c\upsilon|\Gamma|_\infty^2$ and ${\rm diam}(U_i \cap U_{i+1}) \ge c|\Gamma|_\infty$ for $c$ small enough, \eqref{rig-eq: D2.3.2}, the triangle inequality and Lemma \ref{lemma: diff1}(ii) for $D_1 = U_i \cap U_{i+1}$ and $D_2 = N$  yield for all $Z \subset N$, $Z \in  {\cal Y}'$,
$$
|N|\dashint_{Z}  |(A_{U_i}\,x +  b_{U_i}) - (A_{U_{i+1}}\,x +  b_{U_{i+1}})|^2 \, dx  \le C(1+rC_*)  \upsilon^{-3}\vert \Gamma \vert^3_\infty  \eps. 
$$
As $\# \hat{\cal C}$ is uniformly bounded, the triangle inequality yields that the last estimate holds for two mappings $A_{U_i}, b_{U_i}$ and $A_{U_j}, b_{U_j}$ for $1 \le i,j \le n$ for a larger constant $C$. Then setting $A = A_{U_1}$, $b=b_{U_1}$,  using  \eqref{rig-eq: D2.3.2}  and again the triangle inequality we get 
\begin{align}\label{rig-eq: D2.3_subset} 
|N|\dashint_{Z}  |\bar{u}(x)- (A\,x +  b)|^2 \, dx  \le C(1+rC_*)  \upsilon^{-3}\vert \Gamma \vert^3_\infty  \eps 
\end{align}
for all $Z \subset N$, $Z \in  {\cal Y}'$.

It remains to establish \eqref{rig-eq: D2.3.2} for $U \in \hat{\cal C}$. By assumption and Lemma \ref{rig-lemma: D2} the assertion is clear if $U \cap E_{\pm,\pm} = \emptyset$ as then particularly $U \in {\cal C}$. Therefore, we first let $l_2 \ge \frac{l_1}{2}$ and assume that e.g.  $U \cap E_{+,+} \neq \emptyset$. The necessary changes for the case $l_2 \le \frac{l_1}{2}$ are indicated at the end of the proof. 

As in \eqref{rig-eq: D2.6} we find  $\frac{T}{4} +  2 \le t\le  \frac{T}{2}-1$ such that $\sum_{Q \subset U} \#{\cal S}_t(Q) \le \upsilon^{-3}r^{2t+3}$. Again let $\hat{\upsilon} = \upsilon^2 r^{- \frac{2t +1 }{2}}$. As before, the main strategy will be to construct a suitable path in $U$. Let $(\Gamma_l)_l$ be the boundary components  such that the corresponding rectangles $(\partial R_l)_l$ given by \eqref{rig-eq: W prop}(i) and \eqref{rig-eq: W prop}(v), respectively,  satisfy $\partial R_l \cap U \neq \emptyset$ and  $|\partial R_l|_\pi \neq |\partial R_l|_\infty$. Let $V_l \subset \overline{N}$ be the smallest rectangle containing $R_l \cap N$ and $(l_1 +\tau,l_2+\tau)$. We partition $(V_l)_l$ into ${\cal V}_1$ and ${\cal V}_2$ depending on whether $|\pi_1 V_l| \le |\pi_2 V_l| $ or $|\pi_1 V_l| > |\pi_2 V_l|$. Recalling \eqref{rig-eq: VertP} it is not hard to see that $|\pi_j V_l|  = \vert R_l \vert_\pi$ for $V_l \in {\cal V}_j$ for $j=1,2$. Let $d_j = \inf\lbrace s \in \R: s \in  \pi_jV_l \text{ for a } V_l \in {\cal V}_j \rbrace$  and define the stripes 
$$D_1 = (-\infty,d_1) \times (-\infty,d_2) \cap N_{1,+} \cap U, \ \ \  D_2 = (-\infty,d_1) \times (-\infty, d_2)\cap  N_{2,+} \cap U.$$

 \begin{figure}[H]
\centering
\begin{overpic}[width=.6\linewidth,clip]{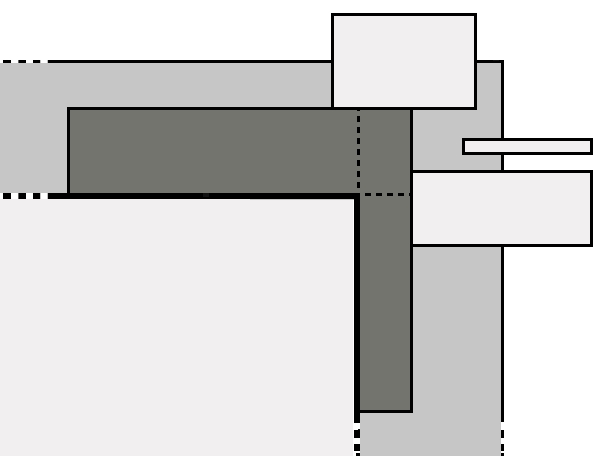}

\put(75,44){\small{$X$}}
\put(175,44){\small{$N(\Gamma)$}}
\put(88,118){\small{$D_2$}}

\put(251,94){\small{$\partial R_{l_1}$}}
\put(242,97){\line(1,0){6}}

\put(251,125){\small{$\partial R_{l_2}$}}
\put(242,127){\line(1,0){6}}

\put(203,175){\small{$\partial R_{l_3}$}}
\put(194,177){\line(1,0){6}}

\put(151,68){\small{$D_1$}}

\end{overpic}
\caption{\small Sketch of a part of $N(\Gamma)$ containing $U^1_{+,+} \cup U^2_{+,+}$. The sets $D_1$, $D_2$ are highlighted in grey.} \label{I}
\end{figure}
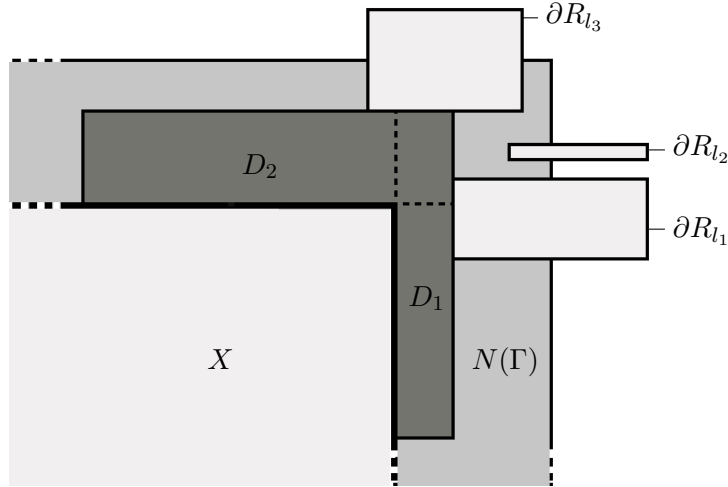

As by assumption $\Vert U \Vert_\pi \le \frac{19}{20} \tau$ for all $U \in {\cal C}$, we find sets $S_j \subset (l_j, d_j)$ with $|S_j| \ge \frac{\tau}{20}$ such that the stripes $\hat{D}_1 = D_1 \cap (S_1 \times \R)  \in {\cal U}^s$ and $\hat{D}_2 = D_2 \cap (\R \times S_2)  \in {\cal U}^s$ satisfy $\partial W_i \cap \hat{D}_j = \emptyset$ for $j=1,2$.  Moreover, observe that $|d_j - l_j| \ge \frac{\tau}{20}$ for $j=1,2$. We cover $D_1$ by vertical paths ${\cal P}_1 = (P^1_i)_i$, $i=1,\ldots,k_1$, and $D_2$ by horizontal paths ${\cal P}_2 = (P^2_i)_i$, $i=1,\ldots,k_2$, consisting of squares $\hat{Q}^{\hat{\upsilon}|\Gamma|_\infty}(\xi) = \hat{Q}(\xi)$, i.e. $k_j = \lceil( 2\hat{\upsilon}\vert \Gamma\vert_\infty)^{-1}(d_j - l_j)\rceil$.  As in \eqref{rig-eq: D2.2} we find that there are subsets $\hat{\cal P}_j \subset {\cal P}_j$ with $\# \hat{\cal P}_j \ge ck_j \ge c\upsilon\hat{\upsilon}^{-1}$ for $c>0$ sufficiently small such that \eqref{rig-eq: PS} and \eqref{rig-eq: D2.2} hold for all $P^j_i \in \hat{\cal P}_j$, $j=1,2$. We can now choose $P^j \in \hat{\cal P}_j$, $j=1,2$, such that \eqref{rig-eq: D1.4} is satisfied possibly passing to a larger constant. Moreover, this can be done in such a way that $ Q^* := P^1 \cap P^2$  satisfies 
\begin{align}\label{rig-eq: D3.1}
\begin{split}
&\sum\nolimits_{\Gamma_k \cap Q^* \neq \emptyset} \vert \Gamma_k \vert_\infty \le \tilde{C}\upsilon^{-1} \hat{\upsilon}^{2} \vert\Gamma\vert_\infty,  \\ &Q^* \cap \Gamma_k = \emptyset \ \text{ for all } \Gamma_k: |\Gamma_k| \ge \tilde{C}\hat{\upsilon}^2\upsilon^{-1}|\Gamma|_\infty, 
\end{split}
\end{align}
for $\tilde{C}>0$ sufficiently large. To see the latter, note that we have $\sim \bar{\tau}^2 (\hat{\upsilon}\vert\Gamma\vert_\infty)^{-2} = \upsilon^2 \hat{\upsilon}^{-2}$ possibilities to combine paths in $\hat{\cal P}_1$, $\hat{\cal P}_2$ such that \eqref{rig-eq: D1.4} holds.  Moreover, we also have $\sim \bar{\tau}^2 (\hat{\upsilon}\vert\Gamma\vert_\infty)^{-2} = \upsilon^2 \hat{\upsilon}^{-2}$ possibilities to combine paths in $\hat{\cal P}_1$, $\hat{\cal P}_2$ such that $Q^*$  additionally has empty intersection with all $\Gamma_k$ satisfying $\vert\Gamma_k\vert_\infty \ge  \tilde{C}\hat{\upsilon}^2\upsilon^{-1}\vert\Gamma\vert_\infty$.  This follows from \eqref{rig-eq: PS} and the fact that  by \eqref{rig-eq: property,b} we derive  
$$\# \lbrace \hat{Q}: \exists \Gamma_k: \tilde{C}\hat{\upsilon}^2\upsilon^{-1}|\Gamma|_\infty \le |\Gamma_k|_\infty \le \bar{C}^2\hat{\upsilon}|\Gamma|_\infty, \Gamma_k \cap \hat{Q} \neq \emptyset \rbrace \le C\tilde{C}^{-1}\upsilon^2 \hat{\upsilon}^{-2}$$ 
for $C=C(h_*,q,\bar{C})$ large enough with $\bar{C}$ from \eqref{rig-eq: PS}. Since all other components $\Gamma_k$ intersect at most four adjacent squares $\hat{Q}$, using again \eqref{rig-eq: property,b} we can select $Q^*$ such that also $\sum\nolimits_{\Gamma_k \cap Q^* \neq \emptyset} \vert \Gamma_k \vert_\infty \le C\upsilon\vert\Gamma\vert_\infty \upsilon^{-2} \hat{\upsilon}^{2}$ holds.

Let $P = \hat{P}^1 \cup Q^* \cup \hat{P}^2$, where $\hat{P}^j$, $j=1,2$, is the connected component of $P^j \setminus  Q^*$ not completely contained in $E_{+,+}$. Denote the midpoints of the squares in $P$ by $(\xi_1, \ldots, \xi_n)$. Recall the definition of $\hat{\cal S}^t_>$ in \eqref{rig-eq: D2.7} and let 
$${\cal K} = \lbrace \hat{Q} = \hat{Q}(\xi_j): \hat{Q} \cap \Gamma_l = \emptyset \ \text{ for all } \ \Gamma_l \in \hat{\cal S}^t_>(P) \rbrace \cup \lbrace Q^*\rbrace.$$
Consider two sets $\hat{Q}(\xi_0), \hat{Q}(\xi_m) \in {\cal K
}$ and the path $(\xi_0, \ldots, \xi_m)$. We can repeat the slicing method of the previous proofs and end up with an estimate of the form  (cf. \eqref{rig-eq: D2.8})
$$ \Vert b_0 - b_m + (A_0 - A_m)\, \cdot \Vert_{L^2(\hat{Q}(\xi_0) \cup \hat{Q}(\xi_m))}   \le C\hat{\upsilon}^{-1} \Big({\cal E}(\hat{Q}(\xi_0)) + {\cal E}(\hat{Q}( \xi_m)) + {\cal E}(\hat{Q}^*) + \hat{\alpha}(D)\Big) $$
for suitable $A_0,A_m \in \R^{2 \times 2}_{\rm skew}$, $b_0,b_m \in \R^2$, where $D = \bigcup^m_{j=0} \hat{Q}(\xi_j)$. In fact, if $\hat{Q}(\xi_0), \hat{Q}(\xi_m) \subset \hat{P}^j$ for some $j=1,2$, this follows immediately. Otherwise, we apply the arguments leading to \eqref{rig-eq: D2.1} on each pair $\hat{Q}(\xi_0), Q^*$ and $Q^*, \hat{Q}(\xi_m)$ and employ the triangle inequality.

Defining $ \tilde{P} = \overline{\bigcup_{\hat{Q} \in {\cal K}} \hat{Q}}$ and arguing as in \eqref{rig-eq: D1.1*},  we then obtain   for $Z \subset U$, $Z \in {\cal Y}'$
\begin{align*}
 \Vert \bar{u} -  (A\, \cdot +  b)\Vert^2_{L^2(\tilde{P} \cap Z)}  & \le C  |\tilde{P} \cap Z||\tilde{P}|^{-1} \hat{\upsilon}^{-3} (({\cal E}(\tilde{P}))^2 + (\alpha(P))^2) \\ &
\le C  |\tilde{P} \cap Z||\tilde{P}|^{-1} (1 + C_* r)   \frac{\eps \hat{\upsilon}}{\upsilon^4}  \vert \Gamma \vert^3_\infty + C\frac{1}{\hat{\upsilon}^3} ({\cal E}(\hat{Q}^*))^2 
\end{align*}
for some $A \in \R^{2 \times 2}_{\rm skew}$ and $b \in \R^2$. In the last step we proceeded as in \eqref{rig-eq: D2.14} (see also \eqref{rig-eq: D1.1}), observing that the paths $\tilde{P}$ defined here and in the proof of Lemma \ref{rig-lemma: D1,5} differ essentially by the square $Q^*$. By \eqref{rig-eq: D1.7} we get $\hat{\upsilon}^3 \upsilon^{-3} \le \hat{\upsilon}r$ and using \eqref{rig-eq: D3.1} as well as  \eqref{rig-eq: property,c} we derive  (cf. \eqref{rig-eq: D1.1})
\begin{align}\label{rig-eq: D3.2}
\begin{split}
\hat{\upsilon}^{-3} ({\cal E}(\hat{Q}^*))^2 &\le C\hat{\upsilon}^{-1} |\Gamma|_\infty^3\eps + C\hat{\upsilon}^{-3} C_* \frac{\eps}{\upsilon^4} \big(\sum\nolimits_{\Gamma_l \cap Q^* \neq \emptyset} (\vert  \Gamma_l \vert_\infty)^{3/2} \big)^2 \\
& \le  C\hat{\upsilon}^{-1} |\Gamma|_\infty^3\eps  + C|\Gamma|_\infty\hat{\upsilon}^{-1}\upsilon^{-1} C_* \frac{\eps}{\upsilon^4} \big(\sum\nolimits_{\Gamma_l \cap Q^* \neq \emptyset} \vert \Gamma_l \vert_\infty \big)^2 \\
&\le C(1 + C_* r) \, \hat{\upsilon}  \frac{\eps}{\upsilon^4} \vert \Gamma \vert^3_\infty.
\end{split}
\end{align}
Consequently, we have re-derived \eqref{rig-eq: D2.14}. Proceeding as in Lemma \ref{rig-lemma: D1,5} we obtain a set $U' \subset U$ with $|U \setminus U'|\le Cr|U|$ such that 
$$|U|\dashint_{Z } |\bar{u}(x)- (A\,x + b)|^2 \, dx \le C (1 + C_* r)  \frac{\eps}{\upsilon^3} \vert \Gamma \vert^3_\infty$$
for $Z \subset U'$, $Z \in  {\cal Y}'$. On the other hand, by Corollary \ref{rig-cor: D2} we find $A_j \in \R^{2 \times 2}_{\rm skew}$, $b_j \in \R^2$ for $j=1,2$ such that
\begin{align}\label{rig-eq: new237}
|U \cap N_{j,+}|\dashint_{Z\cap N_{j,+}} |\bar{u}(x)- (A_j\,x + b_j)|^2 \, dx \le C (1 + C_* r)  \frac{\eps}{\upsilon^3} \vert \Gamma \vert^3_\infty
\end{align}
for all $Z \subset U$,  $Z \in  {\cal Y}'$. Now  \eqref{rig-eq: D2.3.2}   follows by applying  the triangle inequality and Lemma \ref{lemma: diff1}(ii) on $D_1 = U' \cap N_{j,+}$ and $D_2 = U \cap N_{j,+}$.

Finally, the essential difference in the treatment of the case $l_2 \le \frac{l_1}{2}$ is that in the construction of the path $P$ one has to choose two sets $Q^*_1$,$Q^*_2$ where the path changes its direction. Following the above arguments we observe that these sets can be selected so that the required conditions,  in particular \eqref{rig-eq: D3.1},  are satisfied.    \eop

\subsection{Step 4: General case}\label{rig-sec: subsub,  gen-case}

We are eventually in a position to give the proof of Theorem \ref{rig-theorem: D}. We briefly remark that the following proof crucially depends on the trace theorem established in Lemma \ref{rig-sec: sub, trace} and the fact that there are at most two large cracks in the neighborhood of $\Gamma$.

\noindent {\em Proof of Theorem \ref{rig-theorem: D}.} 
In the general situation we possibly have $N \neq \tilde{N}$  with $\tilde{N} = N \setminus (X_1 \cup X_2)$ as introduced at the beginning of Section \ref{rig-sec: subsub,  pre-est}.  Let $\hat{\cal C}$ be the covering considered  in Lemma \ref{rig-lemma: D3}. Let $K_1$, $K_2$ with $\dist(K_1,K_2)\ge c|\Gamma|_\infty$ be the sets given by Lemma \ref{rig-lemma: slicing bound2} and let $\tilde{\cal C}$ be the covering of $N\setminus (K_1 \cup K_2)$ consisting of the connected components of the sets $U\setminus (K_1 \cup K_2)$, $U \in \hat{\cal C}$.  To simplify the exposition we prefer to present first a special case where $K_1,K_2$ have the form $K_- := K_1=(-\tau-l_1,-l_1) \times (-\tau,\tau)$ and $K_+ := K_2 = (l_1, l_1 + \tau) \times (-\tau,\tau)$. Moreover, we suppose that the sets $\Psi^\pm$ associated to boundary components larger than $\hat{\tau}$  -- if they exist at all -- have the form $\Psi^\pm = \Psi^\pm_1 \cup \Psi^\pm_2 \cup \Psi^\pm_3$, where  $\Psi^\pm_1 = (\pm l_1, \pm(l_1 + \tau)) \times (\psi^\pm,2\tau)$, $\Psi^\pm_2 = (\pm l_1, \pm(l_1 + \psi^\pm)) \times (-\psi^\pm,\psi^\pm)$ and  $\Psi^\pm_3 = (\pm l_1, \pm(l_1 + \tau)) \times (-2\tau, -\psi^\pm)$. Here $\psi^\pm$ denote the corresponding values to $\Psi^\pm$ (see Section \ref{rig-sec: subsub,  dod-neigh}).

If $\Psi^-$ or $\Psi^+$ do not exist, we set $\Psi^-_2 = K_1$,  $\Psi^+_2 = K_2$, respectively, and let $\Psi^\pm_j$, $j=1,3$, be the adjacent squares. In addition, we then define $\psi^\pm = \tau$. We will treat both cases simultaneously in the following. 

This special case already covers the fundamental ideas of the proof as the arguments essentially rely  on the property that $\dist(K_1,K_2)\ge c|\Gamma|_\infty$  and the fact that the shapes of all sets are comparable (through homeomorphisms with constants depending on $h_*$) to  squares. We will indicate the necessary adaptions for the general case at the end of the proof.

Let $N'_\pm = N \cap \lbrace  \pm x_2 \ge 0 \rbrace \setminus  (K_1 \cup K_2)$. By Lemma \ref{rig-lemma: slicing bound2} the assumptions of Lemma \ref{rig-lemma: D3} are satisfied on each set $N'_+$ and $N'_-$. Consequently, there are $A_\pm \in \R^{2 \times 2}_{\rm skew}$ and $b_{\pm} \in \R^2$ such that for all $V_{\pm} \subset N'_\pm$, $V_{\pm} \in {\cal Y}$, one has
\begin{align}\label{rig-eq: D5.2}
|N'_\pm|\dashint_{V_\pm} |\bar{u}(x) - (A_\pm\,x -b_\pm)|^2 \,dx \le  C(1 + C_* r) \frac{\eps}{\upsilon^3} \vert \Gamma \vert^3_\infty =:G
\end{align}
by \eqref{rig-eq: D2.3_subset}. We let $\Xi_0^+ = (l_1,l_1 + \psi^+) \times \lbrace 0 \rbrace$, $\Xi^-_0 = (-l_1-\psi^-,-l_1) \times \lbrace 0 \rbrace $ and without restriction (possibly after a small translation in $\e_2$-direction) we can assume ${\cal H}^1( \Xi^\pm_0 \cap \partial W_i)= 0$. The goal is to show
\begin{align}\label{rig-eq: D5.5}
\int_{\Xi^\pm_0}|(A_+ - A_-)\,x + (b_+ -b_-)|^2 \, d{\cal H}^1(x) \le C \frac{\psi^{\pm}}{\upsilon \vert \Gamma \vert^2_\infty} G.
\end{align}
We prove this only for $\Xi^-_0$.  As a preparation let $\tilde{\Psi}^-_1 = (-\frac{\tau}{2} -l_1 ,-l_1 ) \times (\psi^-, \frac{3}{2}\tau)$,  $\tilde{\Psi}^-_3 = (-\frac{\tau}{2} -l_1 ,-l_1 ) \times (- \frac{3}{2}\tau,-\psi^-)$ and $\tilde{\Psi}^- = \big(\overline{\tilde{\Psi}^-_1} \cup \overline{\Psi^-_2} \cup \overline{\tilde{\Psi}^-_3}\big)^\circ$. We observe 
\begin{align}\label{rig-eq: D5.20}
\sum\nolimits_{\Gamma_l \cap  \tilde{\Psi}^- \neq \emptyset} \vert \Gamma_l\vert_{\cal H} \le C(1-\omega_{\rm min})^{-1}\psi^-
\end{align}
for  $C = C(h_*,q)$. If $\psi^- \ge c (1-\omega_{\rm min})\tau$ this follows from \eqref{rig-eq: property,b} and the fact that  $\Gamma_l \subset \overline{R_l} \subset N^{2\hat{\tau}}$ for all $\Gamma_l$ with $\Gamma_l \cap \tilde{\Psi}^- \neq \emptyset$ by \eqref{rig-eq: W prop}(i) and \eqref{rig-eq: new9}. If $\psi^- \le c(1-\omega_{\rm min})\tau$, by \eqref{rig-eq: property,a} we obtain  $|\Psi^- \cap \partial W_i|_{\cal H} \le D (1-\omega_{\rm min})^{-1}\psi^-$,  which for $c = c(h_*)> 0$ sufficiently small is small with respect to $\tau$. Thus, we can assume that $\Gamma_l \subset \Psi^-$ if $\Gamma_l \cap \tilde{\Psi}^- \neq \emptyset$. This implies $\sum\nolimits_{\Gamma_l \cap \tilde{\Psi}^-_i \neq \emptyset} \vert \Gamma_l\vert_{\cal H} \le C|\Psi^- \cap \partial W_i|_{\cal H}$ and gives the assertion  by \eqref{rig-eq: property,a}.

 Recall that $\upsilon \le r (1-\omega_{\rm min})^3$  (see beginning of Section \ref{rig-sec: subsub,  pre-est}). Applying Theorem \ref{rig-theo: korn} we obtain by \eqref{rig-eq: property,a},  \eqref{rig-eq: property,c}, \eqref{rig-eq: D5.20} and the fact that $\psi^- \le \upsilon \vert \Gamma\vert_\infty$  (cf. also \eqref{rig-eq: D1.5}  for  a similar estimate)
\begin{align}\label{rig-eq: D5.7}
\begin{split}
\int_{\tilde{\Psi}^-_i}& |\bar{u}(x) - (A_i\,x +b_i)|^2 \,dx  \\
&\le C| \tilde{\Psi}^-_i|\alpha( \tilde{\Psi}^-_i) +  CC_*\eps\upsilon^{-4} |\partial W_i \cap \tilde{\Psi}^-_i|\big( \sum\nolimits_{\Gamma_l \cap  \tilde{\Psi}^-_i \neq \emptyset} |\Gamma_l|_\infty\big)^2  
 \\
& \le C\upsilon^2|\Gamma|^2_\infty  (1-\omega_{\rm min})^{-1} \eps\psi^- + CC_* (1-\omega_{\rm min})^{-3}\eps\upsilon^{-2}|\Gamma|_\infty^2 \psi^- \\ & \le C(1+ C_*r) \upsilon^{-3} |\Gamma|^2_\infty \psi^- \eps 
\end{split}
\end{align}
for $A_i \in \R^{2 \times 2}_{\rm skew}$ and $b_i \in \R^2$, $i=1,3$. Likewise  using particularly \eqref{rig-eq: D5.20}  and $|\Psi_2^-| \le C (\psi^-)^2$ we get
\begin{align}\label{rig-eq: D5.7*}
\int_{\Psi^-_2} |\bar{u}(x) - &(A_2\,x +b_2)|^2 \,dx  \le C(1+ C_*r) \upsilon^{-4} |\Gamma|_\infty (\psi^-)^2 \eps 
\end{align}
for $A_2 \in \R^{2 \times 2}_{\rm skew}$ and $b_2 \in \R^2$. By \eqref{rig-eq: D5.2} for $V_\pm = \Psi^-_1   \setminus K_1, \Psi^-_3  \setminus K_1$ we see that 
$$\Vert \bar{u} - (A_+ \,x + b_+)\Vert^2_{L^2(\Psi^-_1  \setminus K_1) } + \Vert \bar{u} - (A_- \,x + b_-)\Vert^2_{L^2(\Psi^-_3   \setminus K_1)} \le C\upsilon G.$$
Applying  Lemma \ref{lemma: diff1}(i) on $D_1 = \tilde{\Psi}^-_i  \setminus K_1$, $D_2 =  \tilde{\Psi}^-_i$ $i=1,3$, we then derive  by \eqref{rig-eq: D5.7} employing $\psi^- \le \upsilon |\Gamma|_\infty$ 
\begin{align}\label{rig-eq: D5.4}
\begin{split}
&\tau^2|A_+ - A_1|^2 + |b_+ - b_1 + (A_+ - A_1)\,d_+|^2 \le C(\upsilon \vert \Gamma \vert^2_\infty)^{-1}G, \\
&\tau^2|A_- - A_3|^2 + |b_- - b_3 + (A_- - A_3)\,d_-|^2 \le C(\upsilon \vert \Gamma \vert^2_\infty)^{-1}G,
\end{split}
\end{align}
where $d_-=(-l_1,-\tau)^T$ and $d_+=(-l_1,\tau)^T$. Furthermore,  Lemma \ref{rig-lemma: trace}, \eqref{rig-eq: property,a}, \eqref{rig-eq: property,c}, \eqref{rig-eq: new8}, \eqref{rig-eq: D5.20} and  \eqref{rig-eq: D5.7}    yield
\begin{align}\label{rig-eq: D5.3}
\begin{split}
\int_{\partial  \tilde{\Psi}^-_i} & |\bar{u}(x) - (A_i\,x +b_i)|^2 \, d{\cal H}^1(x)   \\
& \le C (\upsilon \vert \Gamma \vert_\infty)^{-1} \Vert \bar{u} - (A_i \, \cdot + b_i)\Vert^2_{L^2(\tilde{\Psi}^-_i)}  +C\upsilon \vert \Gamma \vert_\infty \alpha(\tilde{\Psi}^-_i) 
\\
& \ \ + CC_*\eps \upsilon^{-4}  \sum\nolimits_{\Gamma_l \cap \tilde{\Psi}^-_i \neq \emptyset} \vert \Gamma_l\vert_{\cal H} \sum\nolimits_{\Gamma_l \cap \tilde{\Psi}^-_i \neq \emptyset}   \vert  \Gamma_l \vert_{\cal H}  \\
& \le C(1+ rC_* )  \upsilon^{-4}\vert \Gamma \vert_\infty    \psi^- \eps \le C \psi^-(\upsilon \vert \Gamma \vert^2_\infty)^{-1} G
\end{split}
\end{align}
for $i=1,3$, where we tacitly assumed that all boundary components are rectangular (cf. discussion after \eqref{rig-eq: D1.2.2}). In the  penultimate step we again used $ \psi^- \le \upsilon|\Gamma|_\infty$  and $\upsilon \le r (1-\omega_{\rm min})^3$. Likewise, we get
\begin{align}\label{rig-eq: D5.31}
\int_{\partial \Psi^-_2 \cup \Xi^-_0} |\bar{u}(x) - (A_2\,x +b_2)|^2 \,d{\cal H}^1(x)  \le  C \psi^-(\upsilon \vert \Gamma \vert^2_\infty)^{-1} G ,
\end{align}
where we replaced $(\upsilon \vert \Gamma \vert_\infty)^{-1}$ by $(\psi^-)^{-1}$  in the second line of \eqref{rig-eq: D5.3} and used \eqref{rig-eq: D5.7*} instead of \eqref{rig-eq: D5.7}. Observe that \eqref{rig-eq: D5.31} is well defined in the sense of traces since ${\cal H}^1( \Xi^-_0 \cap \partial W_i)= 0$. Define $\Xi^-_\pm = (-  \frac{\psi^-}{2} - l_1 , -l_1) \times \lbrace \pm \psi^- \rbrace$ and note that $\Xi^-_+ \subset \partial \Psi_2^- \cap \partial\tilde{\Psi}^-_1$ and $\Xi^-_- \subset \partial \Psi_2^- \cap \partial\tilde{\Psi}^-_3$. Again up to a small translation in $\e_2$-direction we may suppose  ${\cal H}^1( \Xi^-_\pm \cap \partial W_i)= 0$. Combining the estimates \eqref{rig-eq: D5.4}, \eqref{rig-eq: D5.3} for $i=1,3$  we obtain 
 $$\int_{\Xi^-_+} |\bar{u} - (A_+\,x +b_+)|^2 \,d{\cal H}^1 + \int_{\Xi^-_-} |\bar{u} - (A_-\,x +b_-)|^2 \,d{\cal H}^1\le C \psi^-(\upsilon \vert \Gamma \vert^2_\infty)^{-1} G.$$
 By Remark \ref{rem: diff1}(ii) we  estimate the difference of $A_\pm$, $A_2$ and $b_\pm$, $b_2$ on the boundaries $\Xi^-_\pm$  and obtain  by \eqref{rig-eq: D5.31} an expression similar to \eqref{rig-eq: D5.4} 
\begin{align}\label{rig-eq: D5.6}
(\psi^-)^2|A_\pm - A_2|^2 + |b_\pm - b_2 + (A_\pm - A_2)\,d_2|^2 \le C(\upsilon \vert \Gamma \vert^2_\infty)^{-1} G,
\end{align}
where $d_2 = (-l_1,0)^T$. Together with \eqref{rig-eq: D5.31} this leads to
$$\int_{\Xi^-_0} |\bar{u}(x) - (A_\pm\,x +b_\pm)|^2 \,d{\cal H}^1(x)  \le  C\psi^-(\upsilon \vert \Gamma \vert^2_\infty)^{-1} G$$
and then by the triangle inequality we derive
$$\int_{ \Xi^-_0} |(A_+ - A_-)\,x +(b_+ - b_-)|^2 \,d{\cal H}^1(x)  \le  C \psi^-(\upsilon \vert \Gamma \vert^2_\infty)^{-1} G.$$
This gives the desired estimate \eqref{rig-eq: D5.5}. From \eqref{rig-eq: D5.5} applied on both sets, $\Xi^-_0$ and $\Xi^+_0$, we deduce
$$- C(l_1\upsilon)^2 |(A_+ - A_-) \, \e_1|^2   + |-(A_+ - A_-) \, l_1\e_1 + (b_+  -b_-)|^2 \le C(\upsilon \vert \Gamma \vert^2_\infty)^{-1}G$$
and
$$- C(l_1\upsilon)^2 |(A_+ - A_-) \, \e_1|^2   + |(A_+ - A_-) \, l_1\e_1 + (b_+  -b_-)|^2 \le C(\upsilon \vert \Gamma \vert^2_\infty)^{-1}G.$$
Combining these two estimates we find for $\upsilon$ sufficiently small $l_1^2|A_+ - A_-|^2 = 2l_1^2|(A_+ - A_-)\e_1|^2  \le C(\upsilon \vert \Gamma \vert^2_\infty)^{-1}G$ and then also $|b_+ - b_-|^2  \le C(\upsilon \vert \Gamma \vert^2_\infty)^{-1}G$. (This is the step where we fundamentally use $\dist(K_1,K_2)\ge c |\Gamma|_\infty$.) We choose $A=A_-$ and 
$b = b_-$. Recalling the definition of $G$ and $|N| \le \upsilon \vert \Gamma\vert^2_\infty$ we obtain by \eqref{rig-eq: D5.2} for $V_\pm = N'_\pm$
$$\int_{N'_+ \cup N'_-} |\bar{u}(x) - (A\,x +b)|^2 \,dx \le  C(1 + C_* r) \frac{\eps}{\upsilon^3} \vert \Gamma \vert^3_\infty,$$
which  together with the estimates \eqref{rig-eq: D5.7} and \eqref{rig-eq: D5.4} gives \eqref{rig-eq: D1.2.2}(i). Finally,  \eqref{rig-eq: D5.6} yields $(\psi^-)^2|A - A_2|^2 \le C(\upsilon \vert \Gamma \vert^2_\infty)^{-1}G$ and $|b - b_2 + (A - A_2)\, (-l_1,0)^T|^2 \le C (\upsilon \vert \Gamma \vert^2_\infty)^{-1}G$. Then by \eqref{rig-eq: D5.7*} and the fact that $|\Psi^-_2| \le C(\psi^-)^2 \le C\psi^- \upsilon \vert \Gamma\vert_\infty$ we conclude
$$\int_{\Psi^-_2} |\bar{u}(x) - (A\,x +b)|^2 \,dx  \le C(1+ rC_* )  \upsilon^{-3}\vert \Gamma \vert^2_\infty    \psi^- \eps  $$
giving \eqref{rig-eq: D1.2.2}(ii). The estimate for $ \Psi^+_2$ follows analogously. 

It remains to briefly indicate the necessary adaptions for the general case. The main differences are (i) the shape of the sets $\Psi^\pm_{i}, i=1,2,3$ and (ii) the position of the sets $K_1, K_2$. For (i) we observe that  $\Psi^\pm_{i}, i=1,3$, are $C(h_*)$-Lipschitz equivalent to a square by Lemma \ref{rig-lemma: two observations1}(ii) and Lemma \ref{rig-lemma: two observations2}(ii) whereby \eqref{rig-eq: D5.7} can still be derived (cf. Remark \ref{rig-rem: 1}(i)). (Note that the sets are even related  by  affine mappings.) Likewise, an estimate of the form \eqref{rig-eq: D5.7*} can be derived for sets $(\Psi^\pm_2)^* \supset\Psi^\pm_2$ which have been constructed in Section \ref{rig-sec: subsub,  dod-neigh}. Moreover, although not stated explicitly in Section \ref{rig-sec: sub, trace}, the trace estimate used in \eqref{rig-eq: D5.3}, \eqref{rig-eq: D5.31} can also be applied for sets being an affine transformation of a square. The rest of the arguments concerning the difference of infinitesimal rigid motions (see \eqref{rig-eq: D5.4}, \eqref{rig-eq: D5.6}) remains unchanged. For (ii) we observe that in the derivation of $|A_+ -A_-|^2 \le Cl_1^{-2}(\upsilon |\Gamma|_\infty^2)^{-1} G$ we fundamentally used that $\dist(K_1,K_2) \sim l_1$, but the exact position of the sets $K_1,K_2$ was not essential. \eop

We close this section with an estimate for the skew symmetric matrices involved in the above  results which will be needed in the derivation of the nonlinear rigidity estimates in \cite{Friedrich-Schmidt:15}. 
\begin{lemma} \label{rig-lemma: A neigh}
Let be given the situation of Theorem \ref{rig-th: derive prop}  for a function $u \in H^1(W)$ and define $y = \bar{R}\, (\id + u)$, where $\id$ denotes the identity function and $\bar{R} \in SO(2)$. Let  $V \subset Q_\mu$ be a rectangle and  let ${\cal F}(V)$  be the boundary components $(\Gamma_l)_l =  (\Gamma_l(U))_l$ satisfying $N^{\hat{\tau}_l}(\partial R_l) \subset V$ and \eqref{rig-eq: D1}. Then there is a $C_3 = C_3(\sigma,h_*)$  independent of $u,W,\bar{R}$ and $V$ such that   
$$\sum\nolimits_{\Gamma_l \in {\cal F}(V)}   |\Gamma_l|_\infty^2 |A_l|^p \le C_3 \big( \Vert \nabla y - \bar{R}\Vert^p_{L^p(V  \cap W)} +  (\eps s^{-1})^{\frac{p}{2}-1} \eps|\partial U \cap V|_{\cal H} \big)$$
for $p=2,4$, where  $A_l \in \R^{2 \times 2}_{\rm skew}$ is given in  \eqref{rig-eq: extend def_new}.
\end{lemma}

 \begin{rem}\label{rig-rem: z2}
 {\normalfont
  Similarly as in Remark \ref{rig-rem: z}(i) we note that the constant $C_3=C_3(\sigma,h_*)$ has polynomial growth in $\sigma$, i.e. $C_3(\sigma,h_*) \le C(h_*) \sigma^{-z}$ for some $z \in \N$.
 }
 \end{rem}

\Proof Let $p =2,4$.   Consider a component $\Gamma = \Gamma_l(U)$ with corresponding rectangle $R$ and $X$ with $\partial X = \Gamma$.  Let $A \in \R^{2 \times 2}_{\rm skew}$ as given by \eqref{rig-eq: extend def_new}. It suffices to show 
$$|\Gamma|_\infty^2 |A|^p \le C_3 \big( \Vert \nabla y - \bar{R}\Vert^p_{L^p({N}')} +  (\eps s^{-1})^{\frac{p}{2}-1}\eps|\Gamma|_{\cal H}\big)$$
 for this component, where  ${N}' = N^{\hat{\tau}}(\partial R) \setminus \bigcup_{\Gamma_l \in {\cal I}(\Gamma)} N^{\hat{\tau}_l}( \partial R_l)$ and ${\cal I}(\Gamma) = \lbrace \Gamma_l: |\Gamma_l|_\infty \le |\Gamma|_\infty \rbrace$. Then the assertion follows by summation over all components.

  As $\Gamma$ satisfies \eqref{rig-eq: D1}, we observe  that we applied Theorem \ref{rig-theorem: D} on  $\partial R$ in some iteration step, in particular \eqref{rig-eq: D1.2.2} is satisfied.  Recall the rectangular neighborhood of $\partial R$ and its coverings as constructed in Section \ref{rig-sec: subsub,  rec-neigh} with $R = (-l_1,l_1) \times (-l_2,l_2)$. Choose $U \in {\cal C}$ with $U \subset N_{2,+}(\partial R)$ as considered in Lemma \ref{rig-lemma: D2}. By assumption we find  a set $S \subset (l_2, l_2 + \tau)$ with $|S| \ge \frac{1}{2}\frac{\tau}{20}$ such that for $T = (\R \times S) \cap U$ we have $T \cap \bigcup_{\Gamma_l \in {\cal I}(\Gamma)} N^{\hat{\tau}_l}( \partial R_l) = \emptyset$  for $h_*$ small enough (see \eqref{rig-eq: hat tau def}). In view of Remark \ref{rig_rem: connect}(ii) we get $T \subset W$. Moreover, $|T| \ge C\upsilon|\partial R|_\infty^2$ by \eqref{rig-eq: covering def}. It is not restrictive to assume that $S$ is connected as otherwise we follow the subsequent arguments for every connected component of $S$. Recall $|\Gamma|_\infty \le |\partial R|_\infty \le 2|\Gamma|_\infty$ by \eqref{rig-eq: new2}(i). The Poincar\'e inequality and a scaling argument imply
$$\int\nolimits_T |u(x) - \hat{b}|^2 \, dx\le C |T|^{1 - \frac{2}{p}}  |\Gamma|_\infty^2 \Vert \nabla y - \bar{R}\Vert^2_{L^p(T)} $$
for a constant $\hat{b} \in \R^2$ and $p=2,4$.  This together with \eqref{rig-eq: D1.2.2} yields
$$\int\nolimits_T |A\,x  + b - \hat{b}|^2\, dx \le C (\upsilon|\Gamma|_\infty^2)^{ 1 - \frac{2}{p}}   |\Gamma|_\infty^2 \Vert \nabla y - \bar{R}\Vert^2_{L^p(T)} +  C \hat{C}\sum^\infty_{n=0}\Big( \frac{2}{3}\Big)^n \upsilon^{-3}|\Gamma|_\infty^3\eps.$$  
For the constant in the latter part see   \eqref{rig-eq: D2}.  By Lemma \ref{lemma: diff1}(i)  we then find
$|T||\Gamma|_\infty^2|A|^2 \le C \int\nolimits_T |A\,x  + b - \hat{b}|^2\, dx$ and thus by $|T| \ge C\upsilon|\Gamma|_\infty^2$ and an elementary calculation we derive
$$ |\Gamma|^2_\infty |A|^4 \le C\upsilon^{-1} \Vert \nabla y - \bar{R}\Vert^4_{L^4(T)} + C\upsilon^{-8} \eps^2.$$
  Since $W \in {\cal V}^s$, we get $|\partial R|_\infty \ge s$ and obtain $|\Gamma|_{\cal H} \ge |\Gamma|_\infty \ge \frac{1}{2}|\partial R|_\infty \ge \frac{s}{2}$ by \eqref{rig-eq: new2}(i). Recalling  $T \subset N^{\bar{\tau}}(\partial R) \setminus \bigcup_{\Gamma_l \in {\cal I}(\Gamma)} N^{\bar{\tau}_l}( \partial R_l) \subset V$ we find  for $C_3 =C(h_*)\upsilon^{-8}$ with $C(h_*)$ large enough 
$$|\Gamma|^2_\infty |A|^4 \le C_3 \Vert \nabla y - \bar{R}\Vert^4_{L^4({N}')} + C_3  \eps s^{-1} \eps |\Gamma|_{\cal H}.$$
giving the claim for $p=4$. Likewise, for $p
=2$ we deduce
\begin{align*}
|\Gamma|^2_\infty |A|^2 &\le C\upsilon^{-1} \Vert \nabla y - \bar{R}\Vert^2_{L^2(T)} + C\upsilon^{-4} \eps |\Gamma|_\infty   \le C_3 \Vert \nabla y - \bar{R}\Vert^2_{L^2({N}')} + C_3 \eps |\Gamma|_{\cal H}.
\end{align*}
 Recalling the definition of $C_1=C_1(h_*,\sigma)$ before \eqref{eq: before} we see $C_3 \le C(h_*) C_1^2$. This yields $C_3=C_3(h_*,\sigma)$ and also Remark \ref{rig-rem: z2} by Remark \ref{rig-rem: z}(i).
 \eop

\begin{appendix}

\section{Appendix}\label{rig-sec: appe}

In this section we recall the definition and basic properties of functions of bounded deformation, derive estimates on the difference of infinitesimal rigid motions and establish a trace theorem.

\subsection{Functions of bounded deformation}\label{rig-sec: sub, bd}

We collect the definitions and some fundamental properties of functions of bounded deformation. Let $\Omega \subset \R^d$ open, bounded with Lipschitz boundary. Recall that the space $BD(\Omega, \R^d)$ of \emph{functions of bounded deformation} consists of functions $u \in L^1(\Omega, \R^d)$ whose symmetrized distributional derivative $Eu := \frac{(Du)^T + Du}{2}$ is a finite $\R^{d \times d}_{\rm sym}$-valued Radon measure. In \cite{Ambrosio-Coscia-Dal Maso:1997} it is shown that it can be decomposed as 
\begin{align}\label{rig-eq: symmeas}
 Eu = e(u) {\cal L}^d +  E^s u = e(u) {\cal L}^d + [u] \odot \xi_u {\cal H}^{d-1}|_{J_u} + {\cal E}^c(u),
 \end{align}
where $e(u)$ is the absolutely continuous part of $Eu$ with respect to the Lebesgue measure, ${\cal E}^c(u)$ denotes the `Cantor part',  $J_u$ (the `crack path') is an ${\cal H}^{d-1}$-rectifiable set in $\Omega$, $\xi_u$ is a normal of $J_u$ and $[u] = u^+ - u^-$ (the `crack opening') with $u^{\pm}$ being the one-sided limits of $u$ at $J_u$. Here ${\cal L}^d$ denotes the $d$-dimensional Lebesgue measure, ${\cal H}^{d-1}$ denotes the $(d-1)$-dimensional Hausdorff measure and $a \odot b = \frac{1}{2}(a \otimes b + b \otimes a)$.

The space $SBD(\Omega,\R^d)$ of \emph{special functions of bounded deformation} consists of all $u \in BD(\Omega, \R^d)$ with ${\cal E}^c(u) = 0$. If in addition $e(u) \in L^2(\Omega)$ and ${\cal H}^{d-1}(J_u) < \infty$, we write $u \in SBD^2(\Omega)$. For basic properties of this function space we refer to \cite{Ambrosio-Coscia-Dal Maso:1997,  Bellettini-Coscia-DalMaso:98}.

We recall a Korn-Poincar\'e inequality and a trace theorem in BD (see \cite{Bredies:13,Temam:85}). 

\begin{theorem}\label{rig-theo: korn}
Let $\Omega \subset \R^d$ open, bounded,  connected with Lipschitz boundary and let $P: L^2(\Omega,\R^d) \to  L^2(\Omega,\R^d)$ be a linear projection onto the space of infinitesimal rigid motions. Then there is a constant $C>0$,  which is invariant under rescaling of the domain, such that for all $u \in BD(\Omega,\R^d)$
$$\Vert u - Pu \Vert_{L^{\frac{d}{d-1}}(\Omega)} \le C |Eu|(\Omega).$$
\end{theorem}

\begin{theorem}\label{rig-th: tracsbv}
 Let $\Omega \subset \R^d$ open, bounded,  connected with Lipschitz boundary. There exists a constant $C > 0$ such
that the trace mapping $\gamma: BD(\Omega,\R^2) \to  L^1(\partial \Omega, \R^2)$ is well defined and
satisfies the estimate
$$\Vert \gamma u \Vert_{L^1(\partial \Omega)}\le C\big(\Vert u\Vert_{L^1(\Omega)} + |Eu|(\Omega) \big)$$
for each $u \in BD(\Omega,\R^2).$
\end{theorem}

For sets which are related through bi-Lipschitzian homeomorphisms with Lipschitz constants of both
the homeomorphism itself and its inverse uniformly bounded the constants in Theorem \ref{rig-theo: korn} can be chosen independently of these sets, see e.g. \cite{FrieseckeJamesMueller:02}.  
We now present a density result in $SBD^2(\Omega)$ (see \cite[Theorem 3, Remark 5.3]{Chambolle:2004})). 

\begin{theorem}\label{th: sbd dense}
Let $\Omega \subset \R^d$ open, bounded with Lipschitz boundary. Then for every $u \in SBD^2(\Omega) \cap L^2(\Omega)$ and $\delta>0$ there is a function $\tilde{u} \in SBD^2(\Omega) \cap L^2(\Omega)$ such that $J_{\tilde{u}}$ is the finite union of closed connected pieces of $C^1$-hypersurfaces and  $u|_{\Omega\setminus J_{\tilde{u}}} \in H^{1}(\Omega \setminus J_{\tilde{u}})$ such that
\begin{align}\label{rig-eq: V3.20}
\begin{split}
\Vert u - \tilde{u} \Vert_{L^2(\Omega)} +\Vert e(u) - e(\tilde{u}) \Vert_{L^2(\Omega)} \le \delta, \ \ \ \ 
{\cal H}^1(J_{\tilde{u}}) \le {\cal H}^1(J_u) + \delta.
\end{split}
\end{align}
Moreover, if $u \in L^\infty(Q_\mu)$, one can ensure  that $\Vert \tilde{u} \Vert_\infty \le \Vert u \Vert_\infty$.
\end{theorem}

\subsection{Korn-Poincar\'e-inequality and infinitesimal rigid motions}\label{rig-sec: sub, pre} 

In this section we provide estimates for infinitesimal rigid motions and analyze how the constant  of the Korn-Poincar\'e inequality in Theorem \ref{rig-theo: korn} depends  on the set shape. For $A \in \R^{2 \times 2}_{\rm skew}$, $b \in \R^2$ we define for shorthand the mapping $a= a_{A,b}: \R^2 \to \R^2$ with $a(x) = A\,x  + b$ for $x \in \R^2$. By ${\rm diam}(D)$  we denote the diameter of a set $D \subset \R^2$.  We start with an elementary property of infinitesimal rigid motions.

\begin{lemma}\label{lemma: diff1}
Let $c>0$. Then there is a constant $C=C(c)>0$ such that for all measurable sets $D_1\subset D_2 \subset \R^2$, satisfying $|D_1 \setminus B_\rho| \ge \frac{1}{2}|D_1|$ for all balls $B_\rho$ with radius $\rho = c \, {\rm diam}(D_1)$,  and for all infinitesimal rigid motions $a=a_{A,b}$ we have
\begin{align*}
(i)& \ \  ({\rm diam}(D_1))^2|A|^2 \le C\tfrac{1}{|D_1|}\Vert a \Vert^2_{L^2(D_1)}, \ \  \Vert a \Vert^2_{L^\infty(D_2)} \le C\tfrac{1}{|D_1|} \big(\tfrac{{\rm diam}(D_2)}{{\rm diam}(D_1)}\big)^2\Vert a\Vert^2_{L^2(D_1)}, \\  (ii)& \ \   \Vert a\Vert^2_{L^2(Z)} \le C|Z| \ \tfrac{1}{|D_1|} \big(\tfrac{{\rm diam}(D_2)}{{\rm diam}(D_1)}\big)^2\Vert a\Vert^2_{L^2(D_1)}
\end{align*}
for all measurable $Z \subset D_2$. 
\end{lemma}

\Proof Without restriction assume that $A\neq 0$ as otherwise the statement is clear. The assumption $A \in \R^{2 \times 2}_{\rm skew}$ implies that $A$ is invertible and that $|Ay|=\frac{\sqrt2}{2}|A| |y|$ for all $y \in \R^{2}$. Setting $z:= - A^{-1}b$ we find for all $x \in D_1$
\begin{align*}
\tfrac{1}{\sqrt2}|A| |x-z| = |a(x)|.
\end{align*}
By assumption we have $|D_1 \setminus B_{\rho}(z)| \ge  \frac{1}{2}|D_1|$ with $\rho = c\, {\rm diam}(D_1)$. This implies $({\rm diam}(D_1))^2|D_1||A|^2 \le C \Vert a \Vert^2_{L^2(D_1)}$.  Observe that there is some $x_0 \in D_1$ with $|A\,x_0 + b|^2 \le |D_1|^{-1} \Vert a\Vert^2_{L^2(D_1)}$. This together with the bound on $|A|$ yields for all $x \in D_2$
$$|A\,x + b|^2 \le C|A\,x_0+b|^2 + C|A|^2({\rm diam}(D_2))^2 \le C|D_1|^{-1} \big(\tfrac{{\rm diam}(D_2)}{{\rm diam}(D_1)}\big)^2\Vert a\Vert^2_{L^2(D_1)}.$$
This shows (i) and (ii) follows directly from the estimate on $\Vert a \Vert_{L^\infty(D_2)}$. \eop

\begin{rem}\label{rem: diff1}

{\normalfont
(i) The assumption on $D_1$ is particularly satisfied for rectangles for a universal constant $c$ small enough.

(ii) Analogous estimates hold on segments. If, e.g., $S= (0,l) \times \lbrace 0\rbrace$ for $l>0$, then for a universal constant $C>0$ we have for all infinitesimal rigid motions $a=a_{A,b}$
$$l^2|A|^2 + \max\nolimits_{x \in S} |a(x)|^2 \le Cl^{-1}\int_{S} |a(x)|^2 \, d{\cal H}^1(x).$$

}
\end{rem}

In order to quantify how the constant in Theorem \ref{rig-theo: korn} depends on the set shape we will estimate the variation of infinitesimal rigid motions on different squares (cf. also \cite{FrieseckeJamesMueller:02}).  For $s>0$ we partition $\R^2$ up to a set of measure zero into squares $Q^s(p) = p + s(-1,1)^2$ for $p \in I^s := s(1, 1) + 2s\Z^2$. We first have the following estimate.

\begin{lemma}\label{rig-lemma: weak rig2}
There is a universal constant $C>0$ such that for all $s>0$, $N\in \N$, squares $Q^s(\xi_1), \ldots, Q^s(\xi_N)$ and infinitesimal rigid motions $a_i = a_{A_i,b_i}$, $i=1,\ldots, N$ one has for all $j=1,\ldots,N$
$$\max_{1 \le i_1,i_2 \le N} \Vert a_{i_1} - a_{i_2}\Vert_{L^2(Q^s(\xi_j))} \le Cds^{-1}\sum\nolimits^{N-1}_{i=1} \Vert a_i - a_{i+1} \Vert_{L^2(Q^s(\xi_i) \cup Q^s(\xi_{i+1}))},$$
where $d = \max_{1 \le i_1,i_2 \le N} |\xi_{i_1} - \xi_{i_2}|$.
\end{lemma}

\Proof For shorthand we set $Q_i = Q^s(\xi_i)$ and $Q_{i,i+1}= (\overline{Q^s(\xi_i) \cup Q^s(\xi_{i+1})})^\circ$ for $i=1,\ldots,N$.  Fix $Q_j$.   First, for each  $i=1,\ldots,N-1 $ by Lemma \ref{lemma: diff1}(ii) for $D_1 = Q_{i,i+1}$ and $D_2 = Q_j \cup Q_{i,i+1}$ we have 
$$\Vert a_i - a_{i+1} \Vert_{L^2(Q_j \cup Q_{i,i+1})} \le C d s^{-1} \Vert a_i - a_{i+1} \Vert_{L^2(Q_{i,i+1})}.$$
Consider $\xi_{i_1}$ and $\xi_{i_2}$ for $1 \le i_1,i_2 \le N$. Then the triangle inequality   yields
\begin{align*} 
 \Vert a_{i_2} - a_{i_1}\Vert_{L^2(Q_j)} &\le \sum\nolimits^{i_2-1}_{i=i_1} \Vert a_{i+1} - a_i \Vert_{L^2(Q_j)} \le Cds^{-1} \sum\nolimits^{N-1}_{i=1} \Vert a_i - a_{i+1} \Vert_{L^2(Q_{i,i+1})}.
\end{align*}
 \eop

We now   derive a Korn-Poincar\'e-inequality for sets in 
$${\cal U}^s := \Big\{ U \subset \R^2: \exists I \subset I^s: U = \Big(\bigcup\nolimits_{p \in I} \overline{Q^s(p)} \Big)^\circ \Big\}.$$
Let $I^s(U) = \lbrace p \in I^s: Q^s(p) \subset U \rbrace$.  In the following $\dashint_{Z}$ stands for $\frac{1}{|Z|}\int_Z$.

\begin{lemma}\label{rig-lemma: weak rig}
There is a constant $C > 0$ such that for all $s>0$, for all sets $U \in {\cal U}^s$ and $u \in SBD^2(U)$ the following holds: \\
(i) For all infinitesimal rigid motions $a(p)$, $p \in I^s(U)$, we find $A \in \R^{2 \times 2}_{\rm skew}$ and $b \in \R^2$ such that for all $Z \subset U$, $Z \in {\cal U}^s$,   we have 
\begin{align*}
|U|\dashint_{Z}  |u(x)- (A\,x + b)|^2 \, dx  & \le C  s^{-2}|U|  \sum\nolimits_{p \in I^s(U)} \Vert u - a(p) \Vert^2_{L^2(Q^s(p))}  \\
& \ \  + Cs^{-2}|U|\max_{p_1,p_2 \in I^s(U)} \Vert a(p_1) - a(p_2)\Vert^2_{L^2(Q^s(p_1) \cup Q^s(p_2))},
\end{align*}
(ii) If $U$ is connected, there are $A \in \R^{2 \times 2}_{\rm skew}$ and $b \in \R^2$  such that  for all $Z \subset U$, $Z \in {\cal U}^s$ 
\begin{align*}
|U|\dashint_{Z}  |u(x)- (A\,x + b)|^2 \, dx \le C  (s^{-2}|U|)^3 (|E u|(U))^2.
\end{align*}
\end{lemma}

\Proof  (i) We fix $p_0 \in I^s(U)$ and set set $A = A(p_0)$, $b=b(p_0)$. Summing over all $p \in I^s(U)$ with $Q^s(p) \subset Z \subset U$ we derive by the triangle inequality
\begin{align*}
 \Vert u - a(p_0)\Vert^2_{L^2(Z)} &\le C\sum\nolimits_{p \in I^s(U)} \Vert u - a(p) \Vert^2_{L^2(Q^s(p))} \\ 
 & \ \ \ \ \ \ \ \  + C|Z|s^{-2}\max_{p \in I^s(U)} \Vert a(p_0) - a(p)\Vert^2_{L^2(Q^s(p))}. 
\end{align*}
The claim follows by multiplication with $\frac{|U|}{|Z|}$ and the observation that $|Z| \ge 4s^2$.

(ii) For each $Q^s(p) \subset U$ we can apply Theorem \ref{rig-theo: korn} to obtain an infinitesimal rigid motion $a(p)$  such that 
\begin{align}\label{rig-eq: A rigXXX}
\Vert u - a(p)\Vert_{L^2(Q^s(p))} \le C|Eu|(Q^s(p))
\end{align}
for a universal constant (independent of $u$,$s$). Fix arbitrary $p_1, p_2 \in I^s(U)$. As $U$ is connected, there is a path $\xi=(\xi_1=p_1,\ldots,\xi_N=p_2)$ with $N \le |U|(2s)^{-2}$ and $\xi_j - \xi_{j-1} = \pm 2s\e_i$ for some $i=1,2$ for all $j=2,\ldots,N$. As in the previous proof we let $Q_i = Q^s(\xi_i)$ and  $Q_{i,i+1} =(\overline{Q^s(\xi_i) \cup Q^s(\xi_{i+1})})^\circ$. Then we derive by Theorem \ref{rig-theo: korn}
\begin{align}\label{rig-eq: A rig}
\Vert u - a_{i,i+1}\Vert_{L^2(Q_{i,i+1})}  \le C|Eu|(Q_{i,i+1})
\end{align}
for a suitable infinitesimal rigid motion $a_{i,i+1}$. Let $j \in \lbrace i,i+1 \rbrace$. By \eqref{rig-eq: A rigXXX}, \eqref{rig-eq: A rig}, the triangle inequality and Lemma \ref{lemma: diff1}(ii) for $D_1 = Q_j$ and $D_2 = Q_{i,i+1}$ we obtain 
\begin{align*}
\Vert a(\xi_j) - a_{i,i+1}\Vert_{L^2(Q_{i,i+1})} \le C\Vert a(\xi_j) - a_{i,i+1}\Vert_{L^2(Q_{j})} \le C|Eu|(Q_{i,i+1})
\end{align*}
for a universal constant $C>0$ and therefore again by the triangle inequality
\begin{align}\label{rig-eq: A difference}
\Vert a(\xi_{i+1}) - a(\xi_{i})\Vert_{L^2(Q_{i,i+1})}  \le C|Eu|(Q_{i,i+1}).
\end{align}
This allows to estimate the difference of the infinitesimal rigid motions on $Q^s(p_1)$ and $Q^s(p_2)$ by Lemma \ref{rig-lemma: weak rig2} and we derive 
\begin{align*}
\Vert a(p_1) - a(p_2)\Vert_{L^2(Q^s(p_1) \cup Q^s(p_2))}  \le Cs^{-2}|U||Eu|(U),
 \end{align*}
 where in the last step we used $d \le CNs$ and  $N \le |U|(2s)^{-2}$. This together with \eqref{rig-eq: A rigXXX} yields the claim by (i).
 \eop

\begin{rem}\label{rig-rem: 1} 
{\normalfont
(i) The fact that we covered the sets $U$ with squares is not essential. Recalling how we derived  \eqref{rig-eq: A rigXXX}, \eqref{rig-eq: A difference} we could as well cover $U$ with rectangles $R_i = t_i + (-a_i,a_i) \times (-b_i,b_i)$, where $c_1 a_i \le b_i \le c_2a_i$ and $c_1 s \le b_i \le c_2 s$ for constants $0<c_1 < c_2$. The constants in \eqref{rig-eq: A rigXXX}, \eqref{rig-eq: A rig} only depend on $c_1, c_2$ as all the possible shapes are related to $(-s,s)^2$ through bi-Lipschitzian homeomorphisms with Lipschitz constants of both the homeomorphism itself and its inverse bounded (see Section \ref{rig-sec: sub, bd}). 

(ii) In view of the proof, in the choice of $A\,x +b$ we have the freedom to select any of the infinitesimal rigid motions which are given on the squares $Q^s(p)\subset U$.

}
\end{rem}

We now show that under additional assumptions on the energies the constants involved in the above estimates may be refined.

\begin{lemma}\label{rig-lemma: weak rig4}
 There is a universal constant $C>0$ such that for all  $\delta \in (0,\frac{1}{2})$ and for all $B_1,B_2 \in {\cal U}^s$ with $B_1 \subset B_2 \subset \R \times (-s,s)$,  $\frac{1}{1-\delta}{\rm diam}(B_1) \ge  {\rm diam}(B_2)\ge \frac{s}{\delta}$ and $|B_2 \setminus B_1| \le \delta |B_2|$ the following holds: For all $u \in L^2(B_2)$ and infinitesimal rigid motions $a_1$, $a_2$ as well as $H_1,H_2 \ge 0$ such that one has
\begin{align}\label{rig-eq: new22}
\Vert u - a_i\Vert^2_{L^2(B_i \cap Z)} \le |B_i \cap Z||B_2|^{-1} H_i
\end{align}
for $i=1,2$ and all $Z \subset B_2$, $Z \in {\cal U}^s$,  we find
$$\Vert u - a_1 \Vert^2_{L^2(Z)} \le (1+C\delta)|B_2|^{-1}|Z|H_1+ \tfrac{C}{\delta}|B_2|^{-1}|Z|H_2  $$
for all $Z \subset B_2$, $Z \in {\cal U}^s$.
\end{lemma}

\Proof Let $a_i = a_{A_i,b_i}$ for $i=1,2$  and $u \in L^2(B_2)$ be given such that \eqref{rig-eq: new22} holds. As ${\rm diam}(B_1) \ge \frac{1}{2} {\rm diam}(B_2)$, we can choose two squares $Q_1,Q_2 \subset B_1$ with ${\rm diam} (Q_1 \cup Q_2) \ge c \, {\rm diam}(B_2)$ for $c$ small enough. Then Lemma \ref{lemma: diff1}(i) with $D_1 = D_2 = Q_1 \cup Q_2$ together with \eqref{rig-eq: new22} and the triangle inequality yields   
\begin{align}\label{eq:AA}
({\rm diam}(B_2))^2|A_1-A_2|^2 \le C|B_2|^{-1} (H_1 + H_2).
\end{align}
For shorthand we write $\bar{a}= a_1 - a_2$ and $\bar{A} = A_1- A_2$. Let $x_0 \in B_2$ arbitrary. As $|B_2 \setminus B_1| \le \delta |B_2|$ and ${\rm diam}(B_1) \ge (1-\delta){\rm diam}(B_2)$, we find some $x'_0 \in B_1$ with $|x_0-x_0'| \le C\delta {\rm diam}(B_2) $ for a universal $C>0$ large enough. Let $Q \subset B_1$ be the square containing $x_0'$ and observe $|x-x_0| \le C(s+ \delta {\rm diam}(B_2)) \le C\delta {\rm diam}(B_2)$ for all $x \in Q$, where in the last step we used $s \le \delta {\rm diam}(B_2)$. Then applying a scaled version of Young's inequality, the triangle inequality and \eqref{eq:AA} we compute
\begin{align*}
|Q||\bar{a}(x_0)|^2 &= \Vert \bar{a}(x_0)\Vert^2_{L^2(Q)} \le (1+ \delta)\Vert \bar{a} \Vert^2_{L^2(Q)} + (1+\tfrac{1}{\delta})\Vert\bar{A}\, (\cdot-x_0) \Vert^2_{L^2(Q)}\\
& \le (1+\delta)^2\Vert u - a_1\Vert^2_{L^2(Q)} + \tfrac{C}{\delta}\Vert u - a_2\Vert^2_{L^2(Q)}  \\ & \ \ \  + \tfrac{C}{\delta}|Q|\max\nolimits_{x \in Q}|x-x_0|^2  |B_2|^{-1} ({\rm diam}(B_2))^{-2}(H_1 + H_2) \\
& \le (1+C\delta)|B_2|^{-1}|Q|H_1 + \tfrac{C}{\delta}|B_2|^{-1}|Q|H_2 +  C|B_2|^{-1}|Q| \delta (H_1 + H_2)\\
&\le (1+C\delta)|B_2|^{-1}|Q|H_1+ \tfrac{C}{\delta}|B_2|^{-1}|Q|H_2.
\end{align*}  
As $x_0 \in B_2$ was arbitrary, we conclude by \eqref{rig-eq: new22} and Young's inequality for all $Z \subset B_2$, $Z \in {\cal U}^s$
\begin{align*}
\Vert u - a_1 \Vert^2_{L^2(Z)} &\le (1+\delta)\Vert \bar{a} \Vert^2_{L^2(Z)} + (1 + \tfrac{1}{\delta})\Vert u - a_2 \Vert^2_{L^2(Z)}\notag \\
& \le (1+C\delta)|B_2|^{-1}|Z|H_1+ \tfrac{C}{\delta}|B_2|^{-1}|Z|H_2.
\end{align*}
\eop

\subsection{A trace theorem in SBD$^2$}\label{rig-sec: sub, trace}

By the trace theorem for BD functions (Theorem \ref{rig-th: tracsbv}) one can control the $L^1$-norm of the function on the boundary. In our framework we may establish a trace theorem in $L^2$ for SBD$^2$ functions if the jump set is sufficiently regular: Let $Q_\mu = (-\mu,\mu)^2$ and recall the definition of $SBD^2(Q_\mu)$ in Section \ref{rig-sec: sub, bd}. We suppose that some $u \in SBD^2(Q_\mu)$ satisfies $J_{u} = \bigcup_j \Gamma_j \cap Q_\mu$, where $\Gamma_i = \partial R_i$ for rectangles $R_i = (a^i_1,a^i_2) \times (b^i_1,b^i_2) \subset \R^2$ (note that for the application we have in mind we do not require that the rectangles are subsets of $Q_\mu$.). Clearly, as $u \in H^1(Q_\mu\setminus J_{u})$ the trace is well defined in $L^2$. More precisely, we have the following statement. 

\begin{lemma}\label{rig-lemma: trace}
Let $\mu >0$. There is a constant $C>0$ such that for all $u \in SBD^2(Q_\mu)$ with $J_{u} = \bigcup^n_{j=1} \Gamma_j \cap Q_\mu$, where $\Gamma_j = \partial R_j$, one has
\begin{align}\label{rig-eq: trace1}
\begin{split}
\int_{\partial Q_\mu}& |u|^2 \, d{\cal H}^1 \le C\mu \Vert e(u)\Vert^2_{L^2(Q_\mu)} + \tfrac{C}{\mu}\Vert u\Vert^2_{L^2(Q_\mu)} \\
&  + C\sum\nolimits^n_{j=1} {\cal H}^1(\Gamma_j)\ \sum\nolimits^n_{j=1} \Big(({\cal H}^1(\Gamma_j))^{-1}\int_{\Gamma_j \cap Q_\mu} |[u]|^2 \, d{\cal H}^1\Big).
\end{split}
\end{align}
\end{lemma} 

\Proof Let $Q_\mu = (-\mu,\mu)^2$ and $u \in SBD^2(Q_\mu)$ with $J_{u} = \bigcup^n_{j=1} \Gamma_j \cap Q_\mu$. In what follows we drop the subscript $\mu$ for notational convenience. First by approximation of Sobolev functions on Lipschitz sets (see, e.g., \cite[Section 4.2]{EvansGariepy92}) we may assume that $u|_{R_j}$ is smooth for $j=0,\ldots,n$, where $R_0 = Q \setminus \bigcup^n_{j=1} R_j$. We only consider the part $\partial' Q = (-\mu,\mu) \times \lbrace  \mu \rbrace$ of the boundary. Let $\pi_x = \lbrace x \rbrace \times  \R$ and compute for the second component $u_2$ by a slicing argument in $\e_2$-direction:
\begin{align*}
\int^\mu_{-\mu}\int^\mu_{-\mu} &|u_2(x,\mu) - u_2(x,y)|^2 \, dx\,dy  = \int^\mu_{-\mu}\int^\mu_{-\mu} \Big| \int^\mu_y D_2 u_2(x,t) \, dt\Big|^2 \, dx\,dy \\
& \le C\int^\mu_{-\mu}\int^\mu_{-\mu} \Big( \mu\int^\mu_{-\mu} |\partial_2 u_2(x,t)|^2 \, dt + \Big(\sum_{z \in J_u \cap \pi_x} |[u](z)|\Big)^2\Big) \, dx\,dy \\
& \le C\mu^2 \Vert  e(u)\Vert^2_{L^2(Q)} + C\mu \int^\mu_{-\mu} \Big(\sum_{z \in J_u \cap \pi_x} |[u](z)|\Big)^2 \, dx. 
\end{align*}
In the second step we have used H\"older's inequality. We now estimate the term on the right side. As $\Gamma_j$ is a rectangle,  except for two $x$-values there are exactly two points $t^{1}_j, t^{2}_j \in \R$ such that $\Gamma_j \cap \pi_x = \lbrace (x,t^{1}_j),(x, t^{2}_j)\rbrace$ if $\Gamma_j \cap \pi_x \neq \emptyset$ . We write $\vert \Gamma_j \vert_{\cal _H} = {\cal H}^1(\Gamma_j)$, $\vert S \vert_{\cal _H} = \sum_j {\cal H}^1(\Gamma_j)$ for shorthand. Letting $z^{k,x}_j = (x,t^k_j) \in \R^2$ and setting $ |[u](z^{k,x}_j)| = 0$ if $z^{k,x}_j \notin Q \cap\Gamma_j$, we then obtain by the discrete version of Jensen's inequality  
\begin{align*}
\int^\mu_{-\mu} \Big(\sum_{z \in J_u \cap \pi_x} |[u](z)|\Big)^2 \, dx & = 4\int^\mu_{-\mu} \Big(\sum\nolimits_{j}\sum\nolimits_{k=1,2} \frac{\vert \Gamma_j \vert_{\cal H}}{2\vert S \vert_{\cal H}} \ |[u](z^{k,x}_j)|\frac{\vert S \vert_{\cal H}}{\vert \Gamma_j \vert_{\cal H}}  \Big)^2 \, dx \\
& \le 4\int^\mu_{-\mu} \sum\nolimits_{j}\sum\nolimits_{k=1,2} \frac{\vert \Gamma_j \vert_{\cal H}}{2\vert S \vert_{\cal H}} \ \Big(|[u](z^{k,x}_j)|\frac{\vert S \vert_{\cal H}}{\vert \Gamma_j \vert_{\cal H}} \Big)^2 \, dx \\
& \le 2\vert S \vert_{\cal H}\sum\nolimits_j \Big(\vert \Gamma_j \vert^{-1}_{\cal H}\int_{\Gamma_j \cap Q} |[u]|^2 \, d{\cal H}^1\Big).
\end{align*}
Consequently, letting $E$ be the right hand side of \eqref{rig-eq: trace1} we derive
\begin{align*}
 \int_{\partial' Q}|u_2|^2 \, d{\cal H}^1 \le \frac{C}{\mu}\Big(\int^\mu_{-\mu}\int^\mu_{-\mu} |u_2(x,\mu) - u_2(x,y)|^2 \, dx\,dy + \Vert u\Vert^2_{L^2(Q)} \Big) \le CE. 
 \end{align*}
The same argument with slicing in the directions $\xi_1 = \frac{1}{\sqrt{2}}(1,-1)$ and $\xi_1 = \frac{1}{\sqrt{2}}(-1,-1)$ yields
\begin{align*}
\int_{\partial'_1 Q}|u \cdot \xi_1|^2 \, d{\cal H}^1 \le CE, \ \ \ \  \int_{\partial'_2 Q}|u \cdot \xi_2|^2 \, d{\cal H}^1 \le CE, 
\end{align*}
where $\partial'_1 Q = (-\mu,0) \times \lbrace  \mu \rbrace$ and $\partial'_2 Q = (0,\mu) \times \lbrace  \mu \rbrace$. The claim now follows by combination of the previous estimates. \eop

\subsection{Notation}

We provide a list of frequently used notation for convenience of the reader.

{\footnotesize

\begin{multicols}{2}

\textbf{Energies}
\begin{itemize}
\item $\alpha(\cdot)$: Before \eqref{rig-eq: hat tau def}
\item $\hat{\alpha}(\cdot)$: \eqref{rig-eq: shorthand}
\item ${\cal E}(\cdot)$: \eqref{rig-eq: shorthand}
\end{itemize}

\textbf{Measures}
\begin{itemize}
\item $\vert \cdot \vert_{\cal H}, \vert \cdot \vert_\infty, \vert \cdot \vert_*$: \eqref{rig-eq: h*}
\item $\Vert \cdot \Vert_{\cal H}, \Vert \cdot \Vert_\infty, \Vert \cdot \Vert_*$: Below \eqref{rig-eq: h*}
\item $\omega(\cdot)$: Before \eqref{rig-eq: h*, omega}
\item $|\cdot|_\omega$, $\Vert \cdot \Vert_\omega$: \eqref{rig-eq: h*, omega}
\item $|\cdot|_\pi$, $\Vert \cdot \Vert_\pi$:  \eqref{rig-eq: VertP}
\end{itemize}

\textbf{Neighborhoods}
\begin{itemize}
\item ${\cal C}^t  = \lbrace Y^t_0, \ldots Y^t_{m-1} \rbrace$: \eqref{rig-eq: covering def}
\item $J^t = \lbrace Q^t_1, \ldots, Q^t_n \rbrace $: Before \eqref{rig-eq: covering def}
\item $M^t$: \eqref{rig-eq: MMhatdef}
\item $M^t_k$, $k=1,2$: \eqref{rig-eq: new7}
\item $N^t$, $N^t_{j,\pm}$: Before \eqref{rig-eq: tau def}
\item $\tilde{N}$: Before \eqref{rig-eq: new9}
\item $\Psi$, $\psi$: Section \ref{rig-sec: subsub,  dod-neigh}
\item $\Phi$: Section \ref{rig-sec: subsub,  dod-neigh}
\item $\bar{\tau} = \upsilon|\Gamma|_\infty$: \eqref{rig-eq: tau def}
\item $\tau$: \eqref{rig-eq: tau-bartau}
\item $\hat{\tau}$: \eqref{rig-eq: hat tau def}
\item ${\cal Y}^t$: Before \eqref{rig-eq: covering def}
\item ${\cal Y}'$: Before Lemma \ref{rig-lemma: D1}
\end{itemize}

\textbf{Parameters}
\begin{itemize}
\item $h_*$: \eqref{rig-eq: h*}
\item $\omega_{\rm min}:$ Before \eqref{rig-eq: h*, omega}
\item $q$: \eqref{rig-eq: MMhatdef2}
\item $r$: \eqref{rig-eq: D1.2.2}
\item $\upsilon$: \eqref{rig-eq: W prop}
\item $\hat{\upsilon}$: \eqref{rig-eq: D1.7}
\end{itemize}

\textbf{Sets and components}
\begin{itemize}
\item  $\Gamma$, $\Theta$: \eqref{rig-eq: Xdef}
\item $H(\cdot)$:  \eqref{rig-eq: no holes}
\item $R$: \eqref{rig-eq: W prop}(i),(v)
\item ${\cal T}_\eta(\cdot,\cdot)$: Before \eqref{eq: etatau2}
\item ${\cal U}^s$: \eqref {rig-eq: calV-s def} 
\item ${\cal V}^s$: \eqref{rig-eq: calV-s def2}
\item ${\cal V}^s_\lambda$: \eqref{rig-eq: W prop}
\item $X$: Before \eqref{rig-eq: no holes}
\end{itemize}

\end{multicols}

}
\end{appendix}

\textbf{Acknowledgements}  I am very grateful to Bernd Schmidt for many stimulating discussions and valuable comments from which the results of this paper and their exposition have benefited a lot. 


 \typeout{References}

\end{document}